\documentclass[11pt]{article}

\usepackage{amsfonts}
\usepackage{amssymb,amsmath,amsthm, amsfonts}

\def\sqr#1#2{{\vcenter{\vbox{\hrule height.#2pt
              \hbox{\vrule width.#2pt height#1pt \kern#1pt \vrule width.#2pt}
          \hrule height.#2pt}}}}
\def\sqr#1#2{{\vcenter{\vbox{\hrule height.#2pt
              \hbox{\vrule width.#2pt height#1pt \kern#1pt \vrule width.#2pt}
              \hrule height.#2pt}}}}
\def\3n{\negthinspace \negthinspace \negthinspace }
\def\2n{\negthinspace \negthinspace }
\def\1n{\negthinspace }

\def\={\buildrel \triangle \over =}

\def\limsup{\mathop{\overline{\rm lim}}}

\def\max{\mathop{\rm max}}
\def\min{\mathop{\rm min}}

\def\sup{\mathop{\rm sup}}
\def\inf{\mathop{\rm inf}}

\def\({\Big (}
\def\){\Big )}
\def\[{\Big[}
\def\]{\Big]}
\def\be{\begin{equation}}

\def\ee{\end{equation}}

\def\square#1{\vbox{\hrule\hbox{\vrule height#1%
     \kern#1\vrule}\hrule}}
\def\rectangle#1#2{\vbox{\hrule\hbox{\vrule height#1%
     \kern#2\vrule}\hrule}}

\font\tenbb=msbm10 \font\sevenbb=msbm7 \font\fivebb=msbm5

\newfam\bbfam
\scriptscriptfont\bbfam=\fivebb \textfont\bbfam=\tenbb
\scriptfont\bbfam=\sevenbb

\newtheorem{lemma}{Lemma}[section]
\newtheorem{remark}{Remark}[section]

\newtheorem{theorem}{Theorem}[section]
\newtheorem{corollary}{Corollary}[section]

\newtheorem{definition}{Definition}[section]
\newtheorem{proposition}{Proposition}[section]

\makeatletter

\newcommand{\Rmnum}[1]{\expandafter\@slowromancap\romannumeral #1@}
   
   \@addtoreset{equation}{section}
\makeatother

\begin{document}

\title{Mean-field BDSDEs and associated nonlocal semi-linear backward stochastic partial differential equations \footnotemark[1]}
\author{Rainer Buckdahn$^{1,2}$,\,\, Juan Li$^{3, \dag}$,\,\, Chuanzhi Xing$^{3, \dag}$ \\
 {$^1$\small Laboratoire de Math\'{e}matiques de Bretagne Atlantique, Univ Brest,}\\
	{\small UMR CNRS 6205, 6 avenue Le Gorgeu, 29200 Brest, France.}\\
	{$^2$\small  School of Mathematics, Shandong University,}
	{\small Jinan 250100, P. R. China.}\\
	{$^3$\small  School of Mathematics and Statistics, Shandong University, Weihai,}
	{\small Weihai 264209, P. R. China.}\\
 {\small{\it E-mails: rainer.buckdahn@univ-brest.fr,\,\ juanli@sdu.edu.cn,\,\ chuanzhixing@mail.sdu.edu.cn.}}
\date{November 15, 2021}
}
\renewcommand{\thefootnote}{\fnsymbol{footnote}}
\footnotetext[1]{Rainer Buckdahn is supported in part by the ``FMJH Program Gaspard Monge in optimization and operation research", and the ANR (Agence Nationale de la Recherche), France project ANR-16-CE40-0015-01. Juan Li is supported by the NSF of P.R. China (NOs. 12031009, 11871037), National Key R and D Program of China (NO. 2018YFA0703900), Shandong Province (No. JQ201202), NSFC-RS (No. 11661130148; NA150344).}
\footnotetext[2]{Corresponding authors.}
\maketitle

\begin{abstract} In this paper we investigate mean-field backward doubly stochastic differential equations (BDSDEs), i.e., BDSDEs whose driving coefficients also depend on the joint law of the solution process as well as the solution of an associated mean-field forward SDE. Unlike the pioneering paper on BDSDEs by Pardoux-Peng \cite{PP1994}, we handle a driving coefficient in the backward integral of the BDSDE for which the Lipschitz assumption with respect to the law of the solution is sufficient, without assuming that this Lipschitz constant is small enough.
Using the splitting method introduced in \cite{BLPR2017} for mean-field SDEs, we study the unique solutions $(Y^{t,\xi},Z^{t,\xi})$ and $(Y^{t,x,P_{\xi}},Z^{t,x,P_{\xi}})$ of our BDSDEs.
Under suitable regularity assumptions on the coefficients we investigate the first and the second
order derivatives of the solution $(Y^{t,x,P_{\xi}},Z^{t,x,P_{\xi}})$  with respect to $x$, the derivative $(\partial_{\mu}Y^{t,x,P_{\xi}}(y),\partial_{\mu}Z^{t,x,P_{\xi}}(y))$ of the
solution process with respect to the measure $P_{\xi}$, and the derivative of
$(\partial_{\mu}Y^{t,x,P_{\xi}}(y),\partial_{\mu}Z^{t,x,P_{\xi}}(y))$  with respect to $y$. However, as the parameters  $(x,P_\xi)$ and $(x,P_\xi,y)$ run an infinite-dimensional space, unlike Pardoux and Peng, we cannot apply Kolmogorov's continuity criterion to the value function $V(t,x,P_{\xi}):=Y_t^{t,x,P_{\xi}}$, while in the classical case studied in \cite{PP1994} the value function $V(t,x)=Y_t^{t,x}$ can be shown to be of class $C^{1,2}([0,T]\times\mathbb{R}^d)$, we have for our value function $V(t,x,P_{\xi})$ and its derivative $\partial_\mu V(t,x,P_{\xi},y)$ only the $L^2$-differentiability  with respect to $x$ and $y$, respectively. However, we have to use the (mean-field) It\^{o} formula. To overcome this problem the characterisation of $V=(V(t,x,P_{\xi}))$ as the unique solution of the associated mean-field backward stochastic PDE uses the $C_b^{1,2,2}$-functions $\Psi(t,x,P_{\xi}):=E[V(t,x,P_{\xi})\cdot\eta]$ for suitable $\eta\in L^{\infty}(\mathcal{F};\mathbb{R})$. Using a similar idea, we extend the classical mean-field It\^{o} formula to smooth functions of solutions of mean-field BDSDEs.
\end{abstract}



\section{Introduction}

Backward doubly stochastic differential equations (BDSDEs, for short) were introduced by Pardoux and Peng in their pioneering paper \cite{PP1994} in 1994. These equations are a generalisation of backward SDEs involving in addition to the
(forward) It\^{o} integral driven by a Brownian motion $W$ a backward one governed by an independent Brownian motion $B$.
Pardoux and Peng showed namely that the BDSDE gives a (doubly) stochastic interpretation for the classical solution of the following backward stochastic partial differential equation (SPDE, for short):
\begin{equation}\label{eq1.1}
  \begin{split}
    u(t,x)=&\ \Phi(x) +\int_t^T\Big(\mathcal{L}u(s,x)+f(x,u(s,x),(\nabla u\sigma)(s,x) )\Big)ds\\
    &\ +\int_t^Tg(x,u(s,x),(\nabla u\sigma)(s,x) )d\overleftarrow{B_s},\ t\in[0,T],
  \end{split}
\end{equation}
where the second-order differential operator $\mathcal{L}$ is the generator of the forward diffusion SDE associated with the BDSDE.

Stimulated by this pioneering work, a great number of researchers have been attracted by this topic and studied different extensions. Without going into details, let us cite here, for instance, the works by  Bally and
Matoussi \cite{BM2001}, Zhang and Zhao \cite{ZZ2013} and Matoussi, Piozin and Popier \cite{MPP2017}. Other works on BDSDEs were motivated, e.g., by Zakai equation in filtering \cite{LS2001a,LS2001b}, by stochastic control with partial observations or pathwise stochastic control theory \cite{LS1998}.
For other recent developments on BDSDEs, we refer, e.g., to Buckdahn and Ma \cite{BM2001a,BM2001b}, Han, Peng and Wu \cite{HPW2010}, but also Li and Xing \cite{LX2021}, Li, Xing and Peng \cite{LXP2021}, Shi, Wen and Xiong \cite{SWX2020}, Shi, Gu and Liu \cite{SGL2005}, should be mentioned.

The study of mean-field stochastic differential equations (SDEs), known also as McKean-Vlasov SDEs, goes back to the work by Kac \cite{K1956} in 1956. Since then, stimulated by numerous applications, the theory of mean-field SDEs has been dynamically developing. In recent years, new impulses to this research were given by the course given by P.L. Lions at \emph{Coll\`{e}ge de France} \cite{L} (We also refer to the notes by Cardaliaguet \cite{C2012}), in which the author introduced the notion of differentiability with respect to (w.r.t.) the probability measure of a function defined over a space of probability laws with finite second order moment.
Inspired by \cite{L, C2012}, in 2017, Buckdahn, Li, Peng and Rainer \cite{BLPR2017} considered a general mean-field SDE whose coefficients depend on both the solution process and its law. They studied the properties of the solution of this SDE, and in particular, the associated non-local PDEs of mean-field type. Hao and Li \cite{HL2016} extended these studies to mean-field SDEs which, in addition to the driving Brownian motion, are governed by a compensated Poisson random measure. In \cite{L2018}, Li extended this work to general mean-field forward-backward SDEs (FBSDEs) with jumps, in which the coefficients of both the forward and also the backward SDE depend on the solution processes but also on their joint law. She showed in particular that such mean-field FBSDEs give a stochastic interpretation to nonlocal integral-partial differential equations.
We emphasize that the value functions in \cite{BLPR2017}, \cite{HL2016} and \cite{L2018} are deterministic functions.

On the other hand, inspired by the seminal paper \cite{LL2007} by Lasry and Lions (2007),
Buckdahn, Djehiche, Li and Peng \cite{BDLP2009} investigated a nonlinear mean-field BSDE. Since then,
motivated by its various applications, the theory of mean-field BSDEs has been studied by numerous researchers, among them Buckdahn, Li and Peng \cite{BLP2009}, Carmona and Delarue \cite{CD2014}, Li, Liang and Zhang \cite{LLZ2018}, Chen, Xing and Zhang \cite{CXZ2018}. Li and Xing \cite{LX2021} extended these investigations to general mean-field BDSDEs

\begin{equation}\label{eq1.2}
  \begin{split}
  Y_t=\xi&+\int_t^Tf(s,P_{(Y_s,Z_s)},Y_s,Z_s)ds+\int_t^Tg(s,P_{(Y_s,Z_s)},Y_s,Z_s)d\overleftarrow{B_s}\\
  &+\int_t^Th(s,P_{(Y_s,Z_s)})d\overleftarrow{B_s}-\int_t^TZ_sdW_s,\ 0\leq t\leq T.
  \end{split}
\end{equation}
We emphasise that the Lipschitz constant of the coefficient $h$ can be arbitrary without necessarily being small enough.
For more details we refer to Appendix A.2.

Inspired by above works, in this paper we study a general mean-field BDSDE associated with a forward diffusion SDE as well as the related backward SPDEs of mean-field type.
More precisely, we consider the solutions $(X^{t,\xi},X^{t,x,P_{\xi}})$ of the split forward SDE (see, \eqref{eq3.1} and \eqref{eq3.2})
and those of the split BDSDE $(Y^{t,\xi},Z^{t,\xi})$ and $(Y^{t,x,\xi},Z^{t,x,\xi})$ (see, \eqref{eq4.1} and \eqref{eq4.2}). Their existence is stated in Theorem \ref{thA.1} (Appendix A.2), and the fact that  $(Y^{t,x,\xi},Z^{t,x,\xi})$ depends on $\xi$ only through its law ($(Y^{t,x,\xi},Z^{t,x,\xi})=(Y_s^{t,x,P_{\xi}},\ Z_s^{t,x,P_{\xi}})$) is proven by Proposition 4.1. Besides of the study of regularity properties of the solution of the mean-field BDSDE, our main objective is the characterization of the value function $V(t,x,P_{\xi}):=Y_t^{t,x,P_{\xi}}$, a stochastic process adapted to the backward filtration generated by the Brownian motion $B$, as unique solution of the associated mean-field semi-linear backward stochastic PDE (SPDE) which is a new type of SPDE
\begin{equation}\label{SPDE}
	\begin{split}
		&\  V(t,x,P_{\xi})\\
		=&\ \Phi(x,P_{\xi})+\int_t^T\Bigg\{\sum_{i=1}^d\partial_{x_i}V(s,x,P_{\xi})b_i(x,P_{\xi})
		+\frac{1}{2}\sum_{i,j,k=1}^d(\partial_{x_ix_j}^2V)(s,x,P_{\xi})(\sigma_{i,k}\sigma_{j,k})(x,P_{\xi})\\
		&\ + f(x,V(s,x,P_{\xi}),\sum_{i=1}^d\partial_{x_i}V(s,x,P_{\xi})\sigma_i(x,P_{\xi}),P_{(\xi,\psi(s,\xi,P_{\xi}))}) \\
		&\ +E\[\sum_{i=1}^d(\partial_{\mu}V)_i(s,x,P_{\xi},\xi)b_i(\xi,P_{\xi})
		+\frac{1}{2}\sum_{i,j,k=1}^d\partial_{y_i}(\partial_{\mu}V)_j(s,x,P_{\xi},\xi)(\sigma_{i,k}\sigma_{j,k})(\xi,P_{\xi}) \] \Bigg\}ds\\
		&\ + \int_t^T \sum_{j=1}^lg_j(x,V(s,x,P_{\xi}),\sum_{i=1}^d\partial_{x_i}V(s,x,P_{\xi})\sigma_i(x,P_{\xi}),P_{(\xi,\psi(s,\xi,P_{\xi}))})d\overleftarrow{B_s^j}\\
		&\ + \int_t^T \sum_{j=1}^lh_j(P_{(\xi,\psi(s,\xi,P_{\xi}))})d\overleftarrow{B_s^j},\ \ (t,x,\xi,P_{\xi})\in[0,T]\times\mathbb{R}^d\times L^2(\mathcal{G}_t;\mathbb{R}^d)\times\mathcal{P}_2(\mathbb{R}^d),
	\end{split}
\end{equation}
where $\psi(s,x,P_{\xi}):=(V(s,x,P_{\xi}),\sum_{i=1}^d\partial_{x_i}V(s,x,P_{\xi})\sigma_i(x,P_{\xi}))$.
Observed that, when the coefficients don't depend on the law, \eqref{SPDE} reduces to \eqref{eq1.1}. Following the scheme given by Pardoux and Peng in their seminal paper \cite{PP1994}, we study the regularity of $(Y^{t,x,P_{\xi}},Z^{t,x,P_{\xi}})$. However, apart from the fact that Pardoux and Peng restricted to the proof of the first order derivative of $(Y^{t,x},Z^{t,x})$ w.r.t. $x$ and only stated the second order differentiability, to make the proofs not only for the first order derivative of $(Y^{t,x,P_{\xi}},Z^{t,x,P_{\xi}})$ w.r.t. $x$ and w.r.t. the measure, but we discuss also all the details of the proof for its second order derivative w.r.t. $x$ and the derivative w.r.t. $y$ of the derivative w.r.t. the measure  $(\partial_{\mu}Y^{t,x,P_{\xi}}(y),\partial_{\mu}Z^{t,x,P_{\xi}}(y))$, and this is related with different subtle technicalities. In particular, Malliavin calculus will be used to prove some crucial estimates for $Z^{t,x,P_{\xi}}$ and its derivatives (see, the Propositions \ref{prop4.2}, \ref{prop8.1} and \ref{prop6.1+1}). Another difficulty comes from the fact that the derivatives of $(Y^{t,x,P_{\xi}},Z^{t,x,P_{\xi}})$ w.r.t. $x$ (of first and second order) and the (first order) derivative of $(\partial_{\mu}Y^{t,x,P_{\xi}}(y),\partial_{\mu}Z^{t,x,P_{\xi}}(y))$ w.r.t. $y$ are $L^2$-derivatives, i.e., we prove
the $L^2$-regularity of the value function $V(t,x,P_{\xi}):=Y_t^{t,x,P_{\xi}}$. However, we have to apply the (mean-field) It\^{o} formula to $V\big(s,X_s^{t,x,P_{\xi}},P_{X_s^{t,\xi}}\big).$ Pardoux and Peng concluded that their value function $V(t,x)$ is a $C^{1,2}$-function by using Kolmogorov's continuity criterion. But we cannot use it, because in our case the parameter $(t,x,P_\xi)$ of $V=V(t,x,P_{\xi})$ runs an infinite-dimensional space. To overcome that we have for $V$ only $L^2$-regularity, we introduce the deterministic function $\Psi(t,x,P_{\xi}):=E[V(t,x,P_{\xi})\cdot\eta]$, for suitable $\eta\in L^{\infty}(\mathcal{F};\mathbb{R})$, and for $\Psi$ we can show that it is a $C_b^{1,2,2}$-function. Applying the (mean-field) It\^{o} formula to $\Psi\big(s,X_s^{t,x,P_{\xi}},P_{X_s^{t,\xi}}\big)$ will be crucial for the proof that $V$ is a solution of SPDE (\ref{SPDE}), a classical one but with the derivatives w.r.t. $x$ and $y$ in $L^2$-sense. By defining a suitable space to which the solution process $V$ belongs, the same technique will allow to prove the uniqueness.

Last but not least, let us mention that we extend the classical mean-field It\^{o} formula (see \cite{BLPR2017}) to $C_b^{1,2,2}$-functions $F(s,U_s,P_{Y_s})$ applied to solutions of BDSDEs $U=(U_s)$ and $Y=(Y_s)$  (Theorem \ref{th2.1}). Our It\^{o} formula extends namely that by Pardoux-Peng \cite{PP1994} who studied the case $F(s,U_s)=|U_s|^2.$ As, on the other hand, a classical SDE (with backward integral) is a special case of a BDSDE, our It\^{o} formula can also be regarded as extension of the classical mean-field It\^{o} SDE. The proof of Theorem \ref{th2.1} is inspired by our approach in the proof that $V$ is the unique solution of SPDE (\ref{SPDE}).

Our paper is organized as follows: Section 2 is on one hand devoted to preliminaries, but on the other hand we state there also our mean-field It\^{o} formula extended to solutions of BDSDEs.
While Section 3 gives a recall on mean-field SDEs,
Section 4 is devoted to the investigation of the existence and the uniqueness of solutions
and the corresponding estimates for our split mean-field BDSDEs. Section 5 gives a recall on the first
order derivatives of the process $X^{t,x,P_{\xi}}$ with respect to $x$ and the measure $P_{\xi}$, and the
corresponding estimates. In Section 6 the first order derivatives
of $(Y^{t,x,P_{\xi}},Z^{t,x,P_{\xi}})$ w.r.t. $x$ and w.r.t. $P_{\xi}$ are investigated.
Section 7 is devoted to the second order derivatives of $X^{t,x,P_{\xi}}$, and
using an additional assumption on $g$ (Assumption (H8.2)), the second order derivatives of $(Y^{t,x,P_{\xi}},Z^{t,x,P_{\xi}})$ are discussed in Section 8.
Basing on our mean-field It\^{o}'s formula, in Section 9 we prove that the value function $V(t,x,P_{\xi})$
is the unique classical solution of the SPDE \eqref{SPDE}.
Finally, in the Appendix we give the proof of Theorem \ref{th2.1}
(Appendix A.1), we recall some basic results on mean-field BDSDEs (Appendix A.2), and we study the special case of BDSDEs with a coefficient $g$ which is affine in $z$ (Appendix A.3).

\section{Preliminaries}
Let $T>0$ be a fixed time horizon and $(\Omega,\mathcal{F},P)$ be a complete probability space. Let $\{W_t,\ 0\leq t\leq T\}$ and $\{B_t,\ 0\leq t\leq T\}$ be two mutually independent standard Brownian motions defined on $(\Omega,\mathcal{F},P)$, with values in $\mathbb{R}^d$ and in $\mathbb{R}^l$, respectively. We assume that there is a sub-$\sigma$-field $\mathcal{F}^0\subset\mathcal{F}$, containing all $P$-null subsets of $\mathcal{F}$, such that\\
\indent\quad(i)\ \ the Brownian motion $(B,W)$ is independent of $\mathcal{F}^0$;

\indent\quad(ii)\ $\mathcal{F}^{0}$ is `rich enough', i.e., $\mathcal{P}_{2}(\mathbb{R}^{k})=\{P_{\xi},\xi\in L^{2}(\mathcal{F}^{0};\mathbb{R}^{k})\},\ k\geq1$.\\
\noindent Here $P_\xi := P\circ [\xi]^{-1}$ denotes the law of the random variable $\xi$ under the probability $P$.

For $0\leq t\leq s\leq T$, we define $\mathcal{F}_{t,s}:=\mathcal{F}_{t,s}^W\vee\mathcal{F}_{s,T}^B\vee\mathcal{F}^0$, $\mathcal{F}_{s}:=\mathcal{F}_{0,s}$,
 $\mathcal{G}_{s}:=\mathcal{F}_{s}^W\vee\mathcal{F}^0$ and $\mathcal{D}_{s}:=\mathcal{F}_{s}^W\vee\mathcal{F}_{T}^B\vee\mathcal{F}^0$, where for any process $\{\eta_t\}$, $\mathcal{F}_{t,s}^{\eta}=\sigma\{\eta_r-\eta_t;t\leq r\leq s\}$ augmented by the $P$-null sets, and $\mathcal{F}_t^{\eta}=\mathcal{F}_{0,t}^{\eta}$. It should be noted that $\{\mathcal{F}_t,t\in [0,T]\}$ is neither increasing nor decreasing, and so it does not constitute a filtration.
 However, $\{\mathcal{G}_t,t\in [0,T]\}$ as well as $\{\mathcal{D}_t,t\in [0,T]\}$ are filtrations.
 For $n\in \mathbb{N}$\ and\ $x,y\in\mathbb{R}^{n}$, we denote the norm and inner product, respectively, by $|x|=\Big(\sum\limits_{i=1}^{n}x_{i}^{2}\Big)^\frac{1}{2}$ and $\langle x,y\rangle=\sum\limits_{i=1}^{n}x_{i}y_{i}.$

\indent Let us introduce some notations and concepts, which are used frequently in what follows. Recall that $\mathcal{P}_{2}(\mathbb{R}^{k})$ is the set of
all probability measures $\mu$ on $(\mathbb{R}^{k},\mathcal{B}(\mathbb{R}^{k}))$ with finite second moment, i.e., $\displaystyle\int_{\mathbb{R}^k}|x|^{2}\mu(dx)<\infty $. Here $\mathcal{B}(\mathbb{R}^{k})$ denotes the Borel $\sigma$-field over $\mathbb{R}^{k}$. $\mathcal{P}_{2}(\mathbb{R}^{k})$ is endowed with the $2$-Wasserstein metric: For $\mu,\nu\in \mathcal{P}_{2}(\mathbb{R}^k)$,
\begin{equation}\label{eq2.1}
W_{2}(\mu,\nu):=\inf\Big\{\Big(\int_{\mathbb{R}^k\times\mathbb{R}^k}|x-y|^{2}\rho(dxdy)\Big)^{\frac{1}{2}}: \rho\in\mathcal{P}_{2}(\mathbb{R}^{2k}), \rho(.\times\mathbb{R}^k)=\mu,\rho(\mathbb{R}^k\times.)=\nu\Big\}.
\end{equation}
We will also need a weighted $2$-Wasserstein metric on $\mathcal{P}_{2}(\mathbb{R}^{k}\times \mathbb{R}^{k\times d})$:
For $\mu,\mu'\in\mathcal{P}_{2}(\mathbb{R}^{k}\times \mathbb{R}^{k\times d})$,
\begin{equation}\label{eq2.2}
\begin{split}
  W_{2,\gamma_1,\gamma_2}(\mu,\mu'):=\inf\Big\{\Big(E[\gamma_1|\xi-\xi'|^{2}+\gamma_2|\eta-\eta'|^{2}]\Big)^{\frac{1}{2}}\Big|
(\xi,\eta),(\xi',\eta')&\in L^2(\mathcal{F};\mathbb{R}^{k}\times\mathbb{R}^{k\times d}):\\
&P_{(\xi,\eta)}=\mu,P_{(\xi',\eta')}=\mu'\Big\},
\end{split}
\end{equation}
for any fixed $\gamma_1,\gamma_2>0$. It is obvious that $W_{2,\gamma_1,\gamma_2}(\cdot,\cdot)$ is not only a metric but also equivalent to $W_{2}(\cdot,\cdot)$.
Indeed, $(\gamma_1\wedge\gamma_2)^{\frac{1}{2}}W_{2}(\mu,\mu')\leq W_{2,\gamma_1,\gamma_2}(\mu,\mu')\leq (\gamma_1\vee\gamma_2)^{\frac{1}{2}}W_{2}(\mu,\mu')$.

We shall also introduce the following spaces of stochastic processes: For $p>1$, $t\in[0,T]$,\\
\noindent$\bullet \ L^p(\Omega,\mathcal{F}_T, P;\mathbb{R}^d)$ is the set of $\mathcal{F}_T$-measurable random variables $ \xi: \Omega\rightarrow\mathbb{R}^d$ such that $\|\xi\|_{L^{p}}:=$\\
\noindent\mbox{ } \ $\left (E\left[|\xi|^{p}\right]\right)^{\frac{1}{p}}<\infty $.\\
\noindent$\bullet\ \mathcal{S}_{\mathcal{F}}^p(t,T;\mathbb{R}^d)$ is the set of $\{\mathcal{F}_{t,s}\}$-adapted measurable continuous processes  $\eta: \Omega\times[t,T]\rightarrow \mathbb{R}^d$ with\\
\noindent\mbox{ } \ $\|\eta\|_{\mathcal{S}^{p}}:= \big(E\big[\sup\limits_{t\leq s\leq T}|\eta(s)|^{p}\big]\big)^{\frac{1}{p}}<\infty$.\\
\noindent$\bullet\ \mathcal{H}_{\mathcal{F}}^p(t,T;\mathbb{R}^d)$ is the set of $\{\mathcal{F}_{t,s}\}$-adapted measurable processes $\eta: \Omega\times[t,T]\rightarrow \mathbb{R}^d$ with $\|\eta\|_{\mathcal{H}^{p}}:=$\\
\noindent\mbox{ } \ $\displaystyle{\big( E\big[(\int_{t}^{T}|\eta(s)|^{2}ds)^{\frac{p}{2}}\big]\big)^{\frac{1}{p}}<\infty}$.\\
\noindent$\bullet\ C_b^k(\mathbb{R}^p,\mathbb{R}^q)$ is the set of functions of class $C^k$ from $\mathbb{R}^p$ into $\mathbb{R}^q$ whose partial derivatives of all \\
\noindent\mbox{ } \  order less than or equal to $k$ are bounded.\\
\noindent$\bullet\ \mathcal{S}$ is the set of smooth random variables $\xi$ of the form $\xi=\varphi(W(h_1),\cdots,W(h_n);B(k_1),\cdots,B(k_p))$, \\
\noindent\mbox{ } \ $n,\ p\geq0$, with $\varphi\in C_b^{\infty}(\mathbb{R}^{n+p},\mathbb{R})$, $h_1,\cdots,h_n\in L^2(0,T;\mathbb{R}^d)$, $k_1,\cdots,k_p\in L^2(0,T;\mathbb{R}^l)$,  $W(h_i):=$ \\
\noindent\mbox{ } \ $\displaystyle\int_0^T\langle h_i(t),dW_t\rangle$, $\displaystyle B(k_j):=\int_0^T\langle k_j(t),dB_t\rangle$.

 Moreover, if $\xi\in S$ is of the above form, its Malliavin derivative w.r.t. $W$, denoted by $D[\cdot]:=D^{W}[\cdot]$, is given by
 $ \displaystyle D_t\xi=\sum_{i=1}^n\frac{\partial\varphi}{\partial x_i}(W(h_1),\cdots,W(h_n);B(k_1),\cdots,B(k_p))h_i(t),\ 0\leq t\leq T.$\
For $\xi\in S$, $p>1$, we define the norm $\displaystyle ||\xi||_{1,p}=\(E\[|\xi|^p+\big(\int_0^T|D_t\xi|^2dt\big)^{\frac{p}{2}}\]\)^{\frac{1}{p}}.$
From \cite{N1995} we know the operator $D$ has a closed extension to the space $\mathbb{D}^{1,p}$, the closure of $S$ with respect to the norm $||\cdot||_{1,p}$.
Furthermore, we put $\mathbb{D}^{1,\infty}=\displaystyle{\bigcap_{p\geq 2}\mathbb{D}^{1,p}}$.
Observe that if $\xi\in\mathbb{D}^{1,2}$ is $\mathcal{F}_t$-measurable, $D_{\theta}\xi= 0$, $d\theta d P$-a.e., $\theta\in (t,T]$.
We shall denote by $D_{\theta}^i\xi$, the $i$-th component of $D_{\theta}\xi$, $1\leq i\leq d$.

We now recall the notion of differentiability with respect to probability measure of a function defined on $\mathcal{P}_{2}(\mathbb{R}^{d})$.
There are different definitions, we use that introduced by Lions \cite{L}, also refer to the notes by Cardaliaguet \cite{C2012}.
For more details the readers also may refer to Buckdahn, Li, Peng and Rainer \cite{BLPR2017}, Hao and Li \cite{HL2016}, Li \cite{L2018}.
Given a function $\varphi:\mathcal{P}_{2}(\mathbb{R}^{d})\rightarrow\mathbb{R}$, we consider the lifted function
$\widetilde{\varphi}(\xi):=\varphi(P_{\xi}),\ \xi\in L^{2}(\mathcal{F};\mathbb{R}^d) (:=L^{2}(\Omega,\mathcal{F},P;\mathbb{R}^d))$.
If for a given $\mu_0\in\mathcal{P}_{2}(\mathbb{R}^{d})$ there exists a random variable $\xi_0\in L^{2}(\mathcal{F};\mathbb{R}^d)$ satisfying $P_{\xi_0}=\mu_0$, such
that $\widetilde{\varphi}:L^{2}(\mathcal{F};\mathbb{R}^d)\rightarrow\mathbb{R}$ is Fr\'{e}chet differentiable in $\xi_{0}$, $\varphi$ is said to be
differentiable with respect to $\mu_0$. This is equivalent with the existence of a continuous linear mapping
$D\widetilde{\varphi}(\xi_{0}):L^{2}(\mathcal{F};\mathbb{R}^d)\rightarrow\mathbb{R}^d$
(i.e., $D\widetilde{\varphi}(\xi_{0})\in L(L^{2}(\mathcal{F};\mathbb{R}^d),\mathbb{R}^d))$ such that
\begin{equation}\label{eq2.3}
\widetilde{\varphi}(\xi_{0}+\zeta)-\widetilde{\varphi}(\xi_{0})=D\widetilde{\varphi}(\xi_{0})(\zeta)+o(|\zeta|_{L^{2}}),
\end{equation}
for $\zeta\in L^{2}(\mathcal{F};\mathbb{R}^d)$ with $|\zeta|_{L^{2}}\rightarrow0$. Riesz's Representation Theorem allows to show that there
exists a unique $\eta\in L^{2}(\mathcal{F};\mathbb{R}^d)$ such that $D\widetilde{\varphi}(\xi_{0})(\zeta)=E[\eta\cdot\zeta]$,
$\zeta\in L^{2}(\mathcal{F};\mathbb{R}^d)$. It was shown by Lions that $\eta$ is a Borel measurable function of $\xi_0$, refer to Cardaliaguet \cite{C2012},
 $\eta=\psi(\xi_0)$, where $\psi$ is a Borel measurable function depending on $\xi_0$ only through its
law. Combining \eqref{eq2.3} and the above argument, we have
\begin{equation}\label{eq2.4}
\varphi(P_{\xi_0 +\zeta})-\varphi(P_{\xi_{0}})=E[\psi(\xi_{0})\cdot\zeta]+o(|\zeta|_{L^{2}}).
\end{equation}
In the spirit of Lions and Cardaliaguet, the derivative of $\varphi:\mathcal{P}_{2}(\mathbb{R}^{d})\rightarrow\mathbb{R}$ with respect to the
measure $P_{\xi_0}$ is denoted by $\partial_{\mu}\varphi(P_{\xi_{0}},y):=\psi(y),\ y\in\mathbb{R}^d$. Observe that $\partial_{\mu}\varphi(P_{\xi_{0}},y)$ is only
$P_{\xi_0}(dy)$-a.e. uniquely determined; see also Definition 2.1 in Buckdahn, Li, Peng and Rainer \cite{BLPR2017}. We also mention that an equivalent approach for the derivative was developed by Cardaliaguet, Delarue, Lasry and P.L. Lions \cite{C2019}, interested readers are referred to this work.

Now we introduce several spaces which will be used frequently.
\begin{definition} \label{def2.1}
\emph{(1)} We say that $\varphi$ belongs to $C^{1}_{b}(\mathcal{P}_{2}(\mathbb{R}^{d}))$, if $\varphi:\mathcal{P}_{2}(\mathbb{R}^{d})\rightarrow\mathbb{R}$
is differentiable on $\mathcal{P}_{2}(\mathbb{R}^{d})$ and
$\partial_{\mu}\varphi(\cdot,\cdot):\mathcal{P}_{2}(\mathbb{R}^{d})\times \mathbb{R}^{d}\rightarrow\mathbb{R}^{d}$
is bounded and Lipschitz continuous, i.e., there exists some positive constant $L$ such that\\
\indent \emph{(i)} $|\partial_{\mu}\varphi(\mu,y)|\leq L,\ \mu\in\mathcal{P}_{2}(\mathbb{R}^{d}),\ y\in\mathbb{R}$,\\
\indent \emph{(ii)} $|\partial_{\mu}\varphi(\mu,y)-\partial_{\mu}\varphi(\mu',y')|\leq L(W_2(\mu,\mu')+|y-y'|),\
 \mu,\mu'\in\mathcal{P}_{2}(\mathbb{R}^{d}),\ y,y'\in\mathbb{R}$.\\
\emph{(2)} By $C^{2}_{b}(\mathcal{P}_{2}(\mathbb{R}^{d}))$ we denote the space of all functions $\varphi\in C^{1}_{b}(\mathcal{P}_{2}(\mathbb{R}^{d}))$
such that $(\partial_{\mu}\varphi)_j(\mu,\cdot): \mathbb{R}^d\rightarrow\mathbb{R}$ is differentiable, for every $\mu\in\mathcal{P}_{2}(\mathbb{R}^{d})$,
and the derivative $\partial_{y}\partial_{\mu}\varphi:\mathcal{P}_{2}(\mathbb{R}^{d})\times\mathbb{R}^{d}\rightarrow\mathbb{R}^{d}\otimes \mathbb{R}^{d}$
is bounded and Lipschitz continuous.\\
We will use the notation $\partial_{\mu}\varphi(\mu,y):=\((\partial_{\mu}\varphi)_j(\mu,y)\)_{1\leq j\leq d},\ (\mu,y)\in\mathcal{P}_{2}(\mathbb{R}^{d})\times\mathbb{R}^{d}$.
\end{definition}

\begin{definition} \label{def2.2}
\emph{(1)} We say that $\varphi$ belongs to $C_b^{0,2,2}([0,T]\times\mathbb{R}^{d}\times\mathcal{P}_{2}(\mathbb{R}^{d}))$, if $\varphi:[0,T]\times\mathbb{R}^{d}\times\mathcal{P}_{2}(\mathbb{R}^{d})\rightarrow\mathbb{R}$ satisfies\\
\indent \emph{(i)} $\varphi(\cdot,\cdot,\mu)\in C_b^{0,2}([0,T]\times\mathbb{R}^{d})$, for all $\mu\in\mathcal{P}_{2}(\mathbb{R}^{d})$;\\
\indent \emph{(ii)} $\varphi(t,x,\cdot)\in C^{2}_{b}(\mathcal{P}_{2}(\mathbb{R}^{d}))$, for all $(t,x)\in[0,T]\times\mathbb{R}^{d}$;\\
\indent \emph{(iii)} All derivatives of order $1$ and $2$ are continuous on $[0,T]\times\mathbb{R}^{d}\times\mathcal{P}_{2}(\mathbb{R}^{d})\!\times\mathbb{R}^{d}$,
and uniformly bounded over $[0,T]\times\mathbb{R}^{d}\times\mathcal{P}_{2}(\mathbb{R}^{d})\times\mathbb{R}^{d}$,
$\partial_{\mu}\varphi$ and $\partial_{y}(\partial_{\mu}\varphi)$ are Lipschitz continuous w.r.t. $(x,\mu,y)$, uniformly w.r.t. $t$.\\
\emph{(2)}  We say that $\varphi$ belongs to $C_b^{1,2,2}([0,T]\times\mathbb{R}^{d}\times\mathcal{P}_{2}(\mathbb{R}^{d}))$, if $\varphi:[0,T]\times\mathbb{R}^{d}\times\mathcal{P}_{2}(\mathbb{R}^{d})\rightarrow\mathbb{R}$ satisfies\\
\indent \emph{(i)} $\varphi(\cdot,\cdot,\mu)\in C_b^{1,2}([0,T]\times\mathbb{R}^{d})$, for all $\mu\in\mathcal{P}_{2}(\mathbb{R}^{d})$;\\
\indent \emph{(ii)} $\varphi(t,x,\cdot)\in C^{2}_{b}(\mathcal{P}_{2}(\mathbb{R}^{d}))$, for all $(t,x)\in[0,T]\times\mathbb{R}^{d}$;\\
\indent \emph{(iii)} All derivatives of order $1$ and $2$ are continuous on $[0,T]\times\mathbb{R}^{d}\times\mathcal{P}_{2}(\mathbb{R}^{d})\times\mathbb{R}^{d}$,
and uniformly bounded over $[0,T]\times\mathbb{R}^{d}\times\mathcal{P}_{2}(\mathbb{R}^{d})\times\mathbb{R}^{d}$,
$\partial_{\mu}\varphi$ and $\partial_{y}(\partial_{\mu}\varphi)$ are Lipschitz continuous w.r.t. $(x,\mu,y)$, uniformly w.r.t. $t$.\\
\emph{(3)}  We say that $\varphi$ belongs to $C^{1,2}_b([0,T]\times\mathcal{P}_{2}(\mathbb{R}^{d}))$, if $\varphi:[0,T]\times\mathcal{P}_{2}(\mathbb{R}^{d})\rightarrow\mathbb{R}$ satisfies\\
\indent \emph{(i)} $\varphi(\cdot,\mu)\in C_b^{1}([0,T])$, for all $\mu\in\mathcal{P}_{2}(\mathbb{R}^{d})$;\\
\indent \emph{(ii)} $\varphi(t,\cdot)\in C^{2}_{b}(\mathcal{P}_{2}(\mathbb{R}^{d}))$, for all $t\in[0,T]$;\\
\indent \emph{(iii)} All derivatives of order $1$ and $2$ are continuous on $[0,T]\times\mathcal{P}_{2}(\mathbb{R}^{d})\times\mathbb{R}^{d}$,
and uniformly bounded over $[0,T]\times\mathcal{P}_{2}(\mathbb{R}^{d})\times\mathbb{R}^{d}$,
$\partial_{\mu}\varphi$ and $\partial_{y}(\partial_{\mu}\varphi)$ are Lipschitz continuous w.r.t. $(\mu,y)$, uniformly w.r.t. $t$.
\end{definition}

Now we briefly review the notion of a solution of a backward doubly stochastic differential equations (BDSDEs, for short) from Pardoux and Peng \cite{PP1994}.

Let $ f:\ [0,T]\times\Omega\times\mathbb{R}^k\times\mathbb{R}^{k\times d}\rightarrow\mathbb{R}^{k},\
  g:\ [0,T]\times\Omega\times\mathbb{R}^k\times\mathbb{R}^{k\times d}\rightarrow\mathbb{R}^{k\times l}$
be jointly measurable and satisfy:\\
\noindent\textbf{Assumption (H2.1)} (i) $(g(t,\cdot,0,0))_{t\in[0,T]}\in\mathcal{H}_{\mathcal{F}}^2(0,T;\mathbb{R}^{k\times l})$;\\
(ii) $g$ is Lipschitz in $(y,z)$, i.e., there exist constants $C>0$, and $0<\alpha<1$ such that for all\\
 \mbox{ } \ \ \    $t\in[0,T]$, $y_1,y_2\in\mathbb{R}^k$, $z_1,z_2\in\mathbb{R}^{k\times d}$, $P$-a.s.,
 $$|g(t,y_1,z_1)-g(t,y_2,z_2)|^2\leq C|y_1-y_2|^2+\alpha|z_1-z_2|^2;$$
(iii) $(f(t,\cdot,0,0))_{t\in[0,T]}\in\mathcal{H}_{\mathcal{F}}^2(0,T;\mathbb{R}^{k})$;\\
(iv) $f$ is Lipschitz in $(y,z)$, i.e., there exists a constant $C>0$ such that for all $t\in[0,T]$, $y_1,y_2\in\mathbb{R}^k$,\\
 \mbox{ } \ \ \  $z_1,z_2\in\mathbb{R}^{k\times d}$, $P$-a.s.,
      $$|f(t,y_1,z_1)-f(t,y_2,z_2)|\leq C(|y_1-y_2|+|z_1-z_2|).$$

\noindent Given $\xi\in L^2(\Omega,\mathcal{F}_T, P;\mathbb{R}^k)$, Pardoux and Peng \cite{PP1994} studied the following BDSDE:
\begin{equation}\label{eq2.4+1}
  \begin{split}
  Y_t=\xi&+\int_t^Tf(s,Y_s,Z_s)ds+\int_t^Tg(s,Y_s,Z_s)d\overleftarrow{B_s}-\int_t^TZ_sdW_s,\ 0\leq t\leq T,
  \end{split}
\end{equation}
where $W$ and $B$ are two independent Brownian motions, the integral with respect to $B$ is the It\^{o} backward one, denoted by $d\overleftarrow{B}$.
\begin{proposition} \label{prop2.1}
Under Assumption (H2.1), equation (\ref{eq2.4+1}) has a unique solution $(Y,Z)\in\mathcal{S}_{\mathcal{F}}^2(0,T;$ $\mathbb{R}^k)\times\mathcal{H}_{\mathcal{F}}^2(0,T;\mathbb{R}^{k\times d})$.
\end{proposition}
For more details, please, we refer to Theorem 1.1 in Pardoux and Peng \cite{PP1994} or Appendix A.2 of this paper.

 Next we give a general It\^{o}'s formula which will be used later.

\begin{theorem} \label{th2.1}
(It\^{o}'s formula). Let $F\in C_b^{1,2,2}([0,T]\times\mathbb{R}^{d}\times\mathcal{P}_{2}(\mathbb{R}^{d}))$.
Given $f\in\mathcal{H}^2_{\mathcal{F}}(0,T;\mathbb{R}^d)$, $g\in\mathcal{H}^2_{\mathcal{F}}(0,T;\mathbb{R}^{d\times l})$,
 $\xi\in L^{2}(\mathcal{F}_T;\mathbb{R}^d)$, as well as $u\in \mathcal{H}^2_{\mathcal{F}}(0,T;\mathbb{R}^d)$,
 $v\in \mathcal{H}^2_{\mathcal{F}}(0,T;\mathbb{R}^{d\times l})$, $\eta\in L^{2}(\mathcal{F}_T;\mathbb{R}^d)$.
  We consider the solution $(Y,Z)$, $(U,V)\in\mathcal{S}_{\mathcal{F}}^2(0,T;\mathbb{R}^d)\times\mathcal{H}_{\mathcal{F}}^2(0,T;\mathbb{R}^{d\times d})$
  of the following BDSDE, respectively:
\begin{equation}\label{eqA.3.1}
 Y_t=\xi+\int_t^Tf_sds+\int_t^Tg_sd\overleftarrow{B_s}-\int_t^TZ_sdW_s,\ t\in[0,T],
\end{equation}
 and
\begin{equation}\label{eqA.3.2}
U_t=\eta+\int_t^Tu_sds+\int_t^Tv_sd\overleftarrow{B_s}-\int_t^TV_sdW_s,\ t\in[0,T].
\end{equation}
Then, for all $t\in[0,T]$, we have
\begin{equation}\label{eqA.3.3}
\begin{split}
F(t,&\ U_t,P_{Y_t}) = F(T,\eta,P_{\xi})+\int_t^T \Big\{ -(\partial_s F)(s,U_s,P_{Y_s})+\sum_{i=1}^d(\partial_{x_i} F)(s,U_s,P_{Y_s})u_s^i\\
 & +\frac{1}{2}\sum_{i,j,k=1}^d(\partial_{x_ix_j}^2 F)(s,U_s,P_{Y_s})v_s^{ik}v_s^{jk}
    -\frac{1}{2}\sum_{i,j=1}^d\sum_{k=1}^l(\partial_{x_ix_j}^2 F)(s,U_s,P_{Y_s})V_s^{ik}V_s^{jk}\Big\}ds\\
 & +\int_t^T\widehat{E}\[ \sum_{i=1}^d(\partial_{\mu} F)_i(s,U_s,P_{Y_s},\widehat{Y_s})\widehat{f_s^i}
 -\frac{1}{2}\sum_{i,j,k=1}^d\partial_{y_i}(\partial_{\mu}F)_j(s,U_s,P_{Y_s},\widehat{Y_s})\widehat{Z_s^{ik}}\widehat{Z_s^{jk}}\\
 & +\frac{1}{2}\sum_{i,j=1}^d\sum_{k=1}^l\partial_{y_i}(\partial_{\mu}F)_j(s,U_s,P_{Y_s},\widehat{Y_s})\widehat{g_s^{ik}}\widehat{g_s^{jk}}\]ds
  +\int_t^T\sum_{i=1}^d\sum_{j=1}^l (\partial_{x_i} F)(s,U_s,P_{Y_s})v_s^{ij}d\overleftarrow{B_s^j}\\
 & -\int_t^T\sum_{i,j=1}^d (\partial_{x_i} F)(s,U_s,P_{Y_s})V_s^{ij}dW_s^j.
\end{split}
\end{equation}
\end{theorem}
 Here $(\widehat{Y},\widehat{Z},\widehat{f},\widehat{g})$ denotes an independent copy of $(Y,Z,f,g)$, defined on another probability space
$(\widehat{\Omega},\widehat{\mathcal{F}},\widehat{P})$. The expectation $\widehat{E}[\cdot]$ on $(\widehat{\Omega},\widehat{\mathcal{F}},\widehat{P})$
 concerns only random variables endowed with the superscript ``$\ \widehat{}\ $''. For a better readability of the manuscript,
 the proof of this theorem is postponed to Appendix A.1.
\begin{remark} \label{re2.1}
We note that like in \cite{BLPR2017} (there for the ``classcal'' mean-field It\^{o} formula) we do not need the existence of the second order
mixed derivatives $\partial_x\partial_{\mu}F$, $\partial_{\mu}\partial_x F$, $\partial_{\mu}^2F$ for our It\^{o} formula.
This is why they are not introduced in the definition of the space $C_b^{1,2,2}([0,T]\times\mathbb{R}^{d}\times\mathcal{P}_{2}(\mathbb{R}^{d}))$.
\end{remark}
In particular, for dimensions $d=1$, $l=1$, we have the following corollary.
\begin{corollary} \label{corA.3.1}
Let $F\in C_b^{1,2,2}([0,T]\times\mathbb{R}^{2}\times\mathcal{P}_{2}(\mathbb{R}^{2}))$.
Given the solution $(Y,Z)\in\mathcal{S}_{\mathcal{F}}^2(0,T;\mathbb{R})\times\mathcal{H}_{\mathcal{F}}^2(0,T;\mathbb{R})$
of the BDSDE
\begin{equation}\label{eqA.3.38}
 Y_t=\xi+\int_t^Tf_sds+\int_t^Tg_sd\overleftarrow{B_s}-\int_t^TZ_sdW_s,\ t\in[0,T],
\end{equation}
where $f\in\mathcal{H}^2_{\mathcal{F}}(0,T;\mathbb{R})$, $g\in\mathcal{H}^2_{\mathcal{F}}(0,T;\mathbb{R})$ and
 $\xi\in L^{2}(\mathcal{F}_T;\mathbb{R})$, and the solution $X\in\mathcal{S}_{\mathcal{G}}^2(0,T;\mathbb{R})$ of the SDE
\begin{equation}\label{eqA.3.39}
X_t=X_0+\int_0^tb_sds+\int_0^t\sigma_sdW_s,\ t\in[0,T],
\end{equation}
where $b\in\mathcal{H}^2_{\mathcal{G}}(0,T;\mathbb{R})$, $\sigma\in\mathcal{H}^2_{\mathcal{G}}(0,T;\mathbb{R})$,
 $\xi\in L^{2}(\mathcal{G}_0;\mathbb{R})$.
We have, for all $t\in[0,T]$,
\begin{equation}\label{eqA.3.40}
\begin{split}
&\ dF(t,X_t,Y_t,P_{(X_t,Y_t)}) \\
=&\  \Bigg\{ (\partial_t F)(t,X_t,Y_t,P_{(X_t,Y_t)})
+\Big\langle(\partial_{(x,y)} F)(t,X_t,Y_t,P_{(X_t,Y_t)}),
\left(   \begin{matrix}    b_t\\  -f_t \end{matrix}\right)    \Big\rangle\\
 &\ -\frac{1}{2}(\partial_{yy}^2 F)(t,X_t,Y_t,P_{(X_t,Y_t)})|g_t|^2
    +\frac{1}{2}\Big\langle(\partial_{(x,y)}^2 F)(t,X_t,Y_t,P_{(X_t,Y_t)})\cdot\left(   \begin{matrix}
    \sigma_t\\ Z_t \end{matrix}\right),\left(   \begin{matrix}
    \sigma_t\\  Z_t \end{matrix}\right)\Big\rangle\Bigg\}dt\\
&\ +\Bigg\{\widehat{E}\[ \Big\langle(\partial_{\mu} F)(t,X_t,Y_t,P_{(X_t,Y_t)},\widehat{X}_t,\widehat{Y}_t),
 \left(   \begin{matrix} \widehat{b}_t\\  -\widehat{f}_t \end{matrix}\right)\Big\rangle\]\\
 \end{split}
\end{equation}
\begin{equation*}
\begin{split}
 &\ +  \widehat{E}\[ \frac{1}{2}\Big\langle\partial_{(x,y)}(\partial_{\mu}F)(t,X_t,Y_t,P_{(X_t,Y_t)},\widehat{X}_t,\widehat{Y}_t)
  \cdot\left(   \begin{matrix} \widehat{\sigma}_t\\ \widehat{Z}_t \end{matrix}\right),\left(   \begin{matrix}
    \widehat{\sigma}_t\\  \widehat{Z}_t \end{matrix}\right)\Big\rangle\]\\
 &\ -\widehat{E}\[\frac{1}{2}\partial_{y}(\partial_{\mu}F)_2(t,X_t,Y_t,P_{(X_t,Y_t)},\widehat{X}_t,\widehat{Y}_t)\cdot|\widehat{g}_t|^2\]\Bigg\}dt\\
 &\ - (\partial_{y} F)(t,X_t,Y_t,P_{(X_t,Y_t)})f_td\overleftarrow{B_t}
  + \Big\langle(\partial_{(x,y)} F)(t,X_t,Y_t,P_{(X_t,Y_t)}),\left(   \begin{matrix}
    \sigma_t\\  Z_t \end{matrix}\right)\Big\rangle dW_t.
\end{split}
\end{equation*}
\end{corollary}
\begin{proof}
Obviously, due to our assumptions, the process
$((X,Y),(\sigma,Z))\in \mathcal{S}_{\mathcal{F}}^2(0,T;\mathbb{R}^2)\times\mathcal{H}_{\mathcal{F}}^2(0,T;\mathbb{R}^2)$
is the unique solution of the BDSDE
\begin{equation}\label{eqA.3.42}
  d\left(   \begin{matrix}    X_t\\  Y_t \end{matrix}\right)
  =\left(   \begin{matrix}    b_t\\  -f_t \end{matrix}\right)dt
  +\left(   \begin{matrix}    0\\  -g_t \end{matrix}\right)d\overleftarrow{B_t}
  +\left(   \begin{matrix}    \sigma_t\\  Z_t \end{matrix}\right)dW_t,\quad
  \left(   \begin{matrix}    X_T\\  Y_T \end{matrix}\right)=\left(   \begin{matrix}    X_T\\  \xi \end{matrix}\right).
\end{equation}
Hence, Corollary \ref{corA.3.1} follows directly from Theorem \ref{th2.1}.
\end{proof}

\begin{remark}\label{re2.2}
We also observe that the It\^{o} formula studied in \cite{BLPR2017}
is a special case of Theorem \ref{th2.1} and Corollary \ref{corA.3.1}.
\end{remark}

\section{Mean-field stochastic differential equations }
In this section we consider mean-field stochastic differential equations (SDEs). From now on let be given deterministic Lipschitz functions $b:\mathbb{R}^{d}\times\mathcal{P}_{2}(\mathbb{R}^{d})\rightarrow\mathbb{R}^{d}$,
$\sigma:\mathbb{R}^{d}\times\mathcal{P}_{2}(\mathbb{R}^{d})\rightarrow\mathbb{R}^{d\times d}$ satisfying

\noindent\textbf{Assumption (H3.1)} $b$ and $\sigma$ are bounded and Lipschitz continuous on $\mathbb{R}^{d}\times\mathcal{P}_{2}(\mathbb{R}^{d})$.

We consider for the initial data $(t,x)\in[0,T]\times \mathbb{R}^{d}$ and $\xi\in L^{2}(\mathcal{G}_t;\mathbb{R}^d)$ the following both
stochastic differential equations (SDEs):
\begin{equation}\label{eq3.1}
  X_s^{t,\xi}=\xi+\int_t^sb(X_r^{t,\xi},P_{X_r^{t,\xi}})dr+\int_t^s\sigma(X_r^{t,\xi},P_{X_r^{t,\xi}})dW_r,
\end{equation}
and
\begin{equation}\label{eq3.2}
  X_s^{t,x,\xi}=x+\int_t^sb(X_r^{t,x,\xi},P_{X_r^{t,\xi}})dr+\int_t^s\sigma(X_r^{t,x,\xi},P_{X_r^{t,\xi}})dW_r,\ s\in[t,T].
\end{equation}

We recall that under Assumption (H3.1) the both SDEs have a unique solution in $\mathcal{S}^2_{\mathcal{G}}(t,T;\mathbb{R}^d)$
(see, e.g., Buckdahn, Li, Peng and Rainer \cite{BLPR2017}). In particular, the solution $X^{t,\xi}$ of equation \eqref{eq3.1} allows
to determine that of \eqref{eq3.2}. Observe that $X^{t,x,\xi}\in\mathcal{S}^2_{\mathcal{G}}(t,T;\mathbb{R}^d)$ is independent of $\mathcal{G}_t$. As SDE standard
estimates show, we have for some $C\in\mathbb{R}_{+}$ depending only on the Lipschitz constants of $b$ and $\sigma$,
\begin{equation}\label{eq3.3}
  E[\sup_{s\in[t,T]}|X_s^{t,x,\xi}-X_s^{t,x',\xi}|^2]\leq |x-x'|^2,
\end{equation}
for all $t\in[0,T]$, $x,x'\in\mathbb{R}^d$, $\xi\in L^{2}(\mathcal{G}_t;\mathbb{R}^d)$. This allows to substitute the random variable
$\xi$ for $x$ in \eqref{eq3.2} and shows that $X^{t,x,\xi}|_{x=\xi}$ solves the same SDE as $X^{t,\xi}$. From the uniqueness of
the solution we conclude
\begin{equation}\label{eq3.4}
X_s^{t,\xi}=X_s^{t,x,\xi}\big|_{x=\xi}=X_s^{t,\xi,\xi},\ s\in[t,T].
\end{equation}
Moreover, we deduce the following flow property: For all $0\leq t\leq s \leq T,\ x\in\mathbb{R}^d,\ \xi\in L^{2}(\mathcal{G}_t;\mathbb{R}^d)$,
\begin{equation}\label{eq3.5}
(X_r^{s,X_s^{t,x,\xi},X_s^{t,\xi}},X_r^{s,X_s^{t,\xi}})=(X_r^{t,x,\xi},X_r^{t,\xi}),\ r\in[s,T].
\end{equation}
We have to show that the solution $X^{t,x,\xi}$ does not depend on $\xi$ itself but only on its law $P_{\xi}$.
For this, the following lemma is very useful; please refer to Buckdahn, Li, Peng and Rainer \cite{BLPR2017}, or
Proposition 3.1 in Hao and Li \cite{HL2016}, which is formulated for mean-field SDEs with jump.

\begin{lemma} \label{le3.1}
 Suppose Assumption (H3.1) holds true. Then, for all $p\geq 2$ there exists a constant $C_p>0$ only depending on the Lipschitz constants
of $b$ and $\sigma$, such that for all $t\in[0,T]$, $x,\widehat{x}\in\mathbb{R}^d$, $\xi,\widehat{\xi}\in L^{2}(\mathcal{G}_t;\mathbb{R}^d)$, we have the following estimates
\begin{equation}\label{eq3.6}
  \begin{split}
  &\ \emph{(i)}\ E\[\sup_{s\in[t,T]}|X_s^{t,x,\xi}-X_s^{t,\widehat{x},\widehat{\xi}}\big|^p\big|\mathcal{G}_t\]
       \leq C_p\(|x-\widehat{x}|^p+W_2(P_{\xi},P_{\widehat{\xi}})^p\),\\
&\ \emph{(ii)}\ E\[\sup_{s\in[t,T]}|X_s^{t,x,\xi}|^p\big|\mathcal{G}_t\]\leq C_p\(1+|x|^p\),\\
&\ \emph{(iii)}\ \sup_{s\in[t,T]}W_2(P_{X_s^{t,\xi}},P_{X_s^{t,\widehat{\xi}}})\leq C_2W_2(P_{\xi},P_{\widehat{\xi}}),\\
&\ \emph{(iv)}\ E\[\sup_{s\in[t,t+h]}|X_s^{t,x,\xi}-x|^p\big|\mathcal{G}_t\]\leq C_ph^{\frac{p}{2}},
  \end{split}
\end{equation}

\end{lemma}

\begin{remark} \label{re3.1}
An immediate consequence of Lemma \ref{le3.1}-\rm{(i)} is that, given $(t,x)\in[0,T]\times\mathbb{R}^{d}$, the processes
$X^{t,x,\xi_1}$ and $X^{t,x,\xi_2}$ are indistinguishable, whenever the laws of $\xi_1,\xi_2\in L^{2}(\mathcal{G}_t;\mathbb{R}^d)$
 are the same. But this means that we can define
\begin{equation}\label{eq3.7}
X^{t,x,P_{\xi}}:=X^{t,x,\xi},\ (t,x)\in[0,T]\times\mathbb{R}^{d},\ \xi\in L^{2}(\mathcal{G}_t;\mathbb{R}^d).
\end{equation}
Extending the notation introduced in the preceding section for functions to random variables
and processes, we shall consider the lifted process $\widetilde{X}_s^{t,x,\xi}:=X_s^{t,x,P_{\xi}}=X_s^{t,x,\xi}$, $s\in[t,T]$,
$(t,x)\in[0,T]\times\mathbb{R}^{d}$, $\xi\in L^{2}(\mathcal{G}_t;\mathbb{R}^d)$. However, we prefer to continue to write
$X^{t,x,\xi}$ and reserve the notation $\widetilde{X}^{t,x,\xi}$ for an independent copy of $X^{t,x,P_{\xi}}$, which we will introduce later.
On the other hand, note that $X_s^{t,x,P_{\xi}}$ is $\mathcal{F}_s^W$-measurable, $s\in[t,T]$, $(t,x)\in[0,T]\times\mathbb{R}^{d}$, $\xi\in L^{2}(\mathcal{G}_t;\mathbb{R}^d)$.
\end{remark}

We now recall that under natural conditions the solution of a SDE is differentiable in Malliavin's sense and that the derivative is a solution of a linear SDE.

\begin{proposition} \label{prop3.1}
Let $b$ and $\sigma$ satisfy Assumption (H3.1). Moreover, we suppose the coefficients $b_j$, $\sigma_{i,j}$, $1\leq i,j\leq d$, to satisfy the following assumptions:\\
\indent \emph{(i)} $b_j(\cdot,\mu)$, $\sigma_{i,j}(\cdot,\mu)$ belong to $C_b^{1}(\mathbb{R}^{d})$, for all $\mu\in\mathcal{P}_{2}(\mathbb{R}^{d})$;\\
\indent \emph{(ii)} The derivatives $\partial_xb_j$, $\partial_x\sigma_{i,j}$: $\mathbb{R}^{d}\times\mathcal{P}_{2}(\mathbb{R}^{d})\rightarrow\mathbb{R}^{d}$
are bounded and Lipschitz continuous.

Then for all $(t,x)\in[0,T]\times \mathbb{R}^{d}$, $\xi\in L^{2}(\mathcal{G}_t;\mathbb{R}^d)$, $s\in[t,T]$,
$X_s^{t,x,P_{\xi}}\in(\mathbb{D}^{1,2})^d$ and a version of $\{D_{\theta}X_s^{t,x,P_{\xi}}:\theta,s\in[t,T]\}$ is given by:\\
\indent \emph{(i)} $D_{\theta}X_s^{t,x,P_{\xi}}=0$, if $\theta> s$;\\
\indent \emph{(ii)} $\{D_{\theta}X_s^{t,x,P_{\xi}}=(D_{\theta}^iX_s^{t,x,P_{\xi},j})_{1\leq i,j\leq d}:s\in[\theta,T]\}$
is the unique solution of the linear SDE: $s\in[t,T]$, $1\leq i,j\leq d$,
\begin{equation}\label{eq3.8}
  \begin{split}
 D_{\theta}^iX_s^{t,x,P_{\xi},j}=&\ \sigma_{i,j}(X_{\theta}^{t,x,P_{\xi}},P_{X_{\theta}^{t,\xi}})
 +\sum_{k=1}^d\int_{\theta}^s\partial_{x_k}b_j(X_r^{t,x,P_{\xi}},P_{X_r^{t,\xi}}) D_{\theta}^iX_r^{t,x,P_{\xi},k}dr\\
 &\ + \sum_{k,l=1}^d\int_{\theta}^s\partial_{x_k}\sigma_{j,l}(X_r^{t,x,P_{\xi}},P_{X_r^{t,\xi}}) D_{\theta}^iX_r^{t,x,P_{\xi},k}dW_r^l, \ s\in[\theta,T].
  \end{split}
\end{equation}
Furthermore, for all $p\geq 2$ there exists a constant $C_p>0$ only depending on the Lipschitz
constants of $\partial_xb$ and $\partial_x\sigma$, such that, for all $t\in[0,T]$, $x,x'\in\mathbb{R}^d$,
$\xi,\xi'\in L^{2}(\mathcal{G}_t;\mathbb{R}^d)$, $P$-a.s.,
\begin{equation}\label{eq3.9}
  \begin{split}
  &\ \emph{(i)}\ E\[\sup_{s\in[t,T]}|D_{\theta}X_s^{t,x,P_{\xi}}|^p\]\leq C_p;\\
  &\ \emph{(ii)}\ E\[\sup_{s\in[t,T]}|D_{\theta}X_s^{t,x,P_{\xi}}-D_{\theta}X_s^{t,x',P_{\xi'}}|^p\]
       \leq C_p\(|x-x'|^p+W_2(P_{\xi},P_{\xi'})^p\).
  \end{split}
\end{equation}
\end{proposition}
\begin{proof}
It is standard to prove that $X_s^{t,x,P_{\xi}}$ has a Malliavin derivative under our assumptions (see \cite{N1995}).
Moreover, from the assumptions, Lemma \ref{le3.1} and the standard estimates of classical SDEs, it can easily be verified that
\eqref{eq3.9} holds true.
\end{proof}

\section{Mean-field BDSDEs }
In this section we consider mean-field backward doubly stochastic differential
equations (BDSDEs for short) driven by a standard (forward) stochastic integral $dW_t$
and a backward stochastic integral $dB_t$. The existence and the uniqueness of the solution for this type of BDSDEs is proved;
for more details, please, we refer to Appendix A.2.

Let $f: \mathbb{R}^d\times\mathbb{R}\times\mathbb{R}^d\times\mathcal{P}_{2}(\mathbb{R}^d\times\mathbb{R}\times\mathbb{R}^d)\rightarrow\mathbb{R}$,
$g: \mathbb{R}^d\times\mathbb{R}\times\mathbb{R}^d\times\mathcal{P}_{2}(\mathbb{R}^d\times\mathbb{R}\times\mathbb{R}^d)\rightarrow\mathbb{R}^{l}$,
$h: \mathcal{P}_{2}(\mathbb{R}^d\times\mathbb{R}\times\mathbb{R}^d)\rightarrow\mathbb{R}^{l}$
and $\Phi :\mathbb{R}^d\times\mathcal{P}_{2}(\mathbb{R}^d)\rightarrow\mathbb{R}$ be deterministic and satisfy:

\noindent\textbf{Assumption (H4.1)} The functions $f$, $g$, $h$ and $\Phi$ are bounded and Lipschitz,
 i.e., there exist constants $C>0$, and $\alpha_1,\alpha_2>0$ with $0<\alpha_1+\alpha_2<1$
 such that, for all $x,x'\in\mathbb{R}^{d}$, $y,y'\in\mathbb{R}$, $z,z'\in\mathbb{R}^{d}$, $\mu,\mu'\in\mathcal{P}_{2}(\mathbb{R}^{d+1}\times\mathbb{R}^d)$,
\begin{equation*}
\begin{split}
&\mbox{(i) }|f(x,y,z,\mu)-f(x',y',z',\mu')|+|h(\mu)-h(\mu')|+|\Phi(x,\mu)-\Phi(x',\mu')|\\
&\qquad \leq C(|x-x'|+|y-y'|+|z-z'|+W_2(\mu,\mu')),\\
&\mbox{(ii) } |g(x,y,z,\mu)-g(x',y',z',\mu')|^2\leq C(|x-x'|^2+|y-y'|^2)+\alpha_1|z-z'|^2+W_{2,C,\alpha_2}(\mu,\mu')^2.
\end{split}
\end{equation*}
\indent Given $x\in\mathbb{R}^{d}$ and $\xi\in L^{2}(\mathcal{G}_t;\mathbb{R}^d)$ we consider the following both BDSDEs:
\begin{equation}\label{eq4.1}
  \left\{
   \begin{array}{l}
   dY_s^{t,\xi}=-f(\Pi_s^{t,\xi},P_{\Pi_s^{t,\xi}})ds-
   \big(g(\Pi_s^{t,\xi},P_{\Pi_s^{t,\xi}})+h(P_{\Pi_s^{t,\xi}})\big)d\overleftarrow{B_s}+Z_s^{t,\xi}dW_s,\ s\in[t,T],\\
   Y_T^{t,\xi}=\Phi (X_T^{t,\xi},P_{X_T^{t,\xi}}),\\
   \end{array}
  \right.
\end{equation}
and
\begin{equation}\label{eq4.2}
  \left\{
   \begin{array}{l}
   dY_s^{t,x,\xi}=-f(\Pi_s^{t,x,\xi},P_{\Pi_s^{t,\xi}})ds-
   \big(g(\Pi_s^{t,x,\xi},P_{\Pi_s^{t,\xi}})+h(P_{\Pi_s^{t,\xi}})\big)d\overleftarrow{B_s}+Z_s^{t,x,\xi}dW_s,\ s\in[t,T],\\
   Y_T^{t,x,\xi}=\Phi (X_T^{t,x,\xi},P_{X_T^{t,\xi}}),\\
   \end{array}
  \right.
\end{equation}
where $\Pi_s^{t,\xi}:=(X_s^{t,\xi},Y_s^{t,\xi},Z_s^{t,\xi})$, $\Pi_s^{t,x,\xi}:=(X_s^{t,x,\xi},Y_s^{t,x,\xi},Z_s^{t,x,\xi})$.
Recall that the processes $X^{t,\xi}$ and $X^{t,x,\xi}$ are the solution of SDEs \eqref{eq3.1} and \eqref{eq3.2}, respectively.

Under Assumption (H4.1) we know that from Theorem \ref{thA.1} in Appendix A.2 that Eq.\eqref{eq4.1} has a unique solution
 $(Y^{t,\xi},Z^{t,\xi})\in\mathcal{S}_{\mathcal{F}}^2(t,T;\mathbb{R})\times\mathcal{H}_{\mathcal{F}}^2(t,T;\mathbb{R}^{d})$.
On the other hand, once having the solution of \eqref{eq4.1}, under Assumption (H4.1) BDSDE \eqref{eq4.2} becomes
classical and possesses a unique solution $(Y^{t,x,\xi},Z^{t,x,\xi})\in\mathcal{S}_{\mathcal{F}}^2(t,T;\mathbb{R})\times\mathcal{H}_{\mathcal{F}}^2(t,T;\mathbb{R}^{d})$.
Indeed, once we have got $\Pi_s^{t,\xi}:=(X_s^{t,\xi},Y_s^{t,\xi},Z_s^{t,\xi})$, we define
$\widetilde{f}(s,y,z)=f(X_s^{t,x,\xi},y,z,P_{\Pi_s^{t,\xi}})$, $\widetilde{g}(s,y,z)=g(X_s^{t,x,\xi},y,z,P_{\Pi_s^{t,\xi}})$,
$\widetilde{h}(s)=h(P_{\Pi_s^{t,\xi}})$, $\widetilde{\xi}=\Phi(X_T^{t,x,\xi},P_{\Pi_T^{t,\xi}})$.
Obviously, $\widetilde{f}$, $\widetilde{g}$, $\widetilde{h}$ and $\widetilde{\xi}$
satisfy all assumptions of classical BDSDEs, and hence, BDSDE \eqref{eq4.2} has a unique solution
$(Y^{t,x,\xi},Z^{t,x,\xi})\in\mathcal{S}_{\mathcal{F}}^2(t,T;\mathbb{R})\times\mathcal{H}_{\mathcal{F}}^2(t,T;\mathbb{R}^{d})$.

From the flow property \eqref{eq3.5} and the uniqueness of the solution of \eqref{eq4.1} and  \eqref{eq4.2} we have the
following properties: For all $0\leq t\leq s\leq T$, $x\in\mathbb{R}^d$, $\xi\in L^{2}(\mathcal{G}_t;\mathbb{R}^d)$,
\begin{equation}\label{eq4.3}
\begin{split}
&\mbox{(i) } (Y_r^{s,X_s^{t,x,\xi},X_s^{t,\xi}},Y_r^{s,X_s^{t,\xi}})=(Y_r^{t,x,\xi},Y_r^{t,\xi}),\ r\in[s,T],\ P\mbox{-}a.s.;\\
&\mbox{(ii) } (Z_r^{s,X_s^{t,x,\xi},X_s^{t,\xi}},Z_r^{s,X_s^{t,\xi}})=(Z_r^{t,x,\xi},Z_r^{t,\xi}),\ drdP\mbox{-}a.e.\ \mbox{on}\ [s,T]\times\Omega.
\end{split}
\end{equation}
\indent Moreover, we have to show that the solution $(Y^{t,x,\xi},Z^{t,x,\xi})$ does not depend on $\xi$ itself
but only on its law $P_{\xi}$.
\begin{proposition} \label{prop4.0}
Suppose the Assumptions (H3.1) and (H4.1) hold true. Then,
 there is a constant $C>0$ only depending on the Lipschitz constants
of $b$, $\sigma$, $f$, $g$, $h$ and $\Phi$, such that, for $t\in[0,T]$, $x,\widehat{x}\in\mathbb{R}^d$, $\xi,\widehat{\xi}\in L^{2}(\mathcal{G}_t;\mathbb{R}^d)$,
\begin{equation}\label{eq4.3+1+1}
  \begin{split}
  &\ \emph{(i)}\ E[\sup_{s\in[t,T]}|Y_s^{t,x,\xi}|^2+\int_t^T|Z_s^{t,x,\xi}|^2ds\big|\mathcal{G}_t]\leq C;\\
  &\ \emph{(ii)}\ E[\sup_{s\in[t,T]}|Y_s^{t,x,\xi}-Y_s^{t,\widehat{x},\widehat{\xi}}|^2
  +\int_t^T|Z_s^{t,x,\xi}-Z_s^{t,\widehat{x},\widehat{\xi}}|^2ds\big|\mathcal{G}_t]
       \leq C_p\(|x-\widehat{x}|^2+W_2(P_{\xi},P_{\widehat{\xi}})^2\);\\
&\ \emph{(iii)}\ \int_t^TW_2(P_{\Pi_s^{t,\xi}},P_{\Pi_s^{t,\widehat{\xi}}})^2ds\leq CW_2(P_{\xi},P_{\widehat{\xi}})^2.
  \end{split}
\end{equation}
\end{proposition}
\begin{proof}
From Lemma \ref{leA.1}-(1) we get (i) directly. Now we prove (ii) and (iii).

Note that $\Pi^{t,x,\xi}$ is independent of $\mathcal{G}_t$ and, hence, of $\xi\in L^{2}(\mathcal{G}_t;\mathbb{R}^d)$.
This allows to consider $\Pi^{t,x,\xi}\big|_{x=\xi}$, and from the uniqueness of the solution of \eqref{eq4.1} and \eqref{eq4.2},
it follows from \eqref{eq3.4} that $\Pi^{t,\xi}=\Pi^{t,x,\xi}\big|_{x=\xi}$. On the other hand, it also follows that, if
$\xi'\in L^{2}(\mathcal{G}_t;\mathbb{R}^d)$ has the same law as $\xi$, then also $\Pi^{t,\xi',\xi}:=\Pi^{t,x,\xi}\big|_{x=\xi'}$
and $\Pi^{t,\xi}$ are of the same law. Hence, $P_{\Pi_s^{t,\xi}}=P_{\Pi_s^{t,\xi',\xi}}$, $ds$-a.e. Then, for given
$\xi_i\in L^{2}(\mathcal{G}_t;\mathbb{R}^d)$, $i=1,2$, and $\xi_i'\in L^{2}(\mathcal{G}_t;\mathbb{R}^d)$ of the same law as $\xi_i$, we consider the
following BDSDE:
\begin{equation}\label{eq4.4+1}
  \left\{
   \begin{array}{l}
   dY_s^{t,\xi_i',\xi_i}\!=\!-\!f(\Pi_s^{t,\xi_i',\xi_i}\!,P_{\Pi_s^{t,\xi_i}})ds\!-\!
   \big(g(\Pi_s^{t,\xi_i',\xi_i}\!,P_{\Pi_s^{t,\xi_i}})\!+\!h(P_{\Pi_s^{t,\xi_i}})\big)d\overleftarrow{B_s}
   \!+\!Z_s^{t,\xi_i',\xi_i}dW_s,\ \!s\in[t,T],\\
   Y_T^{t,\xi_i',\xi_i}=\Phi (X_T^{t,\xi_i',\xi_i},P_{X_T^{t,\xi_i}}).\\
   \end{array}
  \right.
\end{equation}
Then, taking into account that $P_{\Pi_s^{t,\xi_i}}=P_{\Pi_s^{t,\xi_i',\xi_i}}$, $ds$-a.e., for $i=1,2$,
we get from Theorem \ref{thA.2} and (H4.1) that, for $0\leq t\leq T$,
\begin{equation}\label{eq4.5+1+1}
  \begin{split}
     &\ E\[\int_t^T\(|Y_s^{t,\xi_1',\xi_1}-Y_s^{t,\xi_2',\xi_2}|^{2}+|Z_s^{t,\xi_1',\xi_1}-Z_s^{t,\xi_2',\xi_2}|^{2}\)ds\]\\
\leq&\ C^1E[|X_T^{t,\xi_1',\xi_1}-X_T^{t,\xi_2',\xi_2}|^{2}+W_2(P_{X_T^{t,\xi_1}},P_{X_T^{t,\xi_2}})^2]
 + C^1 E\[\int_t^T|X_s^{t,\xi_1',\xi_1}-X_s^{t,\xi_2',\xi_2}|^{2}ds\],
   \end{split}
\end{equation}
where $C^1$ depends only on the Lipschitz constants of $f$, $g$, $h$ and $\Phi$. From Lemma \ref{le3.1} we have
\begin{equation}\label{eq4.5+1+2}
  \begin{split}
   \mbox{(i)}&\ W_2(P_{X_T^{t,\xi_1}},P_{X_T^{t,\xi_2}})\leq CW_2(P_{\xi_1},P_{\xi_2});\\
   \mbox{(ii)}&\ E[\sup_{s\in[t,T]}|X_s^{t,\xi_1',\xi_1}-X_s^{t,\xi_2',\xi_2}|^{2}]
  =E[E[\sup_{s\in[t,T]}|X_s^{t,x_1,\xi_1}-X_s^{t,x_2,\xi_2}|^{2}\big|\mathcal{G}_t]\big|_{\substack{x_1=\xi_1',\\ x_2=\xi_2'}}]\\
   &\    \leq CE\[|\xi_1'-\xi_2'|^2+W_2(P_{\xi_1},P_{\xi_2})^2\].
  \end{split}
\end{equation}
Then, by combining \eqref{eq4.5+1+1} and \eqref{eq4.5+1+2} we obtain
\begin{equation}\label{eq4.6+1}
  \begin{split}
    &\ E\[\int_t^T\(|Y_s^{t,\xi_1',\xi_1}-Y_s^{t,\xi_2',\xi_2}|^{2}+|Z_s^{t,\xi_1',\xi_1}-Z_s^{t,\xi_2',\xi_2}|^{2}\)ds\]
\leq CE\[|\xi_1'-\xi_2'|^2+W_2(P_{\xi_1},P_{\xi_2})^2\].
  \end{split}
\end{equation}
Furthermore, from the properties of $W_2$, \eqref{eq4.5+1+2}-(ii) and \eqref{eq4.6+1} we get
\begin{equation*}
  \begin{split}
    &\ \int_t^TW_2(P_{\Pi_s^{t,\xi_1}},P_{\Pi_s^{t,\xi_2}})^2ds=\int_t^TW_2(P_{\Pi_s^{t,\xi_1',\xi_1}},P_{\Pi_s^{t,\xi_2',\xi_2}})^2ds\\
\leq&\  E\[\int_t^T\(|X_s^{t,\xi_1',\xi_1}-X_s^{t,\xi_2',\xi_2}|^{2}+|Y_s^{t,\xi_1',\xi_1}-Y_s^{t,\xi_2',\xi_2}|^{2}
+|Z_s^{t,\xi_1',\xi_1}-Z_s^{t,\xi_2',\xi_2}|^{2}\)ds\]\\
\leq&\ CE\[|\xi_1'-\xi_2'|^2+W_2(P_{\xi_1},P_{\xi_2})^2\].
  \end{split}
\end{equation*}
Hence, taking the infimum over all $\xi_1',\xi_2'\in L^{2}(\mathcal{G}_t;\mathbb{R}^d)$ with $P_{\xi_i}=P_{\xi_i'}$, $i=1,2$, we get
\begin{equation}\label{eq4.7+1}
     \int_t^TW_2(P_{\Pi_s^{t,\xi_1}},P_{\Pi_s^{t,\xi_2}})^2ds\leq CW_2(P_{\xi_1},P_{\xi_2})^2,\ \xi_1,\ \xi_2\in L^{2}(\mathcal{G}_t;\mathbb{R}^d).
\end{equation}
This allows now to apply Lemma \ref{leA.1}-(1) to BDSDE \eqref{eq4.2} with
$f_i(s,y,z):=f(X_s^{t,x_i,\xi_i},y,z,P_{\Pi_s^{t,\xi_i}})$,
$g_i(s,y,z):=g(X_s^{t,x_i,\xi_i},y,z,P_{\Pi_s^{t,\xi_i}})+h(P_{\Pi_s^{t,\xi_i}})$,
$\theta_i:=\Phi(X_T^{t,x_i,\xi_i},P_{\Pi_T^{t,\xi_i}})$. Then, thanks to Lemma \ref{le3.1} and \eqref{eq4.7+1},
there is a constant $C>0$ only depending on the Lipschitz constants of $b$, $\sigma$, $f$, $g$, $h$ and $\Phi$, such that,
for $t\in[0,T]$, $x_1,x_2\in\mathbb{R}^d$, $\xi_1,\xi_2\in L^{2}(\mathcal{G}_t;\mathbb{R}^d)$,
\begin{equation*}\label{eq4.8+1}
  \begin{split}
   &\  E[\sup_{s\in[t,T]}|Y_s^{t,x_1,\xi_1}-Y_s^{t,x_2,\xi_2}|^{2}+\int_t^T|Z_s^{t,x_1,\xi_1}-Z_s^{t,x_2,\xi_2}|^2ds\big|\mathcal{G}_t]\\
  \leq&\  \!CE\[|X_T^{t,x_1,\xi_1}\!\!-\!\!X_T^{t,x_2',\xi_2}|^{2}\!\!+\!\!W_2(P_{X_T^{t,\xi_1}}\!,\!P_{X_T^{t,\xi_2}})^{2}
\!\! + \!\!\int_t^T\!\!\!(|X_s^{t,x_1,\xi_1}\!\!-\!\!X_s^{t,x_2,\xi_2}|^2
   \!\!+\!\!W_2(P_{\Pi_s^{t,\xi_1}},\!P_{\Pi_s^{t,\xi_2}})^2)ds\big|\mathcal{G}_t\]\\
  \leq&\ C\(|x_1-x_2|^{2}+W_2(P_{\xi_1},P_{\xi_2})^{2}\).
  \end{split}
\end{equation*}
The proof is complete.
\end{proof}
For higher order moment estimates, we need an additional assumption on $\alpha_1$ and $\alpha_2$.
\medskip

\noindent\textbf{Assumption (H4.2)} For some $ p_0\geq 16$, $\overline{C}_{p}(\alpha_1+\alpha_2)^{\frac{p}{2}}<1$, when $p=p_0,\frac{p_0}{2},\frac{p_0}{8}$. Here $\overline{C}_{p}:=2^{p-1}C^{\ast}_{p}((\frac{p}{p-1})^{p}+1)C_{p}'$, $C^{\ast}_{p}:=2^{-p-2}3^{p}{p}^{3p}+2^{\frac{p}{2}}$,
$C_{p}':=(\frac{p}{p-1})^{p}3^{p-1}\(2C^{p}5^{p-1}\vee (6{p}^3)^{p} 5^{\frac{p}{2}-1}\)$,
$C$ is the Lipschitz constant in Assumption (H4.1).
\begin{proposition} \label{prop4.1}
Suppose the Assumptions (H3.1), (H4.1) and (H4.2) hold true. Then, for all $ p\in[2,p_0]$,
 there is a constant $C_p>0$ only depending on the Lipschitz constants
of $b$, $\sigma$, $f$, $g$, $h$ and $\Phi$, such that, for $t\in[0,T]$, $x,\widehat{x}\in\mathbb{R}^d$, $\xi,\widehat{\xi}\in L^{2}(\mathcal{G}_t;\mathbb{R}^d)$,
\begin{equation}\label{eq4.3+1}
  \begin{split}
  &\ \emph{(i)}\ E[\sup_{s\in[t,T]}|Y_s^{t,x,\xi}|^p+(\int_t^T|Z_s^{t,x,\xi}|^2ds)^{\frac{p}{2}}\big|\mathcal{G}_t]\leq C_p;\\
  &\ \emph{(ii)}\ E[\sup_{s\in[t,T]}|Y_s^{t,x,\xi}\!-\!Y_s^{t,\widehat{x},\widehat{\xi}}|^p
  \!+\!(\int_t^T|Z_s^{t,x,\xi}\!-\!Z_s^{t,\widehat{x},\widehat{\xi}}|^2ds)^{\frac{p}{2}}\big|\mathcal{G}_t]
       \!\leq \!C_p\(|x\!-\!\widehat{x}|^p\!+\!W_2(P_{\xi},P_{\widehat{\xi}})^p\).
  \end{split}
\end{equation}
\end{proposition}
\begin{proof}
From Lemma \ref{leA.1}-(2) we get
$\displaystyle E[\sup_{s\in[t,T]}|Y_s^{t,x,\xi}|^{p_0}+(\int_t^T|Z_s^{t,x,\xi}|^2ds)^{\frac{p_0}{2}}\big|\mathcal{G}_t]\leq C_{p_0},$
then from H\"{o}lder's inequality, for all $ p\in[2,p_0]$ we have $\displaystyle E[\sup_{s\in[t,T]}|Y_s^{t,x,\xi}|^p+(\int_t^T|Z_s^{t,x,\xi}|^2ds)^{\frac{p}{2}}\big|\mathcal{G}_t]\leq C_p$.
 Now we prove \eqref{eq4.3+1}-(ii).

From our assumptions, we can apply Lemma \ref{leA.1}-(2) to BDSDE \eqref{eq4.2} with
$f_i(s,y,z):=f(X_s^{t,x_i,\xi_i},y,z,P_{\Pi_s^{t,\xi_i}})$,
$g_i(s,y,z):=g(X_s^{t,x_i,\xi_i},y,z,P_{\Pi_s^{t,\xi_i}})+h(P_{\Pi_s^{t,\xi_i}})$,
$\theta_i:=\Phi(X_T^{t,x_i,\xi_i},P_{\Pi_T^{t,\xi_i}})$. Then, thanks to Lemma \ref{le3.1} and \eqref{eq4.3+1+1}-(iii),
there is a constant $C_{p_0}>0$ only depending on the Lipschitz constants of $b$, $\sigma$, $f$, $g$, $h$ and $\Phi$, such that,
for $t\in[0,T]$, $x_1,x_2\in\mathbb{R}^d$, $\xi_1,\xi_2\in L^{2}(\mathcal{G}_t;\mathbb{R}^d)$,
\begin{equation*}\label{eq4.8}
  \begin{split}
   &\  E[\sup_{s\in[t,T]}|Y_s^{t,x_1,\xi_1}-Y_s^{t,x_2,\xi_2}|^{p_0}+(\int_t^T|Z_s^{t,x_1,\xi_1}-Z_s^{t,x_2,\xi_2}|^2ds)^{\frac{p_0}{2}}\big|\mathcal{G}_t]\\
  \leq&\  \!C_{p_0}E\[\!|X_T^{t,x_1,\xi_1}\!\!-\!\!X_T^{t,x_2',\xi_2}|^{p_0}\!\!+\!\!W_2(P_{X_T^{t,\xi_1}}\!,\!P_{X_T^{t,\xi_2}})^{p_0}
\!\! + \!\!\(\!\!\int_t^T\!\!\!(|X_s^{t,x_1,\xi_1}\!\!-\!\!X_s^{t,x_2,\xi_2}|^2
   \!\!+\!\!W_2(P_{\Pi_s^{t,\xi_1}},\!P_{\Pi_s^{t,\xi_2}})^2)ds\!\)^{\frac{p_0}{2}}\big|\mathcal{G}_t\!\]\\
  \leq&\ C_{p_0}\(|x_1-x_2|^{p_0}+W_2(P_{\xi_1},P_{\xi_2})^{p_0}\).
  \end{split}
\end{equation*}
Again from H\"{o}lder's inequality, for all $p\in[2,p_0]$ we have
$$E[\sup_{s\in[t,T]}|Y_s^{t,x,\xi}\!-\!Y_s^{t,\widehat{x},\widehat{\xi}}|^p
  \!+\!(\int_t^T|Z_s^{t,x,\xi}\!-\!Z_s^{t,\widehat{x},\widehat{\xi}}|^2ds)^{\frac{p}{2}}\big|\mathcal{G}_t]
       \!\leq \!C_p\(|x\!-\!\widehat{x}|^p\!+\!W_2(P_{\xi},P_{\widehat{\xi}})^p\).$$
The proof is complete.
\end{proof}
\begin{remark} \label{re4.1}
As a matter of fact, in Assumption (H4.1), it is necessary to assume that the functions $f$, $g$, $h$ and $\Phi$ are Lipschitz.
However, the boundedness condition can be replace by a linear growth condition. In other word, if the functions $f$, $g$, $h$ and $\Phi$ are Lipschitz
and of linear growth, analogous arguments allow to show that Propositions \ref{prop4.0} and \ref{prop4.1} still hold true, respectively. On the other hand, notice that in fact $Y_s^{t,x,\xi}$ is $\mathcal{F}_{t,s}^{W}\vee\mathcal{F}_{s,T}^{B}$-measurable, $s\in[t,T]$.
\end{remark}
Recalling that $(Y^{t,\xi},Z^{t,\xi})=(Y^{t,x,\xi},Z^{t,x,\xi})|_{x=\xi}$ we have the following result.
\begin{corollary} \label{cor4.1}
Suppose the Assumptions (H3.1), (H4.1) and (H4.2) hold true. Then, for all $ p\in[2,p_0]$,
 there is a constant $C_p>0$ only depending on the Lipschitz constants
of $b$, $\sigma$, $f$, $g$, $h$ and $\Phi$, such that, for $t\in[0,T]$, $\xi_1,\xi_2\in L^{2}(\mathcal{G}_t;\mathbb{R}^d)$,
\begin{equation}\label{eq4.9}
  \begin{split}
&\ E[\sup_{s\in[t,T]}|Y_s^{t,\xi_1}-Y_s^{t,\xi_2}|^p+(\int_t^T|Z_s^{t,\xi_1}-Z_s^{t,\xi_2}|^2ds)^{\frac{p}{2}}]\\
\leq&\  C_pE\[|\xi_1-\xi_2|^p+W_2(P_{\xi_1},P_{\xi_2})^p\]\leq C_pE[|\xi_1-\xi_2|^p].
  \end{split}
\end{equation}
\end{corollary}
\indent Due to Proposition \ref{prop4.0} the processes $Y^{t,x,\xi}=\{Y_s^{t,x,\xi}\}_{s\in[t,T]}$, and $Z^{t,x,\xi}=\{Z_s^{t,x,\xi}\}_{s\in[t,T]}$
depend on $\xi$ only through its distribution, which means $(Y_s^{t,x,\xi},Z_s^{t,x,\xi})$ and $(Y_s^{t,x,\overline{\xi}},Z_s^{t,x,\overline{\xi}})$
 are indistinguishable as long as $\xi$ and $\overline{\xi}$ have the same distribution. Hence we can define $Y_s^{t,x,P_{\xi}}$, $Z_s^{t,x,P_{\xi}}$ by putting
\begin{equation*}
Y_s^{t,x,P_{\xi}}:=Y_s^{t,x,\xi},\ Z_s^{t,x,P_{\xi}}:=Z_s^{t,x,\xi}.
\end{equation*}
Moreover, from the uniqueness of the solution of the BDSDEs \eqref{eq4.1} and \eqref{eq4.2} it follows that
\begin{equation}\label{eq4.10}
  Y^{t,\xi}=Y^{t,x,P_{\xi}}\big|_{x=\xi},\ Z^{t,\xi}=Z^{t,x,P_{\xi}}\big|_{x=\xi}.
\end{equation}
Finally, from \eqref{eq4.3}, for $0\leq t\leq s\leq T$, $x\in\mathbb{R}^d$, $\xi\in L^{2}(\mathcal{G}_t;\mathbb{R}^d)$, it holds
\begin{equation}\label{eq4.11}
\begin{split}
&\mbox{(i) } (Y_r^{s,X_s^{t,x,P_{\xi}},P_{X_s^{t,\xi}}},Y_r^{s,X_s^{t,\xi}})=(Y_r^{t,x,P_{\xi}},Y_r^{t,\xi}),\ r\in[s,T], P\mbox{-}a.s.,\\
&\mbox{(ii) } (Z_r^{s,X_s^{t,x,P_{\xi}},P_{X_s^{t,\xi}}},Z_r^{s,X_s^{t,\xi}})=(Z_r^{t,x,P_{\xi}},Z_r^{t,\xi}),\ drdP\mbox{-}a.e.\ \mbox{on}\ [s,T]\times\Omega.
\end{split}
\end{equation}

Next, we recall that under natural conditions the solution of a BDSDE is differentiable w.r.t. the Brownian motion $W$ in Malliavin's sense
and that the derivative is a solution of a linear BDSDE.
\begin{proposition} \label{prop4.2}
Let $b$ and $\sigma$ satisfy the assumption in Proposition \ref{prop3.1}, and let Assumption (H4.1) hold true.
 Moreover, we suppose $\Phi$, $f$ and $g$ satisfy: \\
\indent \emph{(i)} For all $\mu\in\mathcal{P}_{2}(\mathbb{R}^{d})$, $\Phi(\cdot,\mu)\in C_b^{1}(\mathbb{R}^{d})$,
$f(\cdot,\mu)\in C_b^{1}(\mathbb{R}^{d+1+d})$, $g(\cdot,\mu)\in C_b^{1}(\mathbb{R}^{d+1+d})$;
\indent \emph{(ii)} The derivatives $\partial_x\Phi:\mathbb{R}^{d}\times\mathcal{P}_{2}(\mathbb{R}^{d})\rightarrow\mathbb{R}^{d}$,
$(\partial_x,\partial_y,\partial_z)f:\mathbb{R}^{d+1+d}\times\mathcal{P}_{2}(\mathbb{R}^{d+1+d})\rightarrow\mathbb{R}^{d+1+d}$,
and $(\partial_x,\partial_y,\partial_z)g_i:\mathbb{R}^{d+1+d}\times\mathcal{P}_{2}(\mathbb{R}^{d+1+d})\rightarrow\mathbb{R}^{d+1+d}$, $1\leq i\leq l$,
are bounded and Lipschitz continuous.\\
\noindent Then, for all $(t,x)\!\in\![0,T]\times \mathbb{R}^{d}$, $\xi\in L^{2}(\mathcal{G}_t;\mathbb{R}^d)$, $s\in[t,T]$,
$(Y_s^{t,x,P_{\xi}},Z_s^{t,x,P_{\xi}})_{t\leq s\leq T}\!\in\! L^2(t,T;(\mathbb{D}^{1,2})^{1+d})$
and a version of $\{D_{\theta}Y_s^{t,x,P_{\xi}},D_{\theta}Z_s^{t,x,P_{\xi}}:\theta,s\in[t,T]\}$ is given by:\\
\indent \emph{(i)} $D_{\theta}Y_s^{t,x,P_{\xi}}=0$, $D_{\theta}Z_s^{t,x,P_{\xi}}=0$, $t\leq s<\theta\leq T$;\\
\indent \emph{(ii)} $\{D_{\theta}Y^{t,x,P_{\xi}}=(D_{\theta}^iY^{t,x,P_{\xi}})_{1\leq i\leq d},
D_{\theta}Z^{t,x,P_{\xi}}=(D_{\theta}^iZ^{t,x,P_{\xi},j})_{1\leq i,j\leq d}:s\in[\theta,T]\}$
is the unique solution of the linear BDSDE: $s\in[t,T]$, $1\leq i,j\leq d$,
\begin{equation}\label{eq4.12}
  \begin{split}
 &D_{\theta}^i Y_s^{t,x,P_{\xi}}= \sum_{j=1}^d\partial_{x_j}\Phi (X_T^{t,x,P_{\xi}},P_{X_T^{t,\xi}})D_{\theta}^iX_T^{t,x,P_{\xi},j}
+\int_s^T\Bigg\{\sum_{j=1}^d\partial_{x_j}f(\Pi_r^{t,x,P_{\xi}},P_{\Pi_r^{t,\xi}})D_{\theta}^iX_r^{t,x,P_{\xi},j}\\
&\ +\partial_{y}f(\Pi_r^{t,x,P_{\xi}},P_{\Pi_r^{t,\xi}})D_{\theta}^iY_r^{t,x,P_{\xi}}
 +\sum_{j=1}^d\partial_{z_j}f(\Pi_r^{t,x,P_{\xi}},P_{\Pi_r^{t,\xi}})D_{\theta}^iZ_r^{t,x,P_{\xi},j}\Bigg\}dr\\
 &\ +\int_s^T\sum_{n=1}^l\Bigg\{\sum_{j=1}^d\partial_{x_j}g_n(\Pi_r^{t,x,P_{\xi}},P_{\Pi_r^{t,\xi}})D_{\theta}^iX_r^{t,x,P_{\xi},j}
+\partial_{y}g_n(\Pi_r^{t,x,P_{\xi}},P_{\Pi_r^{t,\xi}})D_{\theta}^iY_r^{t,x,P_{\xi}}\\
 &\ +\!\sum_{j=1}^d\!\partial_{z_j}g_n(\Pi_r^{t,x,P_{\xi}},P_{\Pi_r^{t,\xi}})D_{\theta}^iZ_r^{t,x,P_{\xi},j}\!\Bigg\}d\overleftarrow{B_r^n}
 \!-\!\int_s^T\!\sum_{j=1}^d\!D_{\theta}^iZ_r^{t,x,P_{\xi},j}dW_r^j,\ d \theta dP\mbox{-a.e.},\ t\leq\theta\leq s,
  \end{split}
\end{equation}
where $\Pi_r^{t,x,P_{\xi}}=(X_r^{t,x,P_{\xi}},Y_r^{t,x,P_{\xi}},Z_r^{t,x,P_{\xi}})$,
$\Pi_r^{t,\xi}=\Pi_r^{t,x,P_{\xi}}\big|_{x=\xi}=(X_r^{t,\xi},Y_r^{t,\xi},Z_r^{t,\xi})$.
Moreover,
\begin{equation}\label{eq4.12+1}
Z_s^{t,x,P_{\xi}}=P\mbox{-}\displaystyle{\lim_{s<u\downarrow s}}D_{s}Y_u^{t,x,P_{\xi}},\ d s dP\mbox{-a.e.}
\end{equation}
Furthermore, if in addition Assumption (H4.2) is satisfied, then for all $ p\in[2,p_0]$, there exists a constant $C_p>0$ only depending on the Lipschitz constants of the coefficients,
 such that for all $t\in[0,T]$, $x,\ x'\in\mathbb{R}^d$,
$\xi,\ \xi'\in L^{2}(\mathcal{G}_t;\mathbb{R}^d)$, $M\geq1$, $P$-a.s.,
\begin{equation}\label{eq4.13}
  \begin{split}
    &\ \mbox{\emph{(i)}}\ E\[\sup_{s\in[t,T]}|D_{\theta}Y_s^{t,x,P_{\xi}}|^p+(\int_t^T|D_{\theta}Z_s^{t,x,P_{\xi}}|^2ds)^{\frac{p}{2}}\]\leq C_p;\\
  \end{split}
\end{equation}
\begin{equation}\nonumber
  \begin{split}
  &\ \mbox{\emph{(ii)}}\ E\[\sup_{s\in[t,T]}|D_{\theta}Y_s^{t,x,P_{\xi}}-D_{\theta}Y_s^{t,x',P_{\xi'}}|^p
  +(\int_t^T|D_{\theta}Z_s^{t,x,P_{\xi}}-D_{\theta}Z_s^{t,x',P_{\xi'}}|^2ds)^{\frac{p}{2}}\]\\
 &\ \hspace{30pt} \leq C_pM^p\(|x-x'|^p+W_2(P_{\xi},P_{\xi'})^p\)+\rho_{M,p,\theta}(t,x,P_{\xi}),
  \end{split}
\end{equation}
where $\rho_{M,p,\theta}(t,x,P_{\xi})\rightarrow0$ and (bounded convergence)
$E[\rho_{M,p,\theta}(t,\xi,P_{\xi})]\rightarrow0$, as $M\rightarrow\infty$.
In particular, for all $ p\in[2,p_0]$ there exists a constant $C_p>0$ only depending on the Lipschitz constants of the coefficients,
 such that for all $x,x'\in\mathbb{R}^d$,
$\xi,\xi'\in L^{2}(\mathcal{G}_t;\mathbb{R}^d)$, $d s dP$-a.e., $s\in[t,T]$,
\begin{equation}\label{eq4.14}
  \begin{split}
    &\ \mbox{\emph{(i)}}\ E[|Z_s^{t,x,P_{\xi}}|^p]\leq C_p;\\
  &\ \mbox{\emph{(ii)}}\ E[|Z_s^{t,x,P_{\xi}}-Z_s^{t,x',P_{\xi'}}|^p] \leq C_pM^p\(|x-x'|^p+W_2(P_{\xi},P_{\xi'})^p\)+\rho_{M,p}(t,x,P_{\xi}),
  \end{split}
\end{equation}
with $M\geq 1$, $\rho_{M,p}(t,x,P_{\xi})\rightarrow0$, as $M\rightarrow\infty$,
$E[\rho_{M,p}(t,\xi,P_{\xi})]\rightarrow0$, as $M\rightarrow\infty$.
\end{proposition}
\begin{proof}
For simplicity of redaction but without loss of generality, let us restrict to the dimensions $d = 1$, $l=1$
and to $f(x,y,z,\gamma)=f(z,\gamma(\mathbb{R}\times\mathbb{R}\times\cdot))$,
$g(x,y,z,\gamma)=g(z,\gamma(\mathbb{R}\times\mathbb{R}\times\cdot))$, $h(\gamma)=0$,
$(x,y,z)\in\mathbb{R}\times\mathbb{R}\times\mathbb{R}$, $\gamma\in\mathcal{P}_2(\mathbb{R}\times\mathbb{R}\times\mathbb{R})$
and $\Phi(x,\gamma)=\Phi(x)$, $(x,\gamma)\in\mathbb{R}\times\mathcal{P}_2(\mathbb{R})$.
Remark that as a direct consequence of assumption (H4.1)-(ii), $|\partial_zg(z,\gamma(\mathbb{R}\times\mathbb{R}\times\cdot))|^2\leq \alpha_1$,
for all $ (z,\gamma)\in\mathbb{R}\times\mathcal{P}_2(\mathbb{R}\times\mathbb{R}\times\mathbb{R})$.

It is standard to prove that $Y_s^{t,x,P_{\xi}}$ and $Z_s^{t,x,P_{\xi}}$ are Malliavin differentiable under our assumptions,
and $Z_r^{t,x,P_{\xi}}=P\mbox{-}\displaystyle{\lim_{r<s\downarrow r}}D_rY_s^{t,x,P_{\xi}}$, so we omit its proof here. Thus, it suffices to prove \eqref{eq4.13}.
 For this note first that Lemma \ref{leA.1}-(2) combined with Lemma \ref{le3.1}, Proposition \ref{prop3.1} applied to BDSDE \eqref{eq4.12} yields
 $\displaystyle E\[\sup_{s\in[t,T]}|D_{\theta}Y_s^{t,x,P_{\xi}}|^{p_0}+(\int_t^T|D_{\theta}Z_s^{t,x,P_{\xi}}|^2ds)^{\frac{p_0}{2}}\]\leq C_{p_0},$
then from H\"{o}lder's inequality, for all $ p\in[2,p_0]$ we have $\displaystyle E\[\sup_{s\in[t,T]}|D_{\theta}Y_s^{t,x,P_{\xi}}|^p+(\int_t^T|D_{\theta}Z_s^{t,x,P_{\xi}}|^2ds)^{\frac{p}{2}}\]\leq C_p$. Now we want to prove \eqref{eq4.13}-(ii).

For all $0\leq t\leq T$, $t\leq\theta\leq s\leq T$, $x,x'\in\mathbb{R}$, $\xi,\xi'\in L^{2}(\mathcal{G}_t;\mathbb{R})$, from \eqref{eq4.12} we get the following BDSDE:
\begin{equation}\nonumber
  \begin{split}
 &\ D_{\theta}Y_s^{t,x,P_{\xi}}-D_{\theta} Y_s^{t,x',P_{\xi'}}
 = \ I(t,x,P_{\xi},x',P_{\xi'})+\int_s^TR(r,x,P_{\xi},x',P_{\xi'})dr+\int_s^TH(r,x,P_{\xi},x',P_{\xi'})d\overleftarrow{B_r}\\
 &\ +\int_s^T\partial_{z}g(Z_r^{t,x,P_{\xi}},P_{Z_r^{t,\xi}})\(D_{\theta}Z_r^{t,x,P_{\xi}}-D_{\theta}Z_r^{t,x',P_{\xi'}}\)d\overleftarrow{B_r}
  -\int_s^T\(D_{\theta}Z_r^{t,x,P_{\xi}}-D_{\theta}Z_r^{t,x',P_{\xi'}}\)dW_r\\
  \end{split}
\end{equation}
\begin{equation}\label{eq4.15}
  \begin{split}
 &\ +\int_s^T\partial_{z}f(Z_r^{t,x,P_{\xi}},P_{Z_r^{t,\xi}})\(D_{\theta}Z_r^{t,x,P_{\xi}}-D_{\theta}Z_r^{t,x',P_{\xi'}}\)dr,
  \end{split}
\end{equation}
where
\begin{equation*}
  \begin{split}
 &\ I(t,x,P_{\xi},x',P_{\xi'}):= \partial_{x}\Phi (X_T^{t,x,P_{\xi}})D_{\theta}X_T^{t,x,P_{\xi}}-\partial_{x}\Phi (X_T^{t,x',P_{\xi'}})D_{\theta}X_T^{t,x',P_{\xi'}},\\
 &\ R(r,x,P_{\xi},x',P_{\xi'}):=\(\partial_{z}f(Z_r^{t,x,P_{\xi}},P_{Z_r^{t,\xi}})
         -\partial_{z}f(Z_r^{t,x',P_{\xi'}},P_{Z_r^{t,\xi'}})\)D_{\theta}Z_r^{t,x',P_{\xi'}},\\
 &\ H(r,x,P_{\xi},x',P_{\xi'}):=\(\partial_{z}g(Z_r^{t,x,P_{\xi}},P_{Z_r^{t,\xi}})
         -\partial_{z}g(Z_r^{t,x',P_{\xi'}},P_{Z_r^{t,\xi'}})\)D_{\theta}Z_r^{t,x',P_{\xi'}}.
  \end{split}
\end{equation*}
From the Lemmas \ref{leA.1} and \ref{le3.1} and the Propositions \ref{prop3.1} and \ref{prop4.1} and our assumptions we have,
\begin{equation}\nonumber
  \begin{split}
  &\  E\[\sup_{s\in[\theta,T]}|D_{\theta}Y_s^{t,x,P_{\xi}}-D_{\theta}Y_s^{t,x',P_{\xi'}}|^{p_0}
  +(\int_{\theta}^T|D_{\theta}Z_s^{t,x,P_{\xi}}-D_{\theta}Z_s^{t,x',P_{\xi'}}|^2ds)^{\frac{p_0}{2}}\]\\
 \leq&\  C_{p_0} E\[|I(t,x,P_{\xi},x',P_{\xi'})|^{p_0}\!+\!(\int_{\theta}^{T}\!\!|R(r,x,P_{\xi},x',P_{\xi'})|dr)^{p_0}\!
 \!+\!\!(\int_{\theta}^{T}\!\!|H(r,x,P_{\xi},x',P_{\xi'})|^2dr)^{\frac{p_0}{2}}\]\\
  \end{split}
\end{equation}
\begin{equation}\label{eq4.16}
  \begin{split}
 \leq&\  C_{p_0} E\[\(\int_{\theta}^T|(\partial_{z}f)(Z_r^{t,x,P_{\xi}},P_{Z_r^{t,\xi}})-(\partial_{z}f)(Z_r^{t,x',P_{\xi'}},P_{Z_r^{t,\xi'}})|^2
 \cdot|D_{\theta}Z_r^{t,x',P_{\xi'}}|^2dr\)^{\frac{p_0}{2}}\]\\
 &\ +C_{p_0} E\[\(\int_{\theta}^T|(\partial_{z}g)(Z_r^{t,x,P_{\xi}},P_{Z_r^{t,\xi}})-(\partial_{z}g)(Z_r^{t,x',P_{\xi'}},P_{Z_r^{t,\xi'}})|^2
 \cdot|D_{\theta}Z_r^{t,x',P_{\xi'}}|^2dr\)^{\frac{p_0}{2}}\]\\
&\ + C_p\(|x-x'|^{p_0}+W_2(P_{\xi},P_{\xi'})^{p_0}\)\\
\leq &\ C_{p_0} E\[\(\int_{\theta}^T\min\big\{C,|Z_r^{t,x,P_{\xi}}-Z_r^{t,x',P_{\xi'}}|^2+W_2(P_{Z_r^{t,\xi}},P_{Z_r^{t,\xi'}})^2\big\}
  \cdot|D_{\theta}Z_r^{t,x',P_{\xi'}}|^2dr\)^{\frac{p_0}{2}}\]  \\
 &\ + C_{p_0}\(|x-x'|^{p_0}+W_2(P_{\xi},P_{\xi'})^{p_0}\)\\
\leq&\ C_{p_0}M^{p_0}E\[(\int_{\theta}^T(|Z_r^{t,x,P_{\xi}}-Z_r^{t,x',P_{\xi'}}|^2+W_2(P_{Z_r^{t,\xi}},P_{Z_r^{t,\xi'}})^2)dr)^{\frac{p_0}{2}}\]\\
  &\ +C_{p_0}E\[\(\int_{\theta}^T|D_{\theta}Z_r^{t,x',P_{\xi'}}|^2I_{\{|D_{\theta}Z_r^{t,x',P_{\xi'}}|\geq M\}}dr\)^{\frac{p_0}{2}}\]
  + C_{p_0}\(|x-x'|^{p_0}+W_2(P_{\xi},P_{\xi'})^{p_0}\)\\
  \leq&\ C_{p_0}M^{p_0}\(|x-x'|^{p_0}+W_2(P_{\xi},P_{\xi'})^{p_0}\)+\rho_{M,p_0,\theta}(t,x',P_{\xi'}),
  \end{split}
\end{equation}
where $\rho_{M,p_0,\theta}(t,x',P_{\xi'})\rightarrow0$, $E[\rho_{M,p_0,\theta}(t,\xi',P_{\xi'})]\rightarrow0$ ($M\rightarrow\infty$) thanks to the
dominated convergence theorem (indeed, $\sup_{x'\in\mathbb{R}}E\[\(\int_t^T|D_{\theta}Z_r^{t,x',P_{\xi'}}|^2dr\)^{\frac{p_0}{2}}\]\leq C_{p_0}<\infty$, $\theta\in[0,T]$).

Again from H\"{o}lder's inequality, for all $p\in[2,p_0]$ we have
\begin{equation*}
  \begin{split}
  &\  E\[\sup_{s\in[t,T]}|D_{\theta}Y_s^{t,x,P_{\xi}}-D_{\theta}Y_s^{t,x',P_{\xi'}}|^p
  +(\int_t^T|D_{\theta}Z_s^{t,x,P_{\xi}}-D_{\theta}Z_s^{t,x',P_{\xi'}}|^2ds)^{\frac{p}{2}}\]\\
 \leq &\ C_pM^p\(|x-x'|^p+W_2(P_{\xi},P_{\xi'})^p\)+\rho_{M,p,\theta}(t,x,P_{\xi}),
  \end{split}
\end{equation*}
where $\rho_{M,p,\theta}(t,x,P_{\xi})\rightarrow0$ and
$E[\rho_{M,p,\theta}(t,\xi,P_{\xi})]\rightarrow0$, as $M\rightarrow\infty$.

Finally, observe that, thanks to \eqref{eq4.12+1}, \eqref{eq4.14} is an immediate consequence of \eqref{eq4.13}.
\end{proof}

\begin{remark} \label{re4.1+1}
Let us point out that for $p=2$, \eqref{eq4.13} and \eqref{eq4.14} also remain true without supposing
 (H4.2). The proof can be carried out similarly to that of Proposition \ref{prop4.0}.
\end{remark}

 Now we introduce the value function
\begin{equation}\label{eq4.17}
V(t,x,P_{\xi}):=Y_t^{t,x,P_{\xi}}.
\end{equation}
Notice that $V(t,x,P_{\xi})$ is $\mathcal{F}_{t,T}^{B}$-measurable, for all $(t,x)$. On the other hand,
from Proposition \ref{prop4.1} we get
\begin{equation}\label{eq4.18}
V(s,X_s^{t,x,P_{\xi}},P_{X_s^{t,\xi}})=Y_s^{s,X_s^{t,x,\xi},P_{X_s^{t,\xi}}}=Y_s^{t,x,P_{\xi}},\ s\in[t,T].
\end{equation}
A direct consequence of Proposition \ref{prop4.1} is the following regularity property of the value function.
\begin{proposition} \label{prop4.4}
Let $b$, $\sigma$, $f$, $g$, $h$ and $\Phi$ satisfy the assumptions in Proposition \ref{prop4.2}.
Then, $[0,T]\times\mathbb{R}^d\times\mathcal{P}_2(\mathbb{R}^d)\ni (t,x,\mu)\rightarrow V(t,x,\mu)$ is continuous
in $L^p$, for all $p\in[2,p_0]$. More precisely, for all $p\in[2,p_0]$,
there exists a constant $C_p\in\mathbb{R}_{+}$ such that, for all $t,t'\in[0,T]$, $x,x'\in\mathbb{R}^d$,
$\mu,\mu'\in\mathcal{P}_2(\mathbb{R}^d)$,
\begin{equation}\label{eq4.18+1}
  \begin{split}
  E[|V(t,x,\mu)-V(t',x',\mu')|^p]
\leq C_p \(|t-t'|^{\frac{p}{2}}+|x-x'|^{p}+W_2(\mu,\mu')^p\).
  \end{split}
\end{equation}
\end{proposition}
\begin{proof}
As a direct consequence of Proposition \ref{prop4.1}, we have that, for all $p\in[2,p_0]$,
there is some constant $C_p\in\mathbb{R}_{+}$ such that, for all $t\in[0,T]$, $x,x'\in\mathbb{R}^d$,
$\mu,\mu'\in\mathcal{P}_2(\mathbb{R}^d)$,
\begin{equation*}
   E[|V(t,x,\mu)-V(t,x',\mu')|^p]\leq C_p \(|x-x'|^{p}+W_2(\mu,\mu')^p\).
\end{equation*}
Let us prove now the $\frac{1}{2}$-H\"{o}lder continuity of $V$ in $L^p$.
Without loss of generality we let $0\leq t<t'$. Then, for all $p\in[2,p_0]$, we have
\begin{equation}\label{eq4.19}
  \begin{split}
  E[|V(t,x,P_{\xi})-V(t',x,P_{\xi})|^p]
\leq C E[|Y_t^{t,x,P_{\xi}}-Y_{t'}^{t,x,P_{\xi}}|^p]
+CE[|Y_{t'}^{t,x,P_{\xi}}-Y_{t'}^{t',x,P_{\xi}}|^p].
  \end{split}
\end{equation}
On the one hand, as $f$, $g$ and $h$ are bounded, it follows from \eqref{eq4.14}-(ii) that
\begin{equation}\label{eq4.20}
  \begin{split}
  &\ E[|Y_t^{t,x,P_{\xi}}-Y_{t'}^{t,x,P_{\xi}}|^p]\\
\leq&\ CE[(\int_t^{t'}|f(\Pi_s^{t,x,P_{\xi}},P_{\Pi_s^{t,\xi}})|ds)^p]
+CE[(\int_t^{t'}|g(\Pi_s^{t,x,P_{\xi}},P_{\Pi_s^{t,\xi}})+h(P_{\Pi_s^{t,\xi}})|^2ds)^{\frac{p}{2}}]\\
 &\ +CE[(\int_t^{t'}|Z_s^{t,x,P_{\xi} }|^2ds)^{\frac{p}{2}}]\\
\leq&\ C|t'-t|^{\frac{p}{2}}+C|t'-t|^{\frac{p-2}{2}}\int_t^{t'}E[|Z_s^{t,x,P_{\xi} }|^p]ds
\leq\ C|t'-t|^{\frac{p}{2}}.
  \end{split}
\end{equation}
On the other hand, we deduce from \eqref{eq4.11}, Proposition \ref{prop4.1} and Lemma \ref{le3.1},
\begin{equation}\label{eq4.21}
  \begin{split}
  &\ E[|Y_{t'}^{t,x,P_{\xi}}-Y_{t'}^{t',x,P_{\xi}}|^p]
  =E[E[|Y_{t'}^{t',X_{t'}^{t,x,P_{\xi}},P_{X_{t'}^{t,\xi}}}-Y_{t'}^{t',x,P_{\xi}}|^p\big|\mathcal{G}_{t'}]]\\
\leq&\ CE[|X_{t'}^{t,x,P_{\xi}}-x|^{p}+W_2(P_{X_{t'}^{t,\xi}},P_{\xi})^p]\\
\leq&\ CE[|X_{t'}^{t,x,P_{\xi}}-x|^{p}+|X_{t'}^{t,\xi}-\xi|^p]\leq\ C|t'-t|^{\frac{p}{2}}.
  \end{split}
\end{equation}
Therefore, combining \eqref{eq4.19}, \eqref{eq4.20} and \eqref{eq4.21} with (ii) of Proposition \ref{prop4.1}, and applying Kolmogorov's continuity criterion, we get that
$V(t,x,P_{\xi})$ has a version continuous with respect to $(t,x)$.
\end{proof}

\begin{remark} \label{re4.2}
\eqref{eq4.18+1} of Proposition \ref{prop4.4} implies that for all fixed $\mu\in\mathcal{P}_2(\mathbb{R}^d)$,
$V(\cdot,\cdot,\mu)=\{V(t,x,\mu),$ $(t,x)\in[0,T]\times\mathbb{R}^d\}$ admits a continuous version.
As this conclusion is based on Kolmogorov's continuity criterion, and so on the finite-dimensionality of $[0,T]\times\mathbb{R}^d$,
it cannot be extended to $(t,x,\mu)\rightarrow V(t,x,\mu)$.
\end{remark}

\section{First order derivatives of $X^{t,x,P_{\xi}}$}

In this section we revisit the first order derivatives of $X^{t,x,P_{\xi}}$ with respect to $x$ and the measure
$P_{\xi}$, studied by Buckdahn, Li, Peng and Rainer \cite{BLPR2017}. For the reader's convenience we give the main results here, for
more details the reader is referred to \cite{BLPR2017}, or to Hao and Li \cite{HL2016}, where mean-field SDEs with jump are studied.

\noindent\textbf{Assumption (H5.1)} The couple of coefficients $(b,\sigma)$ belongs to
 $C_b^{1,1}(\mathbb{R}^{d}\times\mathcal{P}_{2}(\mathbb{R}^{d});\mathbb{R}^{d}\times\mathbb{R}^{d\times d})$, that is,
 the components $b_j$, $\sigma_{i,j}$, $1\leq i,j\leq d$,  have the following properties:\\
\indent (i) $b_j(x,\cdot)$, $\sigma_{i,j}(x,\cdot)$ belong to $C_b^{1}(\mathcal{P}_{2}(\mathbb{R}^{d}))$, for all $x\in\mathbb{R}^{d}$;\\
\indent (ii) $b_j(\cdot,\mu)$, $\sigma_{i,j}(\cdot,\mu)$ belong to $C_b^{1}(\mathbb{R}^{d})$, for all $\mu\in\mathcal{P}_{2}(\mathbb{R}^{d})$;\\
\indent (iii) The derivatives $\partial_xb_j$, $\partial_x\sigma_{i,j}$: $\mathbb{R}^{d}\times\mathcal{P}_{2}(\mathbb{R}^{d})\rightarrow\mathbb{R}^{d}$
and $\partial_{\mu}b_j$, $\partial_{\mu}\sigma_{i,j}$: $\mathbb{R}^{d}\times\mathcal{P}_{2}(\mathbb{R}^{d})\times\mathbb{R}^{d}\rightarrow\mathbb{R}^{d}$
are bounded and Lipschitz continuous.

We begin with recalling the first order derivative of $X^{t,x,P_{\xi}}$ with respect to $x$.
\begin{theorem} \label{th5.1}
Suppose Assumption (H5.1) holds true. Then the $L^2$-derivative of $X^{t,x,P_{\xi}}$ with
respect to $x$ exists, it is denoted by $\partial_xX^{t,x,P_{\xi}}=(\partial_xX^{t,x,P_{\xi},j})_{1\leq j\leq d}$, and it satisfies the
following SDE: $s\in[t,T]$, $1\leq i,j\leq d$,
\begin{equation}\label{eq5.1}
  \begin{split}
 \partial_{x_i}X_s^{t,x,P_{\xi},j}=&\ \delta_{ij}+\sum_{k=1}^d\int_t^s\partial_{x_k}b_j(X_r^{t,x,P_{\xi}},P_{X_r^{t,\xi}}) \partial_{x_i}X_s^{t,x,P_{\xi},k}dr\\
 &\ + \sum_{k,l=1}^d\int_t^s\partial_{x_k}\sigma_{j,l}(X_r^{t,x,P_{\xi}},P_{X_r^{t,\xi}}) \partial_{x_i}X_s^{t,x,P_{\xi},k}dW_r^l,
  \end{split}
\end{equation}
where $\delta_{ij}$ denotes the Kronecker symbol: It equals $1$, if $i = j$, and is equal to zero, otherwise.
\end{theorem}
For the proof the reader is referred to Theorem 3.1 in \cite{BLPR2017}, and for the case with jumps
also to Theorem 4.1 in \cite{HL2016}. From standard estimates for classical SDEs we have
\begin{proposition} \label{prop5.1}
For all $p\geq 2$, there exists a constant $C_p>0$ only depending on the Lipschitz
constants of $\partial_xb$ and $\partial_x\sigma$, such that, for all $t\in[0,T]$, $x,x'\in\mathbb{R}^d$, $\xi,\xi'\in L^{2}(\mathcal{G}_t;\mathbb{R}^d)$,
$P$-a.s.,
\begin{equation*}
  \begin{split}
  &\ \emph{(i)}\ E\[\sup_{s\in[t,T]}|\partial_xX_s^{t,x,P_{\xi}}|^p\big|\mathcal{G}_t\]\leq C_p;\\
  &\ \emph{(ii)}\ E\[\sup_{s\in[t,T]}|\partial_xX_s^{t,x,P_{\xi}}-\partial_xX_s^{t,x',P_{\xi'}}|^p\big|\mathcal{G}_t\]
       \leq C_p\(|x-x'|^p+W_2(P_{\xi},P_{\xi'})^p\);\\
&\ \emph{(iii)}\ E\[\sup_{s\in[t,t+h]}|\partial_xX_s^{t,x,P_{\xi}}-I_{d\times d}|^p\big|\mathcal{G}_t\]\leq C_ph^{\frac{p}{2}},\ 0\leq t\leq t+h\leq T.
  \end{split}
\end{equation*}
Here $I_{d\times d}$ denotes the unit matrix in dimension $d$.
\end{proposition}
The following theorem shows that the unique solution $X^{t,x,\xi}$ of Eq. \eqref{eq3.2} interpreted as a
functional of $\xi\in L^{2}(\mathcal{G}_t;\mathbb{R}^d)$ is Fr\'{e}chet differentiable.
\begin{theorem} \label{th5.2}
Let $(b,\sigma)$ satisfy Assumption (H5.1). Then, for all $0\leq t\leq s\leq T$, $x\in\mathbb{R}^d$, the
lifted process $L^{2}(\mathcal{G}_t;\mathbb{R}^d)\ni\xi\rightarrow X_s^{t,x,\xi}\in L^{2}(\mathcal{G}_s;\mathbb{R}^d)$ is Fr\'{e}chet differentiable, and
the Fr\'{e}chet derivative is given by
\begin{equation*}
  DX_s^{t,x,\xi}(\eta)=\widetilde{E}\[U_s^{t,x,P_{\xi}}(\widetilde{\xi})\widetilde{\eta}\]
  =\(\widetilde{E}\[\sum_{j=1}^dU_{s,i,j}^{t,x,P_{\xi}}(\widetilde{\xi})\widetilde{\eta}_j\]\)_{1\leq i\leq d},
\end{equation*}
for all $\eta\!=\!(\eta_1,\eta_2,\cdots,\eta_d)\!\in\! L^{2}(\mathcal{G}_t;\mathbb{R}^d)$, where, for $y\in\mathbb{R}^d$,
$U^{t,x,P_{\xi}}(y)\!=\!((U_{s,i,j}^{t,x,P_{\xi}}(y))_{s\in[t,T]})_{1\leq i,j\leq d}\in\mathcal{S}^2_{\mathcal{G}}(t,T;\mathbb{R}^{d\times d})$
is the unique solution of the following SDE: $s\in[t,T],\ 1\leq i,j\leq d$,
\begin{equation}\label{eq5.2}
  \begin{split}
  U_{s,i,j}^{t,x,P_{\xi}}(y)=&\ \sum_{k=1}^d\int_t^s\partial_{x_k}b_i(X_r^{t,x,P_{\xi}},P_{X_r^{t,\xi}}) U_{r,k,j}^{t,x,P_{\xi}}(y)dr\\
  &\ +\sum_{k,l=1}^d\int_t^s\partial_{x_k}\sigma_{i,l}(X_r^{t,x,P_{\xi}},P_{X_r^{t,\xi}}) U_{r,k,j}^{t,x,P_{\xi}}(y)dW_r^l\\
  \end{split}
\end{equation}
\begin{equation}\nonumber
  \begin{split}
  &\ +\sum_{k=1}^d\int_t^sE[(\partial_{\mu}b_i)(z,P_{X_r^{t,\xi}},X_r^{t,y,P_{\xi}}) \partial_{x_j}X_s^{t,y,P_{\xi},k}\\
 &\ \hspace{70pt}+(\partial_{\mu}b_i)(z,P_{X_r^{t,\xi}},X_r^{t,\xi})U_{r,k,j}^{t,\xi}(y) ]\big|_{z=X_r^{t,x,P_{\xi}}}dr\\
 &\ +\sum_{k,l=1}^d\int_t^sE[(\partial_{\mu}\sigma_{i,l})(z,P_{X_r^{t,\xi}},X_r^{t,y,P_{\xi}}) \partial_{x_j}X_s^{t,y,P_{\xi},k}\\
&\ \hspace{70pt}+(\partial_{\mu}\sigma_{i,l})(z,P_{X_r^{t,\xi}},X_r^{t,\xi})U_{r,k,j}^{t,\xi}(y) ]\big|_{z=X_r^{t,x,P_{\xi}}}dW_r^l,
  \end{split}
\end{equation}
where $U^{t,\xi}(y)=((U_{s,i,j}^{t,\xi}(y))_{s\in[t,T]})_{1\leq i,j\leq d}=U^{t,x,P_{\xi}}(y)\big|_{x=\xi}\in\mathcal{S}^2_{\mathcal{G}}(t,T;\mathbb{R}^{d\times d})$
satisfies \eqref{eq5.2} with $x$ replaced by $\xi$.
\end{theorem}
\begin{proposition} \label{prop5.2}
For all $p\geq 2$, there exists a constant $C_p>0$ only depending on the Lipschitz
constants of $b$ and $\sigma$, such that, for all $t\in[0,T]$, $x,x',y,y'\in\mathbb{R}^d$, $\xi,\xi'\in L^{2}(\mathcal{G}_t;\mathbb{R}^d)$,
$P$-a.s.,
\begin{equation*}
  \begin{split}
  &\ \emph{(i)}\ E\[\sup_{s\in[t,T]}(|U_s^{t,x,P_{\xi}}(y)|^p+|U_s^{t,\xi}(y)|^p)\]\leq C_p;\\
  &\ \emph{(ii)}\ E\[\sup_{s\in[t,T]}|U_s^{t,x,P_{\xi}}(y)-U_s^{t,x',P_{\xi'}}(y')|^p\]
      \leq C_p\(|x-x'|^p+|y-y'|^p+W_2(P_{\xi},P_{\xi'})^p\);\\
&\ \emph{(iii)}\ E\[\sup_{s\in[t,t+h]}|U_s^{t,x,P_{\xi}}(y)|^p\]\leq C_ph^{\frac{p}{2}},\ 0\leq t\leq t+h\leq T.
  \end{split}
\end{equation*}
\end{proposition}
For the proof of Theorem \ref{th5.2} and Proposition \ref{prop5.2} we refer the reader to Section 4 in \cite{BLPR2017},
and for the case with jumps also to Section 4 in \cite{HL2016} (in the case of jumps we have $h$ instead of $h^{\frac{p}{2}}$ in (iii)).

Furthermore, the derivative of $X_s^{t,x,P_{\xi}}$ with respect
to the probability measure can be defined as follows
\begin{equation*}
  \partial_{\mu}X_s^{t,x,P_{\xi}}(y):=U_s^{t,x,P_{\xi}}(y),\ s\in[t,T],\ t\in[0,T], x\in\mathbb{R}^d,\ \xi\in L^{2}(\mathcal{G}_t;\mathbb{R}^d),\ y \in\mathbb{R}^d.
\end{equation*}
Thus, we get $DX_s^{t,x,\xi}=\widetilde{E}\[\partial_{\mu}X_s^{t,x,P_{\xi}}(\widetilde{\xi})\widetilde{\eta}\]$,
for all $\eta\in L^{2}(\mathcal{G}_t;\mathbb{R}^d)$.
\medskip

As an immediate result of Proposition \ref{prop5.2}, we have
\begin{proposition} \label{prop5.3}
For all $p\geq 2$, there exists a constant $C_p>0$ only depending on the Lipschitz
constants of $b$ and $\sigma$, such that, for all $t\in[0,T]$, $x,x',y,y'\in\mathbb{R}^d$, $\xi,\xi'\in L^{2}(\mathcal{G}_t;\mathbb{R}^d)$, $P$-a.s.,
\begin{equation*}
  \begin{split}
  &\ \emph{(i)}\ E\[\sup_{s\in[t,T]}| \partial_{\mu}X_s^{t,x,P_{\xi}}(y)|^p\big|\mathcal{G}_t\]\leq C_p;\\
  &\ \emph{(ii)}\ E\[\sup_{s\in[t,T]}| \partial_{\mu}X_s^{t,x,P_{\xi}}(y)- \partial_{\mu}X_s^{t,x',P_{\xi'}}(y')|^p\big|\mathcal{G}_t\]
 \leq C_p\(|x-x'|^p+|y-y'|^p+W_2(P_{\xi},P_{\xi'})^p\);\\
&\ \emph{(iii)}\ E\[\sup_{s\in[t,t+h]}| \partial_{\mu}X_s^{t,x,P_{\xi}}(y)|^p\big|\mathcal{G}_t\]\leq C_ph^{\frac{p}{2}},\ 0\leq t\leq t+h\leq T.
  \end{split}
\end{equation*}
\end{proposition}

\section{First order derivatives of $(Y^{t,x,P_{\xi}},Z^{t,x,P_{\xi}})$}
We recall from Proposition \ref{prop4.1} that $(Y^{t,x,P_{\xi}},Z^{t,x,P_{\xi}})$ depends on $\xi$ only through its
law, which allows to define $(Y^{t,x,P_{\xi}},Z^{t,x,P_{\xi}}):=(Y^{t,x,\xi},Z^{t,x,\xi})$. This section is
devoted to the study of the first order derivatives of $(Y^{t,x,P_{\xi}},Z^{t,x,P_{\xi}})$ with respect to $x$ and
$P_{\xi}$, respectively.

\noindent\textbf{Assumption (H6.1)} Let $\Phi\!\in\! C_b^{1,1}(\mathbb{R}^{d}\times\mathcal{P}_{2}(\mathbb{R}^{d}))$,
$f\!\in\! C_b^{1,1}(\mathbb{R}^{d+1+d}\times\mathcal{P}_{2}(\mathbb{R}^{d+1+d}))$,
$g\!\in\! C_b^{1,1}(\mathbb{R}^{d+1+d}\times\mathcal{P}_{2}(\mathbb{R}^{d+1+d});\mathbb{R}^{l})$
and $h\!\in\! C_b^{1}(\mathcal{P}_{2}(\mathbb{R}^{d+1+d});\mathbb{R}^{l})$.
In addition we suppose Assumption (H4.1)-(ii).
\begin{theorem} \label{th6.1}
Under the Assumptions (H5.1) and (H6.1), the $L^2$-derivative of the solution of
Eq. \eqref{eq4.2} with respect to $x$, $(\partial_xY^{t,x,P_{\xi}},\partial_xZ^{t,x,P_{\xi}})$, exists and is the unique solution of
the following BDSDE: $ s\in[t,T],\ 1\leq i\leq d$,
\begin{equation}\label{eq6.1}
  \begin{split}
 \partial_{x_i}Y_s^{t,x,P_{\xi}}=&\ \sum_{j=1}^d\partial_{x_j}\Phi (X_T^{t,x,P_{\xi}},P_{X_T^{t,\xi}})\partial_{x_i}X_T^{t,x,P_{\xi},j}
+\int_s^T\Bigg\{\sum_{j=1}^d\partial_{x_j}f(\Pi_r^{t,x,P_{\xi}},P_{\Pi_r^{t,\xi}})\partial_{x_i}X_r^{t,x,P_{\xi},j}\\
 &\ +\partial_{y}f(\Pi_r^{t,x,P_{\xi}},P_{\Pi_r^{t,\xi}})\partial_{x_i}Y_r^{t,x,P_{\xi}}
 +\sum_{j=1}^d\partial_{z_j}f(\Pi_r^{t,x,P_{\xi}},P_{\Pi_r^{t,\xi}})\partial_{x_i}Z_r^{t,x,P_{\xi},j}\Bigg\}dr\\
 &\ +\int_s^T\sum_{n=1}^l\Bigg\{\sum_{j=1}^d\partial_{x_j}g_n(\Pi_r^{t,x,P_{\xi}},P_{\Pi_r^{t,\xi}})\partial_{x_i}X_r^{t,x,P_{\xi},j}
+\partial_{y}g_n(\Pi_r^{t,x,P_{\xi}},P_{\Pi_r^{t,\xi}})\partial_{x_i}Y_r^{t,x,P_{\xi}}\\
 &\ +\sum_{j=1}^d\partial_{z_j}g_n(\Pi_r^{t,x,P_{\xi}},P_{\Pi_r^{t,\xi}})\partial_{x_i}Z_r^{t,x,P_{\xi},j}\Bigg\}d\overleftarrow{B_r^n}
 -\int_s^T\sum_{j=1}^d\partial_{x_i}Z_r^{t,x,P_{\xi},j}dW_r^j,
  \end{split}
\end{equation}
where $\Pi_r^{t,x,P_{\xi}}=(X_r^{t,x,P_{\xi}},Y_r^{t,x,P_{\xi}},Z_r^{t,x,P_{\xi}})$,
$\Pi_r^{t,\xi}=\Pi_r^{t,x,P_{\xi}}\big|_{x=\xi}=(X_r^{t,\xi},Y_r^{t,\xi},Z_r^{t,\xi})$.
\end{theorem}
As the $L^2$-derivative of the coefficients $f(\Pi_s^{t,x,P_{\xi}},P_{\Pi_s^{t,\xi}})$ and $g(\Pi_s^{t,x,P_{\xi}},P_{\Pi_s^{t,\xi}})$
concern only $\Pi_s^{t,x,P_{\xi}}$ but not the law $P_{\Pi_s^{t,\xi}}$,
the arguments of the proof are standard; the reader is referred, for instance, to \cite{PP1994}.

From the standard estimates-Lemma \ref{leA.1} for classical BDSDEs,
combined with Lemma \ref{le3.1}, Proposition \ref{prop4.1}, Corollary \ref{cor4.1} and Proposition \ref{prop5.1},
we have the following result (as for the estimate (ii) for $(\partial_xY_s^{t,x,P_{\xi}},\partial_xZ_s^{t,x,P_{\xi}})$,
the reader may also refer to the proof of Proposition \ref{prop4.2}).
\begin{proposition} \label{prop6.1}
 Suppose the Assumptions (H4.2), (H5.1) and (H6.1) hold true. Then, for all $p\in[2,p_0]$,
 there exists a constant $C_p>0$ only depending on the Lipschitz
constants of the coefficients, such that for all $t\in[0,T]$, $x,x'\in\mathbb{R}^d$, $\xi,\xi'\in L^{2}(\mathcal{G}_t;\mathbb{R}^d)$,
 $P$-a.s.,
\begin{equation}\label{eq6.2}
  \begin{split}
  &\ \emph{(i)}\ E\[\sup_{s\in[t,T]}|\partial_xY_s^{t,x,P_{\xi}}|^p+(\int_t^T|\partial_xZ_s^{t,x,P_{\xi}}|^2ds)^{\frac{p}{2}}ds\big|\mathcal{G}_t\]\leq C_p,\\
  &\ \emph{(ii)}\ E\[\sup_{s\in[t,T]}|\partial_xY_s^{t,x,P_{\xi}}-\partial_xY_s^{t,x',P_{\xi'}}|^p
  +(\int_t^T|\partial_xZ_s^{t,x,P_{\xi}}-\partial_xZ_s^{t,x',P_{\xi'}}|^2 ds)^{\frac{p}{2}}\big|\mathcal{G}_t\]\\
  &\hspace{1cm}     \leq C_pM^p\(|x-x'|^p+W_2(P_{\xi},P_{\xi'})^p\)+\rho_{M,p}(t,x,P_{\xi}),
  \end{split}
\end{equation}
with $M\geq1$, $\rho_{M,p}(t,x,P_{\xi})\rightarrow0$, as $M\rightarrow\infty$, $E[\rho_{M,p}(t,\xi,P_{\xi})]\rightarrow0$, as $M\rightarrow\infty$.
\end{proposition}

\begin{remark} \label{re6.0}
In analogy to the proof in Proposition \ref{prop4.0}, it can easily be checked that
under the Assumptions (H5.1) and (H6.1), there exists a constant $C>0$ only depending on the Lipschitz
constants of the coefficients, such that for all $t\in[0,T]$, $x,x'\in\mathbb{R}^d$, $\xi,\xi'\in L^{2}(\mathcal{G}_t;\mathbb{R}^d)$,
 $P$-a.s.,
\begin{equation}\label{eq6.2+1+1}
  \begin{split}
  &\ \emph{(i)}\ E\[\sup_{s\in[t,T]}|\partial_xY_s^{t,x,P_{\xi}}|^2+\int_t^T|\partial_xZ_s^{t,x,P_{\xi}}|^2dsds\big|\mathcal{G}_t\]\leq C,\\
  \end{split}
\end{equation}
\begin{equation*}
  \begin{split}
  &\ \emph{(ii)}\ E\[\sup_{s\in[t,T]}|\partial_xY_s^{t,x,P_{\xi}}-\partial_xY_s^{t,x',P_{\xi'}}|^2
  +\int_t^T|\partial_xZ_s^{t,x,P_{\xi}}-\partial_xZ_s^{t,x',P_{\xi'}}|^2 ds\big|\mathcal{G}_t\]\\
  &\hspace{1cm}     \leq CM^2\(|x-x'|^2+W_2(P_{\xi},P_{\xi'})^2\)+\rho_{M}(t,x,P_{\xi}),
  \end{split}
\end{equation*}
with $M\geq1$, $\rho_{M}(t,x,P_{\xi})\rightarrow0$, as $M\rightarrow\infty$, $E[\rho_{M}(t,\xi,P_{\xi})]\rightarrow0$, as $M\rightarrow\infty$.
\end{remark}

\begin{theorem} \label{th6.2}
Assume the Assumptions (H5.1) and (H6.1) hold. Then, for all $0\leq t\leq s \leq T$, $x\in\mathbb{R}^d$,
the lifted processes $L^{2}(\mathcal{G}_t;\mathbb{R}^d)\ni\xi\rightarrow Y_s^{t,x,\xi}:=Y_s^{t,x,P_{\xi}}\in L^{2}(\mathcal{F}_s;\mathbb{R}^d)$,
and $L^{2}(\mathcal{G}_t;\mathbb{R}^d)\ni\xi\rightarrow Z_s^{t,x,\xi}:=Z_s^{t,x,P_{\xi}}\in\mathcal{H}_{\mathcal{F}}^2(t,T;\mathbb{R}^{d})$
are Fr\'{e}chet differentiable, with the Fr\'{e}chet derivatives
\begin{equation}\label{eq6.3}
 \begin{split}
   & DY_s^{t,x,\xi}(\eta)=\overline{E}\[O_s^{t,x,P_{\xi}}(\overline{\xi})\overline{\eta}\],\ s\in[t,T],\ P\mbox{-}a.s.,\\
  & DZ_s^{t,x,\xi}(\eta)=\overline{E}\[Q_s^{t,x,P_{\xi}}(\overline{\xi})\overline{\eta}\],\  dsdP\mbox{-}a.e.,
 \end{split}
\end{equation}
for all $\eta=(\eta_1,\eta_2,\cdots,\eta_d)\in L^{2}(\mathcal{G}_t;\mathbb{R}^d)$, where for all $y\in\mathbb{R}^d$,
$(O^{t,x,P_{\xi}}(y),Q^{t,x,P_{\xi}}(y))=$ \\
$\(((O_{s,j}^{t,x,P_{\xi}}(y))_{s\in[t,T]})_{1\leq j\leq d},((Q_{s,i,j}^{t,x,P_{\xi}}(y))_{s\in[t,T]})_{1\leq i,j\leq d}\)\in\mathcal{S}^2_{\mathcal{F}}(t,T;\mathbb{R}^{d})\times \mathcal{H}^2_{\mathcal{F}}(t,T;\mathbb{R}^{d\times d})$
is the unique solution of the following BDSDE:
\begin{equation}\label{eq6.4}
  \begin{split}
 & O_{s,j}^{t,x,P_{\xi}}(y)= \sum_{k=1}^d\partial_{x_k}\Phi (X_T^{t,x,P_{\xi}},P_{X_T^{t,\xi}})\partial_{\mu}X_{T,j}^{t,x,P_{\xi},k}(y)\\
&\ \ +\int_s^T(\nabla_{\Pi}f)(\Pi_r^{t,x,P_{\xi}},P_{\Pi_r^{t,\xi}})\Gamma_{r,j}^{t,x,P_{\xi}}(y) dr
    +\int_s^T(\nabla_{\Pi}g)(\Pi_r^{t,x,P_{\xi}},P_{\Pi_r^{t,\xi}})\Gamma_{r,j}^{t,x,P_{\xi}}(y) d\overleftarrow{B_r}\\
&\ \  +\sum_{k=1}^dE\[(\partial_{\mu}\Phi)_k(z,P_{X_T^{t,\xi}},X_r^{t,y,P_{\xi}}) \partial_{x_j}X_T^{t,y,P_{\xi},k}
    \!+\!(\partial_{\mu}\Phi)_k(z,P_{X_T^{t,\xi}},X_T^{t,\xi})\partial_{\mu}X_{T,j}^{t,\xi,k}(y) \]\Big|_{z=X_T^{t,x,P_{\xi}}}\\
&\ \ +\int_s^TE\[(\partial_{\mu}f)(z,P_{\Pi_r^{t,\xi}},\Pi_r^{t,y,P_{\xi}}) \partial_{x_j}\Pi_r^{t,y,P_{\xi}}
               +(\partial_{\mu}f)(z,P_{\Pi_r^{t,\xi}},\Pi_r^{t,\xi}) \Gamma_{r,j}^{t,\xi}(y)\]\Big|_{z=\Pi_r^{t,x,P_{\xi}}}dr\\
&\ \ +\int_s^TE\[\big((\partial_{\mu}g)(z,P_{\Pi_r^{t,\xi}},\Pi_r^{t,y,P_{\xi}})+(\partial_{\mu}h)(P_{\Pi_r^{t,\xi}},\Pi_r^{t,y,P_{\xi}})\big)
            \partial_{x_j}\Pi_r^{t,y,P_{\xi}}\]\Big|_{z=\Pi_r^{t,x,P_{\xi}}}d\overleftarrow{B_r}\\
&\ \ +\int_s^TE\[\big((\partial_{\mu}g)(z,P_{\Pi_r^{t,\xi}},\Pi_r^{t,\xi})+(\partial_{\mu}h)(P_{\Pi_r^{t,\xi}},\Pi_r^{t,\xi})\big)
          \Gamma_{r,j}^{t,\xi}(y)\]\Big|_{z=\Pi_r^{t,x,P_{\xi}}}d\overleftarrow{B_r}\\
&\ \ -\int_s^T\sum_{k=1}^dQ_{r,k,j}^{t,x,P_{\xi}}(y)dW_r^k,\ s\in[t,T],\ 1\leq j\leq d,
  \end{split}
\end{equation}
where $(O^{t,\xi},Q^{t,\xi})=(O^{t,\xi,P_{\xi}},Q^{t,\xi,P_{\xi}})$ is the unique solution of the above BDSDE
\eqref{eq6.4} with $x$ replaced by $\xi$,
$\Gamma_{r,j}^{t,x,P_{\xi}}(y)=(\partial_{\mu}X_{r,j}^{t,x,P_{\xi}}(y),O_{r,j}^{t,x,P_{\xi}}(y),Q_{r,j}^{t,x,P_{\xi}}(y))$,
and $\Gamma_{r,j}^{t,\xi}(y)=\Gamma_{r,j}^{t,\xi,P_{\xi}}(y)$.
\end{theorem}
In order to prove Theorem \ref{th6.2} we need the following three lemmas. For simplicity of redaction
but w.l.o.g., let us restrict to the dimensions $d = 1$, $l=1$ and to $f(x,y,z,\gamma)=f(z,\gamma(\mathbb{R}\times\mathbb{R}\times\cdot))$,
$g(x,y,z,\gamma)=g(z,\gamma(\mathbb{R}\times\mathbb{R}\times\cdot))$, $h(\gamma)=0$,
$(x,y,z)\in\mathbb{R}\times\mathbb{R}\times\mathbb{R}$, $\gamma\in\mathcal{P}_2(\mathbb{R}\times\mathbb{R}\times\mathbb{R})$,
and let $\Phi(x,\gamma)=\Phi(x)$, $(x,\gamma)\in\mathbb{R}\times\mathcal{P}_2(\mathbb{R})$.
We first consider the following BDSDE, which is obtained by formal differentiation
of the lifted solution $(Y^{t,x,\xi+q\eta},Z^{t,x,\xi+q\eta})$ of BDSDE \eqref{eq4.2} (with $\xi+q\eta$ instead of $\xi$,
where $\xi,\eta\in L^{2}(\mathcal{G}_t;\mathbb{R})$) with respect to $q$ at $q = 0$. This formal $L^2$-differentiation (which will be made
rigorous later) leads to a pair of processes $(\mathcal{O}^{t,x,\xi},\mathcal{Q}^{t,x,\xi})$ solving the BDSDE:
\begin{equation}\nonumber
  \begin{split}
  \mathcal{O}_{s}^{t,x,\xi}&(\eta)= \partial_{x}\Phi (X_T^{t,x,\xi})\mathcal{U}_{T}^{t,x,\xi}(\eta)
 +\int_s^T\partial_{z}f(Z_r^{t,x,\xi},P_{Z_r^{t,\xi}})\mathcal{Q}_{r}^{t,x,\xi}(\eta) dr\\
  \end{split}
\end{equation}
\begin{equation}\label{eq6.5}
  \begin{split}
&\  +\int_s^T\widehat{E}\[(\partial_{\mu}f)(Z_r^{t,x,\xi},P_{Z_r^{t,\xi}},\widehat{Z}_r^{t,\widehat{\xi},P_{\xi}})
 \partial_{x}\widehat{Z}_r^{t,\widehat{\xi},P_{\xi}}\widehat{\eta}
 +(\partial_{\mu}f)(Z_r^{t,x,\xi},P_{Z_r^{t,\xi}},\widehat{Z}_r^{t,\widehat{\xi}}) \widehat{\mathcal{Q}}_{r}^{t,\widehat{\xi}}(\widehat{\eta}) \]dr\\
&\ +\int_s^T\partial_{z}g(Z_r^{t,x,\xi},P_{Z_r^{t,\xi}})\mathcal{Q}_{r}^{t,x,\xi}(\eta) d\overleftarrow{B_r}\\
&\  +\int_s^T\widehat{E}\[(\partial_{\mu}g)(Z_r^{t,x,\xi},P_{Z_r^{t,\xi}},\widehat{Z}_r^{t,\widehat{\xi},P_{\xi}})
 \partial_{x}\widehat{Z}_r^{t,\widehat{\xi},P_{\xi}}\widehat{\eta}
 +(\partial_{\mu}g)(Z_r^{t,x,\xi},P_{Z_r^{t,\xi}},\widehat{Z}_r^{t,\widehat{\xi}}) \widehat{\mathcal{Q}}_{r}^{t,\widehat{\xi}}(\widehat{\eta}) \]d\overleftarrow{B_r}\\
&\  -\int_s^T\mathcal{Q}_{r}^{t,x,\xi}(\eta)dW_r,\ s\in[t,T],
  \end{split}
\end{equation}
where $(\mathcal{O}^{t,\xi}(\eta),\mathcal{Q}^{t,\xi}(\eta))=(\mathcal{O}^{t,x,\xi}(\eta),\mathcal{Q}^{t,x,\xi}(\eta))\big|_{x=\xi}$
is the solution of \eqref{eq6.5} for $x$ replaced by $\xi$, and
$\mathcal{U}_s^{t,x,\xi}(\eta):=DX_s^{t,x,\xi}(\eta)=\widetilde{E}[\partial_{\mu}X_s^{t,x,P_{\xi}}(\widetilde{\xi})\widetilde{\eta}]$
and  $\mathcal{U}_s^{t,\xi}(\eta)=\mathcal{U}_s^{t,x,\xi}(\eta)\big|_{x=\xi}$, $s\in[t,T]$.
Moreover, $(\widehat{\Omega},\widehat{\mathcal{F}},\widehat{P})$ is a probability space carrying with $(\widehat{\xi},\widehat{\eta},\widehat{B},\widehat{W})$
an (independent) copy of $(\xi,\eta,B,W)$ (defined on $(\Omega,\mathcal{F},P)$); $(\widehat{X}^{t,x,P_{\xi}},\widehat{Y}^{t,x,P_{\xi}},\widehat{Z}^{t,x,P_{\xi}})$
(resp., $(\widehat{X}^{t,\widehat{\xi}},\widehat{Y}^{t,\widehat{\xi}},\widehat{Z}^{t,\widehat{\xi}})$) is the solution of the same equation as that for
$(X^{t,x,P_{\xi}},Y^{t,x,P_{\xi}},Z^{t,x,P_{\xi}})$ (resp., $(X^{t,\xi},Y^{t,\xi},Z^{t,\xi})$), but with the data $(\widehat{\xi},\widehat{B},\widehat{W})$
instead of $(\xi,B,W)$.

Thanks to Theorem \ref{thA.1} BDSDE \eqref{eq6.5} (with $x$ replaced by $\xi$) has a unique solution
$(\mathcal{O}^{t,\xi}(\eta),$ $\mathcal{Q}^{t,\xi}(\eta))\in\mathcal{S}_{\mathcal{F}}^2(t,T;\mathbb{R})\times\mathcal{H}_{\mathcal{F}}^2(t,T;\mathbb{R})$.
Moreover, from Theorem \ref{thA.2} we deduce that there is a constant $C>0$ only depending on the Lipschitz constants of the coefficients, such that
\begin{equation}\label{eq6.6}
  E\[\sup_{s\in[t,T]}|\mathcal{O}_s^{t,\xi}(\eta)|^{2}+\int_t^T|\mathcal{Q}_s^{t,\xi}(\eta)|^2ds\]\leq C.
\end{equation}
Once having $(\mathcal{O}^{t,\xi}(\eta),\mathcal{Q}^{t,\xi}(\eta))$, it follows again from the Theorems \ref{thA.1} and \ref{thA.2} that \eqref{eq6.5}
possesses a unique solution
$(\mathcal{O}^{t,x,\xi}(\eta),\mathcal{Q}^{t,x,\xi}(\eta))\in\mathcal{S}_{\mathcal{F}}^2(t,T;\mathbb{R})\times\mathcal{H}_{\mathcal{F}}^2(t,T;\mathbb{R})$,
 and there is a constant $C>0$ only depending on the Lipschitz constants of the coefficients, such that
\begin{equation}\label{eq6.7}
  E\[\sup_{s\in[t,T]}|\mathcal{O}_s^{t,x,\xi}(\eta)|^2+\int_t^T|\mathcal{Q}_s^{t,x,\xi}(\eta)|^2ds\]\leq C.
\end{equation}

\begin{lemma} \label{le6.1}
Suppose (H5.1) and (H6.1) hold true. Then, for all $(t,x)\in[0,T]\times\mathbb{R}$, $\xi\in L^{2}(\mathcal{G}_t;\mathbb{R})$,
there exist two stochastic processes $O^{t,x,P_{\xi}}(y)\in\mathcal{S}_{\mathcal{F}}^2(t,T;\mathbb{R})$,
$Q^{t,x,P_{\xi}}(y)\in\mathcal{H}_{\mathcal{F}}^2(t,T;\mathbb{R})$, depending measurably on $y\in\mathbb{R}$, such that
\begin{equation*}
 \begin{split}
   & \mathcal{O}_s^{t,x,\xi}(\eta)=\overline{E}[O_s^{t,x,P_{\xi}}(\overline{\xi})\cdot\overline{\eta}],\ s\in[t,T],\ P\mbox{-}a.s.,\ \ \mbox{and} \ \
  \mathcal{Q}_s^{t,x,\xi}(\eta)=\overline{E}[Q_s^{t,x,P_{\xi}}(\overline{\xi})\cdot\overline{\eta}],\  dsdP\mbox{-}a.e.,
 \end{split}
\end{equation*}
for all $\eta\in L^{2}(\mathcal{G}_t;\mathbb{R})$. In particular, for all $x\in\mathbb{R}$, $0\leq t\leq s\leq T$, $\xi\in L^{2}(\mathcal{G}_t;\mathbb{R})$, the mappings
\begin{equation*}
 \begin{split}
   & \mathcal{O}_s^{t,x,\xi}(\cdot):L^{2}(\mathcal{G}_t;\mathbb{R})\rightarrow L^{2}(\mathcal{F}_s;\mathbb{R}),\ \ \mbox{and} \ \
    \mathcal{Q}^{t,x,\xi}(\cdot):L^{2}(\mathcal{G}_t;\mathbb{R})\rightarrow\mathcal{H}_{\mathcal{F}}^2(t,T;\mathbb{R}),
 \end{split}
\end{equation*}
are linear and continuous.
\end{lemma}

\begin{remark} \label{re6.1}
As $(\mathcal{O}_s^{t,\xi}(y),\mathcal{Q}_s^{t,\xi}(y))=(\mathcal{O}_s^{t,x,\xi}(y),\mathcal{Q}_s^{t,x,\xi}(y))\big|_{x=\xi}$,
$s\in[t,T]$, $\xi\in L^{2}(\mathcal{G}_t;\mathbb{R})$, $y\in\mathbb{R}$, we see directly from Lemma \ref{le6.1} that
\begin{equation*}
 \begin{split}
   & \mathcal{O}_s^{t,\xi}(\eta)=\overline{E}[O_s^{t,\xi}(\overline{\xi})\cdot\overline{\eta}],\ s\in[t,T],\ P\mbox{-}a.s.,\ \ \mbox{and} \ \
  \mathcal{Q}_s^{t,\xi}(\eta)=\overline{E}[Q_s^{t,\xi}(\overline{\xi})\cdot\overline{\eta}],\  dsdP\mbox{-}a.e.
 \end{split}
\end{equation*}
\end{remark}

\begin{proof}[Proof (of Lemma \ref{le6.1}).]
For $y\in\mathbb{R}$, let
$(O^{t,x,P_{\xi}}(y),Q^{t,x,P_{\xi}}(y))\in\mathcal{S}_{\mathcal{F}}^2(t,T;\mathbb{R})\times\mathcal{H}_{\mathcal{F}}^2(t,T;\mathbb{R})$
be the unique solution of BDSDE \eqref{eq6.4}, which, for our special case ($d=1$, $l=1$ and $f=f(z,\gamma(\mathbb{R}\times\mathbb{R}\times\cdot))$,
$g=g(z,\gamma(\mathbb{R}\times\mathbb{R}\times\cdot))$, $h=0$) writes as follows
\begin{equation}\label{eq6.8}
  \begin{split}
  O_{s}^{t,x,P_{\xi}}&(y)= \partial_{x}\Phi (X_T^{t,x,P_{\xi}})\partial_{\mu}X_{T}^{t,x,P_{\xi}}(y)
  +\int_s^T \partial_{z}f(Z_r^{t,x,\xi},P_{Z_r^{t,\xi}})Q_{r}^{t,x,P_{\xi}}(y) dr\\
&\  +\int_s^T\widehat{E}\[(\partial_{\mu}f)(Z_r^{t,x,\xi},P_{Z_r^{t,\xi}},\widehat{Z}_r^{t,y,P_{\xi}}) \partial_{x}\widehat{Z}_r^{t,y,P_{\xi}}
 +(\partial_{\mu}f)(Z_r^{t,x,\xi},P_{Z_r^{t,\xi}},\widehat{Z}_r^{t,\widehat{\xi}}) \widehat{Q}_{r}^{t,\widehat{\xi}}(y) \]dr\\
&\  +\int_s^T \partial_{z}g(Z_r^{t,x,\xi},P_{Z_r^{t,\xi}})Q_{r}^{t,x,P_{\xi}}(y) d\overleftarrow{B_r}\\
&\  +\int_s^T\widehat{E}\[(\partial_{\mu}g)(Z_r^{t,x,\xi},P_{Z_r^{t,\xi}},\widehat{Z}_r^{t,y,P_{\xi}}) \partial_{x}\widehat{Z}_r^{t,y,P_{\xi}}
 +(\partial_{\mu}g)(Z_r^{t,x,\xi},P_{Z_r^{t,\xi}},\widehat{Z}_r^{t,\widehat{\xi}}) \widehat{Q}_{r}^{t,\widehat{\xi}}(y) \]d\overleftarrow{B_r}\\
&\  -\int_s^TQ_{r}^{t,x,P_{\xi}}(y)dW_r,\ s\in[t,T],
  \end{split}
\end{equation}
where $O^{t,\xi}(y),Q^{t,\xi}(y)):=(O^{t,x,P_{\xi}}(y),Q^{t,x,P_{\xi}}(y))\big|_{x=\xi}
\in\mathcal{S}_{\mathcal{F}}^2(t,T;\mathbb{R})\times\mathcal{H}_{\mathcal{F}}^2(t,T;\mathbb{R})$
is the unique solution of \eqref{eq6.8} with $x$ replaced by $\xi$. It follows from Theorem \ref{thA.2} that,
there is a constant $C>0$ only depending on the Lipschitz constants of the coefficients, such that,
for all $t\in[0,T]$, $x\in\mathbb{R}$, $\xi\in L^{2}(\mathcal{G}_t;\mathbb{R})$, $y\in\mathbb{R}$,
\begin{equation}\label{eq6.9}
  E\[\sup_{s\in[t,T]}|O_s^{t,\xi}(y)|^2+\int_t^T|Q_s^{t,\xi}(y)|^2ds\]\leq C.
\end{equation}
Then again from Theorem \ref{thA.2} we obtain
\begin{equation}\label{eq6.10}
  E\[\sup_{s\in[t,T]}|O_s^{t,x,P_{\xi}}(y)|^2+\int_t^T|Q_s^{t,x,P_{\xi}}(y)|^2ds\]\leq C.
\end{equation}
Let the couple $(\overline{\xi},\overline{\eta})$ defined on some probability space $(\overline{\Omega},\overline{\mathcal{F}},\overline{P})$
be an independent copy of $(\xi,\eta)$ on $(\Omega,\mathcal{F},P)$ and, in particular, also an independent copy of
$(\widehat{\xi},\widehat{\eta})$ on $(\widehat{\Omega},\widehat{\mathcal{F}},\widehat{P})$.
After substituting the random variable $\xi$ and $\overline{\xi}$ for $x$ and $y$ in \eqref{eq6.8}, respectively, and then multiplying by
$\overline{\eta}$ both sides of the equation, taking expectation $\overline{E}[\cdot]$ yields
\begin{equation}\label{eq6.11}
  \begin{split}
  \overline{E}&[O_{s}^{t,\xi}(\overline{\xi})\cdot\overline{\eta}]
  = \partial_{x}\Phi (X_T^{t,\xi})\overline{E}[\partial_{\mu}X_{T}^{t,\xi}(\overline{\xi})\cdot\overline{\eta}]
  +\overline{E}\[\int_s^T \partial_{z}f(Z_r^{t,\xi},P_{Z_r^{t,\xi}})Q_{r}^{t,\xi}(\overline{\xi})\cdot\overline{\eta} dr\]\\
&  +\overline{E}\[\int_s^T\widehat{E}\[(\partial_{\mu}f)(Z_r^{t,\xi},P_{Z_r^{t,\xi}},\widehat{Z}_r^{t,\overline{\xi},P_{\xi}})
      \partial_{x}\widehat{Z}_r^{t,\overline{\xi},P_{\xi}}\cdot\overline{\eta}
 \!+\!(\partial_{\mu}f)(Z_r^{t,\xi},P_{Z_r^{t,\xi}},\widehat{Z}_r^{t,\widehat{\xi}})
      \widehat{Q}_{r}^{t,\widehat{\xi}}(\overline{\xi})\cdot\overline{\eta} \]dr\]\\
& +\overline{E}\[\int_s^T \partial_{z}g(Z_r^{t,\xi},P_{Z_r^{t,\xi}})Q_{r}^{t,\xi}(\overline{\xi})\cdot\overline{\eta} d\overleftarrow{B_r}\]\\
&  +\overline{E}\[\int_s^T\widehat{E}\[(\partial_{\mu}g)(Z_r^{t,\xi},P_{Z_r^{t,\xi}},\widehat{Z}_r^{t,\overline{\xi},P_{\xi}})
      \partial_{x}\widehat{Z}_r^{t,\overline{\xi},P_{\xi}}\cdot\overline{\eta}
 \!+\!(\partial_{\mu}g)(Z_r^{t,\xi},P_{Z_r^{t,\xi}},\widehat{Z}_r^{t,\widehat{\xi}})
      \widehat{Q}_{r}^{t,\widehat{\xi}}(\overline{\xi})\cdot\overline{\eta} \]d\overleftarrow{B_r}\]\\
&  -\overline{E}[\int_s^TQ_{r}^{t,\xi}(\overline{\xi})\cdot\overline{\eta}dW_r],\ s\in[t,T].
  \end{split}
\end{equation}
Since $(\overline{\xi},\overline{\eta})$ is independent of $(\xi,\eta,Z^{t,x,\xi})$ and
$(\widehat{\xi},\widehat{\eta},\widehat{Z}^{t,x,P_{\xi}})$, and of the same law as $(\widehat{\xi},\widehat{\eta})$,
we have
\begin{equation*}
\overline{E}[\widehat{E}[(\partial_{\mu}f)(Z_r^{t,\xi},P_{Z_r^{t,\xi}},\widehat{Z}_r^{t,\overline{\xi},P_{\xi}})
      \partial_{x}\widehat{Z}_r^{t,\overline{\xi},P_{\xi}}\cdot\overline{\eta}]]
 =  \widehat{E}[(\partial_{\mu}f)(Z_r^{t,\xi},P_{Z_r^{t,\xi}},\widehat{Z}_r^{t,\widehat{\xi},P_{\xi}})
      \partial_{x}\widehat{Z}_r^{t,\widehat{\xi},P_{\xi}}\cdot\widehat{\eta}],
\end{equation*}
and similar for the other terms. From the above equalities and the uniqueness of the solution of Eq. \eqref{eq6.5}
with $x$ replaced by $\xi$ it follows
\begin{equation}\label{eq6.12}
 \begin{split}
   &  \mathcal{O}_s^{t,\xi}(\eta)=\overline{E}[O_s^{t,\xi}(\overline{\xi})\cdot\overline{\eta}],\ s\in[t,T],\ P\mbox{-}a.s.,\\
  & \mathcal{Q}_s^{t,\xi}(\eta)=\overline{E}[Q_s^{t,\xi}(\overline{\xi})\cdot\overline{\eta}],\  dsdP\mbox{-}a.e.
 \end{split}
\end{equation}
Furthermore, from \eqref{eq6.9} we get
\begin{equation}\label{eq6.13}
 \begin{split}
   &  E[|\mathcal{O}_s^{t,\xi}(\eta)|^2]=E[|\overline{E}[O_s^{t,\xi}(\overline{\xi})\cdot\overline{\eta}]|^2]\leq
\overline{E}[E[|O_s^{t,\xi}(\overline{\xi})|^2\cdot|\overline{\eta}|^2]]\\
  & =\overline{E}[E[|O_s^{t,\xi}(y)|^2]\big|_{y=\overline{\xi}}\cdot|\overline{\eta}|^2]\leq
  C\overline{E}[|\overline{\eta}|^2]=CE[|\eta|^2].
 \end{split}
\end{equation}
That means $|\mathcal{O}_s^{t,\xi}(\eta)|_{L^2}\leq C|\eta|_{L^2}$, for every $\eta\in L^{2}(\mathcal{G}_t;\mathbb{R})$. Hence,
$\mathcal{O}_s^{t,\xi}(\cdot):L^{2}(\mathcal{G}_t;\mathbb{R})\rightarrow L^{2}(\mathcal{F}_s;\mathbb{R})$
 is a linear and continuous mapping, for all $s\in[t,T]$, and $|\mathcal{O}_s^{t,\xi}(\cdot)|_{L(L^2,L^2)}\leq C$.

Moreover, we also get
\begin{equation}\label{eq6.14}
 \begin{split}
   &  E[\int_t^T|\mathcal{Q}_s^{t,\xi}(\eta)|^2ds]=E[\int_t^T|\overline{E}[Q_s^{t,\xi}(\overline{\xi})\cdot\overline{\eta}]|^2ds]\leq
\overline{E}[E[\int_t^T|Q_s^{t,\xi}(\overline{\xi})|^2ds\cdot|\overline{\eta}|^2]]\\
  & =\overline{E}[E[\int_t^T|Q_s^{t,\xi}(y)|^2ds]\big|_{y=\overline{\xi}}\cdot|\overline{\eta}|^2]\leq
  CE[|\eta|^2].
 \end{split}
\end{equation}
Therefore, also $\mathcal{Q}_s^{t,\xi}(\cdot):L^{2}(\mathcal{G}_t;\mathbb{R})\rightarrow\mathcal{H}_{\mathcal{F}}^2(t,T;\mathbb{R})$
is a continuous linear mapping. Making use of the above argument, but now for $(\mathcal{O}^{t,x,\xi}(\eta),\mathcal{Q}^{t,x,\xi}(\eta))$
instead of $(\mathcal{O}^{t,\xi}(\eta),\mathcal{Q}^{t,\xi}(\eta))$, we have, for all $\eta\in L^{2}(\mathcal{G}_t;\mathbb{R})$,
\begin{equation}\label{eq6.15}
 \begin{split}
   & \mathcal{O}_s^{t,x,\xi}(\eta)=\overline{E}[O_s^{t,x,P_{\xi}}(\overline{\xi})\cdot\overline{\eta}],\ s\in[t,T],\ P\mbox{-}a.s.,\\
  & \mathcal{Q}_s^{t,x,\xi}(\eta)=\overline{E}[Q_s^{t,x,P_{\xi}}(\overline{\xi})\cdot\overline{\eta}],\  dsdP\mbox{-}a.e.
 \end{split}
\end{equation}
Finally, in analogy to \eqref{eq6.13} and \eqref{eq6.14}, we deduce from \eqref{eq6.10} that also
$(\mathcal{O}^{t,x,\xi}(\eta),\mathcal{Q}^{t,x,\xi}(\eta))$ is a linear and continuous mapping
(over the same spaces as $(\mathcal{O}^{t,\xi}(\eta),\mathcal{Q}^{t,\xi}(\eta))$).
\end{proof}

Notice that, it follows from Theorem \ref{thA.3} that, under the Assumptions (H4.2), (H5.1) and (H6.1),
there is a constant $C_{p_0}>0$ only depending on the Lipschitz constants of the coefficients, such that,
for all $t\in[0,T]$, $x\in\mathbb{R}^d$, $\xi\in L^{2}(\mathcal{G}_t;\mathbb{R}^d)$, $y\in\mathbb{R}^d$,
$\displaystyle E\[\sup_{s\in[t,T]}|O_s^{t,\xi}(y)|^{p_0}+(\int_t^T|Q_s^{t,\xi}(y)|^2ds)^{\frac{p_0}{2}}\]\leq C_{p_0}.$
From H\"{o}lder's inequality, for all $p\in[2,p_0]$ we have
\begin{equation}\label{eq6.37+1}
  E\[\sup_{s\in[t,T]}|O_s^{t,\xi}(y)|^p+(\int_t^T|Q_s^{t,\xi}(y)|^2ds)^{\frac{p}{2}}\]\leq C_p.
\end{equation}
Then again from Theorem \ref{thA.3} we obtain
$\displaystyle E\[\sup_{s\in[t,T]}|O_s^{t,x,P_{\xi}}(y)|^{p_0}+(\int_t^T|Q_s^{t,x,P_{\xi}}(y)|^2ds)^{\frac{p_0}{2}}\]\leq C_{p_0}.$
Again from H\"{o}lder's inequality, for all $p\in[2,p_0]$ we have
\begin{equation}\label{eq6.38+1}
  E\[\sup_{s\in[t,T]}|O_s^{t,x,P_{\xi}}(y)|^p+(\int_t^T|Q_s^{t,x,P_{\xi}}(y)|^2ds)^{\frac{p}{2}}\]\leq C_p.
\end{equation}
Now we prove the following estimate for the solution of Eq. \eqref{eq6.8}.

\begin{proposition} \label{prop6.2}
For all $p\in[2,p_0]$, under the Assumptions (H4.2), (H5.1) and (H6.1),
there is a constant $C_p>0$ only depending on the Lipschitz constants of the coefficients,
such that, for all $t\in[0,T]$, $x,x',y,y'\in\mathbb{R}^d$, $\xi,\xi'\in L^{2}(\mathcal{G}_t;\mathbb{R}^d)$,
$P$-a.s.,
\begin{equation}\label{eq6.16}
  \begin{split}
  &\  E\[\sup_{s\in[t,T]}|O_s^{t,x,P_{\xi}}(y)-O_s^{t,x',P_{\xi'}}(y')|^p+(\int_t^T|Q_s^{t,x,P_{\xi}}(y)-Q_s^{t,x',P_{\xi'}}(y')|^2ds)^{\frac{p}{2}}\]\\
 \leq&\ C_pM^p\(|x-x'|^p+|y-y'|^p+W_2(P_{\xi},P_{\xi'})^p\)\\
  &\ +\rho_{M,p}(t,x,y,P_{\xi})+\rho_{M,p}(t,y,P_{\xi})+E[\rho_{M,p}(t,\xi,y,P_{\xi})],
  \end{split}
\end{equation}
for $M\geq1$, $\rho_{M,p}(t,x,y,P_{\xi})\rightarrow0$, $\rho_{M,p}(t,y,P_{\xi})\rightarrow0$,
$E[\rho_{M,p}(t,\xi,y,P_{\xi})]\rightarrow0$, as $M\rightarrow\infty$.
\end{proposition}

\begin{proof}
For simplicity of redaction we restrict ourselves again to the dimensions $l=1$, $d=1$,
let $f\equiv0$, $g(\Pi_r^{t,x,P_{\xi}},P_{\Pi_r^{t,\xi}})=g(Z_r^{t,x,P_{\xi}},P_{Z_r^{t,\xi}})$,
$h(P_{\Pi_r^{t,\xi}})=h(P_{(Y_r^{t,\xi},Z_r^{t,\xi})})$ and $\Phi=\Phi(x)$.

Then, from \eqref{eq6.4} we get the following BDSDE:
\begin{equation}\label{eq6.16+1}
  \begin{split}
  &O_{s}^{t,x,P_{\xi}}(y)= \partial_{x}\Phi (X_T^{t,x,P_{\xi}})\partial_{\mu}X_{T}^{t,x,P_{\xi}}(y)
  +\int_s^T \Bigg\{\partial_{z}g(Z_r^{t,x,P_{\xi}},P_{Z_r^{t,\xi}})Q_{r}^{t,x,P_{\xi}}(y)\\
&\  +\widehat{E}\[(\partial_{\mu}h)_{1}(P_{(Y_r^{t,\xi},Z_r^{t,\xi})},\!(\widehat{Y}_r^{t,y,P_{\xi}},\widehat{Z}_r^{t,y,P_{\xi}}))\partial_{x}\widehat{Y}_r^{t,y,P_{\xi}}
   \! +\!(\partial_{\mu}h)_1(P_{(Y_r^{t,\xi},Z_r^{t,\xi})},\!(\widehat{Y}_r^{t,\widehat{\xi}},\widehat{Z}_r^{t,\widehat{\xi}})) \widehat{O}_{r}^{t,\widehat{\xi}}(y) \] \\
&\  +\widehat{E}\[\((\partial_{\mu}g)(Z_r^{t,x,P_{\xi}},P_{Z_r^{t,\xi}},\widehat{Z}_r^{t,y,P_{\xi}})
     +(\partial_{\mu}h)_{2}(P_{(Y_r^{t,\xi},Z_r^{t,\xi})},(\widehat{Y}_r^{t,y,P_{\xi}},\widehat{Z}_r^{t,y,P_{\xi}}))\)\cdot\partial_{x}\widehat{Z}_r^{t,y,P_{\xi}}\]   \\
&\  +\widehat{E}\[\((\partial_{\mu}g)(Z_r^{t,x,P_{\xi}},P_{Z_r^{t,\xi}},\widehat{Z}_r^{t,\widehat{\xi}})
  +(\partial_{\mu}h)_2(P_{(Y_r^{t,\xi},Z_r^{t,\xi})},(\widehat{Y}_r^{t,\widehat{\xi}},\widehat{Z}_r^{t,\widehat{\xi}}))\)\cdot \widehat{Q}_{r}^{t,\widehat{\xi}}(y) \]
  \Bigg\}d\overleftarrow{B_r}\\
&\  -\int_s^TQ_{r}^{t,x,P_{\xi}}(y)dW_r,\ s\in[t,T],
  \end{split}
\end{equation}
where $(O_s^{t,\xi}(y),Q_s^{t,\xi}(y))=(O_s^{t,x,\xi}(y),Q_s^{t,x,\xi}(y))\big|_{x=\xi}$.

Let us introduce the $\{\mathcal{F}_{t,s}\}$-adapted process
\begin{equation}\label{eq6.16+2}
  \begin{split}
  &R_{s}^{t,x,P_{\xi}}(y):= O_{s}^{t,x,P_{\xi}}(y)
  -\int_s^T\Bigg\{\widehat{E}\[(\partial_{\mu}h)_{1}(P_{(Y_r^{t,\xi},Z_r^{t,\xi})},(\widehat{Y}_r^{t,y,P_{\xi}},\widehat{Z}_r^{t,y,P_{\xi}}))
          \partial_{x}\widehat{Y}_r^{t,y,P_{\xi}}\]\\
&\ +\!\widehat{E}\[\!(\partial_{\mu}h)_1(P_{(Y_r^{t,\xi}\!,Z_r^{t,\xi})},\!(\widehat{Y}_r^{t,\widehat{\xi}}\!,\widehat{Z}_r^{t,\widehat{\xi}}))
       \widehat{O}_{r}^{t,\widehat{\xi}}(y) \!\]
\!+\!\widehat{E}\[\!(\partial_{\mu}h)_{2}(P_{(Y_r^{t,\xi}\!,Z_r^{t,\xi})},\!(\widehat{Y}_r^{t,y,P_{\xi}}\!,\widehat{Z}_r^{t,y,P_{\xi}}))\partial_{x}\widehat{Z}_r^{t,y,P_{\xi}}\!\]\\
&\   +\!\widehat{E}\[(\partial_{\mu}h)_2(P_{(Y_r^{t,\xi},Z_r^{t,\xi})},(\widehat{Y}_r^{t,\widehat{\xi}},\widehat{Z}_r^{t,\widehat{\xi}})) \widehat{Q}_{r}^{t,\widehat{\xi}}(y) \]
  \Bigg\}d\overleftarrow{B_r},\ s\in[t,T].
  \end{split}
\end{equation}
We note that the integrand of the above stochastic integral w.r.t. $d\overleftarrow{B}$ is deterministic.
Furthermore,
\begin{equation}\label{eq6.16+3}
  \begin{split}
  R_{s}^{t,x,P_{\xi}}(y)&-R_{s}^{t,x',P_{\xi'}}(y')
  = \Xi(x,x',y,y',\xi,\xi')+\int_s^TH(r,x,x',y,y',\xi,\xi')d\overleftarrow{B_r}\\
 &\ +\int_s^T\Bigg\{ (\partial_{z}g)(Z_r^{t,x,P_{\xi}},P_{Z_r^{t,\xi}})(Q_{r}^{t,x,P_{\xi}}(y)-Q_{r}^{t,x',P_{\xi'}}(y') )\\
&\ \hspace{40pt} +\widehat{E}\[(\partial_{\mu}g)(Z_r^{t,x,P_{\xi}},P_{Z_r^{t,\xi}},\widehat{Z}_r^{t,\widehat{\vartheta},P_{\xi}})
     (\widehat{Q}_{r}^{t,\widehat{\vartheta},P_{\xi}}(y)-\widehat{Q}_{r}^{t,\widehat{\vartheta'},P_{\xi'}}(y') )\]\Bigg\}d\overleftarrow{B_r}\\
&\  -\int_s^T(Q_{r}^{t,x,P_{\xi}}(y)-Q_{r}^{t,x',P_{\xi'}}(y'))dW_r,\ s\in[t,T]
  \end{split}
\end{equation}
where $\vartheta,\vartheta'\in L^{2}(\mathcal{G}_t;\mathbb{R})$ with $P_{\vartheta}=P_{\xi}$, $P_{\vartheta'}=P_{\xi'}$, and
\begin{equation*}
\begin{split}
&\ \mbox{(i)}\ \Xi(x,x',y,y',\xi,\xi'):=(\partial_{x}\Phi) (X_T^{t,x,P_{\xi}})\partial_{\mu}X_{T}^{t,x,P_{\xi}}(y)-(\partial_{x}\Phi)
                       (X_T^{t,x',P_{\xi'}})\partial_{\mu}X_{T}^{t,x',P_{\xi'}}(y'),\\
 &\ \mbox{(ii)}\ H(r,x,x',y,y',\xi,\xi'): =
    \big((\partial_{z}g)(Z_r^{t,x,P_{\xi}},P_{Z_r^{t,\xi}})-(\partial_{z}g)(Z_r^{t,x',P_{\xi'}},P_{Z_r^{t,\xi'}})\big)Q_{r}^{t,x',P_{\xi'}}(y')\\
 &\qquad\quad +\widehat{E}\[\big((\partial_{\mu}g)(Z_r^{t,x,P_{\xi}},P_{Z_r^{t,\xi}},\widehat{Z}_r^{t,\widehat{\vartheta},P_{\xi}})
 -(\partial_{\mu}g)(Z_r^{t,x',P_{\xi'}},P_{Z_r^{t,\xi'}},\widehat{Z}_r^{t,\widehat{\vartheta'},P_{\xi'}})\big)\widehat{Q}_{r}^{t,\widehat{\vartheta'},P_{\xi'}}(y')\]\\
  &\qquad\quad +\widehat{E}\[(\partial_{\mu}g)(Z_r^{t,x,P_{\xi}},P_{Z_r^{t,\xi}},\widehat{Z}_r^{t,y,P_{\xi}}) \partial_{x}\widehat{Z}_r^{t,y,P_{\xi}}
  -(\partial_{\mu}g)(Z_r^{t,x',P_{\xi'}},P_{Z_r^{t,\xi'}},\widehat{Z}_r^{t,y',P_{\xi'}}) \partial_{x}\widehat{Z}_r^{t,y',P_{\xi'}}\].
\end{split}
\end{equation*}
From the Propositions \ref{prop5.1} and \ref{prop5.3}, for all $p\geq 2$, it holds
\begin{equation}\label{eq6.16+4}
 E[| \Xi(x,x',y,y',\xi,\xi')|^p]\leq  C_p\(|x-x'|^p+|y-y'|^p+W_2(P_{\xi},P_{\xi'})^p\).
\end{equation}
We now give the estimate for $H(r,x,x',y,y',\xi,\xi')$. From the Propositions \ref{prop4.1} and \ref{prop4.2} as well as their proofs, we have
\begin{equation*}
  \begin{split}
\mbox{(i)}&\  E\[\(\int_t^T|(\partial_{z}g)(Z_r^{t,x,P_{\xi}},P_{Z_r^{t,\xi}})
     -(\partial_{z}g)(Z_r^{t,x',P_{\xi'}},P_{Z_r^{t,\xi'}})|^2\cdot|Q_{r}^{t,x',P_{\xi'}}(y')|^2dr\)^{\frac{p_0}{2}}\]\\
  &\ \leq C_{p_0} E\[\(\int_t^T\min\big\{C,|Z_r^{t,x,P_{\xi}}-Z_r^{t,x',P_{\xi'}}|^2+W_2(P_{Z_r^{t,\xi}},P_{Z_r^{t,\xi'}})^2\big\}
  \cdot|Q_{r}^{t,x',P_{\xi'}}(y')|^2dr\)^{\frac{p_0}{2}}\]  \\
  &\ \leq C_{p_0}M^{p_0}E\[(\int_t^T(|Z_r^{t,x,P_{\xi}}-Z_r^{t,x',P_{\xi'}}|^2+W_2(P_{Z_r^{t,\xi}},P_{Z_r^{t,\xi'}})^2)dr)^{\frac{p_0}{2}}\]\\
  &\qquad +C_{p_0}E\[\(\int_t^T|Q_{r}^{t,x',P_{\xi'}}(y')|^2I_{\{|Q_{r}^{t,x',P_{\xi'}}(y')|\geq M\}}dr\)^{\frac{p_0}{2}}\]\\
  &\ \leq C_{p_0}M^{p_0}\(|x-x'|^{p_0}+W_2(P_{\xi},P_{\xi'})^{p_0}\)+\rho_{M,p_0}(t,x',y',P_{\xi'}),
\end{split}
\end{equation*}
with $\displaystyle \rho_{M,p_0}(t,x',y',P_{\xi'}):=C_{p_0}E[(\int_t^T|Q_{r}^{t,x',P_{\xi'}}(y')|^2I_{\{|Q_{r}^{t,x',P_{\xi'}}(y')|\geq M\}}dr)^{\frac{p_0}{2}}]$.

\noindent In the same way we see
\begin{equation*}
  \begin{split}
\mbox{(ii)}&\  E\[\(\int_t^T\Big|\widehat{E}\[\big((\partial_{\mu}g)(Z_r^{t,x,P_{\xi}},P_{Z_r^{t,\xi}},\widehat{Z}_r^{t,\widehat{\vartheta},P_{\xi}})
 \!-\!(\partial_{\mu}g)(Z_r^{t,x',P_{\xi'}},P_{Z_r^{t,\xi'}},\widehat{Z}_r^{t,\widehat{\vartheta'},P_{\xi'}})\big)
 \widehat{Q}_{r}^{t,\widehat{\vartheta'},P_{\xi'}}(y')\]\Big|^2dr\)^{\frac{p_0}{2}}\]\\
  &\  \leq C_{p_0}M^{p_0}E\[(\int_t^T(|Z_r^{t,x,P_{\xi}}-Z_r^{t,x',P_{\xi'}}|^2+W_2(P_{Z_r^{t,\xi}},P_{Z_r^{t,\xi'}})^2
         +\widehat{E}[|\widehat{Z}_r^{t,\widehat{\vartheta},P_{\xi}}-\widehat{Z}_r^{t,\widehat{\vartheta'},P_{\xi'}}|^2])dr)^{\frac{p_0}{2}}\]\\
  &\qquad +C_{p_0}\widehat{E}\[\(\int_t^T|\widehat{Q}_{r}^{t,\widehat{\vartheta'},P_{\xi'}}(y')|^2
     I_{\{\widehat{Q}_{r}^{t,\widehat{\vartheta'},P_{\xi'}}(y')\geq M\}}dr\)^{\frac{p_0}{2}}\]\\
  &\  \leq C_{p_0}M^{p_0}\(|x-x'|^{p_0}+W_2(P_{\xi},P_{\xi'})^{p_0}+(E[|\vartheta-\vartheta'|^2])^{\frac{p_0}{2}}\)+E[\rho_{M,p_0}(t,\xi',y',P_{\xi'})],
\end{split}
\end{equation*}
for all $t\in[0,T]$, $x,x',y'\in\mathbb{R}$, $\xi,\xi',\vartheta,\vartheta'\in L^{2}(\mathcal{G}_t;\mathbb{R})$ with
$P_{\vartheta}=P_{\xi}$, $P_{\vartheta'}=P_{\xi'}$, and
\begin{equation*}
  \begin{split}
\mbox{(iii)}&\  E\[\!\(\!\!\int_t^T\!\Big|\widehat{E}\[\!(\partial_{\mu}g)(Z_r^{t,x,P_{\xi}}\!,\!P_{Z_r^{t,\xi}},\!\widehat{Z}_r^{t,y,P_{\xi}})
         \partial_{x}\widehat{Z}_r^{t,y,P_{\xi}}
  \!\!-\!\!(\partial_{\mu}g)(Z_r^{t,x',P_{\xi'}}\!,\!P_{Z_r^{t,\xi'}},\!\widehat{Z}_r^{t,y',P_{\xi'}})
      \partial_{x}\widehat{Z}_r^{t,y',P_{\xi'}}\]\Big|^2\!dr\!\)^{\frac{p_0}{2}}\!\]\\
  &\  \leq C_{p_0}  \widehat{E}[(\int_t^T|\partial_{x}\widehat{Z}_r^{t,y,P_{\xi}}- \partial_{x}\widehat{Z}_r^{t,y',P_{\xi'}}|^2dr)^{\frac{p_0}{2}}]
  + C_{p_0}\widehat{E}[(\int_t^T|\partial_{x}\widehat{Z}_r^{t,y',P_{\xi'}}|^2I_{\{|\partial_{x}\widehat{Z}_r^{t,y',P_{\xi'}}|\geq M\}}dr)^{\frac{p_0}{2}}]\\
 &\qquad +C_{p_0}M^{p_0}E[(\int_t^T(|Z_r^{t,x,P_{\xi}}-Z_r^{t,x',P_{\xi'}}|^2+W_2(P_{Z_r^{t,\xi}},P_{Z_r^{t,\xi'}})^2
     +\widehat{E}[|\widehat{Z}_r^{t,y,P_{\xi}}-\widehat{Z}_r^{t,y',P_{\xi'}}|^2])dr)^{\frac{p_0}{2}}]\\
 &\  \leq C_{p_0}M^{p_0}\(|x-x'|^{p_0}+|y-y'|^{p_0}+W_2(P_{\xi},P_{\xi'})^{p_0}\)+\rho_{M,p_0}(t,y',P_{\xi'}).
  \end{split}
\end{equation*}
Consequently,
\begin{equation*}
\begin{split}
 E[(\!\int_t^T\!\!|H(r,x,x',y,y',\xi,\xi')|^2dr)^{\frac{p_0}{2}}]
\leq&\ C_{p_0}M^{p_0}\(\!|x\!-\!x'|^{p_0}\!+\!|y\!-\!y'|^{p_0}\!+\!W_2(P_{\xi},P_{\xi'})^{p_0}\!+\!(E[|\vartheta\!-\!\vartheta'|^2])^{\frac{p_0}{2}}\!\)\\
&\ +\!\rho_{M,p_0}(t,x',y',P_{\xi'})\!+\!\rho_{M,p_0}(t,y',P_{\xi'})\!+\!E[\rho_{M,p_0}(t,\xi',y',P_{\xi'})],
\end{split}
\end{equation*}
where $\rho_{M,p_0}(t,x',y',P_{\xi'})\rightarrow0$, $\rho_{M,p_0}(t,y',P_{\xi'})\rightarrow0$,
$E[\rho_{M,p_0}(t,\xi',y',P_{\xi'})]\rightarrow0$, as $M\rightarrow\infty$.
In fact, $E[\rho_{M,p_0}(t,\vartheta',y',P_{\xi'})]=E[\rho_{M,p_0}(t,\xi',y',P_{\xi'})]$,
 as $\rho_{M,p_0}(t,x',y',P_{\xi'})$ is deterministic and $P_{\vartheta'}=P_{\xi'}$.
Then from H\"{o}lder's inequality, for all $p\in[2,p_0]$ we have
\begin{equation*}\label{eq6.16+5}
\begin{split}
 E[(\!\int_t^T\!\!|H(r,x,x',y,y',\xi,\xi')|^2dr)^{\frac{p}{2}}]
\leq&\ C_{p}M^{p}\(\!|x\!-\!x'|^{p}\!+\!|y\!-\!y'|^{p}\!+\!W_2(P_{\xi},P_{\xi'})^{p}\!+\!(E[|\vartheta\!-\!\vartheta'|^2])^{\frac{p}{2}}\!\)\\
&\ +\!\rho_{M,p}(t,x',y',P_{\xi'})\!+\!\rho_{M,p}(t,y',P_{\xi'})\!+\!E[\rho_{M,p}(t,\xi',y',P_{\xi'})],
\end{split}
\end{equation*}
where $\rho_{M,p}(t,x',y',P_{\xi'})\rightarrow0$, $\rho_{M,p}(t,y',P_{\xi'})\rightarrow0$,
$E[\rho_{M,p}(t,\xi',y',P_{\xi'})]\rightarrow0$, as $M\rightarrow\infty$.
Hence, substituting in \eqref{eq6.16+2} $x=\vartheta$, $x'=\vartheta'$,
 and using the above estimates, we get from Corollary \ref{corA.1},
\begin{equation}\label{eq6.16+6}
  \begin{split}
  &\  E\[\sup_{s\in[t,T]}|R_s^{t,\vartheta,P_{\xi}}(y)-R_s^{t,\vartheta',P_{\xi'}}(y')|^2+\int_t^T|Q_s^{t,\vartheta,P_{\xi}}(y)-Q_s^{t,\vartheta',P_{\xi'}}(y')|^2ds\]\\
 \leq&\   CM^2\(|y-y'|^2+W_2(P_{\xi},P_{\xi'})^2+E[|\vartheta-\vartheta'|^2]\)+E[\rho_{M}(t,\xi',y',P_{\xi'})]+\rho_{M}(t,y',P_{\xi'}).
  \end{split}
\end{equation}
By combining \eqref{eq6.16+6} with \eqref{eq6.16+3} and the preceding estimates we obtain from Corollary \ref{corA.2}
\begin{equation}\label{eq6.16+7}
  \begin{split}
  &\  E\[\sup_{s\in[t,T]}|R_s^{t,x,P_{\xi}}(y)-R_s^{t,x',P_{\xi'}}(y')|^{p_0}+(\int_t^T|Q_s^{t,x,P_{\xi}}(y)-Q_s^{t,x',P_{\xi'}}(y')|^2ds)^{\frac{p_0}{2}}\]\\
 \leq&\ C_{p_0}M^{p_0}\(|x-x'|^{p_0}+|y-y'|^{p_0}+W_2(P_{\xi},P_{\xi'})^{p_0}+(E[|\vartheta-\vartheta'|^2])^{\frac{p_0}{2}}\)\\
  &\ +\rho_{M,p_0}(t,x',y',P_{\xi'})+\rho_{M,p_0}(t,y',P_{\xi'})+E[\rho_{M,p_0}(t,\xi',y',P_{\xi'})].
  \end{split}
\end{equation}
On the other hand, from \eqref{eq6.16+2}
\begin{equation}\nonumber
  \begin{split}
  &\ E\[\sup_{s\in[t,T]}|O_s^{t,x,P_{\xi}}(y)-O_s^{t,x',P_{\xi'}}(y')|^{p_0}\]\\
\leq &\  C_{p_0}E\[\!\sup_{s\in[t,T]}\!|R_s^{t,x,P_{\xi}}\!(y)\!-\!R_s^{t,x',P_{\xi'}}\!(y')|^{p_0}\!\]
\!+\!C_{p_0}\!\Bigg(\!\int_t^T\!\!\Bigg\{\!\!\widehat{E}\[\big|(\partial_{\mu}h)_{1}(\!P_{(Y_r^{t,\xi}\!,Z_r^{t,\xi})},\!(\widehat{Y}_r^{t,y,P_{\xi}}\!,\!\widehat{Z}_r^{t,y,P_{\xi}}))
          \partial_{x}\widehat{Y}_r^{t,y,P_{\xi}}\\
   &\  - (\partial_{\mu}h)_{1}(P_{(Y_r^{t,\xi'},Z_r^{t,\xi'})},(\widehat{Y}_r^{t,y',P_{\xi'}},\widehat{Z}_r^{t,y',P_{\xi'}}))
          \partial_{x}\widehat{Y}_r^{t,y',P_{\xi'}}  \big|^2  \]\\
&\ +\widehat{E}\[\big|(\partial_{\mu}h)_1(P_{(Y_r^{t,\xi}\!,Z_r^{t,\xi})},\!(\widehat{Y}_r^{t,\widehat{\xi}}\!,\widehat{Z}_r^{t,\widehat{\xi}}))
       \widehat{O}_{r}^{t,\widehat{\xi}}(y)
       -(\partial_{\mu}h)_1(P_{(Y_r^{t,\xi'}\!,Z_r^{t,\xi'})},\!(\widehat{Y}_r^{t,\widehat{\xi'}}\!,\widehat{Z}_r^{t,\widehat{\xi'}}))
       \widehat{O}_{r}^{t,\widehat{\xi'}}(y') \big|^2\]\\
&\ \!+\!\widehat{E}\[\!\big|\!(\partial_{\mu}h)_{2}(P_{(Y_r^{t,\xi}\!,Z_r^{t,\xi})}\!,\!(\widehat{Y}_r^{t,y,P_{\xi}}\!,\!\widehat{Z}_r^{t,y,P_{\xi}}\!)\!)
       \partial_{x}\widehat{Z}_r^{t,y,P_{\xi}}
       \!-\!(\partial_{\mu}h)_{2}(P_{(Y_r^{t,\xi'}\!,Z_r^{t,\xi'})}\!,\!(\widehat{Y}_r^{t,y'\!,P_{\xi'}}\!,\widehat{Z}_r^{t,y'\!,P_{\xi'}}\!)\!)
       \partial_{x}\widehat{Z}_r^{t,y'\!,P_{\xi'}}\!\big|^2\!\]\\
&\   +\widehat{E}\[\big|(\partial_{\mu}h)_2(P_{(Y_r^{t,\xi},Z_r^{t,\xi})},(\widehat{Y}_r^{t,\widehat{\xi}},\widehat{Z}_r^{t,\widehat{\xi}})) \widehat{Q}_{r}^{t,\widehat{\xi}}(y)
-(\partial_{\mu}h)_2(P_{(Y_r^{t,\xi'},Z_r^{t,\xi'})},(\widehat{Y}_r^{t,\widehat{\xi'}},\widehat{Z}_r^{t,\widehat{\xi'}})) \widehat{Q}_{r}^{t,\widehat{\xi'}}(y')\big|^2 \]
  \!\Bigg\}\!dr\Bigg)^{\frac{p_0}{2}}\\
\leq &\   C_{p_0}E\[\sup_{s\in[t,T]}|R_s^{t,x,P_{\xi}}(y)-R_s^{t,x',P_{\xi'}}(y')|^{p_0}\]
     +C_{p_0}\Bigg(\int_t^T\(\widehat{E}[|\partial_{x}\widehat{Y}_r^{t,y,P_{\xi}}-\partial_{x}\widehat{Y}_r^{t,y',P_{\xi'}} |^2 ]\\
&\      +\widehat{E}[|\partial_{x}\widehat{Z}_r^{t,y,P_{\xi}}\!-\!\partial_{x}\widehat{Z}_r^{t,y',P_{\xi'}} |^2 ]
       \!+\!\widehat{E}[|\widehat{O}_{r}^{t,\widehat{\vartheta},P_{\xi}}(y)\!-\!\widehat{O}_{r}^{t,\widehat{\vartheta'},P_{\xi'}}(y') |^2 ]
       \!+\!\widehat{E}[|\widehat{Q}_{r}^{t,\widehat{\vartheta},P_{\xi}}(y)\!-\!\widehat{Q}_{r}^{t,\widehat{\vartheta'},P_{\xi'}}(y') |^2 ]
     \!\)dr\!\Bigg)^{\frac{p_0}{2}}\\
  \end{split}
\end{equation}
\begin{equation}\label{eq6.16+8}
  \begin{split}
&\  +C_{p_0}M^{p_0} \Bigg(\int_t^T\(W_2(P_{(Y_r^{t,\xi},Z_r^{t,\xi})},P_{(Y_r^{t,\xi'},Z_r^{t,\xi'})})^2
     +\widehat{E}[|(\widehat{Y}_r^{t,y,P_{\xi}},\widehat{Z}_r^{t,y,P_{\xi}})-(\widehat{Y}_r^{t,y',P_{\xi'}},\widehat{Z}_r^{t,y',P_{\xi'}})|^2]\\
 &\     +\widehat{E}[|(\widehat{Y}_r^{t,\vartheta,P_{\xi}},\widehat{Z}_r^{t,\vartheta,P_{\xi}})
         -(\widehat{Y}_r^{t,\vartheta',P_{\xi'}},\widehat{Z}_r^{t,\vartheta',P_{\xi'}})|^2]\)dr\!\Bigg)^{\frac{p_0}{2}}\\
&\    +C_{p_0} \Bigg(\int_t^T\widehat{E}\[|(\partial_{x}\widehat{Y}_{r}^{t,y',P_{\xi'}},\partial_{x}\widehat{Z}_{r}^{t,y',P_{\xi'}})|^2
              I_{\{|(\partial_{x}\widehat{Y}_{r}^{t,y',P_{\xi'}},\partial_{x}\widehat{Z}_{r}^{t,y',P_{\xi'}})|\geq M\}}\]dr\Bigg)^{\frac{p_0}{2}}\\
&\    +C_{p_0} \Bigg(\int_t^T\widehat{E}\[|(\widehat{O}_{r}^{t,\widehat{\vartheta'},P_{\xi'}}(y'),\widehat{Q}_{r}^{t,\widehat{\vartheta'},P_{\xi'}}(y'))|^2
              I_{\{|(\widehat{O}_{r}^{t,\widehat{\vartheta'},P_{\xi'}}(y'),\widehat{Q}_{r}^{t,\widehat{\vartheta'},P_{\xi'}}(y'))|\geq M\}}\]dr\Bigg)^{\frac{p_0}{2}} .
  \end{split}
\end{equation}
Put
\begin{equation*}
  \begin{split}
    &\rho_{M}(t,y',P_{\xi'}):=\int_t^T\widehat{E}\[|(\partial_{x}\widehat{Y}_{r}^{t,y',P_{\xi'}},\partial_{x}\widehat{Z}_{r}^{t,y',P_{\xi'}})|^2
              I_{\{|(\partial_{x}\widehat{Y}_{r}^{t,y',P_{\xi'}},\partial_{x}\widehat{Z}_{r}^{t,y',P_{\xi'}})|\geq M\}}\]dr,\\
 & \rho_{M}(t,x',y',P_{\xi'}):=\int_t^T\widehat{E}\[|(\widehat{O}_{r}^{t,x',P_{\xi'}}(y'),\widehat{Q}_{r}^{t,x',P_{\xi'}}(y'))|^2
              I_{\{|(\widehat{O}_{r}^{t,x',P_{\xi'}}(y'),\widehat{Q}_{r}^{t,x',P_{\xi'}}(y'))|\geq M\}}\]dr.
  \end{split}
\end{equation*}
From the independence of $\vartheta'\in L^{2}(\mathcal{G}_t;\mathbb{R})$ of $(O_s^{t,x',\xi'}(y'),Q_s^{t,x',\xi'}(y'))$:
\begin{equation*}
  \begin{split}
 & \widehat{E}[\rho_{M}(t,\widehat{\vartheta'},y',P_{\xi'})]
 =\int_t^T\widehat{E}\[|(\widehat{O}_{r}^{t,\widehat{\vartheta'},P_{\xi'}}(y'),\widehat{Q}_{r}^{t,\widehat{\vartheta'},P_{\xi'}}(y'))|^2
              I_{\{|(\widehat{O}_{r}^{t,\widehat{\vartheta'},P_{\xi'}}(y'),\widehat{Q}_{r}^{t,\widehat{\vartheta'},P_{\xi'}}(y'))|\geq M\}}\]dr.
  \end{split}
\end{equation*}
Then, from \eqref{eq6.16+6}, \eqref{eq6.16+7}, \eqref{eq6.16+8}, the Propositions \ref{prop4.1} and \ref{prop6.1},
\begin{equation}\nonumber
  \begin{split}
    E\[\sup_{s\in[t,T]}|O_s^{t,x,P_{\xi}}&(y)-O_s^{t,x',P_{\xi'}}(y')|^{p_0}\]
 \leq C_{p_0}\(\int_t^TE[|O_r^{t,\vartheta,P_{\xi}}(y)-O_r^{t,\vartheta',P_{\xi'}}(y')|^2]dr\)^{\frac{p_0}{2}}\\
  &\  +C_{p_0}M^{p_0}\(|x-x'|^{p_0}+|y-y'|^{p_0}+W_2(P_{\xi},P_{\xi'})^{p_0}+(E[|\vartheta-\vartheta'|^2])^{\frac{p_0}{2}}\)\\
  &\ +\rho_{M,p_0}(t,x',y',P_{\xi'})+\rho_{M,p_0}(t,y',P_{\xi'})+E[\rho_{M,p_0}(t,\xi',y',P_{\xi'})].
  \end{split}
\end{equation}
Furthermore, from H\"{o}lder's inequality, for all $p\in[2,p_0]$ we have
\begin{equation}\label{eq6.16+9}
  \begin{split}
    E\[\sup_{s\in[t,T]}|O_s^{t,x,P_{\xi}}&(y)-O_s^{t,x',P_{\xi'}}(y')|^{p}\]
 \leq C_{p}\(\int_t^TE[|O_r^{t,\vartheta,P_{\xi}}(y)-O_r^{t,\vartheta',P_{\xi'}}(y')|^2]dr\)^{\frac{p}{2}}\\
  &\  +C_{p}M^{p}\(|x-x'|^{p}+|y-y'|^{p_0}+W_2(P_{\xi},P_{\xi'})^{p}+(E[|\vartheta-\vartheta'|^2])^{\frac{p}{2}}\)\\
  &\ +\rho_{M,p}(t,x',y',P_{\xi'})+\rho_{M,p}(t,y',P_{\xi'})+E[\rho_{M,p}(t,\xi',y',P_{\xi'})].
  \end{split}
\end{equation}
First we substitute $x=\vartheta$, $x'=\vartheta'$ in \eqref{eq6.16+9}, then we get
\begin{equation}\label{eq6.16+10}
  \begin{split}
  &\  E\[\sup_{s\in[t,T]}|O_s^{t,\vartheta,P_{\xi}}(y)-O_s^{t,\vartheta',P_{\xi'}}(y')|^2\]
  =E\[E\[\sup_{s\in[t,T]}|O_s^{t,x,P_{\xi}}(y)-O_s^{t,x',P_{\xi'}}(y')|^2\]\big|_{\substack{x=\vartheta\\ x'=\vartheta'}}\]\\
 \leq&\   C\int_t^TE[|O_r^{t,\vartheta,P_{\xi}}(y)-O_r^{t,\vartheta',P_{\xi'}}(y')|^2]dr
 +CM^2\(|y-y'|^2+W_2(P_{\xi},P_{\xi'})^2+E[|\vartheta-\vartheta'|^2]\)\\
&\ +E[\rho_{M}(t,\xi',y',P_{\xi'})]+\rho_{M}(t,y',P_{\xi'}).
  \end{split}
\end{equation}
Hence, from Gronwall's Lemma,
\begin{equation}\label{eq6.16+11}
  \begin{split}
  &\  E\[\sup_{s\in[t,T]}|O_s^{t,\vartheta,P_{\xi}}(y)-O_s^{t,\vartheta',P_{\xi'}}(y')|^2\]\\
 \leq&\  CM^2\(|y-y'|^2+W_2(P_{\xi},P_{\xi'})^2+E[|\vartheta-\vartheta'|^2]\)
 +E[\rho_{M}(t,\xi',y',P_{\xi'})]+\rho_{M}(t,y',P_{\xi'}).
  \end{split}
\end{equation}
Now substituting \eqref{eq6.16+11} in \eqref{eq6.16+9}, and combining this with the result \eqref{eq6.16+7} we obtain
\begin{equation}\label{eq6.16+13}
  \begin{split}
  &\  E\[\sup_{s\in[t,T]}|O_s^{t,x,P_{\xi}}(y)-O_s^{t,x',P_{\xi'}}(y')|^p+(\int_t^T|Q_s^{t,x,P_{\xi}}(y)-Q_s^{t,x',P_{\xi'}}(y')|^2ds)^{\frac{p}{2}}\]\\
 \leq&\ C_pM^p\(|x-x'|^p+|y-y'|^p+W_2(P_{\xi},P_{\xi'})^p+(E[|\vartheta-\vartheta'|^2])^{\frac{p}{2}}\)\\
  &\ +\rho_{M,p}(t,x',y',P_{\xi'})+\rho_{M,p}(t,y',P_{\xi'})+E[\rho_{M,p}(t,\xi',y',P_{\xi'})].
  \end{split}
\end{equation}
Together with the fact that
\begin{equation*}
\begin{split}
  W_{2}(P_{\xi},P_{\xi'})=\inf\Big\{\Big(E[|\vartheta-\vartheta'|^{2}]\Big)^{\frac{1}{2}}\Big|\
\vartheta,\vartheta'\in L^2(\mathcal{G}_t;\mathbb{R}),\ P_{\vartheta}=P_{\xi},\ P_{\vartheta'}=P_{\xi'}\Big\},
\end{split}
\end{equation*}
we get the wished estimate.
\end{proof}

\begin{lemma} \label{le6.2}
Suppose (H5.1) and (H6.1) hold true. Then, for all $0\leq t\leq s\leq T$ and $x\in\mathbb{R}$, the lifted processes
$L^{2}(\mathcal{G}_t;\mathbb{R})\ni\xi\rightarrow Y_s^{t,x,\xi}:=Y_s^{t,x,P_{\xi}}\in L^{2}(\mathcal{F}_s;\mathbb{R})$,
$L^{2}(\mathcal{G}_t;\mathbb{R})\ni\xi\rightarrow Z_{\cdot}^{t,x,\xi}:=Z_{\cdot}^{t,x,P_{\xi}}\in\mathcal{H}_{\mathcal{F}}^2(t,T;\mathbb{R})$
as functionals of $\xi$ are G\^{a}teaux differentiable, and the G\^{a}teaux derivatives in direction
$\eta\in L^{2}(\mathcal{G}_t;\mathbb{R})$ are just
 $\mathcal{O}_s^{t,x,\xi}(\eta),\mathcal{Q}_s^{t,x,\xi}(\eta)$, respectively, i.e.,
\begin{equation*}
 \begin{split}
   & \partial_{\xi}Y_s^{t,x,\xi}(\eta)=\mathcal{O}_s^{t,x,\xi}(\eta)=\overline{E}[O_s^{t,x,P_{\xi}}(\overline{\xi})\cdot\overline{\eta}],\ s\in[t,T],\ P\mbox{-}a.s.,\\
  & \partial_{\xi}Z_s^{t,x,\xi}(\eta)=\mathcal{Q}_s^{t,x,\xi}(\eta)=\overline{E}[Q_s^{t,x,P_{\xi}}(\overline{\xi})\cdot\overline{\eta}],\  dsdP\mbox{-}a.e.\
      \mbox{on}\ [t,T]\times\Omega,
 \end{split}
\end{equation*}
where $\mathcal{O}_s^{t,x,\xi}(\eta)$, $\mathcal{Q}_s^{t,x,\xi}(\eta)$, $O_s^{t,x,P_{\xi}}(y)$, $Q_s^{t,x,P_{\xi}}(y)$
are defined in Lemma \ref{le6.1}.
\end{lemma}
\begin{proof}
The proof is split into two steps. For simplicity of redaction we restrict ourselves to the dimensions $d=1$, $l=1$,
and we let $f\equiv0$, $g(\Pi_r^{t,x,P_{\xi}},P_{\Pi_r^{t,\xi}})=g(Z_r^{t,x,P_{\xi}})$,
$h(P_{\Pi_r^{t,\xi}})=h(P_{Z_r^{t,\xi}})$ and $\Phi=\Phi(x)$.

\noindent \textbf{Step 1.} We prove that the directional derivatives $Y_s^{t,x,\xi}$, $Z_s^{t,x,\xi}$ in all direction
$\eta\in L^{2}(\mathcal{G}_t;\mathbb{R})$ exist, and
\begin{equation*}
 \begin{split}
   &\frac{1}{q}(Y_s^{t,x,P_{\xi+q\eta}}-Y_s^{t,x,P_{\xi}})-\mathcal{O}_s^{t,x,\xi}(\eta)\xrightarrow[q\rightarrow 0]{L^2(\mathcal{F}_s;\mathbb{R})}0 ,\\
   &\frac{1}{q}(Z_{\cdot}^{t,x,P_{\xi+q\eta}}-Z_{\cdot}^{t,x,P_{\xi}})-\mathcal{Q}_{\cdot}^{t,x,\xi}(\eta) \xrightarrow[q\rightarrow 0]{\mathcal{H}_{\mathcal{F}}^2(t,T;\mathbb{R})}0 .
 \end{split}
\end{equation*}
In fact, for all $s\in[t,T]$,
\begin{equation}\label{eq6.26}
  \begin{split}
&\ \frac{1}{q} (Y_s^{t,x,P_{\xi+q\eta}}-Y_s^{t,x,P_{\xi}})-\mathcal{O}_{s}^{t,x,\xi}(\eta)\\
=&\ I_1(x,q)+I_2(s,x,q)+I_3(s,x,q)-\int_s^T\(\frac{Z_r^{t,x,P_{\xi+q\eta}}-Z_r^{t,x,P_{\xi}}}{q}-\mathcal{Q}_{r}^{t,x,\xi}(\eta)\)dW_r,
  \end{split}
\end{equation}
where
\begin{equation*}
  \begin{split}
  &I_1(x,q):= \frac{1}{q}(\Phi (X_T^{t,x,P_{\xi+q\eta}})-\Phi (X_T^{t,x,P_{\xi}}))-\partial_{x}\Phi (X_T^{t,x,P_{\xi}})\mathcal{U}_{T}^{t,x,\xi}(\eta),\\
& I_2(s,x,q):=\int_s^T\Bigg\{\frac{1}{q}\big(g(Z_r^{t,x,P_{\xi+q\eta}})-g(Z_r^{t,x,P_{\xi}})\big)
   -\partial_{z}g(Z_r^{t,x,P_{\xi}})\mathcal{Q}_{r}^{t,x,\xi}(\eta)\Bigg\} d\overleftarrow{B_r},\\
&  I_3(s,x,q):=\int_s^T\!\Bigg\{\frac{1}{q}\big(h(P_{Z_r^{t,\xi+q\eta}})\!-\!h(P_{Z_r^{t,\xi}})\big)
 -\widehat{E}\[(\partial_{\mu}h)(P_{Z_r^{t,\xi}},\widehat{Z}_r^{t,\widehat{\xi},P_{\xi}})
 \partial_{x}\widehat{Z}_r^{t,\widehat{\xi},P_{\xi}}\widehat{\eta}\]\\
&\hspace{90pt} -\widehat{E}\[(\partial_{\mu}h)(P_{Z_r^{t,\xi}},\widehat{Z}_r^{t,\widehat{\xi}}) \widehat{\mathcal{Q}}_{r}^{t,\widehat{\xi}}(\widehat{\eta}) \]\Bigg\}d\overleftarrow{B_r}
   \end{split}
\end{equation*}
(Recall that $\widehat{\mathcal{Q}}_{r}^{t,\widehat{\xi}}(\widehat{\eta})=\widehat{\mathcal{Q}}_{r}^{t,\widehat{\xi},P_{\xi}}(\widehat{\eta})$).
 First at all, we remark that
\begin{equation*}
  \begin{split}
 I_1(x,q)\!=\!&\ \! \int_0^1\!(\partial_x\Phi) \big(X_T^{t,x,P_{\xi}}\!+\!\rho(X_T^{t,x,P_{\xi+q\eta}}\!-\!X_T^{t,x,P_{\xi}})\big)d\rho
 \frac{X_T^{t,x,P_{\xi+q\eta}}\!-\!X_T^{t,x,P_{\xi}}}{q}
-(\partial_{x}\Phi) (X_T^{t,x,P_{\xi}})\mathcal{U}_{T}^{t,x,\xi}(\eta)\\
   \end{split}
\end{equation*}
\begin{equation*}
  \begin{split}
=&\ \int_0^1\((\partial_x\Phi) \big(X_T^{t,x,P_{\xi}}+\rho(X_T^{t,x,P_{\xi+q\eta}}-X_T^{t,x,P_{\xi}})\big)-(\partial_{x}\Phi) (X_T^{t,x,P_{\xi}})\)d\rho \frac{X_T^{t,x,P_{\xi+q\eta}}-X_T^{t,x,P_{\xi}}}{q}\\
&\ +(\partial_{x}\Phi) (X_T^{t,x,P_{\xi}})\big(\frac{X_T^{t,x,P_{\xi+q\eta}}-X_T^{t,x,P_{\xi}}}{q}-\mathcal{U}_{T}^{t,x,\xi}(\eta)\big).
   \end{split}
\end{equation*}
Consequently, as $\partial_x\Phi$ is Lipschitz and bounded, and as $I_1(x,q)$ is independent of $\mathcal{G}_t$,
\begin{equation*}
  \begin{split}
& E[|I_1(x,q)|^2\big|\mathcal{G}_t] \!=\!E[|I_1(x,q)|^2]
  \!\leq\! \frac{C}{q^2}E[|X_T^{t,x,P_{\xi+q\eta}}\!-\!X_T^{t,x,P_{\xi}}|^4]
 \!+\!CE[|\frac{X_T^{t,x,P_{\xi+q\eta}}\!-\!X_T^{t,x,P_{\xi}}}{q}\!-\!\mathcal{U}_{T}^{t,x,\xi}(\eta)|^2].
   \end{split}
\end{equation*}
From Lemma \ref{le3.1} we have
\begin{equation*}
E[|X_T^{t,x,P_{\xi+q\eta}}-X_T^{t,x,P_{\xi}}|^4]\leq CW_2(P_{\xi+q\eta},P_{\xi})^4\leq Cq^4(E[|\eta|^2])^2.
\end{equation*}
On the other hand, from Proposition \ref{prop5.3}, we know
\begin{equation*}
  \begin{split}
&\ E[|\frac{X_T^{t,x,P_{\xi+q\eta}}\!-\!X_T^{t,x,P_{\xi}}}{q}-\mathcal{U}_{T}^{t,x,\xi}(\eta)|^2]
=  E[|\widehat{E}[  \int_0^1\big(\partial_{\mu}X_T^{t,x,P_{\xi\!+\!\rho q\eta}}(\widehat{\xi}+\rho q\widehat{\eta})-
 \partial_{\mu}X_T^{t,x,P_{\xi}}(\widehat{\xi})\big)d\rho\cdot \widehat{\eta}]|^2]\\
\leq &\ E\[\(\widehat{E}\[  \int_0^1\big|\partial_{\mu}X_T^{t,x,P_{\xi+\rho q\eta}}(\widehat{\xi}+\rho q\widehat{\eta})-
 \partial_{\mu}X_T^{t,x,P_{\xi}}(\widehat{\xi})\big|^2d\rho\]\cdot \widehat{E}[|\widehat{\eta}|^2]\)\]\\
\leq &\ E[|\eta|^2] \int_0^1\widehat{E}\[E[\big|\partial_{\mu}X_T^{t,x,P_{\xi+\rho q\eta}}(y)-
 \partial_{\mu}X_T^{t,x,P_{\xi}}(y')\big|^2]\big|_{y=\widehat{\xi}+\rho q\widehat{\eta},\ y'=\widehat{\xi}}\]d\rho\\
\leq &\ E[|\eta|^2] \int_0^1\widehat{E}\[\big(W_2(P_{\xi+\rho q\eta},P_{\xi})^2+|y-y'|^2\big)\big|_{y=\widehat{\xi}+\rho q\widehat{\eta},\ y'=\widehat{\xi}}\]d\rho
\leq  Cq^2(E[|\eta|^2])^2.
   \end{split}
\end{equation*}
This shows that
\begin{equation}\label{eq6.27}
  E[|I_1(x,q)|^2\big|\mathcal{G}_t] =E[|I_1(x,q)|^2]\leq Cq^2(E[|\eta|^2])^2.
\end{equation}
We now consider $I_2(s,x,q)$. We define
$Z_r^{t,x,\xi}(\eta,\rho):=Z_r^{t,x,P_{\xi}}+\rho(Z_r^{t,x,P_{\xi+q\eta}}-Z_r^{t,x,P_{\xi}}).$
Then, making use of the fact that $g\in C_b^{1}(\mathbb{R})$, we get
\begin{equation}\nonumber
  \begin{split}
 I_2(s,x,q)=&\ \int_s^T\!\Big\{\int_0^1\frac{1}{q}\partial_{\rho}\big[g(Z_r^{t,x,\xi}(\eta,\rho))\big]d\rho
 - (\partial_{z}g)(Z_r^{t,x,P_{\xi}})\mathcal{Q}_{r}^{t,x,\xi}(\eta) \Big\}d\overleftarrow{B_r}\\
  \end{split}
\end{equation}
\begin{equation}\label{eq6.28+1}
  \begin{split}
=&\  \int_s^T\!\Big\{\!\int_0^1\! (\partial_{z}g)(Z_r^{t,x,\xi}(\eta,\rho))\frac{Z_r^{t,x,P_{\xi+q\eta}}\!-\!Z_r^{t,x,P_{\xi}}}{q}d\rho
 \!-\! (\partial_{z}g)(Z_r^{t,x,P_{\xi}})\mathcal{Q}_{r}^{t,x,\xi}(\eta) \Big\}d\overleftarrow{B_r}\\
=&\  I_{2,1}(s,x,q)+ I_{2,2}(s,x,q),
  \end{split}
\end{equation}
where
\begin{equation*}
  \begin{split}
 I_{2,1}(s,x,q):=&\ \int_s^T\!\(\int_0^1 (\partial_{z}g)(Z_r^{t,x,\xi}(\eta,\rho))d\rho\)\cdot\(\frac{Z_r^{t,x,P_{\xi+q\eta}}-Z_r^{t,x,P_{\xi}}}{q}
 - \mathcal{Q}_{r}^{t,x,\xi}(\eta) \)d\overleftarrow{B_r},\\
 I_{2,2}(s,x,q):=&\ \int_s^T\!\(\int_0^1\big( (\partial_{z}g)(Z_r^{t,x,\xi}(\eta,\rho))- (\partial_{z}g)(Z_r^{t,x,P_{\xi}})\big)d\rho\)
\cdot\mathcal{Q}_{r}^{t,x,\xi}(\eta) d\overleftarrow{B_r}.
  \end{split}
\end{equation*}
Let $\Xi_{r}^{1,q}:=\int_0^1\big( (\partial_{z}g)(Z_r^{t,x,\xi}(\eta,\rho))- (\partial_{z}g)(Z_r^{t,x,P_{\xi}})\big)d\rho$, $r\in[t,T]$.
Notice, as $g\in C_b^{1}(\mathbb{R})$, combining with Proposition \ref{prop4.1}, it follows that
$|\Xi_{r}^{1,q}|\leq C$, $\Xi_{r}^{1,q}\rightarrow0$, $q\downarrow0$, in $drdP$-measure.
On the other hand, from \eqref{eq6.7}, $ E[\int_t^T|\mathcal{Q}_s^{t,x,\xi}(\eta)|^2ds]\leq C$, and, consequently, from
 the dominated convergence theorem,
\begin{equation}\label{eq6.28+2}
  \begin{split}
 E[\sup_{s\in[t,T]}| I_{2,2}(s,x,q)|^2]\leq C E\[\int_t^T|\Xi_{r}^{1,q}|^2
\cdot|\mathcal{Q}_{r}^{t,x,\xi}(\eta)|^2 dr\]\rightarrow0,\ \mbox{as}\ q\downarrow0.
  \end{split}
\end{equation}
Recall that due to Assumption (H4.1), $|\partial_{z}g(z)|^2\leq \alpha_1<1$, $z\in\mathbb{R}$. Therefore
\begin{equation}\label{eq6.28+2'}
  \begin{split}
 E[| I_{2,1}(s,x,q)|^2]\leq \alpha_1 E\[\int_s^T|\frac{Z_r^{t,x,P_{\xi+q\eta}}-Z_r^{t,x,P_{\xi}}}{q}
 - \mathcal{Q}_{r}^{t,x,\xi}(\eta) |^2dr\].
  \end{split}
\end{equation}
 Next, we investigate $I_3(s,x,q)$. For this we observe that
 \begin{equation}\label{eq6.28+3}
  \begin{split}
I_3(s,x,q)=I_{3,1}(s,x,q)+I_{3,2}(s,x,q),
  \end{split}
\end{equation}
where
\begin{equation*}
  \begin{split}
 I_{3,1}(s,x,q):=&\ \int_s^T\!\Bigg\{\frac{1}{q}\big(h(P_{Z_r^{t,\xi+q\eta}})\!-\!h(P_{Z_r^{t,\xi,P_{\xi+q\eta}}})\big)
 -\widehat{E}\[(\partial_{\mu}h)(P_{Z_r^{t,\xi}},\widehat{Z}_r^{t,\widehat{\xi},P_{\xi}})
 \partial_{x}\widehat{Z}_r^{t,\widehat{\xi},P_{\xi}}\widehat{\eta}\]\Bigg\}d\overleftarrow{B_r},\\
 I_{3,2}(s,x,q):=&\ \int_s^T\!\Bigg\{\frac{1}{q}\big(h(P_{Z_r^{t,\xi,P_{\xi+q\eta}}})\!-\!h(P_{Z_r^{t,\xi}})\big)
  -\widehat{E}\[(\partial_{\mu}h)(P_{Z_r^{t,\xi}},\widehat{Z}_r^{t,\widehat{\xi}}) \widehat{\mathcal{Q}}_{r}^{t,\widehat{\xi}}(\widehat{\eta}) \]\Bigg\}d\overleftarrow{B_r}.
  \end{split}
\end{equation*}
(i) Computation for $ I_{3,1}(s,x,q)$: Using Theorem \ref{th6.1} and Remark \ref{re6.0} we obtain
\begin{equation}\label{eq6.28+4}
  \begin{split}
 I_{3,1}(s,x,q)=&\ \int_s^T\!\Big\{\int_0^1\frac{1}{q}\partial_{\rho}\big[h(P_{Z_r^{t,\xi+\rho q\eta,P_{\xi+q\eta}}})\big]d\rho
 -\widehat{E}\[(\partial_{\mu}h)(P_{Z_r^{t,\xi}},\widehat{Z}_r^{t,\widehat{\xi},P_{\xi}})
 \partial_{x}\widehat{Z}_r^{t,\widehat{\xi},P_{\xi}}\widehat{\eta}\] \Big\}d\overleftarrow{B_r}\\
=&\  \int_s^T\!\Big\{\int_0^1 \widehat{E}\[(\partial_{\mu}h)(P_{Z_r^{t,\xi+\rho q\eta,P_{\xi+q\eta}}},\widehat{Z}_r^{t,\widehat{\xi}+\rho q\widehat{\eta},P_{\xi+q\eta}})
 \partial_{x}\widehat{Z}_r^{t,\widehat{\xi}+\rho q\widehat{\eta},P_{\xi+q\eta}}\widehat{\eta}\]d\rho\\
&\ -\widehat{E}\[(\partial_{\mu}h)(P_{Z_r^{t,\xi}},\widehat{Z}_r^{t,\widehat{\xi},P_{\xi}})
 \partial_{x}\widehat{Z}_r^{t,\widehat{\xi},P_{\xi}}\widehat{\eta}\] \Big\}d\overleftarrow{B_r}.
  \end{split}
\end{equation}
Now, due to Proposition \ref{prop4.0}, as $\partial_{\mu}h$ is Lipschitz,
\begin{equation}\label{eq6.28+5}
  \begin{split}
 &\ \int_s^T\int_0^1 \widehat{E}\[|(\partial_{\mu}h)(P_{Z_r^{t,\xi+\rho q\eta,P_{\xi+q\eta}}},\widehat{Z}_r^{t,\widehat{\xi}+\rho q\widehat{\eta},P_{\xi+q\eta}})
 -(\partial_{\mu}h)(P_{Z_r^{t,\xi}},\widehat{Z}_r^{t,\widehat{\xi},P_{\xi}})|^2\]d\rho dr\\
\leq&\ C\int_s^T\int_0^1\( W_2(P_{Z_r^{t,\xi+\rho q\eta,P_{\xi+q\eta}}},P_{Z_r^{t,\xi}})^2
+\widehat{E}[|\widehat{Z}_r^{t,\widehat{\xi}+\rho q\widehat{\eta},P_{\xi+q\eta}}-\widehat{Z}_r^{t,\widehat{\xi},P_{\xi}}|^2]\)d\rho dr\\
\leq&\ C\int_s^T\int_0^1 \widehat{E}[|\widehat{Z}_r^{t,\widehat{\xi}+\rho q\widehat{\eta},P_{\xi+q\eta}}-\widehat{Z}_r^{t,\widehat{\xi},P_{\xi}}|^2]d\rho dr\\
  \end{split}
\end{equation}
\begin{equation}\nonumber
  \begin{split}
= &\ C\int_0^1\widehat{E}\[\widehat{E}\[\int_s^T|\widehat{Z}_r^{t,x,P_{\xi+q\eta}}-\widehat{Z}_r^{t,x',P_{\xi}}|^2 dr\]
     \big|_{x=\widehat{\xi}+\rho q\widehat{\eta},\ x'=\widehat{\xi}}\] d\rho\\
\leq&\ C\int_0^1\widehat{E}\[|x-x'|^2\big|_{x=\widehat{\xi}+\rho q\widehat{\eta},\ x'=\widehat{\xi}}
   +W_2(P_{\xi+q\eta},P_{\xi})^2\]d\rho\\
\leq&\ C q^2E[|\eta|^2], \ q> 0.
  \end{split}
\end{equation}
From \eqref{eq6.2+1+1} in Remark \ref{re6.0},
\begin{equation}\label{eq6.28+6}
  \begin{split}
 &\quad \int_s^T\int_0^1\widehat{E}\[| \partial_{x}\widehat{Z}_r^{t,\widehat{\xi}+\rho q\widehat{\eta},P_{\xi+q\eta}}
 -\partial_{x}\widehat{Z}_r^{t,\widehat{\xi},P_{\xi}}|^2\]d\rho dr\\
&=  \int_0^1\widehat{E}\[\widehat{E}\[\int_s^T|\partial_{x}\widehat{Z}_r^{t,x,P_{\xi+q\eta}}-\partial_{x}\widehat{Z}_r^{t,x',P_{\xi}}|^2 dr\]
     \big|_{x=\widehat{\xi}+\rho q\widehat{\eta},\ x'=\widehat{\xi}}\] d\rho\\
& \leq \int_0^1\widehat{E}\[ \(C M^2\big(|\widehat{\xi}+\rho q\widehat{\eta}-\widehat{\xi}|^2
   +W_2(P_{\xi+q\eta},P_{\xi})^2\big)+\rho_{M}(t,x',P_{\xi})\big|_{x'=\widehat{\xi}}\)\wedge C'\]d\rho\\
& \leq \widehat{E}\[\(C M^2\big(q^2|\widehat{\eta}|^2+q^2\widehat{E}[|\widehat{\eta}|^2]\big)
+\rho_{M}(t,\widehat{\xi},P_{\xi})\)\wedge C'\]\\
& \xrightarrow[]{q\rightarrow 0} \widehat{E}\[\rho_{M}(t,\widehat{\xi},P_{\xi})\wedge C'\]\xrightarrow[]{ M\rightarrow\infty} 0.
  \end{split}
\end{equation}
For simplicity, we put
$\Xi_{r}^{2,q,\rho}\!:=(\partial_{\mu}h)(P_{Z_r^{t,\xi+\rho q\eta,P_{\xi+q\eta}}},\widehat{Z}_r^{t,\widehat{\xi}+\rho q\widehat{\eta},P_{\xi+q\eta}})$,
$\Xi_{r}^{2}\!:=(\partial_{\mu}h)(P_{Z_r^{t,\xi}},\widehat{Z}_r^{t,\widehat{\xi},P_{\xi}})$,
$\Xi_{r}^{3,q,\rho}\!:=\partial_{x}\widehat{Z}_r^{t,\widehat{\xi}+\rho q\widehat{\eta},P_{\xi+q\eta}}$,
$\Xi_{r}^{3}\!:=\partial_{x}\widehat{Z}_r^{t,\widehat{\xi},P_{\xi}}$.
Obviously, as $\partial_{\mu}h$ is Lipschitz, $|\Xi_{r}^{2,q,\rho}|\leq C$, $r\in[t,T]$, $\rho\in[0,1]$, $q> 0$,
and from above it follows that
\begin{equation*}
  \Xi^{2,q,\cdot}\xrightarrow[]{q\rightarrow 0}\Xi^{2},\ \Xi^{3,q,\cdot}\xrightarrow[]{q\rightarrow 0}\Xi^{3}\ \mbox{in }\
  L^2([0,1]\times[t,T]\times\widehat{\Omega},d\rho drd\widehat{P}).
\end{equation*}
Hence, $\{\Xi^{3,q,\cdot},q>0\}$ is uniformly $L^2$-integrable, and as $|\Xi_{r}^{2,q,\rho}|\leq C$, also
$\{\Xi^{2,q,\cdot}\cdot\Xi^{3,q,\cdot},q>0\}$ is uniformly $L^2$-integrable.
Then $\Xi^{2,q,\cdot}\cdot\Xi^{3,q,\cdot}\rightarrow\Xi^{2}\cdot\Xi^{3}$ in $d\rho drd\widehat{P}$-measure,
imply that $\Xi^{2,q,\cdot}\cdot\Xi^{3,q,\cdot}\rightarrow\Xi^{2}\cdot\Xi^{3}$ in $L^2([0,1]\times[t,T]\times\widehat{\Omega},d\rho drd\widehat{P})$.
Consequently,
\begin{equation}\label{eq6.28+7}
  \begin{split}
 E[\sup_{s\in[t,T]}| I_{3,1}(s,x,q)|^2]
 \leq C\widehat{E}\[\int_t^T\int_0^1| \Xi_r^{2,q,\rho}\cdot\Xi_r^{3,q,\rho} -\Xi_r^{2}\cdot\Xi_r^{3}|^2 d\rho dr\]\cdot\widehat{E}[|\widehat{\eta}|^2]
 \rightarrow0,\ q\rightarrow0.
  \end{split}
\end{equation}
(ii) Computation for $ I_{3,2}(s,x,q)$: Letting $Z_r^{t,\xi}(\eta,\rho):=Z_r^{t,\xi}+\rho(Z_r^{t,\xi,P_{\xi+q\eta}}-Z_r^{t,\xi})$,
 and making use of the fact that $h\in C_b^{1}(\mathcal{P}_{2}(\mathbb{R}))$, we get
\begin{equation}\label{eq6.28+8}
  \begin{split}
 I_{3,2}(s,x,q)=&\ \int_s^T\!\Big\{\int_0^1\frac{1}{q}\partial_{\rho}\big[h(P_{Z_r^{t,\xi}(\eta,\rho)})\big]d\rho
  -\widehat{E}\[(\partial_{\mu}h)(P_{Z_r^{t,\xi}},\widehat{Z}_r^{t,\widehat{\xi}}) \widehat{\mathcal{Q}}_{r}^{t,\widehat{\xi}}(\widehat{\eta}) \] \Big\}d\overleftarrow{B_r}\\
=&\  \int_s^T\!\Big\{\int_0^1 \widehat{E}\[(\partial_{\mu}h)(P_{Z_r^{t,\xi}(\eta,\rho)},\widehat{Z}_r^{t,\widehat{\xi}}(\widehat{\eta},\rho))
 \frac{1}{q}(\widehat{Z}_r^{t,\widehat{\xi},P_{\xi+q\eta}}-\widehat{Z}_r^{t,\widehat{\xi}})\]d\rho\\
&\ -\widehat{E}\[(\partial_{\mu}h)(P_{Z_r^{t,\xi}},\widehat{Z}_r^{t,\widehat{\xi}}) \widehat{\mathcal{Q}}_{r}^{t,\widehat{\xi}}(\widehat{\eta}) \] \Big\}d\overleftarrow{B_r}
= I_{3,2,1}(s,x,q)+I_{3,2,2}(s,x,q),
  \end{split}
\end{equation}
where
\begin{equation*}
  \begin{split}
 I_{3,2,1}(s,x,q):=&\ \int_s^T\!\Big\{\int_0^1 \widehat{E}\[(\partial_{\mu}h)(P_{Z_r^{t,\xi}(\eta,\rho)},\widehat{Z}_r^{t,\widehat{\xi}}(\widehat{\eta},\rho))
 \( \frac{\widehat{Z}_r^{t,\widehat{\xi},P_{\xi+q\eta}}-\widehat{Z}_r^{t,\widehat{\xi}}}{q}
 -\widehat{\mathcal{Q}}_{r}^{t,\widehat{\xi}}(\widehat{\eta})\)\]d\rho\Big\}d\overleftarrow{B_r},\\
 I_{3,2,2}(s,x,q):=&\ \int_s^T\!\Big\{\int_0^1 \widehat{E}\[\((\partial_{\mu}h)(P_{Z_r^{t,\xi}(\eta,\rho)},\widehat{Z}_r^{t,\widehat{\xi}}(\widehat{\eta},\rho))
 -(\partial_{\mu}h)(P_{Z_r^{t,\xi}},\widehat{Z}_r^{t,\widehat{\xi}})\)\cdot\widehat{\mathcal{Q}}_{r}^{t,\widehat{\xi}}(\widehat{\eta})\]d\rho\Big\}d\overleftarrow{B_r}.
  \end{split}
\end{equation*}
We define $\Xi_{r}^{4,q,\rho}:=(\partial_{\mu}h)(P_{Z_r^{t,\xi}(\eta,\rho)},\widehat{Z}_r^{t,\widehat{\xi}}(\widehat{\eta},\rho))
 -(\partial_{\mu}h)(P_{Z_r^{t,\xi}},\widehat{Z}_r^{t,\widehat{\xi}})$, and in analogy to the computation for $I_{3,1}(s,x,q)$ we see
 $|\Xi_{r}^{4,q,\rho}|\leq C$, $r\in[t,T]$, $\rho\in[0,1]$, $q> 0$,
and $\Xi^{4,q,\cdot}\xrightarrow[]{q\rightarrow 0}0$ in $L^2([0,1]\times[t,T]\times\widehat{\Omega},d\rho drd\widehat{P})$.
As $ E[\int_t^T|\mathcal{Q}_s^{t,\xi}(\eta)|^2ds]\leq C$, it follows from the dominated
convergence theorem that $\Xi^{4,q,\cdot}\!\cdot\!\widehat{\mathcal{Q}}^{t,\widehat{\xi}}(\widehat{\eta})\xrightarrow[]{q\rightarrow 0}0$ in
$L^2([0,1]\times[t,T]\times\widehat{\Omega},d\rho drd\widehat{P})$. Thus,
\begin{equation}\label{eq6.28+9}
  \begin{split}
 E[\sup_{s\in[t,T]}| I_{3,2,2}(s,x,q)|^2]
 \leq C\widehat{E}\[\int_t^T\int_0^1|\Xi_{r}^{4,q,\rho}|^2\cdot|\widehat{\mathcal{Q}}_{r}^{t,\widehat{\xi}}(\widehat{\eta})|^2 d\rho dr\]
 \rightarrow0,\ q\rightarrow0.
  \end{split}
\end{equation}
Now we set $J(s,x,q):=\frac{1}{q} (Y_s^{t,x,P_{\xi+q\eta}}-Y_s^{t,x,P_{\xi}})-\mathcal{O}_{s}^{t,x,\xi}(\eta)-I_3(s,x,q)$.
Note that $J(s,x,q)$ is $\mathcal{F}_s$-measurable and independent of $\mathcal{G}_t$ (and, hence, of $\xi\in L^{2}(\mathcal{G}_t;\mathbb{R})$).
Thus, from \eqref{eq6.26}, we know
\begin{equation}\label{eq6.28+10}
  \begin{split}
 J(s,x,q)+\int_s^T\(\frac{Z_r^{t,x,P_{\xi+q\eta}}-Z_r^{t,x,P_{\xi}}}{q}-\mathcal{Q}_{r}^{t,x,\xi}(\eta)\)dW_r
= I_1(x,q)+I_2(s,x,q).
  \end{split}
\end{equation}
Hence, for any $\delta>0$,
\begin{equation*}
  \begin{split}
& E\[|J(s,x,q)|^2\!+\!\int_s^T|\frac{Z_r^{t,x,P_{\xi+q\eta}}\!-\!Z_r^{t,x,P_{\xi}}}{q}\!-\!\mathcal{Q}_{r}^{t,x,\xi}(\eta)|^2dr\]
\leq C_{\delta}E[|I_1(x,q)|^2]+(1+\delta)E[|I_2(s,x,q)|^2],
  \end{split}
\end{equation*}
where $C_{\delta}>0$ depends on $\delta>0$.
Furthermore, combining the preceding estimate with \eqref{eq6.27}, \eqref{eq6.28+1}, \eqref{eq6.28+2} and \eqref{eq6.28+2'}, we obtain
for $\delta>0$ such that $\alpha_1(1+\delta)^2<1$, the existence of a constant $ C_{\delta}>0$ depending on $\delta>0$ such that
\begin{equation}\label{eq6.28+11}
  \begin{split}
&\  E[|J(s,x,q)|^2]+(1-\alpha_1(1+\delta)^2)E\[\int_s^T|\frac{Z_r^{t,x,P_{\xi+q\eta}}-Z_r^{t,x,P_{\xi}}}{q}-\mathcal{Q}_{r}^{t,x,\xi}(\eta)|^2dr\]\\
\leq&\ C_{\delta}\(E[|I_1(x,q)|^2]+E[|I_{2,2}(s,x,q)|^2]\)\rightarrow0,\ q\downarrow0.
  \end{split}
\end{equation}
Consequently,
\begin{equation*}
 \begin{split}
   &\frac{1}{q}(Z_{\cdot}^{t,x,P_{\xi+q\eta}}-Z_{\cdot}^{t,x,P_{\xi}})-\mathcal{Q}_{\cdot}^{t,x,\xi}(\eta) \xrightarrow[q\rightarrow 0]{\mathcal{H}_{\mathcal{F}}^2(t,T;\mathbb{R})}0 .
 \end{split}
\end{equation*}
Substituting in \eqref{eq6.28+11} $\xi$ for $x$, we get from the independence of the integrands of $\mathcal{G}_t$ that
\begin{equation}\label{eq6.28+12}
  \begin{split}
&\  E[|J(s,\xi,q)|^2]+(1-\alpha_1(1+\delta)^2)E\[\int_s^T|\frac{Z_r^{t,\xi,P_{\xi+q\eta}}-Z_r^{t,\xi,P_{\xi}}}{q}-\mathcal{Q}_{r}^{t,\xi}(\eta)|^2dr\]\\
\leq&\ C_{\delta}\(E[|I_1(\xi,q)|^2]+E[|I_{2,2}(s,\xi,q)|^2]\),
  \end{split}
\end{equation}
and with the same arguments as those used for \eqref{eq6.28+11}, we show that
\begin{equation}\label{eq6.28+13}
  \begin{split}
E[|I_1(\xi,q)|^2]+E[|I_{2,2}(s,\xi,q)|^2]\rightarrow0,\ q\downarrow0.
  \end{split}
\end{equation}
Then also
\begin{equation}\label{eq6.28+14}
  \begin{split}
\frac{1}{q}(Z_{\cdot}^{t,\xi,P_{\xi+q\eta}}-Z_{\cdot}^{t,\xi,P_{\xi}})-\mathcal{Q}_{\cdot}^{t,\xi}(\eta) \xrightarrow[q\rightarrow 0]{\mathcal{H}_{\mathcal{F}}^2(t,T;\mathbb{R})}0 .
  \end{split}
\end{equation}
Thus, as $\partial_{\mu}h$ is bounded, we get
\begin{equation}\label{eq6.28+15}
  \begin{split}
 E[| I_{3,2,1}(s,x,q)|^2]
 \leq CE\[\int_t^T|\frac{Z_r^{t,\xi,P_{\xi+q\eta}}-Z_r^{t,\xi,P_{\xi}}}{q}-\mathcal{Q}_{r}^{t,\xi}(\eta)|^2 dr\]
 \rightarrow0,\ q\rightarrow0.
  \end{split}
\end{equation}
Recalling the definition of $J(s,x,q)$, \eqref{eq6.28+3} and \eqref{eq6.28+8}, we see
\begin{equation}\label{eq6.28+16}
  \begin{split}
 \frac{1}{q}(Y_s^{t,x,P_{\xi+q\eta}}-Y_s^{t,x,P_{\xi}})-\mathcal{O}_{s}^{t,x,\xi}(\eta)
= J(s,x,q)+ I_{3,1}(s,x,q)+ I_{3,2,1}(s,x,q)+ I_{3,2,2}(s,x,q).
  \end{split}
\end{equation}
Consequently, combining with \eqref{eq6.28+7}, \eqref{eq6.28+9}, \eqref{eq6.28+11}, \eqref{eq6.28+15} and \eqref{eq6.28+16}, we obtain
\begin{equation*}
 \begin{split}
   &\frac{1}{q}(Y_s^{t,x,P_{\xi+q\eta}}-Y_s^{t,x,P_{\xi}})-\mathcal{O}_s^{t,x,\xi}(\eta)\xrightarrow[q\rightarrow 0]{L^2(\mathcal{F}_s;\mathbb{R})}0 .
 \end{split}
\end{equation*}

\noindent\textbf{Step 2.} In Step 1 we have proved that the directional derivatives of $Y^{t,x,\xi}$, $Z^{t,x,\xi}$ in all
direction $\eta\in L^{2}(\mathcal{G}_t;\mathbb{R})$ exist and the directional directives
$\partial_{\xi}Y_s^{t,x,\xi}(\eta)$, $\partial_{\xi}Z_s^{t,x,\xi}(\eta)$ coincide with $\mathcal{O}_{s}^{t,x,\xi}(\eta)$,
$\mathcal{Q}_{s}^{t,x,\xi}(\eta)$.
Recall that $\mathcal{O}_{s}^{t,x,\xi}(\eta)$, $\mathcal{Q}_{s}^{t,x,\xi}(\eta)$ are linear and
continuous mappings. Consequently, $Y_s^{t,x,\xi}$, $Z_s^{t,x,\xi}$ as functionals of $\xi$ are G\^{a}teaux
differentiable, and furthermore, from Lemma \ref{le6.1} the G\^{a}teaux derivatives can be characterized by
\begin{equation*}
 \begin{split}
   & \partial_{\xi}Y_s^{t,x,\xi}(\eta)=\mathcal{O}_s^{t,x,\xi}(\eta)=\overline{E}[O_s^{t,x,P_{\xi}}(\overline{\xi})\cdot\overline{\eta}],\ s\in[t,T],\ P\mbox{-}a.s.,\\
  & \partial_{\xi}Z_s^{t,x,\xi}(\eta)=\mathcal{Q}_s^{t,x,\xi}(\eta)=\overline{E}[Q_s^{t,x,P_{\xi}}(\overline{\xi})\cdot\overline{\eta}],\  dsdP\mbox{-}a.e.
 \end{split}
\end{equation*}
The proof is complete.
\end{proof}
In order to prove $\partial_{\xi}Y_s^{t,x,\xi}(\eta)$, $\partial_{\xi}Z_s^{t,x,\xi}(\eta)$
are Fr\'{e}chet derivatives, we want to show that
\begin{equation*}
 \begin{split}
 &\  L^{2}(\mathcal{G}_t;\mathbb{R})\ni\xi\rightarrow \mathcal{O}_s^{t,x,\xi}\in L(L^{2}(\mathcal{G}_t;\mathbb{R}),L^{2}(\mathcal{F}_s;\mathbb{R})),\\
&\ L^{2}(\mathcal{G}_t;\mathbb{R})\ni\xi\rightarrow \mathcal{Q}_{\cdot}^{t,x,\xi}\in L(L^{2}(\mathcal{G}_t;\mathbb{R}),\mathcal{H}_{\mathcal{F}}^2(t,T;\mathbb{R}))
 \end{split}
\end{equation*}
are continuous.

\begin{lemma} \label{le6.3}
Under the Assumptions (H5.1) and (H6.1), for all $t\in[0,T]$, $x\in\mathbb{R}$, the linear functionals
$\partial_{\xi}Y_s^{t,x,\cdot}=\mathcal{O}_s^{t,x,\cdot}(=(\xi\rightarrow\mathcal{O}_s^{t,x,\xi}))\in L(L^{2}(\mathcal{G}_t;\mathbb{R}),L^{2}(\mathcal{F}_s;\mathbb{R}))$,
$\partial_{\xi}Z_s^{t,x,\cdot}=\mathcal{Q}_s^{t,x,\cdot}(=(\xi\rightarrow\mathcal{Q}_s^{t,x,\xi}))\in L(L^{2}(\mathcal{G}_t;\mathbb{R}),\mathcal{H}_{\mathcal{F}}^2(t,T;\mathbb{R}))$
 are continuous.
\end{lemma}

\begin{proof}
We only prove that $\partial_{\xi}Y_s^{t,x,\xi}\!=\!\mathcal{O}_s^{t,x,\xi}$, $s\!\in\![t,T]$
is continuous with respect to $\xi$. The continuity of $\partial_{\xi}Z_s^{t,x,\cdot}=\mathcal{Q}_s^{t,x,\cdot}$ can be proved with a similar
argument. From \eqref{eq6.15} as well as \eqref{eq6.16} we have
\begin{equation*}
 \begin{split}
 &\ |\partial_{\xi}Y_s^{t,x,P_{\xi}}-\partial_{\xi}Y_s^{t,x,P_{\xi'}}|^2_{L(L^{2}(\mathcal{G}_t;\mathbb{R}),L^{2}(\mathcal{F}_s;\mathbb{R}))}
=   \sup_{\eta\in L^{2}(\mathcal{G}_t;\mathbb{R}),\ |\eta|_{L^2}\leq 1}E[|\partial_{\xi}Y_s^{t,x,P_{\xi}}(\eta)-\partial_{\xi}Y_s^{t,x,P_{\xi'}}(\eta)|^2]\\
= &\  \sup_{\eta\in L^{2}(\mathcal{G}_t;\mathbb{R}),\ |\eta|_{L^2}\leq 1}
     E[|\overline{E}[(O_s^{t,x,P_{\xi}}(\overline{\xi})-O_s^{t,x,P_{\xi'}}(\overline{\xi}'))\cdot\overline{\eta}]|^2]\\
\leq&\  \sup_{\eta\in L^{2}(\mathcal{G}_t;\mathbb{R}),\ |\eta|_{L^2}\leq 1}
     E\[\overline{E}[|O_s^{t,x,P_{\xi}}(\overline{\xi})-O_s^{t,x,P_{\xi'}}(\overline{\xi}')|^2]\cdot\overline{E}[|\overline{\eta}|^2]\]\\
\leq&\  E\[\overline{E}[|O_s^{t,x,P_{\xi}}(\overline{\xi})-O_s^{t,x,P_{\xi'}}(\overline{\xi}')|^2]\]
=    \overline{E}\[E[|O_s^{t,x,P_{\xi}}(y)-O_s^{t,x,P_{\xi'}}(y')|^2]\big|_{y=\overline{\xi},\ y'=\overline{\xi}'}\] \\
\leq&\ CM^2\(\overline{E}[|\overline{\xi}-\overline{\xi}'|^2]+W_2(P_{\xi},P_{\xi'})^2\)
         +E[\rho_{M}(t,x,\xi,P_{\xi})]+E[\rho_{M}(t,\xi,P_{\xi})]+E[\rho_{M}(t,\xi,\xi,P_{\xi})]\\
\leq&\ CM^2E[|\xi-\xi'|^2]+E[\rho_{M}(t,x,\xi,P_{\xi})]+E[\rho_{M}(t,\xi,P_{\xi})]+E[\rho_{M}(t,\xi,\xi,P_{\xi})],\ t\leq s \leq T,
 \end{split}
\end{equation*}
for all $M\geq1$, where $E[\rho_{M}(t,x,\xi,P_{\xi})]\rightarrow0$, $E[\rho_{M}(t,\xi,P_{\xi})]\rightarrow0$,
$E[\rho_{M}(t,\xi,\xi,P_{\xi})]\rightarrow0$, as $M\rightarrow\infty$.
Consequently,
\begin{equation*}
 \begin{split}
 &\ \limsup_{\xi'\xrightarrow[]{L^2(\mathcal{G}_t;\mathbb{R})}\xi} |\partial_{\xi}Y_s^{t,x,P_{\xi}}-\partial_{\xi}Y_s^{t,x,P_{\xi'}}|^2_{L(L^{2}(\mathcal{G}_t;\mathbb{R}),L^{2}(\mathcal{F}_s;\mathbb{R}))}\\
\leq&\ E[\rho_{M}(t,x,\xi,P_{\xi})]+E[\rho_{M}(t,\xi,P_{\xi})]+E[\rho_{M}(t,\xi,\xi,P_{\xi})]\rightarrow 0,\ \mbox{as}\ M\rightarrow\infty.
 \end{split}
\end{equation*}
The proof is complete.
\end{proof}

So far, combining the Lemmas \ref{le6.1}, \ref{le6.2} and \ref{le6.3}, Theorem \ref{th6.2} has been proved. As shown in
Section 5, $(O^{t,x,P_{\xi}},Q^{t,x,P_{\xi}})$ is the derivative of $(Y^{t,x,P_{\xi}},Z^{t,x,P_{\xi}})$ with respect
to the measure $P_{\xi}$, i.e., $\partial_{\mu}Y_s^{t,x,P_{\xi}}(y):=O_s^{t,x,P_{\xi}}(y)$,
$\partial_{\mu}Z_s^{t,x,P_{\xi}}(y):=Q_s^{t,x,P_{\xi}}(y)$.
As a direct result of \eqref{eq6.38+1} and Proposition \ref{prop6.2} we have

\begin{proposition} \label{prop6.3}
For all $p\in[2,p_0]$, under the Assumptions (H4.2), (H5.1) and (H6.1),
there is a constant $C_p>0$ only depending on the Lipschitz constants of the coefficients,
such that, for all $t\in[0,T]$, $x,x',y,y'\in\mathbb{R}^d$, $\xi,\xi'\in L^{2}(\mathcal{G}_t;\mathbb{R}^d)$,
$P$-a.s.,
\begin{equation}\label{eq6.38}
  \begin{split}
\emph{(i)}&\ E\[\sup_{s\in[t,T]}|\partial_{\mu}Y_s^{t,x,P_{\xi}}(y)|^p+(\int_t^T|\partial_{\mu}Z_s^{t,x,P_{\xi}}(y)|^2ds)^{\frac{p}{2}}\]\leq C_p,\\
\emph{(ii)}&\ E\[\sup_{s\in[t,T]}|\partial_{\mu}Y_s^{t,x,P_{\xi}}(y)-\partial_{\mu}Y_s^{t,x',P_{\xi'}}(y')|^p
   +(\int_t^T|\partial_{\mu}Z_s^{t,x,P_{\xi}}(y))-\partial_{\mu}Z_s^{t,x',P_{\xi'}}(y')|^2ds)^{\frac{p}{2}}\]\\
&\ \leq C_pM^p\(|x-x'|^p+|y-y'|^p+W_2(P_{\xi},P_{\xi'})^p\)\\
  &\hspace{20pt}+\rho_{M,p}(t,x,y,P_{\xi})+\rho_{M,p}(t,y,P_{\xi})+E[\rho_{M,p}(t,\xi,y,P_{\xi})],
  \end{split}
\end{equation}
with $M\geq1$, $\rho_{M,p}(t,x,y,P_{\xi})\rightarrow0$, $\rho_{M,p}(t,y,P_{\xi})\rightarrow0$,
$E[\rho_{M,p}(t,\xi,y,P_{\xi})]\rightarrow0$, as $M\rightarrow\infty$.
\end{proposition}

\begin{remark} \label{re6.3}
In analogy to the proof in Proposition \ref{prop4.0}, it can easily be checked that
under the Assumptions (H5.1) and (H6.1), there exists a constant $C>0$ only depending on the Lipschitz
constants of the coefficients, such that for all $t\in[0,T]$, $x,x',y,y'\in\mathbb{R}^d$, $\xi,\xi'\in L^{2}(\mathcal{G}_t;\mathbb{R}^d)$,
 $P$-a.s.,
\begin{equation}\nonumber
  \begin{split}
\emph{(i)}&\ E\[\sup_{s\in[t,T]}|\partial_{\mu}Y_s^{t,x,P_{\xi}}(y)|^2+\int_t^T|\partial_{\mu}Z_s^{t,x,P_{\xi}}(y)|^2ds\]\leq C,\\
\emph{(ii)}&\ E\[\sup_{s\in[t,T]}|\partial_{\mu}Y_s^{t,x,P_{\xi}}(y)-\partial_{\mu}Y_s^{t,x',P_{\xi'}}(y')|^
   +\int_t^T|\partial_{\mu}Z_s^{t,x,P_{\xi}}(y))-\partial_{\mu}Z_s^{t,x',P_{\xi'}}(y')|^2ds\]\\
&\ \leq CM^2\(|x\!-\!x'|^2\!+\!|y\!-\!y'|^2\!+\!W_2(P_{\xi},P_{\xi'})^2\)
\!+\!\rho_{M}(t,x,y,P_{\xi})\!+\!\rho_{M}(t,y,P_{\xi})\!+\!E[\rho_{M}(t,\xi,y,P_{\xi})],
  \end{split}
\end{equation}
with $M\geq1$, $\rho_{M}(t,x,y,P_{\xi})\rightarrow0$, $\rho_{M}(t,y,P_{\xi})\rightarrow0$,
$E[\rho_{M}(t,\xi,y,P_{\xi})]\rightarrow0$, as $M\rightarrow\infty$.
\end{remark}

\section{Second order derivatives of $X^{t,x,P_{\xi}}$}

In this section we investigate the second order derivatives of $X^{t,x,P_{\xi}}$. For this we use the
following hypothesis.

\noindent\textbf{Assumption (H7.1)} The couple of coefficients $(b,\sigma)$ belongs to
 $C_b^{2,2}(\mathbb{R}^{d}\times\mathcal{P}_{2}(\mathbb{R}^{d});\mathbb{R}^{d}\times\mathbb{R}^{d\times d})$, that is,
$(b,\sigma)\in C_b^{1,1}(\mathbb{R}^{d}\times\mathcal{P}_{2}(\mathbb{R}^{d});\mathbb{R}^{d}\times\mathbb{R}^{d\times d})$
[see Assumption (H5.1)] and the derivatives of the components $b_j$, $\sigma_{i,j}$, $1\leq i,j\leq d$, have the following properties:\\
\indent (i) $\partial_{x_k}b_i(\cdot,\mu)$, $\partial_{x_k}\sigma_{i,j}(\cdot,\mu)\in
C_b^{1}(\mathbb{R}^{d})$, for all $\mu\in\mathcal{P}_{2}(\mathbb{R}^{d})$, $1\leq k\leq d$;\\
\indent (ii) $\partial_{\mu}b_j(x,\mu,\cdot)$, $\partial_{\mu}\sigma_{i,j}(x,\mu,\cdot)\in
C_b^{1}(\mathbb{R}^{d}\rightarrow\mathbb{R}^{d})$, for all $x\in\mathbb{R}^{d}$, $\mu\in\mathcal{P}_{2}(\mathbb{R}^{d})$;\\
\indent (iii) All the derivatives of $b_j$, $\sigma_{i,j}$ up to order 2 are bounded and Lipschitz.
\begin{theorem} \label{th7.1}
Under Assumption (H7.1) the first order derivatives $x\rightarrow\partial_{x_i}X^{t,x,P_{\xi}}$,
$\partial_{\mu}X^{t,x,P_{\xi}}(y)\in\mathcal{S}^2_{\mathcal{G}}(t,T;\mathbb{R}^d)$
are differentiable w.r.t. $x$ and $y$, respectively, and for
 \begin{equation*}
   M_{s,i,j}^{t,x,P_{\xi}}(y):=\big(\partial_{x_ix_j}^2X_s^{t,x,P_{\xi}},\partial_{y_i}(\partial_{\mu}X_s^{t,x,P_{\xi}}(y))\big),\ 1\leq i,j\leq d,
 \end{equation*}
we have that, for all $p\geq 2$, there exists a constant $C_p\in\mathbb{R}_{+}$ such that, for all
$t\in[0,T]$, $x,x',y,y'\in\mathbb{R}^d$, $\xi,\xi'\in L^{2}(\mathcal{G}_t;\mathbb{R}^d)$, $1\leq i,j\leq d$,
\begin{equation}\label{eq7.1}
   \begin{split}
 &\ \emph{(i)}\ E\[\sup_{s\in[t,T]}| M_{s,i,j}^{t,x,P_{\xi}}(y)|^p\]\leq C_p;\\
  &\ \emph{(ii)}\ E\[\sup_{s\in[t,T]}| M_{s,i,j}^{t,x,P_{\xi}}(y)- M_{s,i,j}^{t,x',P_{\xi'}}(y')|^p\]
 \leq C_p\(|x-x'|^p+|y-y'|^p+W_2(P_{\xi},P_{\xi'})^p\);\\
&\ \emph{(iii)}\ E\[\sup_{s\in[t,t+h]}| M_{s,i,j}^{t,x,P_{\xi}}(y)|^p\]\leq C_ph^{\frac{p}{2}},\ 0\leq t\leq t+h\leq T.
  \end{split}
\end{equation}
\end{theorem}
For the proof we refer to Section 5 in \cite{BLPR2017}, or Proposition 5.1 in \cite{HL2016}.

We now recall that the first order derivative $\partial_{x}X^{t,x,P_{\xi}}$ is differentiable in Malliavin's sense and that the derivative is a solution of a linear SDE.
\begin{proposition} \label{prop7.1}
Let Assumption (H7.1) hold true.
Then for all $(t,x)\in[0,T]\times \mathbb{R}^{d}$, $\xi\in L^{2}(\mathcal{G}_t;\mathbb{R}^d)$, $s\in[t,T]$, $1\leq k\leq d$,
the first order derivative $\partial_{x_k}X^{t,x,P_{\xi}}\in(\mathbb{D}^{1,2})^d$
 and a version of $\{D_{\theta}[\partial_{x_k}X_s^{t,x,P_{\xi}}]:\theta,s\in[t,T]\}$ is given by:\\
\indent \emph{(i)} $D_{\theta}[\partial_{x_k}X_s^{t,x,P_{\xi}}]=0$, if $\theta> s$;\\
\indent \emph{(ii)} $\{D_{\theta}[\partial_{x_k}X^{t,x,P_{\xi}}]=(D_{\theta}^i[\partial_{x_k}X^{t,x,P_{\xi},j}])_{1\leq i,j\leq d}:s\in[\theta,T]\}$
is the unique solution of the linear SDE
\begin{equation}\label{eq7.1+1}
\begin{split}
& D_{\theta}[\partial_{x}X_s^{t,x,P_{\xi}}]\!=\! \partial_{x}\sigma(X_{\theta}^{t,x,P_{\xi}}\!,P_{X_{\theta}^{t,\xi}}) \partial_{x}X_{\theta}^{t,x,P_{\xi}}
\!+\!\int_{\theta}^s\!\partial_{xx}^2b(X_r^{t,x,P_{\xi}}\!,P_{X_r^{t,\xi}}) D_{\theta}[X_r^{t,x,P_{\xi}}]\partial_{x}X_r^{t,x,P_{\xi}}dr\\
&\quad\!+\!\int_{\theta}^s\!\partial_{x}b(X_r^{t,x,P_{\xi}}\!,P_{X_r^{t,\xi}}) D_{\theta}[\partial_{x}X_r^{t,x,P_{\xi}}]dr
\!+\!\! \int_{\theta}^s\!\partial_{xx}^2\sigma(X_r^{t,x,P_{\xi}}\!,P_{X_r^{t,\xi}}) D_{\theta}[X_r^{t,x,P_{\xi}}]\partial_{x}X_r^{t,x,P_{\xi}}dW_r\\
 &\quad\!+\int_{\theta}^s\partial_{x}\sigma(X_r^{t,x,P_{\xi}},P_{X_r^{t,\xi}}) D_{\theta}[\partial_{x}X_r^{t,x,P_{\xi}}]dW_r, \ \theta\leq s\leq T.
\end{split}
\end{equation}
Furthermore, for all $p\geq 2$ there exists a constant $C_p>0$ only depending on the Lipschitz
constants of the coefficients $b$, $\sigma$ and their first and second order derivatives, such that, for all $t\in[0,T]$, $x,x'\in\mathbb{R}^d$,
$\xi,\xi'\in L^{2}(\mathcal{G}_t;\mathbb{R}^d)$, $1\leq k\leq d$, $P$-a.s.,
\begin{equation}\label{eq7.3}
  \begin{split}
  &\ \emph{(i)}\ E\[\sup_{s\in[t,T]}|D_{\theta}[\partial_{x_k}X_s^{t,x,P_{\xi}}]|^p\]\leq C_p;\\
  &\ \emph{(ii)}\ E\[\sup_{s\in[t,T]}|D_{\theta}[\partial_{x_k}X_s^{t,x,P_{\xi}}]-D_{\theta}[\partial_{x_k}X_s^{t,x',P_{\xi'}}]|^p\]
       \leq C_p\(|x-x'|^p+W_2(P_{\xi},P_{\xi'})^p\).
  \end{split}
\end{equation}
\end{proposition}
\begin{proof}
It is standard to prove that $\partial_{x}X_s^{t,x,P_{\xi}}$ has a Malliavin derivative $D_{\theta}[\partial_{x}X_s^{t,x,P_{\xi}}]$ under our assumptions,
and $D_{\theta}[\partial_{x}X_s^{t,x,P_{\xi}}]$ satisfies SDE \eqref{eq7.1+1}.
So we omit proving this here. Moreover, with the help of Lemma \ref{le3.1}, the Propositions \ref{prop3.1} and \ref{prop5.1}
as well as standard estimates for classical SDEs
it can easily be verified that \eqref{eq7.3} holds true.
\end{proof}

\section{Second order derivatives of $(Y^{t,x,P_{\xi}},Z^{t,x,P_{\xi}})$}
This section is devoted to the study of second order derivatives of $(Y^{t,x,P_{\xi}},Z^{t,x,P_{\xi}})$.

\noindent\textbf{Assumption (H8.1)} Let $\Phi\in C_b^{2,2}(\mathbb{R}^{d}\times\mathcal{P}_{2}(\mathbb{R}^{d}))$,
$f\in C_b^{2,2}(\mathbb{R}^{d+1+d}\times\mathcal{P}_{2}(\mathbb{R}^{d+1+d}))$
and $h\!\in\! C_b^{2}(\mathcal{P}_{2}(\mathbb{R}^{d+1+d});\mathbb{R}^{l})$.

\medskip
First we recall that under natural conditions the first order derivative $(\partial_xY^{t,x,P_{\xi}},\partial_xZ^{t,x,P_{\xi}})$
is differentiable w.r.t. the Brownian motion $W$ in Malliavin's sense and that the derivative is a solution of a linear BDSDE.
For this we need an additional assumption on $g$.

\noindent\textbf{Assumption (H8.2)} The coefficient $g$ is affine in $z$: For all $x\in\mathbb{R}^{d}$, $y\in\mathbb{R}$, $z\in\mathbb{R}^{d}$, $\mu\in\mathcal{P}_{2}(\mathbb{R}^{d+1}\times\mathbb{R}^d)$,
$$g(x,y,z,\mu)=g^1(x,y,\mu)+g^2(\mu(\cdot\times\mathbb{R}\times\mathbb{R}^{d}))z,$$
where
$g^1\in C_b^{2,2}(\mathbb{R}^{d+1}\times\mathcal{P}_{2}(\mathbb{R}^{d+1+d});\mathbb{R}^{l})$
and $g^2\in C_b^{2}(\mathcal{P}_{2}(\mathbb{R}^{d});\mathbb{R}^{l})$.
In addition we suppose $|g^2|^2\leq \alpha_1$,
$\sum_{k=1}^{d}\sum_{i=1}^l|(\partial_{\mu}g^1_{i})_{d+1+k}|^2\leq \alpha_2$,
where constants $\alpha_1,\alpha_2>0$ with $0\!<\!\alpha_1\!+\!\alpha_2\!<\!1$.

\begin{remark}\label{re8.1}
In order to understand better why we need Assumption (H8.2), consider for $d=l=1$ and
for the functions $\Phi,g\in C^2_b(\mathbb{R})$ with $|\partial_zg|^2\leq\alpha_1$, the BDSDE
\begin{equation*}
   \begin{split}
   Y_s^{t,x,P_{\xi}}=\Phi (X_T^{t,x,P_{\xi}})+\int_s^Tg(Z_r^{t,x,P_{\xi}})d\overleftarrow{B_r}-\int_s^TZ_r^{t,x,P_{\xi}}dW_r,\ s\in[t,T].
   \end{split}
\end{equation*}
Then, as we have seen, $(\partial_xY^{t,x,P_{\xi}},\partial_xZ^{t,x,P_{\xi}})$ is the solution of the linear BDSDE
\begin{equation}\label{eq8.0.0}
  \begin{split}
 \partial_{x}Y_s^{t,x,P_{\xi}}\!\!=\!\! \partial_{x}\Phi (X_T^{t,x,P_{\xi}})\partial_{x}X_T^{t,x,P_{\xi}}
\! \!+\!\!\int_s^T\!\partial_{z}g(Z_r^{t,x,P_{\xi}})\partial_{x}Z_r^{t,x,P_{\xi}}d\overleftarrow{B_r}
 \!\!-\!\!\int_s^T\!\partial_{x}Z_r^{t,x,P_{\xi}}dW_r,\ s\in[t,T],
  \end{split}
\end{equation}
and the formal second derivative $(\partial_{xx}^2Y_s^{t,x,P_{\xi}},\partial_{xx}^2Z_s^{t,x,P_{\xi}})$ should solve the BDSDE
\begin{equation*}
  \begin{split}
\partial_{xx}^2&Y_s^{t,x,P_{\xi}}= \partial_{xx}^2\Phi(X_T^{t,x,P_{\xi}}) (\partial_{x}X_T^{t,x,P_{\xi}})^2
+\partial_{x}\Phi(X_T^{t,x,P_{\xi}}) \partial_{xx}^2X_T^{t,x,P_{\xi}}\\
&\ +\!\!\int_{s}^T\!\!\(\partial_{zz}^2g(Z_r^{t,x,P_{\xi}}) (\partial_{x}Z_r^{t,x,P_{\xi}})^2
   \!\!+\!\!\partial_{z}g(Z_r^{t,x,P_{\xi}}) \partial_{xx}^2Z_r^{t,x,P_{\xi}}\) d\overleftarrow{B_r}
  \!\!-\!\!\int_s^T\!\!\partial_{xx}^2Z_r^{t,x,P_{\xi}}dW_r, \  s\in[t,T].
  \end{split}
\end{equation*}
However, in order to make sense to the integral $\displaystyle \int_{s}^T\!\!\partial_{zz}^2g(Z_r^{t,x,P_{\xi}})
(\partial_{x}Z_r^{t,x,P_{\xi}})^2 d\overleftarrow{B_r}$, $s\!\in\![t,T]$,
we need $\displaystyle P($ $\int_{t}^T|\partial_{x}Z_r^{t,x,P_{\xi}}|^4dr<\infty)\!=\!1$, while equation \eqref{eq8.0.0} only
allows to conclude that $\displaystyle E[(\!\int_{t}^T\!\!\!|\partial_{x}Z_r^{t,x,P_{\xi}}|^2\!dr)^p]$\\
$<\infty$, $p\geq1$. This is the reason why we have to avoid $(\partial_{x}Z^{t,x,P_{\xi}})^2$ in the stochastic
integral driven by $\overleftarrow{B_r}$, and this is why we suppose (H8.2).
\end{remark}

\begin{proposition} \label{prop8.1}
Let the Assumptions (H4.2), (H7.1), (H8.1) and (H8.2) hold true.
Then for all $(t,x)\in[0,T]\times \mathbb{R}^{d}$, $\xi\in L^{2}(\mathcal{G}_t;\mathbb{R}^d)$, $s\in[t,T]$, $1\leq k\leq d$,
$(\partial_{x_k}Y_s^{t,x,P_{\xi}},\partial_{x_k}Z_s^{t,x,P_{\xi}})\!\in\! L^2(t,T;(\mathbb{D}^{1,2})^{1+d})$,
and a version of $\{D_{\theta}[\partial_{x_k}Y_s^{t,x,P_{\xi}}],D_{\theta}[\partial_{x_k}Z_s^{t,x,P_{\xi}}],\theta,s\in[t,T]\}$ is given by:\\
\indent \emph{(i)} $D_{\theta}[\partial_{x}Y_s^{t,x,P_{\xi}}]=0$, $D_{\theta}[\partial_{x}Z_s^{t,x,P_{\xi}}]=0$, $t\leq s<\theta\leq T$;\\
\indent \emph{(ii)} $\{D_{\theta}[\partial_{x}Y^{t,x,P_{\xi}}]=(D_{\theta}^i[\partial_{x}Y^{t,x,P_{\xi}}])_{1\leq i\leq d},
D_{\theta}[\partial_{x}Z^{t,x,P_{\xi}}]=(D_{\theta}^i[\partial_{x}Z^{t,x,P_{\xi},j}])_{1\leq i,j\leq d},s\in[\theta,T]\}$
is the unique solution of the linear BDSDE (Here for simplicity written for $d = 1$, $l=1$
and $f(\Pi_s^{t,x,\xi},P_{\Pi_s^{t,\xi}})=f(Z_s^{t,x,P_{\xi}})$,
$g(\Pi_s^{t,x,\xi},P_{\Pi_s^{t,\xi}})=g^2(P_{X_s^{t,\xi}})Z_s^{t,x,P_{\xi}}$, $h(P_{\Pi_s^{t,\xi}})=h(P_{(Y_s^{t,\xi},Z_s^{t,\xi})})$,
 $\Phi(X_T^{t,x,P_{\xi}},$\\
 $P_{X_T^{t,\xi}})=\Phi(X_T^{t,x,P_{\xi}})$)
\begin{equation}\label{eq8.0+2}
\begin{split}
& D_{\theta}[\partial_{x}Y_s^{t,x,P_{\xi}}]= \partial_{xx}^2\Phi(X_T^{t,x,P_{\xi}}) D_{\theta}[X_T^{t,x,P_{\xi}}]\partial_{x}X_T^{t,x,P_{\xi}}
+\partial_{x}\Phi(X_T^{t,x,P_{\xi}}) D_{\theta}[\partial_{x}X_T^{t,x,P_{\xi}}]\\
&\qquad+\int_{s}^T\(\partial_{zz}^2f(Z_r^{t,x,P_{\xi}}) D_{\theta}[Z_r^{t,x,P_{\xi}}]\partial_{x}Z_r^{t,x,P_{\xi}}
   +\partial_{z}f(Z_r^{t,x,P_{\xi}}) D_{\theta}[\partial_{x}Z_r^{t,x,P_{\xi}}]\)dr\\
 &\qquad+\int_{s}^Tg^2(P_{X_r^{t,\xi}}) D_{\theta}[\partial_{x}Z_r^{t,x,P_{\xi}}]d\overleftarrow{B_r}
 -\int_s^TD_{\theta}[\partial_{x}Z_r^{t,x,P_{\xi}}]dW_r, \ \theta\leq s\leq T.
\end{split}
\end{equation}
Moreover, $\partial_{x}Z_s^{t,x,P_{\xi}}=P\mbox{-}\displaystyle{\lim_{s<u\downarrow s}}D_{s}[\partial_{x}Y_u^{t,x,P_{\xi}}]$, $d s dP$-a.e.
Furthermore, there exists a constant $C_p>0$ only depending on the Lipschitz constants of the coefficients,
 such that for all $t\in[0,T]$, $x,x'\in\mathbb{R}^d$,
$\xi,\xi'\in L^{2}(\mathcal{G}_t;\mathbb{R}^d)$, $1\leq k\leq d$, $P$-a.s.,
\begin{equation}\label{eq8.0}
  \begin{split}
    &\ \mbox{\emph{(i)}}\ E\[\sup_{s\in[t,T]}|D_{\theta}[\partial_{x_k}Y_s^{t,x,P_{\xi}}]|^p+(\int_t^T|[D_{\theta}\partial_{x_k}Z_s^{t,x,P_{\xi}}]|^2ds)^{\frac{p}{2}}\]\leq C_p,\
   \mbox{for all}\ p\in[2,\frac{p_0}{2}] ;\\
    &\ \mbox{\emph{(ii)}}\ E\[\!\sup_{s\in[t,T]}\!|D_{\theta}[\partial_{x_k}Y_s^{t,x,P_{\xi}}]\!-\!D_{\theta}[\partial_{x_k}Y_s^{t,x',P_{\xi'}}]|^p
  \!+\!(\!\int_t^T\!|D_{\theta}[\partial_{x_k}Z_s^{t,x,P_{\xi}}]\!-\!D_{\theta}[\partial_{x_k}Z_s^{t,x',P_{\xi'}}]|^2ds)^{\frac{p}{2}}\!\]\\
 &\ \hspace{30pt} \leq C_p\(|x-x'|^p+W_2(P_{\xi},P_{\xi'})^p\),\ \mbox{for all}\ p\in[2,\frac{p_0}{8}].
  \end{split}
\end{equation}
In particular, there exists a constant $C_p>0$ only depending on the Lipschitz constants of the coefficients,
 such that for all $x,x'\in\mathbb{R}^d$,
$\xi,\xi'\in L^{2}(\mathcal{G}_t;\mathbb{R}^d)$, $d s dP$-a.e., $s\in[t,T]$,
\begin{equation}\label{eq8.0+1}
  \begin{split}
    &\ \mbox{\emph{(i)}}\ E[|\partial_{x_k}Z_s^{t,x,P_{\xi}}|^p]\leq C_p,\ \mbox{for all}\ p\in[2,\frac{p_0}{2}];\\
  &\ \mbox{\emph{(ii)}}\ E[|\partial_{x_k}Z_s^{t,x,P_{\xi}}-\partial_{x_k}Z_s^{t,x',P_{\xi'}}|^p] \leq C_p\(|x-x'|^p+W_2(P_{\xi},P_{\xi'})^p\),\
  \mbox{for all}\ p\in[2,\frac{p_0}{8}].
  \end{split}
\end{equation}
\end{proposition}
\begin{proof}
Here for simplicity written for $d = 1$, $l=1$
and $f(\Pi_s^{t,x,\xi},P_{\Pi_s^{t,\xi}})=f(Z_s^{t,x,P_{\xi}})$,
$g(\Pi_s^{t,x,\xi},P_{\Pi_s^{t,\xi}})=g^2(P_{X_s^{t,\xi}})Z_s^{t,x,P_{\xi}}$, $h(P_{\Pi_s^{t,\xi}})=h(P_{(Y_s^{t,\xi},Z_s^{t,\xi})})$
and $\Phi(X_T^{t,x,P_{\xi}},P_{X_T^{t,\xi}})=\Phi(X_T^{t,x,P_{\xi}})$.

It is standard to prove that $\partial_{x}Y_s^{t,x,P_{\xi}}$ and $\partial_{x}Z_s^{t,x,P_{\xi}}$ are Malliavin differentiable under our assumptions,
and that they satisfy BDSDE \eqref{eq8.0+2}, and
\begin{equation}\label{eq8.0+2+1}
  \partial_{x}Z_r^{t,x,P_{\xi}}=P\mbox{-}\displaystyle{\lim_{r<s\downarrow r}}D_r[\partial_{x}Y_s^{t,x,P_{\xi}}],
\end{equation}
so we omit proving this here. Thus, it suffices to prove \eqref{eq8.0}.
 For this, note that, from Lemma \ref{leA.1}-(2) and the Propositions \ref{prop3.1}, \ref{prop5.1}, \ref{prop6.1+1} and \ref{prop6.1+2},
 as well as from Theorem \ref{thA.3} and the H\"{o}lder inequality, we get
 $\displaystyle E\[\sup_{s\in[t,T]}|D_{\theta}[\partial_{x}Y_s^{t,x,P_{\xi}}]|^{\frac{p_0}{2}}+(\int_t^T|[D_{\theta}\partial_{x}Z_s^{t,x,P_{\xi}}]|^2ds)^{\frac{p_0}{4}}\]\leq C_{\frac{p_0}{2}}.$ Again from H\"{o}lder's inequality, for all $p\in[2,\frac{p_0}{2}]$ we have
\begin{equation}\label{eq8.0+1+1+1}
  \begin{split}
   E\[\sup_{s\in[t,T]}|D_{\theta}[\partial_{x}Y_s^{t,x,P_{\xi}}]|^p+(\int_t^T|[D_{\theta}\partial_{x}Z_s^{t,x,P_{\xi}}]|^2ds)^{\frac{p}{2}}\]\leq C_p.
  \end{split}
\end{equation}
 Moreover, from \eqref{eq8.0+1+1+1} and \eqref{eq8.0+2+1} we obtain \eqref{eq8.0+1}-(i). In analogy to estimate \eqref{eq8.0+1}-(i), we also show that
\begin{equation}\label{eq8.0+1+1}
  \begin{split}
 E[|D_{\theta}[Z_s^{t,x,P_{\xi}}]|^p]\leq C_p,\ s\in[\theta,T],\ \theta\leq T,\ x\in\mathbb{R},\
 \xi\in L^{2}(\mathcal{G}_t;\mathbb{R}),\ p\in[2,\frac{p_0}{2}] .
  \end{split}
\end{equation}
Now we prove estimate \eqref{eq8.0}-(ii).
For all $0\leq t\leq T$, $t\leq\theta\leq s\leq T$, $x,x'\in\mathbb{R}$, $\xi,\xi'\in L^{2}(\mathcal{G}_t;\mathbb{R})$,
we get from \eqref{eq8.0+2} the following BDSDE:
\begin{equation}\label{eq8.0+3}
  \begin{split}
D_{\theta}&[\partial_{x}Y_s^{t,x,P_{\xi}}]-D_{\theta}[\partial_{x}Y_s^{t,x',P_{\xi'}}]
 =\int_s^T\partial_{z}f(Z_r^{t,x,P_{\xi}})\(D_{\theta}[\partial_{x}Z_r^{t,x,P_{\xi}}]-D_{\theta}[\partial_{x}Z_r^{t,x',P_{\xi'}}]\)dr\\
 & \!+\!\int_s^T\!\!\!g^2(P_{X_r^{t,\xi}})\(\!D_{\theta}[\partial_{x}Z_r^{t,x,P_{\xi}}]\!-\!D_{\theta}[\partial_{x}Z_r^{t,x',P_{\xi'}}]\!\)d\overleftarrow{B_r}\!
  \!-\!\!\int_s^T\!\!\!\(\!D_{\theta}[\partial_{x}Z_r^{t,x,P_{\xi}}]\!-\!D_{\theta}[\partial_{x}Z_r^{t,x',P_{\xi'}}]\!\)dW_r\\
& +I(t,x,P_{\xi},x',P_{\xi'})+\int_s^TR(r,x,P_{\xi},x',P_{\xi'})dr+\int_s^TH(r,x,P_{\xi},x',P_{\xi'})d\overleftarrow{B_r},
  \end{split}
\end{equation}
where
\begin{equation*}
  \begin{split}
 & I(t,x,P_{\xi},x',P_{\xi'}):= \(\partial_{xx}^2\Phi(X_T^{t,x,P_{\xi}}) D_{\theta}[X_T^{t,x,P_{\xi}}]\partial_{x}X_T^{t,x,P_{\xi}}
+\partial_{x}\Phi(X_T^{t,x,P_{\xi}}) D_{\theta}[\partial_{x}X_T^{t,x,P_{\xi}}]\)\\
&\hspace{100pt}-\(\partial_{xx}^2\Phi(X_T^{t,x',P_{\xi'}}) D_{\theta}[X_T^{t,x',P_{\xi'}}]\partial_{x}X_T^{t,x',P_{\xi'}}
+\partial_{x}\Phi(X_T^{t,x',P_{\xi'}}) D_{\theta}[\partial_{x}X_T^{t,x',P_{\xi'}}]\),\\
 & R(r,x,P_{\xi},x',P_{\xi'}):=\partial_{zz}^2f(Z_r^{t,x,P_{\xi}}) D_{\theta}[Z_r^{t,x,P_{\xi}}]\partial_{x}Z_r^{t,x,P_{\xi}}
 -\partial_{zz}^2f(Z_r^{t,x',P_{\xi'}}) D_{\theta}[Z_r^{t,x',P_{\xi'}}]\partial_{x}Z_r^{t,x',P_{\xi'}}\\
 &\hspace{100pt}+\(\partial_{z}f(Z_r^{t,x,P_{\xi}})-\partial_{z}f(Z_r^{t,x',P_{\xi'}})\)D_{\theta}[\partial_{x}Z_r^{t,x',P_{\xi'}}],\\
 & H(r,x,P_{\xi},x',P_{\xi'}):=\(g^2(P_{X_r^{t,\xi}})-g^2(P_{X_r^{t,\xi'}})\)D_{\theta}[\partial_{x}Z_r^{t,x',P_{\xi'}}].
  \end{split}
\end{equation*}
From \eqref{eq8.0}-(i), \eqref{eq8.0+1+1}, Lemmas \ref{leA.1}, Remark \ref{re4.1} and \ref{le3.1}, and the Propositions \ref{prop3.1},
\ref{prop5.1}, \ref{prop6.1+1}, \ref{prop6.1+2} and \ref{prop7.1} as well as our assumptions, it follows that
\begin{equation}\nonumber
  \begin{split}
  &\  E\[\sup_{s\in[\theta,T]}|D_{\theta}[\partial_{x}Y_s^{t,x,P_{\xi}}]-D_{\theta}[\partial_{x}Y_s^{t,x',P_{\xi'}}]|^{\frac{p_0}{8}}
  +(\int_{\theta}^T|D_{\theta}[\partial_{x}Z_s^{t,x,P_{\xi}}]-D_{\theta}[\partial_{x}Z_s^{t,x',P_{\xi'}}]|^2ds)^{\frac{p_0}{16}}\]\\
 \leq&\  C_{\frac{p_0}{8}} E\[|I(t,x,P_{\xi},x',P_{\xi'})|^{\frac{p_0}{8}}\!+\!(\int_{\theta}^{T}\!|R(r,x,P_{\xi},x',P_{\xi'})|dr)^{\frac{p_0}{8}}\!
     +\!(\int_{\theta}^{T}\!|H(r,x,P_{\xi},x',P_{\xi'})|^2dr)^{\frac{p_0}{16}}\]\\
 \leq&\ \! C_{\frac{p_0}{8}}\!\!\int_{\theta}^T\!\! (E[|(\partial_{zz}^2f)(Z_r^{t,x,P_{\xi}}\!)\!-\!(\partial_{zz}^2f)(Z_r^{t,x',P_{\xi'}}\!)|^{\frac{p_0}{4}}])^{\frac{1}{2}}
         (E[|\partial_{x}Z_r^{t,x',P_{\xi'}}|^{\frac{p_0}{2}}])^{\frac{1}{4}} (E[|D_{\theta}[Z_r^{t,x',P_{\xi'}}]|^{\frac{p_0}{2}}])^{\frac{1}{4}}dr\\
&\ +C_{\frac{p_0}{8}} \(E\[\(\int_{\theta}^T|\partial_{x}Z_r^{t,x,P_{\xi}}-\partial_{x}Z_r^{t,x',P_{\xi'}}|^2dr\)^{\frac{p_0}{8}}\]\)^{\frac{1}{2}}
    \(E\[\(\int_{\theta}^T|D_{\theta}[Z_r^{t,x',P_{\xi'}}]|^2dr\)^{\frac{p_0}{8}}\]\)^{\frac{1}{2}}\\
&\ +C_{\frac{p_0}{8}} \(E\[\(\int_{\theta}^T|D_{\theta}[Z_r^{t,x,P_{\xi}}]-D_{\theta}[Z_r^{t,x',P_{\xi'}}]|^2dr\)^{\frac{p_0}{8}}\]\)^{\frac{1}{2}}
    \(E\[\(\int_{\theta}^T|\partial_{x}Z_r^{t,x,P_{\xi}}|^2dr\)^{\frac{p_0}{8}}\]\)^{\frac{1}{2}}\\
&\ +C_{\frac{p_0}{8}} \!\(\!E\[\!\(\int_{\theta}^T\!|(\partial_{z}f)(Z_r^{t,x,P_{\xi}})\!-\!(\partial_{z}f)(Z_r^{t,x',P_{\xi'}})|^2dr\)^{\frac{p_0}{8}}\]\)^{\frac{1}{2}}
    \!\(\!E\[\!\(\!\int_{\theta}^T\!|D_{\theta}[\partial_{x}Z_r^{t,x',P_{\xi'}}]|^2dr\!\)^{\frac{p_0}{8}}\]\)^{\frac{1}{2}}\\
  \end{split}
\end{equation}
\begin{equation}\label{eq8.0+4}
  \begin{split}
 &\ +C_{\frac{p_0}{8}} E\[\(\int_{\theta}^T|g^2(P_{X_r^{t,\xi}})\!-\!g^2(P_{X_r^{t,\xi'}})|^2
 |D_{\theta}[\partial_{x}Z_r^{t,x',P_{\xi'}}]|^2dr\)^{\frac{p_0}{16}}\]
   \!+\! C_{\frac{p_0}{8}}\(\!|x\!-\!x'|^{\frac{p_0}{8}}\!+\!W_2(P_{\xi},P_{\xi'})^{\frac{p_0}{8}}\!\)\\
\leq&\  C_{\frac{p_0}{8}}\!\int_{\theta}^T\! (E[|Z_r^{t,x,P_{\xi}}\!-\!Z_r^{t,x',P_{\xi'}}|^{\frac{p_0}{4}}])^{\frac{1}{2}}dr
       +C_{\frac{p_0}{8}} \(E\[\(\int_{\theta}^T|\partial_{x}Z_r^{t,x,P_{\xi}}-\partial_{x}Z_r^{t,x',P_{\xi'}}|^2dr\)^{\frac{p_0}{8}}\]\)^{\frac{1}{2}}\\
&\ +\!C_{\frac{p_0}{8}} \!\(\!E\[\!\(\!\int_{\theta}^T\!|D_{\theta}[Z_r^{t,x,P_{\xi}}]\!-\!D_{\theta}[Z_r^{t,x',P_{\xi'}}]|^2dr\!\)^{\frac{p_0}{8}}\]\)^{\frac{1}{2}}
     \!+\!C_{\frac{p_0}{8}} \!\(\!E\[\!\(\int_{\theta}^T\!|Z_r^{t,x,P_{\xi}}\!-\!Z_r^{t,x',P_{\xi'}}|^2dr\)^{\frac{p_0}{8}}\]\)^{\frac{1}{2}}\\
 &\ +C_{\frac{p_0}{8}} \sup_{r\in[\theta,T]}W_2(P_{X_r^{t,\xi}},P_{X_r^{t,\xi'}})^{\frac{p_0}{8}}
      E\[\!\(\!\int_{\theta}^T \!|D_{\theta}[\partial_{x}Z_r^{t,x',P_{\xi'}}]|^2dr\)^{\frac{p_0}{16}}\]
   \!+\! C_p\(\!|x\!-\!x'|^{\frac{p_0}{8}}\!+\!W_2(P_{\xi},P_{\xi'})^{\frac{p_0}{8}}\!\)\\
\leq &\ C_{\frac{p_0}{8}}\(|x-x'|^{\frac{p_0}{8}}+W_2(P_{\xi},P_{\xi'})^{\frac{p_0}{8}}\).
  \end{split}
\end{equation}
Then from H\"{o}lder's inequality, for all $p\in[2,\frac{p_0}{8}]$ we have
\begin{equation*}
  \begin{split}
    &\  E\[\!\sup_{s\in[t,T]}\!|D_{\theta}[\partial_{x_k}Y_s^{t,x,P_{\xi}}]\!-\!D_{\theta}[\partial_{x_k}Y_s^{t,x',P_{\xi'}}]|^p
  \!+\!(\!\int_t^T\!|D_{\theta}[\partial_{x_k}Z_s^{t,x,P_{\xi}}]\!-\!D_{\theta}[\partial_{x_k}Z_s^{t,x',P_{\xi'}}]|^2ds)^{\frac{p}{2}}\!\]\\
 \leq&\ C_p\(|x-x'|^p+W_2(P_{\xi},P_{\xi'})^p\).
  \end{split}
\end{equation*}
Moreover, statement \eqref{eq8.0+1}-(ii) follows immediately from the above estimate and \eqref{eq8.0+2+1}.
\end{proof}
Now we investigate the $L^2$-derivative of $(\partial_xY_s^{t,x,P_{\xi}},\partial_xZ_s^{t,x,P_{\xi}})$ with
respect to $x$. We notice that, similar to Theorem \ref{th6.1}, with standard arguments we obtain the following Theorem.
\begin{theorem} \label{th8.1}
Suppose the Assumptions (H4.2), (H7.1), (H8.1) and (H8.2) hold true.
Then, for all $t\in[0,T]$, $x\in\mathbb{R}^d$, $\xi\in L^{2}(\mathcal{G}_t;\mathbb{R}^d)$,
the mapping $(\partial_xY_s^{t,x,P_{\xi}},\partial_xZ_s^{t,x,P_{\xi}})\in
 \mathcal{S}_{\mathcal{F}}^2(t,T;\mathbb{R}^d)\times\mathcal{H}_{\mathcal{F}}^2(t,T;\mathbb{R}^{d\times d})$
 is differentiable in $x$, and its derivative
$(\partial_{xx}^2Y_s^{t,x,P_{\xi}},\partial_{xx}^2Z_s^{t,x,P_{\xi}})$ is the unique solution of a BDSDE.
 Moreover, there is some constant $C_p>0$ only depending on the bounds and
the Lipschitz constants of the coefficients $b$, $\sigma$, $f$, $g$, $h$, $\Phi$ and their first and second order derivatives,
such that, for all $t\in[0,T]$, $x,x'\in\mathbb{R}^d$, $\xi,\xi'\in L^{2}(\mathcal{G}_t;\mathbb{R}^d)$,
\begin{equation*}
  \begin{split}
  &\ \emph{(i)}\ E\[\sup_{s\in[t,T]}|\partial_{xx}^2Y_s^{t,x,P_{\xi}}|^p+(\int_t^T|\partial_{xx}^2Z_s^{t,x,P_{\xi}}|^2ds)^{\frac{p}{2}}\]\leq C_p,\
  \mbox{for all}\ p\in[2,\frac{p_0}{2}];\\
  \end{split}
\end{equation*}
\begin{equation*}
  \begin{split}
  &\ \emph{(ii)}\ E\[\sup_{s\in[t,T]}|\partial_{xx}^2Y_s^{t,x,P_{\xi}}-\partial_{xx}^2Y_s^{t,x',P_{\xi'}}|^p
  +(\int_t^T|\partial_{xx}^2Z_s^{t,x,P_{\xi}}-\partial_{xx}^2Z_s^{t,x',P_{\xi'}}|^2ds)^{\frac{p}{2}}\]\\
 &\ \hspace{20pt} \leq C_p\(|x-x'|^p+W_2(P_{\xi},P_{\xi'})^p\),\ \mbox{for all}\ p\in[2,\frac{p_0}{8}] .
  \end{split}
\end{equation*}
\end{theorem}
\begin{proof}
Following the approach for the first-order derivatives, we restrict here ourselves
to the dimensions $d = 1$, $l=1$ and to $f(\Pi_s^{t,x,\xi},P_{\Pi_s^{t,\xi}})=f(Z_s^{t,x,P_{\xi}},P_{Z_s^{t,\xi}})$,
$g(\Pi_s^{t,x,\xi},P_{\Pi_s^{t,\xi}})=g^2(P_{X_s^{t,\xi}})Z_s^{t,x,P_{\xi}}$, $h(P_{\Pi_s^{t,\xi}})=0$
and $\Phi(X_T^{t,x,P_{\xi}},P_{X_T^{t,\xi}})=\Phi(X_T^{t,x,P_{\xi}})$.

First we notice that in analogy to Proposition \ref{prop6.1+2}, as the $L^2$-derivative of $(\partial_xY_s^{t,x,P_{\xi}},\partial_xZ_s^{t,x,P_{\xi}})$ with
respect to $x$ does not concern the law $P_{Z^{t,\xi}}$, the arguments of the proof are standard. Therefore, it can easily be checked that its derivative solves
the equation
\begin{equation*}
\begin{split}
\partial_{xx}^2Y_s^{t,x,P_{\xi}}=&\ \partial_{xx}^2\Phi(X_T^{t,x,P_{\xi}}) (\partial_{x}X_T^{t,x,P_{\xi}})^2
+\partial_{x}\Phi(X_T^{t,x,P_{\xi}}) \partial_{xx}^2X_T^{t,x,P_{\xi}}\\
&\ +\int_{s}^T\(\partial_{zz}^2f(Z_r^{t,x,P_{\xi}},P_{Z_r^{t,\xi}}) (\partial_{x}Z_r^{t,x,P_{\xi}})^2
   +\partial_{z}f(Z_r^{t,x,P_{\xi}},P_{Z_r^{t,\xi}}) \partial_{xx}^2Z_r^{t,x,P_{\xi}}\)dr\\
\end{split}
\end{equation*}
\begin{equation*}
\begin{split}
 &\ +\int_{s}^Tg^2(P_{X_r^{t,\xi}}) \partial_{xx}^2Z_r^{t,x,P_{\xi}}d\overleftarrow{B_r}
 -\int_s^T\partial_{xx}^2Z_r^{t,x,P_{\xi}}dW_r, \ t\leq s\leq T.
\end{split}
\end{equation*}
From Lemma \ref{leA.1}-(2), Theorem \ref{th7.1} and the Propositions \ref{prop5.1} and \ref{prop6.1+2}, we now get (i) directly.
Finally, using the similar techniques as the proof of Proposition \ref{prop8.1}, we obtain estimate (ii).
\end{proof}
Next, we consider the $L^2$-differentiability of
$ y\rightarrow (\partial_{\mu}Y^{t,x,P_{\xi}}(y),\partial_{\mu}Z^{t,x,P_{\xi}}(y)).$

\begin{theorem} \label{th8.2}
Suppose the Assumptions (H4.2), (H7.1), (H8.1) and (H8.2) hold true.
Then, for all $t\in[0,T]$, $x\in\mathbb{R}^d$ and $\xi\in L^{2}(\mathcal{G}_t;\mathbb{R}^d)$,
the process $(\partial_{\mu}Y^{t,x,P_{\xi}}(y),\partial_{\mu}Z^{t,x,P_{\xi}}(y))\in
 \mathcal{S}_{\mathcal{F}}^2(t,T;\mathbb{R}^d)\times\mathcal{H}_{\mathcal{F}}^2(t,T;\mathbb{R}^{d\times d})$
 is $L^2$-differentiable with respect to $y$, and its derivative
$(\partial_y\partial_{\mu}Y_s^{t,x,P_{\xi}}\!(y),\partial_y\partial_{\mu}Z_s^{t,x,P_{\xi}}$
$(y))$ is the unique solution of a BDSDE (see \eqref{eq8.2} in the proof).

 Moreover, there is some constant $C_p>0$ only depending on the bounds and
the Lipschitz constants of the coefficients $b$, $\sigma$, $f$, $g$, $h$, $\Phi$ and their first and second order derivatives,
such that, for all $t\in[0,T]$, $x,x',y,y'\in\mathbb{R}^d$, $\xi,\xi'\in L^{2}(\mathcal{G}_t;\mathbb{R}^d)$,
\begin{equation*}
  \begin{split}
  &\ \emph{(i)}\ E\[\sup_{s\in[t,T]}|\partial_y\partial_{\mu}Y_s^{t,x,P_{\xi}}(y)|^p+(\int_t^T|\partial_y\partial_{\mu}Z_s^{t,x,P_{\xi}}(y)|^2ds)^{\frac{p}{2}}\]\leq C_p,\
   \mbox{for all}\ p\in[2,p_0];\\
  &\ \emph{(ii)}\ E\[\sup_{s\in[t,T]}|\partial_y\partial_{\mu}Y_s^{t,x,P_{\xi}}(y)-\partial_y\partial_{\mu}Y_s^{t,x',P_{\xi'}}(y')|^p
  +(\int_t^T|\partial_y\partial_{\mu}Z_s^{t,x,P_{\xi}}(y)-\partial_y\partial_{\mu}Z_s^{t,x',P_{\xi'}}(y')|^2ds)^{\frac{p}{2}}\]\\
 &\ \hspace{20pt} \leq C_p\(|x-x'|^p+|y-y'|^p+W_2(P_{\xi},P_{\xi'})^p\),\  \mbox{for all}\ p\in[2,\frac{p_0}{2}].
  \end{split}
\end{equation*}
\end{theorem}
\begin{proof}
For simplicity of the computation we use again the special case considered in Proposition \ref{prop8.1} with,
in addition, $f=0$. Then from \eqref{eq6.40} it follows that, for $\xi\in L^{2}(\mathcal{G}_t;\mathbb{R})$, $x,y\in\mathbb{R}$ and $t\in[0,T]$,
$(\partial_{\mu}Y_{\cdot}^{t,x,P_{\xi}}(y),\partial_{\mu}Z_{\cdot}^{t,x,P_{\xi}}(y))$
is a solution of the following BDSDE: For $ s\in[t,T]$,
\begin{equation}\label{eq8.1}
  \begin{split}
 \partial_{\mu}&Y_s^{t,x,P_{\xi}}(y)
= \partial_{x}\Phi (X_T^{t,x,P_{\xi}})\partial_{\mu}X_{T}^{t,x,P_{\xi}}(y)
  +\int_s^Tg^2(P_{X_r^{t,\xi}})\partial_{\mu}Z_r^{t,x,P_{\xi}}(y) d\overleftarrow{B_r}\\
&\   +\int_s^T\widehat{E}\[\!(\partial_{\mu}g^2)(P_{X_r^{t,\xi}},\widehat{X}_r^{t,y,P_{\xi}}) \partial_{x}\widehat{X}_r^{t,y,P_{\xi}}
 \!+\!(\partial_{\mu}g^2)(P_{X_r^{t,\xi}},\widehat{X}_r^{t,\widehat{\xi}}) \partial_{\mu}\widehat{X}_r^{t,\widehat{\xi},P_{\xi}}(y) \]Z_r^{t,x,P_{\xi}}d\overleftarrow{B_r}\\
&\   +\int_s^T\widehat{E}\[\Big\langle(\partial_{\mu}h)(P_{(Y_r^{t,\xi},Z_r^{t,\xi})},(\widehat{Y}_r^{t,y,P_{\xi}},\widehat{Z}_r^{t,y,P_{\xi}})),
\left(   \begin{matrix}    \partial_{x}\widehat{Y}_r^{t,y,P_{\xi}}\\ \partial_{x}\widehat{Z}_r^{t,y,P_{\xi}} \end{matrix}\right)\Big\rangle\]d\overleftarrow{B_r}\\
&\ +\int_s^T\widehat{E}\[\Big\langle(\partial_{\mu}h)(P_{(Y_r^{t,\xi},Z_r^{t,\xi})},(\widehat{Y}_r^{t,\widehat{\xi}},\widehat{Z}_r^{t,\widehat{\xi}})),
 \left(   \begin{matrix}   \partial_{\mu}\widehat{Y}_r^{t,\widehat{\xi},P_{\xi}}(y)\\ \partial_{\mu}\widehat{Z}_r^{t,\widehat{\xi},P_{\xi}}(y) \end{matrix}\right) \Big\rangle\]d\overleftarrow{B_r} -\int_s^T\partial_{\mu}Z_r^{t,x,P_{\xi}}(y)dW_r.
  \end{split}
\end{equation}
 Differentiating formally the above BDSDE \eqref{eq8.1} with respect to $y$, we get the following BDSDE:
\begin{equation}\label{eq8.2}
  \begin{split}
   \partial_{y}&(\partial_{\mu}Y_s^{t,x,P_{\xi}}(y))
  \!= \! \partial_{x}\Phi (X_T^{t,x,P_{\xi}} ) \partial_{y}(\partial_{\mu}X_{T}^{t,x,P_{\xi}}(y))
  \!+\!\int_s^T\! g^2(P_{X_r^{t,\xi}})\partial_{y}(\partial_{\mu}Z_r^{t,x,P_{\xi}}(y))  d\overleftarrow{B_r}\\
&\  \!+\!\int_s^T\!\widehat{E}\[(\partial_{\mu}g^2)(P_{X_r^{t,\xi}},\widehat{X}_r^{t,y,P_{\xi}}) \partial_{xx}^2\widehat{X}_r^{t,y,P_{\xi}}
  \! +\!\partial_{y}((\partial_{\mu}g^2)(P_{X_r^{t,\xi}},\widehat{X}_r^{t,y,P_{\xi}}))
   (\partial_{x}\widehat{X}_r^{t,y,P_{\xi}})^2\]Z_r^{t,x,P_{\xi}}d\overleftarrow{B_r}\\
&\  \!+\!\!\int_s^T\!\!\widehat{E}\[(\partial_{\mu}g^2)(P_{X_r^{t,\xi}},\widehat{X}_r^{t,\widehat{\xi}}) \partial_{y}
      (\partial_{\mu}\widehat{X}_{r}^{t,\widehat{\xi}}(y)) \]Z_r^{t,x,P_{\xi}}d\overleftarrow{B_r}
 \! -\!\int_s^T\!\!\partial_{y}(\partial_{\mu}Z_r^{t,x,P_{\xi}}(y))dW_r\!+\!R_s^{t,y,P_{\xi}}\!,\ \!s\!\in\![t,T],
  \end{split}
\end{equation}
where
\begin{equation*}
  \begin{split}
  R_s^{t,y,P_{\xi}}:= &\ \int_s^T\Bigg\{\widehat{E}\[\Big\langle
  \partial_{(y,z)}(\partial_{\mu}h)(P_{(Y_r^{t,\xi},Z_r^{t,\xi})},(\widehat{Y}_r^{t,y,P_{\xi}},\widehat{Z}_r^{t,y,P_{\xi}}))\cdot
  \left(   \begin{matrix}    \partial_{x}\widehat{Y}_r^{t,y,P_{\xi}}\\ \partial_{x}\widehat{Z}_r^{t,y,P_{\xi}} \end{matrix}\right),
  \left(   \begin{matrix}    \partial_{x}\widehat{Y}_r^{t,y,P_{\xi}}\\ \partial_{x}\widehat{Z}_r^{t,y,P_{\xi}} \end{matrix}\right)\Big\rangle\]\\
&\ + \widehat{E}\[\Big\langle
(\partial_{\mu}h)(P_{(Y_r^{t,\xi},Z_r^{t,\xi})},(\widehat{Y}_r^{t,y,P_{\xi}},\widehat{Z}_r^{t,y,P_{\xi}})),
  \left(   \begin{matrix}    \partial_{xx}^2\widehat{Y}_r^{t,y,P_{\xi}}\\ \partial_{xx}^2\widehat{Z}_r^{t,y,P_{\xi}} \end{matrix}\right)\Big\rangle\]\\
&\ +\widehat{E}\[\Big\langle(\partial_{\mu}h)(P_{(Y_r^{t,\xi},Z_r^{t,\xi})},(\widehat{Y}_r^{t,\widehat{\xi}},\widehat{Z}_r^{t,\widehat{\xi}})),
\left(   \begin{matrix}    \partial_{y}(\partial_{\mu}\widehat{Y}_r^{t,\widehat{\xi},P_{\xi}}(y))\\ \partial_{y}(\partial_{\mu}\widehat{Z}_r^{t,\widehat{\xi},P_{\xi}}(y)) \end{matrix}\right)
    \Big\rangle\]\Bigg\}d\overleftarrow{B_r} .
   \end{split}
\end{equation*}
Observe that $(R_s^{t,y,P_{\xi}})$ is $(\mathcal{F}_{s,T}^B\vee\mathcal{F}^0)$-adapted.

We first prove estimate (i). For this we consider Eq. \eqref{eq8.2} with $x$ replaced by $\xi$, then from Theorem \ref{thA.1} this equation has a unique solution
$(\partial_y(\partial_{\mu}Y^{t,\xi,P_{\xi}}(y)),\partial_y(\partial_{\mu}Z^{t,\xi,P_{\xi}}(y)))
\in\mathcal{S}_{\mathcal{F}}^2(t,T;\mathbb{R})\times\mathcal{H}_{\mathcal{F}}^2(t,T;\mathbb{R})$, and, furthermore, from Theorem \ref{thA.3}
combined with Proposition \ref{prop8.1} and the Theorems \ref{th7.1} and \ref{th8.1}, we get that there is
some $C_{p_0}>0$ only depending on the bounds and the Lipschitz constants of the coefficients and
their derivatives of order $1$ and $2$ such that
$$ E\[\sup_{s\in[t,T]}|\partial_y(\partial_{\mu}Y_s^{t,\xi,P_{\xi}}(y))|^{p_0}
+(\int_t^T|\partial_y(\partial_{\mu}Z_s^{t,\xi,P_{\xi}}(y))|^2ds)^{\frac{p_0}{2}}\]\leq C_{p_0},$$
for all $t\in[0,T]$, $y\in\mathbb{R}$, $\xi\in L^2(\mathcal{G}_t;\mathbb{R}).$
Then from H\"{o}lder's inequality, for all $p\in[2,p_0]$ we have
\begin{equation}\label{eq8.3}
  E\[\sup_{s\in[t,T]}|\partial_y(\partial_{\mu}Y_s^{t,\xi,P_{\xi}}(y))|^p+(\int_t^T|\partial_y(\partial_{\mu}Z_s^{t,\xi,P_{\xi}}(y))|^2ds)^{\frac{p}{2}}\]\leq C_p.
\end{equation}
Then it follows from Theorem \ref{thA.1} that equation \eqref{eq8.2} has a unique solution
$(\partial_y(\partial_{\mu}Y^{t,x,P_{\xi}}(y)),$ $\partial_y(\partial_{\mu}Z^{t,x,P_{\xi}}(y)))
\in\mathcal{S}_{\mathcal{F}}^2(t,T;\mathbb{R})\times\mathcal{H}_{\mathcal{F}}^2(t,T;\mathbb{R})$, and making use of
Theorem \ref{thA.3} we see that
$$  E\[\sup_{s\in[t,T]}|\partial_y(\partial_{\mu}Y_s^{t,x,P_{\xi}}(y))|^{p_0}
 +(\int_t^T|\partial_y(\partial_{\mu}Z_s^{t,x,P_{\xi}}(y))|^2ds)^{\frac{p_0}{2}}\]\leq C_{p_0},$$
for all $t\in[0,T]$, $y\in\mathbb{R}$, $\xi\in L^2(\mathcal{G}_t;\mathbb{R})$. Again from H\"{o}lder's inequality, for all $p\in[2,p_0]$ we have
\begin{equation}\label{eq8.4}
  E\[\sup_{s\in[t,T]}|\partial_y(\partial_{\mu}Y_s^{t,x,P_{\xi}}(y))|^p+(\int_t^T|\partial_y(\partial_{\mu}Z_s^{t,x,P_{\xi}}(y))|^2ds)^{\frac{p}{2}}\]\leq C_p,
\end{equation}

Next we prove estimate (ii). Let $x,x',y,y'\in\mathbb{R}$, and let $\xi,\xi',\vartheta,\vartheta'\in L^2(\mathcal{G}_t;\mathbb{R})$
be such that $P_{\vartheta}\!=\!P_{\xi}$, $P_{\vartheta'}\!=\!P_{\xi'}$.
Note that $(X_s^{t,x,P_{\xi}},Y_s^{t,x,P_{\xi}},Z_s^{t,x,P_{\xi}})$ and $(O_s^{t,x,P_{\xi}}(y),Q_s^{t,x,P_{\xi}}(y))$, $t\leq s\leq T$,
are independent of $\mathcal{G}_t$. Hence, from \eqref{eq8.2} we get the following BDSDE:
\begin{equation}\label{eq8.5}
  \begin{split}
 &\  \big(\partial_{y}(\partial_{\mu}Y_s^{t,x,P_{\xi}}(y))-R_s^{t,y,P_{\xi}}\big)-\big(\partial_{y}(\partial_{\mu}Y_s^{t,x',P_{\xi'}}(y'))-R_s^{t,y',P_{\xi'}}\big)\\
  =&\ I_1(x,y,P_{\xi},x',y',P_{\xi'})+\int_s^TH(r,x,y,P_{\xi},x',y',P_{\xi'})d\overleftarrow{B_r}\\
&\  +\int_s^T g^2(P_{X_r^{t,\xi}})
    \Big(\partial_{y}(\partial_{\mu}Z_r^{t,x,P_{\xi}}(y))-\partial_{y}(\partial_{\mu}Z_r^{t,x',P_{\xi'}}(y'))\Big)d\overleftarrow{B_r}\\
&\  +\int_s^T\widehat{E}\[(\partial_{\mu}g^2)(P_{X_r^{t,\xi}},\widehat{X}_r^{t,\widehat{\vartheta},P_{\xi}})
     \Big(\partial_{y}(\partial_{\mu}\widehat{X}_{r}^{t,\widehat{\vartheta},P_{\xi}}(y))
     -\partial_{y}(\partial_{\mu}\widehat{X}_{r}^{t,\widehat{\vartheta}',P_{\xi'}}(y'))\Big) \]Z_r^{t,x,P_{\xi}}d\overleftarrow{B_r}\\
&\  -\int_s^T\Big(\partial_{y}(\partial_{\mu}Z_r^{t,x,P_{\xi}}(y))-\partial_{y}(\partial_{\mu}Z_r^{t,x',P_{\xi'}}(y'))\Big)dW_r,\ s\in[t,T],
  \end{split}
\end{equation}
where
\begin{equation*}
  \begin{split}
 & I_1(x,y,P_{\xi},x',y',P_{\xi'}):= (\partial_{x}\Phi) (X_T^{t,x,P_{\xi}}) \partial_{y}(\partial_{\mu}X_{T}^{t,x,P_{\xi}}(y))
         -(\partial_{x}\Phi) (X_T^{t,x',P_{\xi'}}) \partial_{y}(\partial_{\mu}X_{T}^{t,x',P_{\xi'}}(y')),\\
  & H(r,x,y,P_{\xi},x',y',P_{\xi'}):=
   \Big(g^2(P_{X_r^{t,\xi}})-g^2(P_{X_r^{t,\xi'}})\Big) \cdot\partial_{y}(\partial_{\mu}Z_r^{t,x',P_{\xi'}}(y'))\\
  &\quad + \widehat{E}\[\Big((\partial_{\mu}g^2)(P_{X_r^{t,\xi}},\widehat{X}_r^{t,\widehat{\vartheta},P_{\xi}})Z_r^{t,x,P_{\xi}}
 -(\partial_{\mu}g^2)(P_{X_r^{t,\xi'}},\widehat{X}_r^{t,\widehat{\vartheta}',P_{\xi'}})Z_r^{t,x',P_{\xi'}}\Big) \cdot\partial_{y}(\partial_{\mu}\widehat{X}_{r}^{t,\widehat{\vartheta}',P_{\xi'}}\!(y')) \]\\
&\quad +\(\widehat{E}\[(\partial_{\mu}g^2)(P_{X_r^{t,\xi}},\widehat{X}_r^{t,y,P_{\xi}}) \partial_{xx}^2\widehat{X}_r^{t,y,P_{\xi}}\]
 -\widehat{E}\[(\partial_{\mu}g^2)(P_{X_r^{t,\xi'}},\widehat{X}_r^{t,y',P_{\xi'}}) \partial_{xx}^2\widehat{X}_r^{t,y',P_{\xi'}}\]\)\\
&\quad +\(\widehat{E}\[\partial_{y}((\partial_{\mu}g^2)(P_{X_r^{t,\xi}},\widehat{X}_r^{t,y,P_{\xi}})) (\partial_{x}\widehat{X}_r^{t,y,P_{\xi}})^2\]
-\widehat{E}\[\partial_{y}((\partial_{\mu}g^2)(P_{X_r^{t,\xi'}},\widehat{X}_r^{t,y',P_{\xi'}})) (\partial_{x}\widehat{X}_r^{t,y',P_{\xi'}})^2\]\).
  \end{split}
\end{equation*}
From Lemma \ref{le3.1}, the Theorems \ref{th7.1} and \ref{th8.1} as well as Propositions \ref{prop4.1} and \ref{prop6.1},
 the same argument as in the proof of the Propositions \ref{prop6.2} and \ref{prop8.1} allow to deduce with the help of Corollary \ref{corA.2}
\begin{equation}\label{eq8.6}
  \begin{split}
  &\  E\[\sup_{s\in[t,T]}|\big(\partial_{y}(\partial_{\mu}Y_s^{t,x,P_{\xi}}(y))-R_s^{t,y,P_{\xi}}\big)
  -\big(\partial_{y}(\partial_{\mu}Y_s^{t,x',P_{\xi'}}(y'))-R_s^{t,y',P_{\xi'}}\big)|^p\\
&\  +(\int_t^T|\partial_{y}(\partial_{\mu}Z_s^{t,x,P_{\xi}}(y))\!-\!\partial_{y}(\partial_{\mu}Z_s^{t,x',P_{\xi'}}(y'))|^2ds)^{\frac{p}{2}}\]\\
 \leq&\ C_p\(|x-x'|^p+|y-y'|^p+W_2(P_{\xi},P_{\xi'})^p+(E[|\vartheta-\vartheta'|^2])^{\frac{p}{2}}\),\ \mbox{for all}\ p\in[2,\frac{p_0}{2}].
  \end{split}
\end{equation}
To get the estimate (ii) of the theorem we have still to estimate $R_s^{t,y,P_{\xi}}-R_s^{t,y',P_{\xi'}}$.
For this we remark that, thanks to Theorem \ref{th8.1}, Remark \ref{re4.1} and the Propositions \ref{prop6.1+1}, \ref{prop6.1+2} and \ref{prop8.1}-\eqref{eq8.0+1},
\begin{equation*}
  \begin{split}
\mbox{(i)}&\  \int_t^T\Big|\widehat{E}\[\partial_{z}(\partial_{\mu}h)_2(P_{(Y_r^{t,\xi},Z_r^{t,\xi})},(\widehat{Y}_r^{t,y,P_{\xi}},\widehat{Z}_r^{t,y,P_{\xi}}))
  (\partial_{x}\widehat{Z}_r^{t,y,P_{\xi}})^2\\
     &\qquad\qquad   - \partial_{z}(\partial_{\mu}h)_2(P_{(Y_r^{t,\xi'},Z_r^{t,\xi'})},(\widehat{Y}_r^{t,y',P_{\xi'}},\widehat{Z}_r^{t,y',P_{\xi'}}))
  (\partial_{x}\widehat{Z}_r^{t,y',P_{\xi'}})^2\]\Big|^2\!dr\\
&\  \leq C\!\int_t^T\!\!\Big|\!\widehat{E}\[\!\Big(\!W_2(P_{(Y_r^{t,\xi}\!,Z_r^{t,\xi})},P_{(Y_r^{t,\xi'}\!,Z_r^{t,\xi'})})
  \!+\!|\widehat{Y}_r^{t,y,P_{\xi}}\!\!-\!\!\widehat{Y}_r^{t,y',P_{\xi'}}|
  \!+\!|\widehat{Z}_r^{t,y,P_{\xi}}\!\!-\!\!\widehat{Z}_r^{t,y',P_{\xi'}}|\!\Big)(\partial_{x}\widehat{Z}_r^{t,y',P_{\xi'}})^2\!\]\!\Big|^2dr\\
     &\qquad  + C\int_t^T\(\widehat{E}\[|(\partial_{x}\widehat{Z}_r^{t,y,P_{\xi}})^2-(\partial_{x}\widehat{Z}_r^{t,y',P_{\xi'}})^2|\]\)^2dr\\
&\  \leq C\!\int_t^T\!\!\widehat{E}[|\partial_{x}\widehat{Z}_r^{t,y',P_{\xi'}}|^4]
     \widehat{E}\[\!W_2(P_{(Y_r^{t,\xi}\!,Z_r^{t,\xi})},P_{(Y_r^{t,\xi'}\!,Z_r^{t,\xi'})})^2
 \! \!+\!\!|\widehat{Y}_r^{t,y,P_{\xi}}\!\!-\!\!\widehat{Y}_r^{t,y',P_{\xi'}}|^2
 \! \!+\!\!|\widehat{Z}_r^{t,y,P_{\xi}}\!\!-\!\!\widehat{Z}_r^{t,y',P_{\xi'}}|^2\!\] dr\\
\end{split}
\end{equation*}
\begin{equation*}
  \begin{split}
     &\qquad  + C\int_t^T\widehat{E}\[|\partial_{x}\widehat{Z}_r^{t,y,P_{\xi}}|^2+|\partial_{x}\widehat{Z}_r^{t,y',P_{\xi'}}|^2\]\cdot
     \widehat{E}\[|\partial_{x}\widehat{Z}_r^{t,y,P_{\xi}}-\partial_{x}\widehat{Z}_r^{t,y',P_{\xi'}}|^2\]dr\\
  &\  \leq C\(|y-y'|^2+W_2(P_{\xi},P_{\xi'})^2\)
  +C\! \int_t^T\!E\[E\[|Y_r^{t,x,P_{\xi}}\!-\!Y_r^{t,x',P_{\xi'}}|^2\!+\!|Z_r^{t,x,P_{\xi}}\!-\!Z_r^{t,x',P_{\xi'}}|^2\]\big|_{\substack{x=\vartheta\\ x'=\vartheta'}}\]dr\\
&\  \leq C\(|y-y'|^2+W_2(P_{\xi},P_{\xi'})^2+E[|\vartheta-\vartheta'|^2]\)  ,
\end{split}
\end{equation*}
and similarly, we see
\begin{equation*}
  \begin{split}
\mbox{(ii)}&\  \int_t^T\Big|\widehat{E}\[(\partial_{\mu}h)_2(P_{(Y_r^{t,\xi},Z_r^{t,\xi})},(\widehat{Y}_r^{t,y,P_{\xi}},\widehat{Z}_r^{t,y,P_{\xi}}))
  \partial_{xx}^2\widehat{Z}_r^{t,y,P_{\xi}}\\
     &\qquad\qquad   - (\partial_{\mu}h)_2(P_{(Y_r^{t,\xi'},Z_r^{t,\xi'})},(\widehat{Y}_r^{t,y',P_{\xi'}},\widehat{Z}_r^{t,y',P_{\xi'}}))
  \partial_{xx}^2\widehat{Z}_r^{t,y',P_{\xi'}}\]\Big|^2\!dr\\
 &\  \leq \!C\!\int_t^T\!\!\widehat{E}[|\partial_{xx}^2\widehat{Z}_r^{t,y',P_{\xi'}}|^2]
     \widehat{E}\[\!W_2(P_{(Y_r^{t,\xi}\!,Z_r^{t,\xi})},P_{(Y_r^{t,\xi'}\!,Z_r^{t,\xi'})})^2
 \! \!+\!\!|\widehat{Y}_r^{t,y,P_{\xi}}\!\!-\!\!\widehat{Y}_r^{t,y',P_{\xi'}}|^2
 \! \!+\!\!|\widehat{Z}_r^{t,y,P_{\xi}}\!\!-\!\!\widehat{Z}_r^{t,y',P_{\xi'}}|^2\!\] \!dr\\
\end{split}
\end{equation*}
\begin{equation*}
  \begin{split}
     &\qquad  + C\int_t^T \widehat{E}\[|\partial_{xx}^2\widehat{Z}_r^{t,y,P_{\xi}}-\partial_{xx}^2\widehat{Z}_r^{t,y',P_{\xi'}}|^2\]dr\\
 &\  \leq C\(|y-y'|^2+W_2(P_{\xi},P_{\xi'})^2+E[|\vartheta-\vartheta'|^2]\)\cdot\(1+\int_t^T\widehat{E}[|\partial_{xx}^2\widehat{Z}_r^{t,y',P_{\xi'}}|^2]dr\)\\
  &\  \leq C\(|y-y'|^2+W_2(P_{\xi},P_{\xi'})^2+E[|\vartheta-\vartheta'|^2]\) .
\end{split}
\end{equation*}
Using \eqref{eq8.4} and previous estimates for our BDSDE with affine $g$, we also have
\begin{equation*}
  \begin{split}
\mbox{(iii)}&\  \int_t^T\Big|\widehat{E}\[(\partial_{\mu}h)_2(P_{(Y_r^{t,\xi},Z_r^{t,\xi})},(\widehat{Y}_r^{t,\widehat{\xi}},\widehat{Z}_r^{t,\widehat{\xi}}))
  \partial_{y}(\partial_{\mu}\widehat{Z}_r^{t,\widehat{\xi},P_{\xi}}(y))\\
     &\qquad\qquad    - (\partial_{\mu}h)_2(P_{(Y_r^{t,\xi'},Z_r^{t,\xi'})},(\widehat{Y}_r^{t,\widehat{\xi'}},\widehat{Z}_r^{t,\widehat{\xi'}}))
  \partial_{y}(\partial_{\mu}\widehat{Z}_r^{t,\widehat{\xi'},P_{\xi'}}(y))\]\Big|^2\!dr\\
&\  \leq C\!\int_t^T\!\!\widehat{E}[|\partial_{y}(\partial_{\mu}\widehat{Z}_r^{t,\widehat{\vartheta'},P_{\xi'}}(y))|^2]
     \widehat{E}\[\!W_2(P_{(Y_r^{t,\xi}\!,Z_r^{t,\xi})},P_{(Y_r^{t,\xi'}\!,Z_r^{t,\xi'})})^2
 \! \!+\!\!|\widehat{Y}_r^{t,\widehat{\xi}}\!\!-\!\!\widehat{Y}_r^{t,\widehat{\xi'}}|^2
 \! \!+\!\!|\widehat{Z}_r^{t,\widehat{\xi}}\!\!-\!\!\widehat{Z}_r^{t,\widehat{\xi'}}|^2\!\] dr\\
     &\qquad  + C\int_t^T \widehat{E}\[|\partial_{y}(\partial_{\mu}\widehat{Z}_r^{t,\widehat{\vartheta},P_{\xi}}(y))
                           -\partial_{y}(\partial_{\mu}\widehat{Z}_r^{t,\widehat{\vartheta'},P_{\xi'}}(y))|^2\]dr\\
&\  \leq C\!\int_t^T\!\!\widehat{E}[|\partial_{y}(\partial_{\mu}\widehat{Z}_r^{t,\widehat{\vartheta'},P_{\xi'}}(y))|^2]
   E\[E\[|Y_r^{t,x,P_{\xi}}\!-\!Y_r^{t,x',P_{\xi'}}|^2\!+\!|Z_r^{t,x,P_{\xi}}\!-\!Z_r^{t,x',P_{\xi'}}|^2\]\big|_{\substack{x=\vartheta\\ x'=\vartheta'}}\]dr\\
     &\qquad  + CE\[E\[\int_t^T|\partial_{y}(\partial_{\mu}Z_r^{t,x,P_{\xi}}(y))
                           -\partial_{y}(\partial_{\mu}Z_r^{t,x',P_{\xi'}}(y))|^2dr\]\big|_{\substack{x=\vartheta\\ x'=\vartheta'}}\]\\
  &\  \leq C\(W_2(P_{\xi},P_{\xi'})^2+E[|\vartheta-\vartheta'|^2]\).
\end{split}
\end{equation*}
With similar estimates for the remaining terms of $R_s^{t,y,P_{\xi}}-R_s^{t,y',P_{\xi'}}$ we get that, for $\vartheta,\vartheta'\in L^2(\mathcal{G}_t;\mathbb{R})$
with $P_{\vartheta}\!=\!P_{\xi}$, $P_{\vartheta'}\!=\!P_{\xi'}$:
\begin{equation}\nonumber
  \begin{split}
  E\[\sup_{r\in[s,T]}|R_r^{t,y,P_{\xi}}-R_r^{t,y',P_{\xi'}}|^{\frac{p_0}{2}}\]
\leq&\   C_{\frac{p_0}{2}}\(|y-y'|^p+W_2(P_{\xi},P_{\xi'})^{\frac{p_0}{2}}+(E[|\vartheta-\vartheta'|^2])^{\frac{p_0}{4}}\)\\
 &\ + C_{\frac{p_0}{2}}\(\!E\[\!\int_s^T\!|\partial_{y}(\partial_{\mu}Y_r^{t,\vartheta,P_{\xi}}(y))
                           \!-\!\partial_{y}(\partial_{\mu}Y_r^{t,\vartheta',P_{\xi'}}(y))|^2dr\]\)^{\frac{p_0}{4}}.
  \end{split}
\end{equation}
Furthermore, from H\"{o}lder's inequality, for all $p\in[2,\frac{p_0}{2}]$ we have
\begin{equation}\label{eq8.7}
  \begin{split}
  E\[\sup_{r\in[s,T]}|R_r^{t,y,P_{\xi}}-R_r^{t,y',P_{\xi'}}|^{p}\]
\leq&\   C_{p}\(|y-y'|^p+W_2(P_{\xi},P_{\xi'})^{\frac{p}{2}}+(E[|\vartheta-\vartheta'|^2])^{\frac{p}{2}}\)\\
 &\ + C_{p}\(\!E\[\!\int_s^T\!|\partial_{y}(\partial_{\mu}Y_r^{t,\vartheta,P_{\xi}}(y))
                           \!-\!\partial_{y}(\partial_{\mu}Y_r^{t,\vartheta',P_{\xi'}}(y))|^2dr\]\)^{\frac{p}{2}}.
  \end{split}
\end{equation}
Then, substituting $x=\vartheta$, $x'=\vartheta'$ in \eqref{eq8.6}, using the above estimate, we obtain, for $s\in[t,T]$,
\begin{equation*}
  \begin{split}
 &\ E\[\sup_{r\in[s,T]}|\partial_{y}(\partial_{\mu}Y_r^{t,\vartheta,P_{\xi}}(y))
                           \!-\!\partial_{y}(\partial_{\mu}Y_r^{t,\vartheta',P_{\xi'}}(y))|^2\]\\
=&\  E\[E\[\sup_{r\in[s,T]}|\partial_{y}(\partial_{\mu}Y_r^{t,x,P_{\xi}}(y))
                           \!-\!\partial_{y}(\partial_{\mu}Y_r^{t,x',P_{\xi'}}(y))|^2\]\big|_{\substack{x=\vartheta\\ x'=\vartheta'}}\]\\
\leq&\ CE\[\sup_{r\in[s,T]}|R_r^{t,y,P_{\xi}}-R_r^{t,y',P_{\xi'}}|^2\]+ C\(|y-y'|^2+W_2(P_{\xi},P_{\xi'})^2+E[|\vartheta-\vartheta'|^2]\)\\
\leq&\ C\(|y-y'|^2+W_2(P_{\xi},P_{\xi'})^2+E[|\vartheta-\vartheta'|^2]\)
        +C E\[\!\int_s^T\!|\partial_{y}(\partial_{\mu}Y_r^{t,\vartheta,P_{\xi}}(y))
                           \!-\!\partial_{y}(\partial_{\mu}Y_r^{t,\vartheta',P_{\xi'}}(y))|^2dr\],
  \end{split}
\end{equation*}
and thus, from Gronwall's Lemma,
\begin{equation}\label{eq8.8}
  \begin{split}
 & E\[\sup_{r\in[s,T]}\!|\partial_{y}(\partial_{\mu}Y_r^{t,\vartheta,P_{\xi}}(y))
                           \!-\!\partial_{y}(\partial_{\mu}Y_r^{t,\vartheta',P_{\xi'}}(y))|^2\]
\leq C\(|y\!-\!\!y'|^2\!+\!W_2(P_{\xi},P_{\xi'})^2\!+\!E[|\vartheta\!-\!\vartheta'|^2]\).
  \end{split}
\end{equation}
Hence, \eqref{eq8.7} yields
\begin{equation}\label{eq8.9}
  \begin{split}
  E\[\sup_{r\in[s,T]}|R_r^{t,y,P_{\xi}}-R_r^{t,y',P_{\xi'}}|^p\]
\leq&\   C_p\(|y-y'|^p+W_2(P_{\xi},P_{\xi'})^p+(E[|\vartheta-\vartheta'|^2])^{\frac{p}{2}}\).
  \end{split}
\end{equation}
Then, from \eqref{eq8.6} and \eqref{eq8.9} we conclude
\begin{equation}\label{eq8.10}
  \begin{split}
   & E\[\!\sup_{s\in[t,T]}|\partial_{y}(\partial_{\mu}Y_s^{t,x,P_{\xi}}(y))\!-\!\partial_{y}(\partial_{\mu}Y_s^{t,x',P_{\xi'}}(y'))|^p
  \!+\!(\int_t^T|\partial_{y}(\partial_{\mu}Z_s^{t,x,P_{\xi}}(y))\!-\!\partial_{y}(\partial_{\mu}Z_s^{t,x',P_{\xi'}}(y'))|^2ds)^{\frac{p}{2}}\!\]\\
& \leq C_p\(|x-x'|^p+|y-y'|^p+W_2(P_{\xi},P_{\xi'})^p+(E[|\vartheta-\vartheta'|^2])^{\frac{p}{2}}\),\ \mbox{for all}\ p\in[2,\frac{p_0}{2}].
  \end{split}
\end{equation}
Finally, from \eqref{eq8.10}, by taking the infimum of $(E[|\vartheta-\vartheta'|^2])^{\frac{p}{2}}$ over all
$\vartheta,\vartheta'\in L^2(\mathcal{G}_t;\mathbb{R})$
with $P_{\vartheta}\!=\!P_{\xi}$, $P_{\vartheta'}\!=\!P_{\xi'}$, we get estimate (ii).

Now it still remains to show that the formal derivative $(\partial_y(\partial_{\mu}Y^{t,x,P_{\xi}}(y)),\partial_y(\partial_{\mu}Z^{t,x,P_{\xi}}(y)))$
is really the $L^2$-derivative of $(\partial_{\mu}Y^{t,x,P_{\xi}}(y),\partial_{\mu}Z^{t,x,P_{\xi}}(y))$.
For this we estimate
\begin{equation*}
 \begin{split}
   & \mbox{(i) }\frac{1}{q}(\partial_{\mu}Y^{t,x,P_{\xi}}(y+q)-\partial_{\mu}Y^{t,x,P_{\xi}}(y))-\partial_y(\partial_{\mu}Y^{t,x,P_{\xi}}(y)),\\
   & \mbox{(ii) }\frac{1}{q}(\partial_{\mu}Z^{t,x,P_{\xi}}(y+q)-\partial_{\mu}Z^{t,x,P_{\xi}}(y))-\partial_y(\partial_{\mu}Z^{t,x,P_{\xi}}(y)).
 \end{split}
\end{equation*}
We note that, from \eqref{eq8.1} and \eqref{eq8.2}, for all $s\in[t,T]$,
\begin{equation}\label{eq8.11}
  \begin{split}
&\ \frac{1}{q}(\partial_{\mu}Y_s^{t,x,P_{\xi}}(y+q)-\partial_{\mu}Y_s^{t,x,P_{\xi}}(y))-\partial_y(\partial_{\mu}Y_s^{t,x,P_{\xi}}(y))\\
=&\ I^{q}(x,y)+I_s^{1,q}(x,y)+I_s^{2,q}(x,y)+I_s^{3,q}(y)+I_s^{4,q}(y)\\
&\ +\int_s^T g^2(P_{X_r^{t,\xi}})\big(\frac{1}{q}(\partial_{\mu}Z_r^{t,x,P_{\xi}}(y+q)-\partial_{\mu}Z_r^{t,x,P_{\xi}}(y))
   -\partial_{y}(\partial_{\mu}Z_r^{t,x,P_{\xi}}(y)) \big)d\overleftarrow{B_r}\\
   &\  -\int_s^T\(\frac{1}{q}(\partial_{\mu}Z_r^{t,x,P_{\xi}}(y+q)-\partial_{\mu}Z_r^{t,x,P_{\xi}}(y))-\partial_y(\partial_{\mu}Z_r^{t,x,P_{\xi}}(y))\)dW_r,
  \end{split}
\end{equation}
where
\begin{equation*}
  \begin{split}
 & I^{q}(x,y):=(\partial_{x}\Phi) (X_T^{t,x,P_{\xi}})\big(\frac{1}{q}(\partial_{\mu}X_{T}^{t,x,P_{\xi}}(y+q)-\partial_{\mu}X_{T}^{t,x,P_{\xi}}(y))
     -\partial_{y}(\partial_{\mu}X_{T}^{t,x,P_{\xi}}(y))\big),\\
 & I_s^{1,q}(x,y)\! :=\!\!\int_s^T\!\!\widehat{E}\[\!(\partial_{\mu}g^2)(P_{X_r^{t,\xi}},\widehat{X}_r^{t,\widehat{\xi}})
   \big(\frac{1}{q} (\partial_{\mu}\widehat{X}_r^{t,\widehat{\xi},P_{\xi}}(y\!+\!q)\!-\!\partial_{\mu}\widehat{X}_r^{t,\widehat{\xi},P_{\xi}}(y))
      \!-\!\partial_{y}(\partial_{\mu}\widehat{X}_{r}^{t,\widehat{\xi}}(y))\big)\! \]Z_r^{t,x,P_{\xi}}d\overleftarrow{B_r},\\
 & I_s^{2,q}(x,y) :=\int_s^T\!\!\widehat{E}\[\frac{1}{q}\big((\partial_{\mu}g^2)(P_{X_r^{t,\xi}},\widehat{X}_r^{t,y+q,P_{\xi}})\partial_{x}\widehat{X}_r^{t,y+q,P_{\xi}}
      -(\partial_{\mu}g^2)(P_{X_r^{t,\xi}},\widehat{X}_r^{t,y,P_{\xi}})\partial_{x}\widehat{X}_r^{t,y,P_{\xi}}\big) \\
 &\hspace{55pt}   \!-\!\partial_{y}((\partial_{\mu}g^2)(P_{X_r^{t,\xi}},\widehat{X}_r^{t,y,P_{\xi}}))(\partial_{x}\widehat{X}_r^{t,y,P_{\xi}})^2
    \!-\!(\partial_{\mu}g^2)(P_{X_r^{t,\xi}},\widehat{X}_r^{t,y,P_{\xi}})\partial_{xx}^2\widehat{X}_r^{t,y,P_{\xi}}\!\]Z_r^{t,x,P_{\xi}}d\overleftarrow{B_r},\\
 & I_s^{3,q}(y) := \int_s^T\!\! \Bigg\{ \frac{1}{q}\(\widehat{E}\[\Big\langle
    (\partial_{\mu}h)(P_{(Y_r^{t,\xi},Z_r^{t,\xi})},(\widehat{Y}_r^{t,y+q,P_{\xi}},\widehat{Z}_r^{t,y+q,P_{\xi}})),
    \left(   \begin{matrix}    \partial_{x}\widehat{Y}_r^{t,y+q,P_{\xi}}\\ \partial_{x}\widehat{Z}_r^{t,y+q,P_{\xi}} \end{matrix}\right)\Big\rangle\]\\
 &\hspace{55pt} -\widehat{E}\[\Big\langle(\partial_{\mu}h)(P_{(Y_r^{t,\xi},Z_r^{t,\xi})},(\widehat{Y}_r^{t,y,P_{\xi}},\widehat{Z}_r^{t,y,P_{\xi}})),
  \left(   \begin{matrix}    \partial_{x}\widehat{Y}_r^{t,y,P_{\xi}}\\ \partial_{x}\widehat{Z}_r^{t,y,P_{\xi}} \end{matrix}\right)\Big\rangle\] \) \\
  \end{split}
\end{equation*}
\begin{equation*}
  \begin{split}
 &\hspace{55pt}  -\widehat{E}\[\Big\langle\partial_{y}(\partial_{\mu}h)(P_{(Y_r^{t,\xi},Z_r^{t,\xi})},(\widehat{Y}_r^{t,y,P_{\xi}},\widehat{Z}_r^{t,y,P_{\xi}}))\cdot
  \left(   \begin{matrix}    \partial_{x}\widehat{Y}_r^{t,y,P_{\xi}}\\ \partial_{x}\widehat{Z}_r^{t,y,P_{\xi}} \end{matrix}\right),
  \left(   \begin{matrix}    \partial_{x}\widehat{Y}_r^{t,y,P_{\xi}}\\ \partial_{x}\widehat{Z}_r^{t,y,P_{\xi}} \end{matrix}\right)\Big\rangle\]\\
&\hspace{55pt} - \widehat{E}\[\Big\langle(\partial_{\mu}h)(P_{(Y_r^{t,\xi},Z_r^{t,\xi})},(\widehat{Y}_r^{t,y,P_{\xi}},\widehat{Z}_r^{t,y,P_{\xi}})),
  \left(   \begin{matrix}    \partial_{xx}^2\widehat{Y}_r^{t,y,P_{\xi}}\\ \partial_{xx}^2\widehat{Z}_r^{t,y,P_{\xi}} \end{matrix}\right)\Big\rangle\]\Bigg\}d\overleftarrow{B_r},\\
 & I_s^{4,q}(y) := \int_s^T\!\! \Bigg\{
   \widehat{E}\[\Big\langle(\partial_{\mu}h)(P_{(Y_r^{t,\xi},Z_r^{t,\xi})},(\widehat{Y}_r^{t,\widehat{\xi}},\widehat{Z}_r^{t,\widehat{\xi}})),
 \left(   \begin{matrix}    \frac{1}{q}\big(\partial_{\mu}\widehat{Y}_r^{t,\widehat{\xi},P_{\xi}}(y+q)
                              -\partial_{\mu}\widehat{Y}_r^{t,\widehat{\xi},P_{\xi}}(y)\big)\\
 \frac{1}{q}\big(\partial_{\mu}\widehat{Z}_r^{t,\widehat{\xi},P_{\xi}}(y+q)
                              -\partial_{\mu}\widehat{Z}_r^{t,\widehat{\xi},P_{\xi}}(y)\big) \end{matrix}\right)
    \Big\rangle\]\\
&\hspace{55pt}    -\widehat{E}\[\Big\langle(\partial_{\mu}h)(P_{(Y_r^{t,\xi},Z_r^{t,\xi})},(\widehat{Y}_r^{t,\widehat{\xi}},\widehat{Z}_r^{t,\widehat{\xi}})),
\left(   \begin{matrix}     \partial_{y}(\partial_{\mu}\widehat{Y}_r^{t,\widehat{\xi},P_{\xi}}(y))\\
     \partial_{y}(\partial_{\mu}\widehat{Z}_r^{t,\widehat{\xi},P_{\xi}}(y)) \end{matrix}\right)
    \Big\rangle\]\Bigg\}d\overleftarrow{B_r}.
  \end{split}
\end{equation*}
From Theorem \ref{th7.1} we know
\begin{equation*}\label{eq8.12}
  \begin{split}
& \mbox{(i)}\ E[|I^{q}(x,y)|^2]\leq C E\[|\big(\frac{1}{q}(\partial_{\mu}X_{T}^{t,x,P_{\xi}}(y\!+\!q)\!-\!\partial_{\mu}X_{T}^{t,x,P_{\xi}}(y))
     \!-\!\partial_{y}(\partial_{\mu}X_{T}^{t,x,P_{\xi}}(y))\big)|^2\]\rightarrow0,\ \mbox{as}\ q\rightarrow0;\\
& \mbox{(ii)}\ E[\sup_{s\in[t,T]}|I_s^{1,q}(x,y)|^2]\!\leq\! C E[\int_t^T\!|Z_r^{t,x,P_{\xi}}|^2dr]
\widehat{E}\[\!\sup_{r\in[t,T]}\!
   \big|\frac{\partial_{\mu}\widehat{X}_r^{t,\widehat{\xi}}(y\!+\!q)\!-\!\partial_{\mu}\widehat{X}_r^{t,\widehat{\xi}}(y)}{q}
      \!-\!\partial_{y}(\partial_{\mu}\widehat{X}_{r}^{t,\widehat{\xi}}(y))\big|^2\! \]\\
&\hspace{120pt}      \rightarrow0,\ \mbox{as}\ q\rightarrow0.
  \end{split}
\end{equation*}
Moreover, Theorem \ref{th7.1} combined with Lemma \ref{le3.1} and Proposition \ref{prop5.1} also allows to show
\begin{equation*}\label{eq8.13}
  \begin{split}
 E[\sup_{s\in[t,T]}|I_s^{2,q}(x,y)|^2]  \rightarrow0,\ \mbox{as}\ q\rightarrow0.
  \end{split}
\end{equation*}
Furthermore, as $I_s^{3,q}(y)$, $I_s^{4,q}(y)$ and so also $J_s^{q}(y):=I_s^{3,q}(y)+I_s^{4,q}(y)$ are $\mathcal{F}^{B}_{s,T}$-measurable, $s\in[t,T)$,
it follows that from \eqref{eq8.11}
\begin{equation}\label{eq8.14}
  \begin{split}
   & E\[\!\sup_{s\in[t,T]}|\frac{1}{q}(\partial_{\mu}Y_s^{t,x,P_{\xi}}(y+q)-\partial_{\mu}Y_s^{t,x,P_{\xi}}(y))-\partial_y(\partial_{\mu}Y_s^{t,x,P_{\xi}}(y))-J_s^{q}(y)|^2\]\\
 & +E\[\int_t^T|\frac{1}{q}(\partial_{\mu}Z_r^{t,x,P_{\xi}}(y+q)-\partial_{\mu}Z_r^{t,x,P_{\xi}}(y))-\partial_y(\partial_{\mu}Z_r^{t,x,P_{\xi}}(y))|^2ds\]\\
\leq&\ C\( E[|I^{q}(x,y)|^2]+E[\sup_{s\in[t,T]}|I_s^{1,q}(x,y)|^2]+E[\sup_{s\in[t,T]}|I_s^{2,q}(x,y)|^2]\)\rightarrow0,\ \mbox{as}\ q\rightarrow0.
  \end{split}
\end{equation}
Substituting in \eqref{eq8.14} $\xi$ for $x$, we get from the independence of the integrands of $\mathcal{G}_t$,
\begin{equation*}
  \begin{split}
  &\ E\[\int_t^T|\frac{1}{q}(\partial_{\mu}Z_r^{t,\xi,P_{\xi}}(y+q)-\partial_{\mu}Z_r^{t,\xi,P_{\xi}}(y))-\partial_y(\partial_{\mu}Z_r^{t,\xi,P_{\xi}}(y))|^2ds\]\\
\leq&\ C\( E[|I^{q}(\xi,y)|^2]+E[\sup_{s\in[t,T]}|I_s^{1,q}(\xi,y)|^2]+E[\sup_{s\in[t,T]}|I_s^{2,q}(\xi,y)|^2]\),
  \end{split}
\end{equation*}
and with the same arguments as those used for \eqref{eq8.14}, we show that
\begin{equation*}
  \begin{split}
 E[|I^{q}(\xi,y)|^2]+E[\sup_{s\in[t,T]}|I_s^{1,q}(\xi,y)|^2]+E[\sup_{s\in[t,T]}|I_s^{2,q}(\xi,y)|^2]\rightarrow0,\ \mbox{as}\ q\rightarrow0.
  \end{split}
\end{equation*}
Then also
\begin{equation*}
  \begin{split}
  E\[\int_t^T|\frac{1}{q}(\partial_{\mu}Z_r^{t,\xi,P_{\xi}}(y+q)-\partial_{\mu}Z_r^{t,\xi,P_{\xi}}(y))-\partial_y(\partial_{\mu}Z_r^{t,\xi,P_{\xi}}(y))|^2ds\]
\rightarrow0,\ \mbox{as}\ q\rightarrow0.
  \end{split}
\end{equation*}
On the other hand, from Theorem \ref{th8.1} we have that
\begin{equation}\label{eq8.15}
  \begin{split}
   & E\[\!\sup_{s\in[t,T]}\!|\frac{\partial_{x}Y_s^{t,x+q,P_{\xi}}\!-\!\partial_{x}Y_s^{t,x,P_{\xi}}}{q}\!-\!\partial_{xx}^2Y_s^{t,x,P_{\xi}}|^2
  \!+\!\int_t^T\!\!|\frac{\partial_{x}Z_s^{t,x+q,P_{\xi}}\!-\!\partial_{x}Z_s^{t,x,P_{\xi}}}{q}\!-\!\partial_{xx}^2Z_s^{t,x,P_{\xi}}|^2ds\]\\
  &\rightarrow0,\ \mbox{as}\ q\rightarrow0.
  \end{split}
\end{equation}
Then, using Theorem \ref{th8.1} and \eqref{eq8.0+1}, we get with the help of standard estimates that
\begin{equation*}\label{eq8.16}
  \begin{split}
 E[\sup_{s\in[t,T]}|I_s^{3,q}(y)|^2]  \rightarrow0,\ \mbox{as}\ q\rightarrow0.
  \end{split}
\end{equation*}
Consequently, as
\begin{equation*}\label{eq8.17}
  \begin{split}
 E[\sup_{s\in[t,T]}|I_s^{4,q}(y)|^2]
 \leq&\ C \widehat{E}\[\int_t^T|\frac{1}{q}\big(\partial_{\mu}\widehat{Y}_r^{t,\widehat{\xi},P_{\xi}}(y+q)
                              -\partial_{\mu}\widehat{Y}_r^{t,\widehat{\xi},P_{\xi}}(y)\big)
                             - \partial_{y}(\partial_{\mu}\widehat{Y}_r^{t,\widehat{\xi},P_{\xi}}(y))|^2dr\]\\
 &\  +C \widehat{E}\[\int_t^T|\frac{1}{q}\big(\partial_{\mu}\widehat{Z}_r^{t,\widehat{\xi},P_{\xi}}(y+q)
                              -\partial_{\mu}\widehat{Z}_r^{t,\widehat{\xi},P_{\xi}}(y)\big)
                             - \partial_{y}(\partial_{\mu}\widehat{Z}_r^{t,\widehat{\xi},P_{\xi}}(y))|^2dr\] ,
  \end{split}
\end{equation*}
it follows that
\begin{equation*}\label{eq8.18}
  \begin{split}
 E[\sup_{s\in[t,T]}|J_s^{q}(y)|^2]
 \leq C \widehat{E}\[\int_t^T|\frac{1}{q}\big(\partial_{\mu}\widehat{Y}_r^{t,\widehat{\xi},P_{\xi}}(y+q)
                              -\partial_{\mu}\widehat{Y}_r^{t,\widehat{\xi},P_{\xi}}(y)\big)
                             - \partial_{y}(\partial_{\mu}\widehat{Y}_r^{t,\widehat{\xi},P_{\xi}}(y))|^2dr\]+R^{q}(y),
  \end{split}
\end{equation*}
where
\begin{equation*}\label{eq8.19}
  \begin{split}
R^{q}(y):=CE[\sup_{s\in[t,T]}|I_s^{3,q}(y)|^2]\!+\!C \widehat{E}\[\int_t^T|\frac{1}{q}\big(\partial_{\mu}\widehat{Z}_r^{t,\widehat{\xi},P_{\xi}}(y+q)
                              \!-\!\partial_{\mu}\widehat{Z}_r^{t,\widehat{\xi},P_{\xi}}(y)\big)
                             \!- \!\partial_{y}(\partial_{\mu}\widehat{Z}_r^{t,\widehat{\xi},P_{\xi}}(y))|^2dr\],
  \end{split}
\end{equation*}
with $R^{q}(y)\rightarrow0,\ \mbox{as}\ q\rightarrow0.$

Combining this with \eqref{eq8.14} and applying Gronwall's Lemma, we see that also
\begin{equation*}\label{eq8.20}
  \begin{split}
  \widehat{E}\[\sup_{s\in[t,T]}|\frac{1}{q}(\partial_{\mu}\widehat{Y}_s^{t,\widehat{\xi},P_{\xi}}(y+q)-\partial_{\mu}\widehat{Y}_s^{t,\widehat{\xi},P_{\xi}}(y))
  -\partial_y(\partial_{\mu}\widehat{Y}_s^{t,\widehat{\xi},P_{\xi}}(y))|^2\]
\rightarrow0,\ \mbox{as}\ q\rightarrow0.
  \end{split}
\end{equation*}
Consequently, from \eqref{eq8.14}, we get the wished result: \begin{equation*}
 \begin{split}
   & \frac{1}{q}(\partial_{\mu}Y_{\cdot}^{t,x,P_{\xi}}(y+q)-\partial_{\mu}Y_{\cdot}^{t,x,P_{\xi}}(y))
   -\partial_y(\partial_{\mu}Y_{\cdot}^{t,x,P_{\xi}}(y))\xrightarrow[q\rightarrow 0]{\mathcal{S}_{\mathcal{F}}^2(t,T;\mathbb{R})}0,\\
   & \frac{1}{q}(\partial_{\mu}Z_{\cdot}^{t,x,P_{\xi}}(y+q)-\partial_{\mu}Z_{\cdot}^{t,x,P_{\xi}}(y))
   -\partial_y(\partial_{\mu}Z_{\cdot}^{t,x,P_{\xi}}(y))\xrightarrow[q\rightarrow 0]{\mathcal{H}_{\mathcal{F}}^2(t,T;\mathbb{R})}0.
 \end{split}
\end{equation*}
\end{proof}

\section{Related backward SPDEs of mean-field type}
The objective of this section is to study the related backward SPDEs of mean-field type. We will
prove that $V(t,x,P_{\xi})$ defined by \eqref{eq4.17} is the unique classical solution of the following
nonlocal
semi-linear backward SPDE of mean-field type:
\begin{equation}\nonumber
  \begin{split}
  &\  V(t,x,P_{\xi})\\
  =&\ \Phi(x,P_{\xi})+\int_t^T\Bigg\{\sum_{i=1}^d\partial_{x_i}V(s,x,P_{\xi})b_i(x,P_{\xi})
 +\frac{1}{2}\sum_{i,j,k=1}^d(\partial_{x_ix_j}^2V)(s,x,P_{\xi})(\sigma_{i,k}\sigma_{j,k})(x,P_{\xi})\\
  \end{split}
\end{equation}
\begin{equation}\label{eq9.1}
  \begin{split}
 &\ + f(x,V(s,x,P_{\xi}),\sum_{i=1}^d\partial_{x_i}V(s,x,P_{\xi})\sigma_i(x,P_{\xi}),P_{(\xi,\psi(s,\xi,P_{\xi}))}) \\
 &\ +E\[\sum_{i=1}^d(\partial_{\mu}V)_i(s,x,P_{\xi},\xi)b_i(\xi,P_{\xi})
  +\frac{1}{2}\sum_{i,j,k=1}^d\partial_{y_i}(\partial_{\mu}V)_j(s,x,P_{\xi},\xi)(\sigma_{i,k}\sigma_{j,k})(\xi,P_{\xi}) \] \Bigg\}ds\\
 &\ + \int_t^T \sum_{j=1}^lg_j(x,V(s,x,P_{\xi}),\sum_{i=1}^d\partial_{x_i}V(s,x,P_{\xi})\sigma_i(x,P_{\xi}),P_{(\xi,\psi(s,\xi,P_{\xi}))})d\overleftarrow{B_s^j},\\
 &\ + \int_t^T \sum_{j=1}^lh_j(P_{(\xi,\psi(s,\xi,P_{\xi}))})d\overleftarrow{B_s^j},\ (t,x,\xi,P_{\xi})\in[0,T]\times\mathbb{R}^d\times L^2(\mathcal{G}_t;\mathbb{R}^d)\times\mathcal{P}_2(\mathbb{R}^d),
  \end{split}
\end{equation}
where $\psi(s,x,P_{\xi}):=(V(s,x,P_{\xi}),\sum_{i=1}^d\partial_{x_i}V(s,x,P_{\xi})\sigma_i(x,P_{\xi}))$,
and the derivatives $\partial_{x_i}V$, $\partial_{x_ix_j}^2V$ and $\partial_{y_i}(\partial_{\mu}V)$ are in $L^2$-sense.

\begin{proposition} \label{prop9.1}
(Representation Formulas). Under the Assumptions (H4.2), (H7.1), (H8.1) and (H8.2) we have
the following representation formulas:
\begin{equation}\label{eq9.2}
  \begin{split}
&\ Y_s^{t,x,P_{\xi}}=V(s,X_s^{t,x,P_{\xi}},P_{X_s^{t,\xi}}),\ P\mbox{-}a.s.,\ s\in[t,T];\\
&\ Z_s^{t,x,P_{\xi}}=(\partial_xV)(s,X_s^{t,x,P_{\xi}},P_{X_s^{t,\xi}})\sigma(X_s^{t,x,P_{\xi}},P_{X_s^{t,\xi}}),\ dsdP\mbox{-}a.e.
  \end{split}
\end{equation}
Moreover,
\begin{equation}\label{eq9.3}
  E\[|Z_s^{t,x,P_{\xi}}-(\partial_xV)(s,x,P_{\xi})\sigma(x,P_{\xi})|^2\]\leq C(s-t),\ 0\leq t\leq s\leq T.
\end{equation}
\end{proposition}

\begin{remark} \label{re9.1}
Due to \eqref{eq4.10} the solution of BDSDE \eqref{eq4.1} has the following representation formulas:
\begin{equation}\label{eq9.4}
  \begin{split}
&\ Y_s^{t,\xi}=V(s,X_s^{t,\xi},P_{X_s^{t,\xi}}),\ P\mbox{-}a.s.,\ s\in[t,T];\\
&\ Z_s^{t,\xi}=(\partial_xV)(s,X_s^{t,\xi},P_{X_s^{t,\xi}})\sigma(X_s^{t,\xi},P_{X_s^{t,\xi}}),\ dsdP\mbox{-}a.e.
  \end{split}
\end{equation}
\end{remark}

\begin{proof}
The representation formulas for $Y^{t,x,P_{\xi}}$ is an immediate consequence of \eqref{eq4.11}.
Let us prove that for $Z_s^{t,x,P_{\xi}}$.
For simplicity, we suppose $d=l=1$, $b=0$, $f=0$, $g^1=0$ and $\Phi(x,\mu)=\Phi(x)$.

From Proposition \ref{prop3.1} and Theorem \ref{th5.1} we have
\begin{equation}\label{eq9.5}
  \begin{split}
  \partial_{x}X_s^{t,x,P_{\xi}}=& 1+ \int_t^s\partial_{x}\sigma(X_r^{t,x,P_{\xi}},P_{X_r^{t,\xi}}) \partial_{x}X_s^{t,x,P_{\xi}}dW_r,\ t\leq s\leq T,
  \end{split}
\end{equation}
and
\begin{equation}\label{eq9.6}
  \begin{split}
 D_{\theta}[X_s^{t,x,P_{\xi}}]=& \sigma(X_{\theta}^{t,x,P_{\xi}},P_{X_{\theta}^{t,\xi}})
  + \int_{\theta}^s\partial_{x}\sigma(X_r^{t,x,P_{\xi}},P_{X_r^{t,\xi}}) D_{\theta}[X_r^{t,x,P_{\xi}}]dW_r, \ t\leq \theta\leq s\leq T.
  \end{split}
\end{equation}
It follows from the uniqueness of the solution of SDE \eqref{eq9.6} that
\begin{equation}\label{eq9.7}
  \begin{split}
D_{\theta}[X_s^{t,x,P_{\xi}}]= \partial_{x}X_s^{t,x,P_{\xi}}(\partial_{x}X_{\theta}^{t,x,P_{\xi}})^{-1}\sigma(X_{\theta}^{t,x,P_{\xi}},P_{X_{\theta}^{t,\xi}}),\
t\leq \theta\leq s\leq T.
  \end{split}
\end{equation}
Furthermore, from Propositions \ref{prop6.1+1} and \ref{prop6.1+2}, we get for $t\leq\theta\leq s\leq T$:
\begin{equation}\label{eq9.8}
  \begin{split}
 \partial_{x}Y_s^{t,x,P_{\xi}}=& \partial_{x}\Phi (X_T^{t,x,P_{\xi}})\partial_{x}X_T^{t,x,P_{\xi}}
+\int_s^Tg^2(P_{X_r^{t,\xi}})\partial_{x}Z_r^{t,x,P_{\xi}}d\overleftarrow{B_r}
 -\int_s^T\partial_{x}Z_r^{t,x,P_{\xi}}dW_r;
  \end{split}
\end{equation}
\begin{equation}\label{eq9.9}
  \begin{split}
D_{\theta} [Y_s^{t,x,P_{\xi}}]=&\ \partial_{x}\Phi (X_T^{t,x,P_{\xi}})D_{\theta}X_T^{t,x,P_{\xi}}
  +\int_s^Tg^2(P_{X_r^{t,\xi}})D_{\theta}[Z_r^{t,x,P_{\xi}}]d\overleftarrow{B_r}
 \!-\!\int_s^T\!D_{\theta}[Z_r^{t,x,P_{\xi}}]dW_r.
  \end{split}
\end{equation}
Then, as $(\partial_{x}X_{\theta}^{t,x,P_{\xi}})^{-1}\sigma(X_{\theta}^{t,x,P_{\xi}},P_{X_{\theta}^{t,\xi}})$
is $\mathcal{F}_{\theta}^{W}$-measurable (and, hence independent of $\mathcal{F}_{0,T}^{B}$), it follows from \eqref{eq9.7}
that $(\partial_{x}Y_s^{t,x,P_{\xi}}(\partial_{x}X_{\theta}^{t,x,P_{\xi}})^{-1}\sigma(X_{\theta}^{t,x,P_{\xi}}\!,\!P_{X_{\theta}^{t,\xi}}),
\partial_{x}Z_s^{t,x,P_{\xi}}(\partial_{x}X_{\theta}^{t,x,P_{\xi}})^{-1}\sigma(X_{\theta}^{t,x,P_{\xi}}\!,\!P_{X_{\theta}^{t,\xi}}))$
is a solution of \eqref{eq9.9}. From the uniqueness of the solution of equation \eqref{eq9.9} we get
\begin{equation}\label{eq9.10}
  \begin{split}
D_{\theta}[Y_s^{t,x,P_{\xi}}]= \partial_{x}Y_s^{t,x,P_{\xi}}(\partial_{x}X_{\theta}^{t,x,P_{\xi}})^{-1}\sigma(X_{\theta}^{t,x,P_{\xi}},P_{X_{\theta}^{t,\xi}}),\
s\in[t,T],\ P\mbox{-}a.s.
  \end{split}
\end{equation}
Hence,
\begin{equation*}
  \begin{split}
&\ Z_{\theta}^{t,x,P_{\xi}} =P\mbox{-}\displaystyle{\lim_{\theta<s\downarrow \theta}}D_{\theta}[Y_s^{t,x,P_{\xi}}]
=\partial_{x}Y_{\theta}^{t,x,P_{\xi}}(\partial_{x}X_{\theta}^{t,x,P_{\xi}})^{-1}\sigma(X_{\theta}^{t,x,P_{\xi}},P_{X_{\theta}^{t,\xi}})\\
=&\ (\partial_{x}Y)_{\theta}^{\theta,X_{\theta}^{t,x,P_{\xi}},P_{X_{\theta}^{t,\xi}}}\sigma(X_{\theta}^{t,x,P_{\xi}},P_{X_{\theta}^{t,\xi}})
=(\partial_xV)(\theta,X_{\theta}^{t,x,P_{\xi}},P_{X_{\theta}^{t,\xi}})\sigma(X_\theta{}^{t,x,P_{\xi}},P_{X_{\theta}^{t,\xi}}),\ d\theta dP\mbox{-}a.e.
  \end{split}
\end{equation*}
Moreover,
\begin{equation*}
  \begin{split}
&\ E\[|Z_{\theta}^{t,x,P_{\xi}} -(\partial_xV)(\theta,x,P_{\xi})\sigma(x,P_{\xi})|^2\]\\
=&\ E\[|(\partial_{x}Y)_{\theta}^{\theta,X_{\theta}^{t,x,P_{\xi}},P_{X_{\theta}^{t,\xi}}}\sigma(X_{\theta}^{t,x,P_{\xi}},P_{X_{\theta}^{t,\xi}})
-\partial_{x}Y_{\theta}^{\theta,x,P_{\xi}}\sigma(x,P_{\xi})|^2\]\\
\leq&\ C E\[|(\partial_{x}Y)_{\theta}^{\theta,X_{\theta}^{t,x,P_{\xi}},P_{X_{\theta}^{t,\xi}}}-\partial_{x}Y_{\theta}^{\theta,x,P_{\xi}}|^2\]
       + C \(E\[|\sigma(X_{\theta}^{t,x,P_{\xi}},P_{X_{\theta}^{t,\xi}})-\sigma(x,P_{\xi})|^4\]\)^{\frac{1}{2}}\\
\leq&\ C(\theta-t),\ 0\leq t\leq \theta\leq T.
  \end{split}
\end{equation*}
\end{proof}

If we tried now to translate the approach of Pardoux and Peng \cite{PP1994} to our setting,
in order to prove that $V$ is a classical solution of SPDE \eqref{eq9.1}, we should
show that $V(t,\cdot,\cdot)\in C_{b}^{2,2}(\mathbb{R}^d\times\mathcal{P}_2(\mathbb{R}^d))$.
However, in \cite{PP1994} Kolmogorov's continuity criterion played a crucial role for the
proof. Here, in our mean-field context we meet the problem that $(t,x,P_{\xi})\in[0,T]\times\mathbb{R}^d \times\mathcal{P}_2(\mathbb{R}^d)$
runs an infinite dimensional space which excludes the application of Kolmogorov's continuity criterion.
The consequence is that we have to content with the continuity and differentiability of first and second order derivatives of $V(t,\cdot,\cdot)$
in the only $L^2$-sense. However, we can make the following observation.
\begin{lemma} \label{le9.1}
Let $\varphi:\Omega\times\mathbb{R}^d\rightarrow\mathbb{R}$ be $\mathcal{F}\otimes\mathcal{B}(\mathbb{R}^d)$-measurable and
$x\rightarrow \varphi(\cdot,x)$ $L^2$-differentiable. Then, for all $\eta\in L^{\infty}(\mathcal{F};\mathbb{R})$,
the deterministic function $\Psi(x):=E[\varphi(\cdot,x)\cdot\eta]$,
$x\in\mathbb{R}^d$, is differentiable w.r.t. $x$ on $\mathbb{R}^d$, and
$\partial_x\Psi(x)=E[\partial_x\varphi(\cdot,x)\cdot\eta]$, $x\in\mathbb{R}^d$,
where $\partial_x\varphi(x)$ denotes the $L^2$-derivative of $\varphi(\cdot,\cdot)$ at $x$.
\end{lemma}

\begin{proof}
For simplicity, let $d=1$. For all $\eta\in L^{2}(\mathcal{F})$, $x\in\mathbb{R}$:
\begin{equation*}
\begin{split}
 \frac{1}{q}(\Psi(x+q)-\Psi(x))&\ =E\[\(\frac{1}{q}(\varphi(\cdot,x+q)-\varphi(\cdot,x))-\partial_x\varphi(\cdot,x)\)\cdot\eta\]
  +E[\partial_x\varphi(\cdot,x)\cdot\eta]\\
  &\ \rightarrow E[\partial_x\varphi(\cdot,x)\cdot\eta],\ \mbox{as}\ q\rightarrow0.
\end{split}
\end{equation*}
\end{proof}

Combining the Propositions \ref{prop6.1+2} and  \ref{prop6.6} in Appendix A.3 and the Theorems \ref{th8.1} and \ref{th8.2} in Section 8,
we know that, for all $t\in[0,T]$,  (i) $x\rightarrow V(t,x,P_{\xi})$ is twice $L^2$-differentiable,
 (ii) $V(t,x,\cdot):\mathcal{P}_2(\mathbb{R}^d)\rightarrow \mathbb{R}$ is differentiable,
 (iii) $y\rightarrow (\partial_{\mu}V)(t,x,P_{\xi},y)$ is $L^2$-differentiable, and we have
(iv) the Lipschitz property in $L^2$ of all these derivatives (with Lipschitz constants independent of $t$).

Now, for $\eta\in L^{\infty}(\mathcal{F};\mathbb{R})$, we define $\Psi(t,x,P_{\xi}):=\Psi_{\eta}(t,x,P_{\xi}):=E[V(t,x,P_{\xi})\cdot\eta]$,
$(t,x,P_{\xi})\in[0,T]\times\mathbb{R}^d\times \mathcal{P}_2(\mathbb{R}^d)$. It can easily be verified that $\Psi(t,\cdot,\cdot)\in C_b^{2,2}(\mathbb{R}^d\times\mathcal{P}_2(\mathbb{R}^d))$. Moreover,
the following proposition studies the regularity properties with respect to $t$ of $\Psi(t,x,P_{\xi})$.
For this we make the following additional assumption on $\eta$:\\
\noindent\textbf{Assumption (H9.1)} The random variable $\eta\in L^{2}(\Omega,\mathcal{F}_{0,T}^{B},P;\mathbb{R})$ is such that,
for the ($\mathbb{F}^{B}=(\mathcal{F}_{s,T}^{B})_{0\leq s\leq T}$)-adapted process $\theta^{\eta}\in \mathcal{H}_{\mathcal{F}_{\cdot,T}^{B}}^2(0,T;\mathbb{R})$ with
$\displaystyle\eta=E[\eta]+\int_0^T\theta^{\eta}_sd\overleftarrow{B_s}$, there exists a constant $C_{\eta}\in \mathbb{R}_{+}$, such that
$|\theta^{\eta}_s|\leq C_{\eta}$, $dsdP$-a.e.

\begin{proposition} \label{prop9.2}
Under assumptions (H4.2), (H7.1), (H8.1) and (H8.2), for all $\eta\in L^{\infty}(\mathcal{F};\mathbb{R})$ with (H9.1),
$\Psi\in C_b^{0,2,2}([0,T]\times\mathbb{R}^d\times\mathcal{P}_2(\mathbb{R}^d))$, and for all
$\varphi\in\{\Psi,\partial_{x}\Psi,\partial_{xx}^2\Psi,\partial_{\mu}\Psi,\partial_{y}(\partial_{\mu}\Psi)\}$ it holds
\begin{equation}\label{eq9.11}
  |\varphi(t+q,x,P_{\xi},y)-\varphi(t,x,P_{\xi},y)|\leq C_{\eta}'\sqrt{q},\ 0\leq t\leq t+q\leq T,\
  (x,P_{\xi},y)\in\mathbb{R}^d \times\mathcal{P}_2(\mathbb{R}^d)\times \mathbb{R}^d,
\end{equation}
where $\partial_{\mu}\Psi(t,x,P_{\xi},y)=E[\partial_{\mu}V(t,x,P_{\xi},y)\cdot\eta]$,
and the constant $C_{\eta}'\in \mathbb{R}_{+}$ depends on $\eta$ and $C_{\eta}$.
\end{proposition}
\begin{proof}
Let us prove \eqref{eq9.11} in four steps.
For simplicity, we suppose $d=l=1$, $b=0$, $f=0$, $g^1=0$, and $\Phi(x,\mu)=\Phi(x)$.

\noindent\textbf{Step 1.} From Proposition \ref{prop4.4}, we know
\begin{equation}\label{eq9.12}
  E[|V(t+q,x,P_{\xi})-V(t,x,P_{\xi})|^p]\leq C q^{\frac{p}{2}},\ (t,x,P_{\xi})\in[0,T]\times\mathbb{R} \times\mathcal{P}_2(\mathbb{R}),\ p\in[2,p_0].
\end{equation}
It follows that $\Psi(\cdot,x,\mu)\in C^{0}([0,T])$, for all $(x,P_{\xi})\in\mathbb{R} \times\mathcal{P}_2(\mathbb{R})$, and
\begin{equation}\label{eq9.13}
  |\Psi(t+q,x,P_{\xi})-\Psi(t,x,P_{\xi})|\leq  C_{\eta}'\sqrt{q},\ 0\leq t\leq t+q\leq T,\ (x,P_{\xi})\in\mathbb{R} \times\mathcal{P}_2(\mathbb{R}).
\end{equation}

\noindent\textbf{Step 2.}
From \eqref{eq4.17} and the flow property of $Y_{\cdot}^{t,x,P_{\xi}}$ we have
\begin{equation}\label{eq9.14}
  \begin{split}
 &\   \partial_xV(t+q,x,P_{\xi})-\partial_xV(t,x,P_{\xi})=\partial_xY_{t+q}^{t+q,x,P_{\xi}}-\partial_xY_t^{t,x,P_{\xi}}\\
=&\ \(\partial_xY_{t+q}^{t+q,x,P_{\xi}}-(\partial_xY)_{t+q}^{t+q,X_{t+q}^{t,x,P_{\xi}},P_{X_{t+q}^{t,\xi}}}\)
    +(\partial_xY)_{t+q}^{t+q,X_{t+q}^{t,x,P_{\xi}},P_{X_{t+q}^{t,\xi}}}\(1-\partial_xX_{t+q}^{t,x,P_{\xi}}\)\\
  &\  +\(\partial_xY_{t+q}^{t,x,P_{\xi}}-\partial_xY_{t}^{t,x,P_{\xi}}\).
  \end{split}
\end{equation}
On the other hand, from Lemma \ref{le3.1} and the Propositions \ref{prop6.1+2}, \ref{prop5.1} and \ref{prop8.1}, we obtain, for all $p\in[2,\frac{p_0}{2}]$,
\begin{equation*}
  \begin{split}
  &\ \mbox{(i)}\ E\[|\partial_xY_{t+q}^{t+q,x,P_{\xi}}-(\partial_xY)_{t+q}^{t+q,X_{t+q}^{t,x,P_{\xi}},P_{X_{t+q}^{t,\xi}}}|^p\]
     \leq C_p\(E[|x-X_{t+q}^{t,x,P_{\xi}}|^p]+W_2(P_{\xi},P_{X_{t+q}^{t,\xi}})^p\)\leq C_pq^{\frac{p}{2}},\\
  &\ \mbox{(ii)}\ E\[|(\partial_xY)_{t+q}^{t+q,X_{t+q}^{t,x,P_{\xi}},P_{X_{t+q}^{t,\xi}}}\(1-\partial_xX_{t+q}^{t,x,P_{\xi}}\)|^p\]
      \leq C_p(E[|1-\partial_xX_{t+q}^{t,x,P_{\xi}}|^{2p}])^{\frac{1}{2}}\leq C_pq^{\frac{p}{2}},\\
  &\ \mbox{(iii)}\ E\[|\partial_xY_{t+q}^{t,x,P_{\xi}}-\partial_xY_{t}^{t,x,P_{\xi}}|^p\]
       \leq C_pq^{\frac{p-2}{2}}\int_t^{t+q}E[|\partial_xZ_{s}^{t,x,P_{\xi}}|^p]ds\leq C_pq^{\frac{p}{2}}.
  \end{split}
\end{equation*}
Hence, for all $p\in[2,\frac{p_0}{2}]$,
\begin{equation*}
  \begin{split}
   E\[|\partial_xV(t+q,x,P_{\xi})-\partial_xV(t,x,P_{\xi})|^p\]
     \leq  C_pq^{\frac{p}{2}},\ 0\leq t\leq t+q\leq T,\ (x,P_{\xi})\in\mathbb{R} \times\mathcal{P}_2(\mathbb{R}).
  \end{split}
\end{equation*}
It follows that
\begin{equation}\label{eq9.15}
  |\partial_x\Psi(t+q,x,P_{\xi})-\partial_x\Psi(t,x,P_{\xi})|\leq  C_{\eta}'\sqrt{q},\ 0\leq t\leq t+q\leq T,\ (x,P_{\xi})\in\mathbb{R} \times\mathcal{P}_2(\mathbb{R}).
\end{equation}

\noindent\textbf{Step 3.} To prove $|\partial_{xx}^2\Psi(t+q,x,P_{\xi})-\partial_{xx}^2\Psi(t,x,P_{\xi})|\leq  C_{\eta}'\sqrt{q}$.
We observe that,
\begin{equation}\label{eq9.16}
   \partial_{xx}^2V(t+q,x,P_{\xi})-\partial_{xx}^2V(t,x,P_{\xi})
   =\(\partial_{xx}^2[Y_{t+q}^{t+q,x,P_{\xi}}]-\partial_{xx}^2[Y_{t+q}^{t,x,P_{\xi}}]\)+\(\partial_{xx}^2[Y_{t+q}^{t,x,P_{\xi}}]-\partial_{xx}^2[Y_{t}^{t,x,P_{\xi}}]\).
\end{equation}
On one hand, as $Y_{t+q}^{t,x,P_{\xi}}=Y_{t+q}^{t+q,X_{t+q}^{t,x,P_{\xi}},P_{X_{t+q}^{t,\xi}}}$,
and $\partial_x[Y_{t+q}^{t,x,P_{\xi}}]=(\partial_xY)_{t+q}^{t+q,X_{t+q}^{t,x,P_{\xi}},P_{X_{t+q}^{t,\xi}}}\cdot\partial_xX_{t+q}^{t,x,P_{\xi}}$,
we have $\partial_{xx}^2[Y_{t+q}^{t,x,P_{\xi}}]=(\partial_xY)_{t+q}^{t+q,X_{t+q}^{t,x,P_{\xi}},P_{X_{t+q}^{t,\xi}}}\cdot\partial_{xx}^2X_{t+q}^{t,x,P_{\xi}}
+(\partial_{xx}^2Y)_{t+q}^{t+q,X_{t+q}^{t,x,P_{\xi}},P_{X_{t+q}^{t,\xi}}}\cdot(\partial_xX_{t+q}^{t,x,P_{\xi}})^2$.
From Lemma \ref{le3.1}, the Theorems \ref{th7.1} and \ref{th8.1} as well as the Propositions \ref{prop6.1+2} and \ref{prop5.1}, we deduce that
\begin{equation}\label{eq9.17}
  \begin{split}
 &\   E\[|\partial_{xx}^2[Y_{t+q}^{t+q,x,P_{\xi}}]-\partial_{xx}^2[Y_{t+q}^{t,x,P_{\xi}}]|^p\]\\
 \leq &\ C_p E\[|\partial_{xx}^2[Y_{t+q}^{t+q,x,P_{\xi}}]-(\partial_{xx}^2Y)_{t+q}^{t+q,X_{t+q}^{t,x,P_{\xi}},P_{X_{t+q}^{t,\xi}}}|^p\]
 +C_p \(E\[|(\partial_xX_{t+q}^{t,x,P_{\xi}})^2-1|^{2p}\]\)^{\frac{1}{2}}\\
&\  +C_p \(E\[|\partial_{xx}^2X_{t+q}^{t,x,P_{\xi}}|^{2p}\]\)^{\frac{1}{2}}\\
 \leq &\ C_pq^{\frac{p}{2}}+C_p \(E\[|\partial_xX_{t+q}^{t,x,P_{\xi}}-1|^{4p}\]\)^{\frac{1}{4}}\cdot \(E\[|\partial_xX_{t+q}^{t,x,P_{\xi}}+1|^{4p}\]\)^{\frac{1}{4}}\\
 \leq &\ C_pq^{\frac{p}{2}},\ 0\leq t\leq t+q\leq T,\ (x,P_{\xi})\in\mathbb{R} \times\mathcal{P}_2(\mathbb{R}),\ p\in[2,\frac{p_0}{8}].
  \end{split}
\end{equation}
On the other hand,
\begin{equation*}
\begin{split}
\partial_{xx}^2Y_{t+q}^{t,x,P_{\xi}}-\partial_{xx}^2Y_{t}^{t,x,P_{\xi}}=
 -\int_{t}^{t+q}g^2(P_{X_r^{t,\xi}}) \partial_{xx}^2Z_r^{t,x,P_{\xi}}d\overleftarrow{B_r}
 +\int_{t}^{t+q}\partial_{xx}^2Z_r^{t,x,P_{\xi}}dW_r,
\end{split}
\end{equation*}
and, combining with assumption (H9.1), we get
\begin{equation*}
\begin{split}
&\ E\[\(\partial_{xx}^2Y_{t+q}^{t,x,P_{\xi}}-\partial_{xx}^2Y_{t}^{t,x,P_{\xi}}\)\cdot\eta\]\\
=&\ -E\[\(\int_{t}^{t+q}g^2(P_{X_r^{t,\xi}}) \partial_{xx}^2Z_r^{t,x,P_{\xi}}d\overleftarrow{B_r}\)\cdot\eta\]
 +E\[E\[\int_{t}^{t+q}\partial_{xx}^2Z_r^{t,x,P_{\xi}}dW_r\big|\mathcal{F}_{0,T}^{B}\]\eta\]\\
=&\ -E\[\int_{t}^{t+q}g^2(P_{X_r^{t,\xi}}) \partial_{xx}^2Z_r^{t,x,P_{\xi}}\cdot\theta^{\eta}_rdr\].
\end{split}
\end{equation*}
Therefore, we deduce from Theorem \ref{th8.1} that
\begin{equation}\label{eq9.18}
\begin{split}
&\ \Big|E\[\(\partial_{xx}^2Y_{t+q}^{t,x,P_{\xi}}-\partial_{xx}^2Y_{t}^{t,x,P_{\xi}}\)\eta\]\Big|
\leq C_{\eta}'E\[\int_{t}^{t+q}| \partial_{xx}^2Z_r^{t,x,P_{\xi}}|dr\]\\
\leq&\ C_{\eta}'\sqrt{q}\(\!E\[\!\int_{t}^{t+q}| \partial_{xx}^2Z_r^{t,x,P_{\xi}}|^2dr\]\)^{\frac{1}{2}}
\leq  C_{\eta}'\sqrt{q}.
\end{split}
\end{equation}
Consequently, making use of \eqref{eq9.16}, \eqref{eq9.17} and \eqref{eq9.18}, we obtain
\begin{equation}\label{eq9.19}
\begin{split}
  &\ \Big| \partial_{xx}^2\Psi(t+q,x,P_{\xi})-\partial_{xx}^2\Psi(t,x,P_{\xi})\Big|\\
  \leq &\  E\[\Big|\partial_{xx}^2[Y_{t+q}^{t+q,x,P_{\xi}}]-\partial_{xx}^2[Y_{t+q}^{t,x,P_{\xi}}]\Big|\cdot|\eta|\]
  +\Big|E\[\!\(\!\partial_{xx}^2Y_{t+q}^{t,x,P_{\xi}}\!-\!\partial_{xx}^2Y_{t}^{t,x,P_{\xi}}\)\eta\]\Big|\\
\leq&\  C_{\eta}'\sqrt{q},\ 0\leq t\leq t+q\leq T,\ (x,P_{\xi})\in\mathbb{R} \times\mathcal{P}_2(\mathbb{R}).
\end{split}
\end{equation}

\noindent\textbf{Step 4.}
As the expectation $E[\cdot]:L^2(\mathcal{F})\rightarrow \mathbb{R}$ is a bounded linear operator, it follows from Proposition \ref{prop6.6} that
$ L^{2}(\mathcal{G}_t;\mathbb{R})\ni\xi\rightarrow\widetilde{\Psi}(t,x,\xi):=\Psi(t,x,P_{\xi})=E[V(t,x,P_{\xi})\cdot\eta]$
 is Fr\'{e}chet differentiable and, for all $\zeta\in L^{2}(\mathcal{G}_t;\mathbb{R})$,
\begin{equation*}
\begin{split}
&\  D_{\xi}[\widetilde{\Psi}(t,x,\xi)](\zeta)
  =E\[D_{\xi}[V(t,x,P_{\xi})](\zeta)\cdot\eta\]=E\[\overline{E}\[\partial_{\mu}V(t,x,P_{\xi})(\overline{\xi})\cdot\overline{\zeta}\]\cdot\eta\]\\
=&\ \overline{E}\[E\[\partial_{\mu}V(t,x,P_{\xi})(\overline{\xi})\cdot\eta\]\cdot\overline{\zeta}\]
=\overline{E}\[E\[\partial_{\mu}Y_t^{t,x,P_{\xi}}(\overline{\xi})\cdot\eta\]\cdot\overline{\zeta}\],
\end{split}
\end{equation*}
that is, $\partial_{\mu}\Psi(t,x,P_{\xi},y)=E[\partial_{\mu}Y_t^{t,x,P_{\xi}}(y)\cdot\eta].$
Furthermore, with similar arguments as above we show that
\begin{equation}\label{eq9.20}
\begin{split}
  &\ | \partial_{\mu}\Psi(t+q,x,P_{\xi},y)-\partial_{\mu}\Psi(t,x,P_{\xi},y)|
+ | \partial_{y}(\partial_{\mu}\Psi)(t+q,x,P_{\xi},y)-\partial_{y}(\partial_{\mu} \Psi)(t,x,P_{\xi},y)|\\
\leq&\  C_{\eta}'\sqrt{q},\ 0\leq t\leq t+q\leq T,\ (x,P_{\xi},y)\in\mathbb{R} \times\mathcal{P}_2(\mathbb{R})\times\mathbb{R}.
\end{split}
\end{equation}
\end{proof}

\begin{definition} \label{def9.1}
 We say that random field $\varphi$ belongs to $\mathfrak{C}^{0,2,2}(\Omega\times[0,T]\times\mathbb{R}^{d}\times\mathcal{P}_{2}(\mathbb{R}^{d}))$, if $\varphi:\Omega\times[0,T]\times\mathbb{R}^{d}\times\mathcal{P}_{2}(\mathbb{R}^{d})\rightarrow\mathbb{R}$ satisfies:\\
\indent \emph{(i)} $\varphi(t,x,\mu)$ is $\mathcal{F}_{t,T}^B$-measurable, $(t,x,\mu)\in[0,T]\times\mathbb{R}^{d}\times\mathcal{P}_{2}(\mathbb{R}^{d})$;\\
\indent \emph{(ii)} $x\rightarrow \varphi(t,x,\mu)$ is twice continuously $L^2$-differentiable;\\
\indent \emph{(iii)} $\mu\rightarrow \varphi(t,x,\mu)$ is differentiable;\\
\indent \emph{(iv)} $y\rightarrow \partial_{\mu}\varphi(t,x,\mu,y)$ is continuously $L^2$-differentiable;\\
\indent \emph{(v)} The first and second order derivatives are $L^2$-continuous on $[0,T]\times\mathbb{R}^{d}\times\mathcal{P}_{2}(\mathbb{R}^{d})\times\mathbb{R}^{d}$;\\
\indent \emph{(vi)} $\Gamma\in C_b^{0,2,2}([0,T]\times\mathbb{R}^{d}\times\mathcal{P}_{2}(\mathbb{R}^{d}))$ [see Definition \ref{def2.2}],
         where $\Gamma(t,x,\mu):=E[\varphi(t,x,\mu)\cdot\eta]$, for all $\eta\in L^{\infty}(\mathcal{F}_{0,T}^B;\mathbb{R})$ satisfying (H9.1), and all $(t,x,\mu)\in[0,T]\times\mathbb{R}^{d}\times\mathcal{P}_{2}(\mathbb{R}^{d})$.
\end{definition}

Now we are able to establish and to prove our main result.

\begin{theorem} \label{th9.1}
Under the Assumptions (H4.2), (H7.1), (H8.1) and (H8.2), $V\in \mathfrak{C}^{0,2,2}(\Omega\times[0,T]\times\mathbb{R}^{d}\times\mathcal{P}_{2}(\mathbb{R}^{d}))$
is a classical solution of backward SPDE \eqref{eq9.1}, and it is
unique in $\mathfrak{C}^{0,2,2}(\Omega\times[0,T]\times\mathbb{R}^{d}\times\mathcal{P}_{2}(\mathbb{R}^{d}))$.
\end{theorem}
\begin{proof}
For simplicity of redaction but without loss of generality, let us restrict to the dimensions
$d = 1$, $l=1$, and to the coefficients $b=0$, $f=0$, $g^1=0$, $h(P_{\Pi_s^{t,\xi}})=h(P_{X_s^{t,\xi}})$, $\Phi(x,\mu)=\Phi(x)$.

We prove that $V(t,x,P_{\xi})$ is a solution of \eqref{eq9.1}. Let $\eta\in L^{\infty}(\mathcal{F}_{0,T}^B;\mathbb{R})$ be such that (H9.1) holds true, and
$s\rightarrow\theta_s^{\eta}$ is continuous. It follows from the Propositions \ref{prop6.1+2}, \ref{prop6.6} and \ref{prop9.2}
and the Theorems \ref{th8.1} and \ref{th8.2} that $V\in \mathfrak{C}^{0,2,2}(\Omega\times[0,T]\times\mathbb{R}\times\mathcal{P}_{2}(\mathbb{R}))$.
 For $0\leq t< t+q\leq T$, we have
\begin{equation}\label{eq9.21}
  \begin{split}
 &\   \Psi(t+q,x,P_{\xi})-\Psi(t,x,P_{\xi})
 =E[(V(t+q,x,P_{\xi})-V(t,x,P_{\xi}))\cdot\eta]\\
 =&\ -E[(V(t+q,X_{t+q}^{t,x,P_{\xi}},P_{X_{t+q}^{t,\xi}})-V(t+q,x,P_{\xi}))\cdot\eta]+E[(Y_{t+q}^{t,x,P_{\xi}}-Y_t^{t,x,P_{\xi}})\cdot\eta].
  \end{split}
\end{equation}
Note that, as $X_{t+q}^{t,x,P_{\xi}}$ is $\mathcal{F}_{t+q}^W$-measurable, and, hence, independent of $\mathcal{F}_{0,T}^B$, while
$V(t+q,x,P_{\xi})$ is $\mathcal{F}_{t+q,T}^B$-measurable, $x\in\mathbb{R}$,
$E[V(t+q,X_{t+q}^{t,x,P_{\xi}},P_{X_{t+q}^{t,\xi}})\cdot\eta]=E[\Psi(t+q,X_{t+q}^{t,x,P_{\xi}},P_{X_{t+q}^{t,\xi}})]$.
As $\Psi\in C_b^{0,2,2}([0,T]\times\mathbb{R}\times\mathcal{P}_{2}(\mathbb{R}))$, we can apply Theorem \ref{th2.1} (It\^{o}'s formula) as follows:
\begin{equation}\label{eq9.22}
  \begin{split}
 &\   \Psi(t+q,x,P_{\xi})-\Psi(t,x,P_{\xi}) =E[(V(t+q,x,P_{\xi})-V(t,x,P_{\xi}))\cdot\eta]\\
 =&\ E[(Y_{t+q}^{t,x,P_{\xi}}-Y_t^{t,x,P_{\xi}})\cdot\eta]-E[\Psi(t+q,X_{t+q}^{t,x,P_{\xi}},P_{X_{t+q}^{t,\xi}})-\Psi(t+q,x,P_{\xi})]\\
 =&\ E[(Y_{t+q}^{t,x,P_{\xi}}-Y_t^{t,x,P_{\xi}})\cdot\eta]
      -\int_t^{t+q} E\[\frac{1}{2}(\partial_{xx}^2\Psi)(t+q,X_{s}^{t,x,P_{\xi}},P_{X_{s}^{t,\xi}})\sigma(X_s^{t,x,P_{\xi}},P_{X_s^{t,\xi}})^2\\
 &\ +\widehat{E}\[\frac{1}{2}\partial_{y}(\partial_{\mu}\Psi)(t+q,X_{s}^{t,x,P_{\xi}},P_{X_{s}^{t,\xi}},\widehat{X}_{s}^{t,\widehat{\xi}})
     \sigma(\widehat{X}_{s}^{t,\widehat{\xi}},P_{X_{s}^{t,\xi}})^2 \]\]ds.
  \end{split}
\end{equation}
But, from the BDSDE for $(Y^{t,x,P_{\xi}},Z^{t,x,P_{\xi}})$,
\begin{equation}\label{eq9.23}
  \begin{split}
&\ E[(Y_{t+q}^{t,x,P_{\xi}}-Y_t^{t,x,P_{\xi}})\cdot\eta]\\
 =&\  -\!E\[\eta\int_t^{t+q}\!g^2(P_{X_r^{t,\xi}})Z_r^{t,x,P_{\xi}}d\overleftarrow{B_r}\]
    \!-\!E\[\eta\int_t^{t+q}\!h(P_{X_r^{t,\xi}})d\overleftarrow{B_r}\]
    \!+\!E\[\eta\int_t^{t+q}\!Z_r^{t,x,P_{\xi}}dW_r\]\\
 =&\  -E\[\int_t^{t+q}\!g^2(P_{X_r^{t,\xi}})Z_r^{t,x,P_{\xi}}\cdot\theta_r^{\eta}dr\]
    \!-\!E\[\int_t^{t+q}\!h(P_{X_r^{t,\xi}})\cdot\theta_r^{\eta}dr\] \\
 =&\  -E\[\int_t^{t+q}\!g^2(P_{\xi})(\partial_xV)(r,x,P_{\xi})\sigma(x,P_{\xi})\cdot\theta_r^{\eta}dr\]
    \!-\!E\[\int_t^{t+q}\!h(P_{\xi})\cdot\theta_r^{\eta}dr\] \!-\!R_{t,t+q},
  \end{split}
\end{equation}
where
\begin{equation*}
  \begin{split}
R_{t,t+q}:=&\ E\[\int_t^{t+q}\!\(g^2(P_{X_r^{t,\xi}})Z_r^{t,x,P_{\xi}}-g^2(P_{\xi})(\partial_xV)(r,x,P_{\xi})\sigma(x,P_{\xi})\)\cdot\theta_r^{\eta}dr\]\\
    &\ +\!E\[\int_t^{t+q}\!\(h(P_{X_r^{t,\xi}})-h(P_{\xi})\)\cdot\theta_r^{\eta}dr\].
  \end{split}
\end{equation*}
Thanks to \eqref{eq9.3} in Proposition \ref{prop9.1}, we obtain $ |R_{t,t+q}| \leq Cq^{\frac{3}{2}}$.
Using our estimates of $X_r^{t,x,P_{\xi}}$, $P_{X_r^{t,\xi}}$ and the $\frac{1}{2}$-H\"{o}lder continuity of $t\rightarrow\varphi(t,x,P_{\xi})$,
for $\varphi=\partial_{xx}^2\Psi,\partial_y(\partial_{\mu}\Psi)$ (see Lemma \ref{le3.1} and Proposition \ref{prop9.2}), we get now from \eqref{eq9.22}
\begin{equation}\label{eq9.24}
  \begin{split}
&\ \Psi(t+q,x,P_{\xi})-\Psi(t,x,P_{\xi})\\
 =&\  -E\[\int_t^{t+q}\!g^2(P_{\xi})(\partial_xV)(r,x,P_{\xi})\sigma(x,P_{\xi})\cdot\theta_r^{\eta}dr\]
    \!-\!E\[\int_t^{t+q}\!h(P_{\xi})\cdot\theta_r^{\eta}dr\] \\
 &\   -\int_t^{t+q}\(\frac{1}{2}(\partial_{xx}^2\Psi)(r,x,P_{\xi})\sigma(x,P_{\xi})^2
+\widehat{E}\[\frac{1}{2}\partial_{y}(\partial_{\mu}\Psi)(r,x,P_{\xi},\widehat{\xi})\sigma(\widehat{\xi},P_{\xi})^2 \]\)dr \!-\!R_{t,t+q}',
  \end{split}
\end{equation}
with $|R_{t,t+q}'| \leq Cq^{\frac{3}{2}}$. Recall the boundedness and the continuity of $s\rightarrow\theta_s^{\eta}$.
Therefore, it can easily be verified that $\Psi(t,x,P_{\xi})$ is differentiable w.r.t. $t$, and
\begin{equation}\label{eq9.25}
  \begin{split}
 \partial_t\Psi(t,x,P_{\xi}) =&\  -\frac{1}{2}(\partial_{xx}^2\Psi)(t,x,P_{\xi})\sigma(x,P_{\xi})^2
-\frac{1}{2}\widehat{E}\[\partial_{y}(\partial_{\mu}\Psi)(t,x,P_{\xi},\widehat{\xi})\sigma(\widehat{\xi},P_{\xi})^2 \]\\
 &\  -E\[g^2(P_{\xi})(\partial_xV)(t,x,P_{\xi})\sigma(x,P_{\xi})\cdot\theta_t^{\eta}\]
    \!-\!E\[h(P_{\xi})\cdot\theta_t^{\eta}\],\ t\in[0,T].
  \end{split}
\end{equation}
Hence, for all $\eta\in L^{\infty}(\mathcal{F}_{0,T}^B;\mathbb{R})$ that satisfy (H9.1) and are such that
$s\rightarrow\theta_s^{\eta}$ continuous,
\begin{equation}\label{eq9.26}
  \begin{split}
  &\  E[(V(t,x,P_{\xi})-\Phi(x,P_{\xi}))\cdot\eta]= \Psi(t,x,P_{\xi})-\Psi(T,x,P_{\xi}) \\
=&\  \int_t^T\(\frac{1}{2}(\partial_{xx}^2\Psi)(r,x,P_{\xi})\sigma(x,P_{\xi})^2
+\frac{1}{2}\widehat{E}\[\partial_{y}(\partial_{\mu}\Psi)(r,x,P_{\xi},\widehat{\xi})\sigma(\widehat{\xi},P_{\xi})^2 \]\)dr\\
 &\  +E\[\int_t^Tg^2(P_{\xi})(\partial_xV)(r,x,P_{\xi})\sigma(x,P_{\xi})\cdot\theta_r^{\eta}dr\]
    \!+\!E\[\int_t^Th(P_{\xi})\cdot\theta_r^{\eta}dr\]\\
=&\  \int_t^T\(\frac{1}{2}(\partial_{xx}^2\Psi)(r,x,P_{\xi})\sigma(x,P_{\xi})^2
+\frac{1}{2}\widehat{E}\[\partial_{y}(\partial_{\mu}\Psi)(r,x,P_{\xi},\widehat{\xi})\sigma(\widehat{\xi},P_{\xi})^2 \]\)dr\\
 &\  +E\[\int_t^Tg^2(P_{\xi})(\partial_xV)(r,x,P_{\xi})\sigma(x,P_{\xi})d\overleftarrow{B_r}\cdot\eta\]
    \!+\!E\[\int_t^Th(P_{\xi})d\overleftarrow{B_r}\cdot\eta\]\\
=&\  E\[\eta\cdot I(t,x, P_{\xi})\] ,
  \end{split}
\end{equation}
where
\begin{equation*}
  \begin{split}
I(t,x, P_{\xi}):=&\ \int_t^T\(\frac{1}{2}(\partial_{xx}^2V)(r,x,P_{\xi})\sigma(x,P_{\xi})^2
+\frac{1}{2}\widehat{E}\[\partial_{y}(\partial_{\mu}V)(r,x,P_{\xi},\widehat{\xi})\sigma(\widehat{\xi},P_{\xi})^2 \]\)dr\\
 &\  +\int_t^Tg^2(P_{\xi})(\partial_xV)(r,x,P_{\xi})\sigma(x,P_{\xi})d\overleftarrow{B_r}
    \!+\!\int_t^Th(P_{\xi})d\overleftarrow{B_r}.
  \end{split}
\end{equation*}
But, since these $\eta$ satisfying (H9.1) such that $s\rightarrow\theta_s^{\eta}$ is continuous
form a dense subset of $L^2(\mathcal{F}_{0,T}^B;\mathbb{R})$, it follows that
\begin{equation*}
  \begin{split}
  V(t,x,P_{\xi})=&\ \Phi(x,P_{\xi})+\int_t^T\(\frac{1}{2}(\partial_{xx}^2V)(r,x,P_{\xi})\sigma(x,P_{\xi})^2
+\frac{1}{2}\widehat{E}\[\partial_{y}(\partial_{\mu}V)(r,x,P_{\xi},\widehat{\xi})\sigma(\widehat{\xi},P_{\xi})^2 \]\)dr\\
 &\  +\int_t^Tg^2(P_{\xi})(\partial_xV)(r,x,P_{\xi})\sigma(x,P_{\xi})d\overleftarrow{B_r}
    \!+\!\int_t^Th(P_{\xi})d\overleftarrow{B_r},\ t\in[0,T],\ P\mbox{-a.s.}
  \end{split}
\end{equation*}

Let us now prove the uniqueness of solution of backward SPDE \eqref{eq9.1} in $\mathfrak{C}^{0,2,2}(\Omega\times[0,T]\times\mathbb{R}\times\mathcal{P}_{2}(\mathbb{R}))$.
Let $U\in\mathfrak{C}^{0,2,2}(\Omega\times[0,T]\times\mathbb{R}\times\mathcal{P}_{2}(\mathbb{R}))$ be a solution of our backward SPDE \eqref{eq9.1},
\begin{equation}\label{eq9.27}
  \begin{split}
 &\  U(t+q,x,P_{\xi})-U(t,x,P_{\xi})\\
=&\ -\int_t^{t+q}\(\frac{1}{2}(\partial_{xx}^2U)(r,x,P_{\xi})\sigma(x,P_{\xi})^2
+\frac{1}{2}\widehat{E}\[\partial_{y}(\partial_{\mu}U)(r,x,P_{\xi},\widehat{\xi})\sigma(\widehat{\xi},P_{\xi})^2 \]\)dr\\
 &\  -\int_t^{t+q}\(g^2(P_{\xi})(\partial_xU)(r,x,P_{\xi})\sigma(x,P_{\xi})+h(P_{\xi})\)d\overleftarrow{B_r},\ 0\leq t< t+q\leq T.
  \end{split}
\end{equation}
Hence, for all $\eta\!\in\! L^{\infty}(\mathcal{F}_{0,T}^B;\mathbb{R})$ that satisfies (H9.1) such that
$s\!\rightarrow\!\theta_s^{\eta}$ is continuous, we have
$\Gamma(t,x,P_{\xi})\!:=\!\Gamma_{\eta}(t,x,P_{\xi}):=E[U(t,x,P_{\xi})\cdot\eta]$.
From \eqref{eq9.27} it follows that $\Gamma(t,x,P_{\xi})$ is differentiable w.r.t. $t$, and
\begin{equation}\label{eq9.28}
  \begin{split}
  &\ \partial_t\Gamma(t,x,P_{\xi}) =\partial_tE[U(t,x,P_{\xi})\cdot\eta]\\
=&\  -\frac{1}{2}E\[(\partial_{xx}^2U)(t,x,P_{\xi})\sigma(x,P_{\xi})^2\cdot\eta\]
-\frac{1}{2}\widehat{E}\[E\[\partial_{y}(\partial_{\mu}U)(t,x,P_{\xi},\widehat{\xi})\cdot\eta\]\sigma(\widehat{\xi},P_{\xi})^2 \]\\
 &\  -E\[\(g^2(P_{\xi})(\partial_xU)(t,x,P_{\xi})\sigma(x,P_{\xi})+h(P_{\xi})\)\cdot\theta_t^{\eta}\]\\
=&\ -\frac{1}{2}(\partial_{xx}^2\Gamma)(t,x,P_{\xi})\sigma(x,P_{\xi})^2
-\frac{1}{2}\widehat{E}\[\partial_{y}(\partial_{\mu}\Gamma)(t,x,P_{\xi},\widehat{\xi})\sigma(\widehat{\xi},P_{\xi})^2 \]\\
 &\  -E\[\(g^2(P_{\xi})(\partial_xU)(t,x,P_{\xi})\sigma(x,P_{\xi})+h(P_{\xi})\)\cdot\theta_t^{\eta}\],\ t\in[0,T].
  \end{split}
\end{equation}
Hence, as $X_s^{t,x,P_{\xi}}$ is independent of $\mathcal{F}_{0,T}^B$,
\begin{equation}\label{eq9.29}
  \begin{split}
  (\partial_t\Gamma)(s,&X_s^{t,x,P_{\xi}},P_{X_s^{t,\xi}})
= -\frac{1}{2}(\partial_{xx}^2\Gamma)(s,X_s^{t,x,P_{\xi}},P_{X_s^{t,\xi}}) \sigma(X_s^{t,x,P_{\xi}},P_{X_s^{t,\xi}})^2\\
 &\     -\frac{1}{2}\widehat{E}\[\partial_{y}(\partial_{\mu}\Gamma)(s,X_s^{t,x,P_{\xi}},P_{X_s^{t,\xi}},\widehat{X}_s^{t,\widehat{\xi}})
     \sigma(\widehat{X}_s^{t,\widehat{\xi}},P_{X_s^{t,\xi}})^2 \]\\
 &\  -E\[\(g^2(P_{X_s^{t,\xi}})(\partial_xU)(s,X_s^{t,x,P_{\xi}},P_{X_s^{t,\xi}}) \sigma(X_s^{t,x,P_{\xi}},P_{X_s^{t,\xi}})+h(P_{X_s^{t,\xi}})\)
     \cdot\theta_s^{\eta}\Big|\mathcal{F}_{T}^W\].
  \end{split}
\end{equation}
On the other hand, as $\Gamma\in C_{b}^{0,2,2}([0,T]\times\mathbb{R}\times\mathcal{P}_{2}(\mathbb{R}))$,
 from the It\^{o} formula (It's only the classical It\^{o} formula from \cite{BLPR2017}):
\begin{equation*}\label{eq9.30}
  \begin{split}
  d[\Gamma(s,X_s^{t,x,P_{\xi}},P_{X_s^{t,\xi}})]
\!=\!&\ \((\partial_t\Gamma)(s,X_s^{t,x,P_{\xi}},P_{X_s^{t,\xi}})
  \!+\!\frac{1}{2}(\partial_{xx}^2\Gamma)(s,X_s^{t,x,P_{\xi}},P_{X_s^{t,\xi}}) \sigma(X_s^{t,x,P_{\xi}},P_{X_s^{t,\xi}})^2\\
 &\     +\frac{1}{2}\widehat{E}\[\partial_{y}(\partial_{\mu}\Gamma)(s,X_s^{t,x,P_{\xi}},P_{X_s^{t,\xi}},\widehat{X}_s^{t,\widehat{\xi}})
     \sigma(\widehat{X}_s^{t,\widehat{\xi}},P_{X_s^{t,\xi}})^2 \]\)ds\\
 &\  +(\partial_{x}\Gamma)(s,X_s^{t,x,P_{\xi}},P_{X_s^{t,\xi}}) \sigma(X_s^{t,x,P_{\xi}},P_{X_s^{t,\xi}})dW_s.
  \end{split}
\end{equation*}
Consequently, using \eqref{eq9.29}, we have
\begin{equation}\label{eq9.31}
  \begin{split}
 &\ d[\Gamma(s,X_s^{t,x,P_{\xi}},P_{X_s^{t,\xi}})]\\
=&\ -E\[\(g^2(P_{X_s^{t,\xi}})(\partial_xU)(s,X_s^{t,x,P_{\xi}},P_{X_s^{t,\xi}}) \sigma(X_s^{t,x,P_{\xi}},P_{X_s^{t,\xi}})+h(P_{X_s^{t,\xi}})\)
     \cdot\theta_s^{\eta}\Big|\mathcal{F}_{T}^W\]ds\\
 &\  +(\partial_{x}\Gamma)(s,X_s^{t,x,P_{\xi}},P_{X_s^{t,\xi}}) \sigma(X_s^{t,x,P_{\xi}},P_{X_s^{t,\xi}})dW_s, \ s\in[t,T].
  \end{split}
\end{equation}
Then, for all $\zeta\in L^{\infty}(\mathcal{F}_{T}^W;\mathbb{R})$, and for all $s\in[t,T]$,
as $(\zeta,X_{\cdot}^{t,x,P_{\xi}})$ is independent of $\mathcal{F}_{0,T}^B$,
\begin{equation*}\label{eq9.32}
  \begin{split}
 &\ E\[\!\(\!U(s,\!X_s^{t,x,P_{\xi}}\!,\!P_{X_s^{t,\xi}})\!\!-\!\!\Phi(X_T^{t,x,P_{\xi}}\!,\!P_{X_T^{t,\xi}})\!\)\eta\zeta\!\]
\!\!=\!\! E\[\!\(\!\Gamma(s,\!X_s^{t,x,P_{\xi}}\!,\!P_{X_s^{t,\xi}})\!\!-\!\!\Gamma(T,\!X_T^{t,x,P_{\xi}}\!,\!P_{X_T^{t,\xi}})\!\)\zeta\!\]\\
=&\ E\[\int_s^T\(g^2(P_{X_r^{t,\xi}})(\partial_xU)(r,X_r^{t,x,P_{\xi}},P_{X_r^{t,\xi}}) \sigma(X_r^{t,x,P_{\xi}},P_{X_r^{t,\xi}})+h(P_{X_r^{t,\xi}})\)
     \cdot\theta_r^{\eta}\cdot\zeta dr\]\\
  \end{split}
\end{equation*}
\begin{equation*}
  \begin{split}
 &\  -E\[\int_s^T(\partial_{x}\Gamma)(r,X_r^{t,x,P_{\xi}},P_{X_r^{t,\xi}}) \sigma(X_r^{t,x,P_{\xi}},P_{X_r^{t,\xi}})dW_r\cdot\zeta\]\\
=&\  E\[\int_s^T\(g^2(P_{X_r^{t,\xi}})(\partial_xU)(r,X_r^{t,x,P_{\xi}},P_{X_r^{t,\xi}}) \sigma(X_r^{t,x,P_{\xi}},P_{X_r^{t,\xi}})+h(P_{X_r^{t,\xi}})\)d\overleftarrow{B_r}
     \cdot\eta\cdot\zeta \]\\
 &\  -E\[\int_s^T(\partial_{x}U)(r,X_r^{t,x,P_{\xi}},P_{X_r^{t,\xi}}) \sigma(X_r^{t,x,P_{\xi}},P_{X_r^{t,\xi}})dW_r\cdot\eta\cdot\zeta\],
  \end{split}
\end{equation*}
and so
\begin{equation}\label{eq9.33}
  \begin{split}
 E\!\[\!\Big\{&U(s,\!X_s^{t,x,P_{\xi}}\!,\!P_{X_s^{t,\xi}}\!)\!-\!\(\!\Phi(X_T^{t,x,P_{\xi}}\!,\!P_{X_T^{t,\xi}}\!)
\!+\! \! \int_s^T\!\!\!\(g^2(P_{X_r^{t,\xi}})(\partial_xU)(r,\!X_r^{t,x,P_{\xi}}\!,\!P_{X_r^{t,\xi}}) \sigma(X_r^{t,x,P_{\xi}}\!\!,\!P_{X_r^{t,\xi}}\!)\\
 &   +h(P_{X_r^{t,\xi}})\)d\overleftarrow{B_r}
 -\int_s^T(\partial_{x}U)(r,X_r^{t,x,P_{\xi}},P_{X_r^{t,\xi}}) \sigma(X_r^{t,x,P_{\xi}},P_{X_r^{t,\xi}})dW_r\)\Big\}\cdot\eta\cdot\zeta\]=0,
  \end{split}
\end{equation}
for all $\zeta\in L^{\infty}(\mathcal{F}_{T}^W;\mathbb{R})$, and for all $\eta\in L^{\infty}(\mathcal{F}_{0,T}^B;\mathbb{R})$ that satisfy (H9.1) such that
$s\rightarrow\theta_s^{\eta}$ continuous.
This shows that, $P$-a.s.,
\begin{equation*}
  \begin{split}
  U(s,X_s^{t,x,P_{\xi}},P_{X_s^{t,\xi}})=&\ \Phi(X_T^{t,x,P_{\xi}},P_{X_T^{t,\xi}})
+ \int_s^T\(g^2(P_{X_r^{t,\xi}})(\partial_xU)(r,X_r^{t,x,P_{\xi}},P_{X_r^{t,\xi}}) \sigma(X_r^{t,x,P_{\xi}},P_{X_r^{t,\xi}})\\
 &\   +h(P_{X_r^{t,\xi}})\)d\overleftarrow{B_r}
 \!-\!\int_s^T(\partial_{x}U)(r,X_r^{t,x,P_{\xi}},P_{X_r^{t,\xi}}) \sigma(X_r^{t,x,P_{\xi}},P_{X_r^{t,\xi}})dW_r,\ s\in[t,T],
  \end{split}
\end{equation*}
i.e., $( U(s,X_s^{t,x,P_{\xi}},P_{X_s^{t,\xi}}),(\partial_{x}U)(s,X_s^{t,x,P_{\xi}},P_{X_s^{t,\xi}}) \sigma(X_s^{t,x,P_{\xi}},P_{X_s^{t,\xi}}))$
is a solution of BDSDE \eqref{eq4.2}, and from the uniqueness of the solution it follows:
$$( U(s,X_s^{t,x,P_{\xi}},P_{X_s^{t,\xi}}),(\partial_{x}U)(s,X_s^{t,x,P_{\xi}},P_{X_s^{t,\xi}}) \sigma(X_s^{t,x,P_{\xi}},P_{X_s^{t,\xi}}))
=(Y_s^{t,x,P_{\xi}},Z_s^{t,x,P_{\xi}}),$$
and, in particular, for $s=t$, $U(t,x,P_{\xi})=Y_t^{t,x,P_{\xi}}=V(t,x,P_{\xi})$,
$(t,x,P_{\xi})\in[0,T]\times\mathbb{R}\times\mathcal{P}_{2}(\mathbb{R})$. The proof is complete.
\end{proof}

\appendix
\section{Appendix}
\subsection{The proof of Theorem \ref{th2.1}}
\begin{proof}
Without loss of generality, we may prove the theorem for the case $F\in C_b^{1,2,2}([0,T]\times\mathbb{R}^{d}\times\mathcal{P}_{2}(\mathbb{R}^{d}))$,
$f=0$ and the dimensions $d=1$ and $l=1$.
The general case with $f$ and dimensions $d \geq 1$, $l\geq 1$ can be obtained by a straight-forward extension.

\noindent\textbf{Step 1.}
Firstly, we consider the special case $F(s,x,\mu)=G(\mu)$, $(s,x,\mu)\in[0,T]\times\mathbb{R}\times\mathcal{P}_{2}(\mathbb{R})$.
For this we consider the solution $(Y,Z)$ of the BDSDE $\displaystyle Y_t=\xi+\int_t^Tg_sd\overleftarrow{B_s}-\int_t^TZ_sdW_s$, $t\in[0,T]$,
for $\xi\in L^{2}(\mathcal{F}_T;\mathbb{R})$ and $g\in\mathcal{H}^2_{\mathcal{F}}(0,T;\mathbb{R})$,
where $g$ and $\xi$ are bounded by a positive constant $C$.

Let $n,k\in\mathbb{N}$, $n\geq C$, $1\leq k\leq 2^n$. We put $t_k^n:=k2^{-n}T$, $\varphi^{+}(s):=k2^{-n}T$,
$\varphi^{-}(s):=(k-1)2^{-n}T$, and we define the process $g^n_s:=E[g_s\big|\mathcal{F}^{W}_{\varphi^{-}(s)}\vee\mathcal{F}^{B}_{\varphi^{+}(s),T}]$,
for $s\in(t_{k-1}^n,t_k^n]$. As $g^n\in\mathcal{H}^2_{\mathcal{F}}(0,T;\mathbb{R})$ we can introduce the process $(Y^n,Z^n)$ as solution of the BDSDE
\begin{equation}\label{eqA.3.4}
 Y_t^n=\xi+\int_t^Tg_s^n d\overleftarrow{B_s}-\int_t^TZ_s^ndW_s,\ t\in[0,T].
\end{equation}
Lemma \ref{leA.1}, but also standard estimates allow to verify easily that
\begin{equation}\label{eqA.3.5}
  \begin{split}
    E\[\sup_{s\in[0,T]}|Y_s-Y_s^n|^{2}+\int_{0}^{T}|Z_s-Z^n_s|^{2}ds\]
\leq   C E\[\int_{0}^{T}|g_s-g^n_s|^{2}ds\] \rightarrow0,\ n\rightarrow\infty.
  \end{split}
\end{equation}
We remark that, for $1\leq k\leq 2^n$, $Y_{t_k^n}^n\in L^2(\mathcal{F}_{t_k^n};\mathbb{R})$,
where $\mathcal{F}_{t_k^n}=\mathcal{F}_{t_k^n}^W\vee\mathcal{F}_{t_k^n,T}^B$.
As $W=(W_t)$ is an $(\mathcal{F}_{t}^W\vee\mathcal{F}_{t_k^n,T}^B)$-Brownian motion and has w.r.t. the
same filtration the martingale representation property, there exists a unique square integrable,
$(\mathcal{F}_{t}^W\vee\mathcal{F}_{t_k^n,T}^B)$-adapted process $\{(\widetilde{Y}_t^n,\widetilde{Z}_t^n),t\in[t_{k-1}^n,t_k^n]\}$
such that
\begin{equation*}
  \begin{split}
 \widetilde{Y}_t^n=Y_{t_k^n}^n-\int_t^{t_k^n}\widetilde{Z}_r^ndW_r,\ t\in[t_{k-1}^n,t_k^n].
  \end{split}
\end{equation*}
But as $\{g_r^n,r\in[t_{k-1}^n,t_k^n]\}$ is bounded and
$\mathcal{B}([t_{k-1}^n,t_k^n])\otimes(\mathcal{F}_{t_{k-1}^n}^W\vee\mathcal{F}_{t_k^n,T}^B)$-measurable,
the process $\displaystyle\{\int_t^{t_k^n}g_r^n d\overleftarrow{B_r},\ t\in[t_{k-1}^n,t_k^n]\}$ is in $\mathcal{H}^2_{\mathcal{F}_{t_{k-1}^n}^W\vee\mathcal{F}_{\cdot,T}^B}(t_{k-1}^n,t_k^n;\mathbb{R})$
and the process $\displaystyle\overline{Y}_t^n:=\widetilde{Y}_t^n+\int_t^{t_k^n}g_r^n d\overleftarrow{B_r},\ t\in[t_{k-1}^n,t_k^n]$,
belongs to $\mathcal{S}^2_{\mathcal{F}}(t_{k-1}^n,t_k^n;\mathbb{R})$.
Consequently, $(\overline{Y}_t^n,\widetilde{Z}_t^n),\ t\in[t_{k-1}^n,t_k^n]$,
is a solution of BDSDE \eqref{eqA.3.4},
\begin{equation*}
  \begin{split}
 \overline{Y}_t^n=\widetilde{Y}_t^n+\int_t^{t_k^n}g_r^n d\overleftarrow{B_r}
 =Y_{t_k^n}^n+\int_t^{t_k^n}g_r^n d\overleftarrow{B_r}-\int_t^{t_k^n}\widetilde{Z}_r^ndW_r,\ t\in[t_{k-1}^n,t_k^n],
  \end{split}
\end{equation*}
and from the uniqueness of the solution of \eqref{eqA.3.4} it follows that
\begin{equation*}
  \begin{split}
    & Y_{t}^n=\overline{Y}_t^n=\widetilde{Y}_t^n+\int_t^{t_k^n}g_r^n d\overleftarrow{B_r},\ t\in[t_{k-1}^n,t_k^n],\ P\mbox{-a.s.},\\
    & Z_{t}^n=\widetilde{Z}_t^n,\ dtdP\mbox{-a.e.}\ \mbox{on}\  [t_{k-1}^n,t_k^n]\times\Omega,
  \end{split}
\end{equation*}
and so, in particular, we have that $Z_{t}^n$ is $(\mathcal{F}_{t}^W\vee\mathcal{F}_{t_k^n,T}^B)$-measurable,
$t\in[t_{k-1}^n,t_k^n]$, $1\leq k\leq 2^n$. Moreover, as $g$ and $\xi$ are bounded,
it follows from Lemma \ref{leA.1} that, for all $p\geq2$, there is some $C_p\in\mathbb{R}_{+}$ such that
\begin{equation}\label{eqA.3.6}
  E\[\sup_{s\in[0,T]}|Y_s^n|^p+(\int_0^T|Z_s^n|^2ds)^{\frac{p}{2}}\] \leq C_p,\ n\geq1.
\end{equation}
Then, for $1\leq k\leq 2^n$, $t_{k-1}^n\leq t<t+h\leq t_{k}^n$, with the notation $\Xi_{t,t+h}^{n,\lambda}:=Y_t^n+\lambda(Y_{t+h}^n-Y_t^n)$, it follows from the chain rule that
\begin{equation}\label{eqA.3.8}
  \begin{split}
 &\ G(P_{Y_{t+h}^n})-G(P_{Y_t^n})=\widetilde{G}(Y_{t+h}^n)-\widetilde{G}(Y_t^n)
 = \int_0^1\frac{d}{d\lambda}\widetilde{G}(\Xi_{t,t+h}^{n,\lambda})d\lambda\\
 = &\ \int_0^1\widehat{E}\big[(\partial_{\mu}G)(P_{\Xi_{t,t+h}^{n,\lambda}},\widehat{\Xi}_{t,t+h}^{n,\lambda})\cdot (\widehat{Y}_{t+h}^n-\widehat{Y}_t^n)\big]d\lambda
 =I_{1}-I_{2},
  \end{split}
\end{equation}
where
\begin{equation*}
  \begin{split}
&I_1:= \int_0^1\widehat{E}\big[(\partial_{\mu}G)(P_{\Xi_{t,t+h}^{n,\lambda}},\widehat{\Xi}_{t,t+h}^{n,\lambda})
       \cdot (\int_{t}^{t+h}\widehat{Z}_r^nd\widehat{W}_r)\big]d\lambda ,\\
&I_2:= \int_0^1\widehat{E}\big[(\partial_{\mu}G)(P_{\Xi_{t,t+h}^{n,\lambda}},\widehat{\Xi}_{t,t+h}^{n,\lambda})
       \cdot (\int_{t}^{t+h}\widehat{g}_r^n d\overleftarrow{\widehat{B}_r})\big]d\lambda .
  \end{split}
\end{equation*}
Here $(\widehat{\Omega},\widehat{\mathcal{F}},\widehat{P})$ is a probability space carrying with $(\widehat{\xi},\widehat{g}^n,\widehat{B},\widehat{W})$
an independent copy of $(\xi,g^n,B,\!W)$ (which is defined on $(\Omega,\mathcal{F},P)$); $(\widehat{Y}^{n},\widehat{Z}^{n})$
 is the solution of the same equation as that for $(Y^{n},Z^{n})$, but with the data $(\widehat{\xi},\widehat{g}^n,\widehat{B},\widehat{W})$
instead of $(\xi,g^n,B,W)$.

We first make the computation for $I_1$. It is easy to check that
\begin{equation}\label{eqA.3.9}
  \begin{split}
 I_{1} =&\ \int_0^1\widehat{E}\[\((\partial_{\mu}G)(P_{\Xi_{t,t+h}^{n,\lambda}},\widehat{\Xi}_{t,t+h}^{n,\lambda})
          -(\partial_{\mu}G)(P_{\Xi_{t,t+h}^{n,\lambda}},\widehat{Y}_{t}^{n})\)
       \cdot (\int_{t}^{t+h}\widehat{Z}_r^nd\widehat{W}_r)\]d\lambda\\
=&\ \int_0^1\lambda\int_0^1\widehat{E}\[\partial_{y}(\partial_{\mu}G)(P_{\Xi_{t,t+h}^{n,\lambda}},\widehat{\Xi}_{t,t+h}^{n,\rho\lambda})
      \cdot(\widehat{Y}_{t+h}^n-\widehat{Y}_t^n) \cdot (\int_{t}^{t+h}\widehat{Z}_r^nd\widehat{W}_r)\]d\rho d\lambda   \\
=&\ I_{1,1}+I_{1,2},
  \end{split}
\end{equation}
where
\begin{equation*}
  \begin{split}
&I_{1,1}:= \int_0^1\lambda\int_0^1\widehat{E}\[\partial_{y}(\partial_{\mu}G)(P_{\Xi_{t,t+h}^{n,\lambda}},\widehat{\Xi}_{t,t+h}^{n,\rho\lambda})
           \cdot (\int_{t}^{t+h}\widehat{Z}_r^nd\widehat{W}_r)^2\]d\rho d\lambda ,\\
&I_{1,2}:= -\int_0^1\lambda\int_0^1\widehat{E}\[\partial_{y}(\partial_{\mu}G)(P_{\Xi_{t,t+h}^{n,\lambda}},\widehat{\Xi}_{t,t+h}^{n,\rho\lambda})
      \cdot (\int_{t}^{t+h}\widehat{g}_r^n d\overleftarrow{\widehat{B}_r})\cdot (\int_{t}^{t+h}\widehat{Z}_r^nd\widehat{W}_r)\]d\rho d\lambda .
  \end{split}
\end{equation*}
Let us now consider $I_{1,1}$. Obviously,
\begin{equation}\label{eqA.3.10}
  \begin{split}
 I_{1,1}=&\ \int_0^1\lambda\int_0^1\widehat{E}\[\partial_{y}(\partial_{\mu}G)(P_{\Xi_{t,t+h}^{n,\lambda}},\widehat{Y}_{t}^{n})
           \cdot (\int_{t}^{t+h}\widehat{Z}_r^nd\widehat{W}_r)^2\]d\rho d\lambda+R_{1,1}\\
=&\ \int_0^1\lambda\widehat{E}\[\partial_{y}(\partial_{\mu}G)(P_{\Xi_{t,t+h}^{n,\lambda}},\widehat{Y}_{t}^{n})
           \cdot (\int_{t}^{t+h}|\widehat{Z}_r^n|^2dr)\] d\lambda+R_{1,1}\\
=&\ \frac{1}{2}\widehat{E}\[\partial_{y}(\partial_{\mu}G)(P_{Y_{t}^{n}},\widehat{Y}_{t}^{n})
           \cdot (\int_{t}^{t+h}|\widehat{Z}_r^n|^2dr)\]+R_{1,1}+ R_{1,2} ,
  \end{split}
\end{equation}
where
\begin{equation*}
  \begin{split}
&R_{1,1}:= \int_0^1\lambda\int_0^1\widehat{E}\[\(\partial_{y}(\partial_{\mu}G)(P_{\Xi_{t,t+h}^{n,\lambda}},\widehat{\Xi}_{t,t+h}^{n,\rho\lambda})
           -\partial_{y}(\partial_{\mu}G)(P_{\Xi_{t,t+h}^{n,\lambda}},\widehat{Y}_{t}^{n})\)
           \cdot (\int_{t}^{t+h}\widehat{Z}_r^nd\widehat{W}_r)^2\]d\rho d\lambda ,\\
&R_{1,2}:= \int_0^1\lambda\widehat{E}\[\(\partial_{y}(\partial_{\mu}G)(P_{\Xi_{t,t+h}^{n,\lambda}},\widehat{Y}_{t}^{n})
           -\partial_{y}(\partial_{\mu}G)(P_{Y_{t}^{n}},\widehat{Y}_{t}^{n})\)\cdot (\int_{t}^{t+h}|\widehat{Z}_r^n|^2dr)\] d\lambda.
  \end{split}
\end{equation*}
Making use of the fact that $\partial_y(\partial_{\mu}G)$ is Lipschitz, we get
\begin{equation*}
  \begin{split}
   |R_{1,1}|\leq&\ C\int_0^1\lambda\int_0^1\widehat{E}\[|\widehat{\Xi}_{t,t+h}^{n,\rho\lambda}-\widehat{Y}_t^n|
\cdot(\int_{t}^{t+h}\widehat{Z}_r^nd\widehat{W}_r)^2\]d\rho d\lambda\\
\leq&\ C E\[|Y_{t+h}^n-Y_t^n|\cdot(\int_{t}^{t+h}Z_r^ndW_r)^2\]\\
\leq&\ C E\[|\int_{t}^{t+h}Z_r^ndW_r|^3\]+C E\[|\int_t^{t+h}g_r^n d\overleftarrow{B_r}|\cdot(\int_{t}^{t+h}Z_r^ndW_r)^2\].
  \end{split}
\end{equation*}
But, as $\displaystyle\int_t^{t+h}g_r^n d\overleftarrow{B_r}$ is $\mathcal{F}_{t_{k-1}^n}^W\vee\mathcal{F}_{t,T}^B\vee\mathcal{F}^0$-measurable
and $\displaystyle\int_{t}^{t+h}Z_r^ndW_r$ is $\mathcal{F}_{t+h}^W\vee\mathcal{F}_{t_{k}^n,T}^B$-measurable, both
are conditionally independent knowing $\mathcal{F}_{t_{k-1}^n}^W\vee\mathcal{F}_{t_{k}^n,T}^B$. Thus,
\begin{equation*}
  \begin{split}
   &\ E\[|\int_t^{t+h}g_r^n d\overleftarrow{B_r}|\cdot(\int_{t}^{t+h}Z_r^ndW_r)^2\]\\
=&\ E\[E\[|\int_t^{t+h}g_r^n d\overleftarrow{B_r}|\big|\mathcal{F}_{t_{k-1}^n}^W\vee\mathcal{F}_{t_{k}^n,T}^B\]
       \cdot E\[(\int_{t}^{t+h}Z_r^ndW_r)^2\big|\mathcal{F}_{t_{k-1}^n}^W\vee\mathcal{F}_{t_{k}^n,T}^B\]\]\\
  \end{split}
\end{equation*}
\begin{equation*}
  \begin{split}
\leq&\ C E\[E\[(\int_t^{t+h}|g_r^n|^2 dr)^{\frac{1}{2}}\big|\mathcal{F}_{t_{k-1}^n}^W\vee\mathcal{F}_{t_{k}^n,T}^B\]
      \cdot E\[\int_{t}^{t+h}|Z_r^n|^2dr\big|\mathcal{F}_{t_{k-1}^n}^W\vee\mathcal{F}_{t_{k}^n,T}^B\]\]     \\
\leq&\  C\sqrt{h} E\[\int_{t}^{t+h}|Z_r^n|^2dr\],
  \end{split}
\end{equation*}
where we have used the Burkholder-Davis-Gundy inequality and the boundedness of $g$.
Consequently,
\begin{equation*}
  \begin{split}
   |R_{1,1}|\leq  C E\[\((\int_t^{t+h}|Z_r^n|^2 dr)^{\frac{1}{2}}+\sqrt{h}\)\cdot\int_{t}^{t+h}|Z_r^n|^2dr\].
  \end{split}
\end{equation*}
The estimate of $R_{1,2}$ is similar but easier:
\begin{equation*}
  \begin{split}
 |R_{1,2}|\leq&\ C\int_0^1\lambda W_2(P_{\Xi_{t,t+h}^{n,\lambda}},P_{Y_{t}^{n}})d\lambda\cdot\widehat{E}\[ \int_{t}^{t+h}|\widehat{Z}_r^n|^2dr\]\\
   \leq&\ C(E[|Y_{t+h}^n-Y_t^n|^2])^{\frac{1}{2}}\cdot E[\int_{t}^{t+h}|Z_r^n|^2dr]\\
   \leq&\ C\(\sqrt{h}+(E[\int_t^{t+h}|Z_r^n|^2 dr])^{\frac{1}{2}}\)\cdot E[\int_{t}^{t+h}|Z_r^n|^2dr].
  \end{split}
\end{equation*}
Next we compute $I_{1,2}$. As $\displaystyle (Y_{t+h}^n,\int_{t}^{t+h}Z_r^ndW_r)$ is $\mathcal{F}_{t+h}^W\vee\mathcal{F}_{t+h,T}^B$-measurable, it holds
\begin{equation}\label{eqA.3.17}
  \begin{split}
-I_{1,2}=&\ \int_0^1\lambda\int_0^1\widehat{E}\[\partial_{y}(\partial_{\mu}G)(P_{\Xi_{t,t+h}^{n,\lambda}},\widehat{\Xi}_{t,t+h}^{n,\rho\lambda})
      \cdot (\int_{t}^{t+h}\widehat{g}_r^n d\overleftarrow{\widehat{B}_r})\cdot (\int_{t}^{t+h}\widehat{Z}_r^nd\widehat{W}_r)\]d\rho d\lambda \\
=&\  \int_0^1\lambda\widehat{E}\[\partial_{y}(\partial_{\mu}G)(P_{\Xi_{t,t+h}^{n,\lambda}},\widehat{Y}_{t+h}^{n})
      \cdot (\int_{t}^{t+h}\widehat{g}_r^n d\overleftarrow{\widehat{B}_r})\cdot (\int_{t}^{t+h}\widehat{Z}_r^nd\widehat{W}_r)\] d\lambda+R_{1,3}
=    R_{1,3}  ,
  \end{split}
\end{equation}
where
\begin{equation*}
  \begin{split}
R_{1,3}:=&\ \int_0^1\lambda\int_0^1\widehat{E}\[\(\partial_{y}(\partial_{\mu}G)(P_{\Xi_{t,t+h}^{n,\lambda}},\widehat{\Xi}_{t,t+h}^{n,\rho\lambda})
      -\partial_{y}(\partial_{\mu}G)(P_{\Xi_{t,t+h}^{n,\lambda}},\widehat{Y}_{t+h}^{n})\)\\
      &\    \times (\int_{t}^{t+h}\widehat{g}_r^n d\overleftarrow{\widehat{B}_r})\cdot (\int_{t}^{t+h}\widehat{Z}_r^nd\widehat{W}_r)\]d\rho d\lambda.
  \end{split}
\end{equation*}
Similar to the estimate of $R_{1,1}$, by using again the conditional independence of
$\displaystyle\int_t^{t+h}g_r^n d\overleftarrow{B_r}$ and $\displaystyle\int_{t}^{t+h}Z_r^ndW_r$
knowing $\mathcal{F}_{t_{k-1}^n}^W\vee\mathcal{F}_{t_{k}^n,T}^B$, we have
\begin{equation*}
  \begin{split}
|R_{1,3}|\leq&\ C E\[|Y_{t+h}^n-Y_t^n|\cdot|\int_{t}^{t+h}g_r^n d\overleftarrow{B_r}|\cdot|\int_{t}^{t+h}Z_r^ndW_r|\]\\
 \leq&\  C E\[|\int_{t}^{t+h}g_r^n d\overleftarrow{B_r}|^2\cdot|\int_{t}^{t+h}Z_r^ndW_r|\]
     + C E\[|\int_{t}^{t+h}g_r^n d\overleftarrow{B_r}|\cdot|\int_{t}^{t+h}Z_r^ndW_r|^2\]\\
 =&\   CE\[E\[|\int_{t}^{t+h}g_r^n d\overleftarrow{B_r}|^2\big|\mathcal{F}_{t_{k-1}^n}^W\vee\mathcal{F}_{t_{k}^n,T}^B\]
    \cdot E\[|\int_{t}^{t+h}Z_r^ndW_r|\big|\mathcal{F}_{t_{k-1}^n}^W\vee\mathcal{F}_{t_{k}^n,T}^B\] \]\\
   &\   + C E\[E\[|\int_{t}^{t+h}g_r^n d\overleftarrow{B_r}|\big|\mathcal{F}_{t_{k-1}^n}^W\vee\mathcal{F}_{t_{k}^n,T}^B\]
  \cdot E\[|\int_{t}^{t+h}Z_r^ndW_r|^2\big|\mathcal{F}_{t_{k-1}^n}^W\vee\mathcal{F}_{t_{k}^n,T}^B\]\]\\
 \leq&\ ChE[(\int_{t}^{t+h}|Z_r^n|^2dr)^{\frac{1}{2}}]+C\sqrt{h}E[\int_{t}^{t+h}|Z_r^n|^2dr].
  \end{split}
\end{equation*}
Finally, summarizing the above computation for $I_1$, we get
\begin{equation}\label{eqA.3.23}
  \begin{split}
 I_{1}= \frac{1}{2}\widehat{E}\[\partial_{y}(\partial_{\mu}G)(P_{Y_{t}^{n}},\widehat{Y}_{t}^{n})
           \cdot (\int_{t}^{t+h}|\widehat{Z}_r^n|^2dr)\]+R_1,
  \end{split}
\end{equation}
where $R_1=R_{1,1}+ R_{1,2}-R_{1,3}$, and
\begin{equation*}
  |R_1|\!\leq\! C E\[\((\int_t^{t+h}\!\!|Z_r^n|^2 dr)^{\frac{1}{2}}+(E[\int_{t}^{t+h}\!\!|Z_r^n|^2dr] )^{\frac{1}{2}}+\sqrt{h}\)
  \cdot\(\int_{t}^{t+h}\!\!|Z_r^n|^2dr+h\)\],\ t_{k-1}^n\!\leq\! t\!<\!t+h\!\leq\! t_{k}^n.
\end{equation*}

Next we compute $I_2$. Let $\Gamma_{t,t+h}^{n,\rho,\lambda}:=Y_{t+h}^n-\rho(1-\lambda)(Y_{t+h}^n-Y_t^n)$. Then, as $Y_{t+h}^n$
is $ \mathcal{F}^{W}_{t+h}\vee\mathcal{F}^{B}_{t+h,T}$-measurable,
\begin{equation}\label{eqA.3.24}
  \begin{split}
 I_2=&\ \int_0^1\widehat{E}\big[(\partial_{\mu}G)(P_{\Xi_{t,t+h}^{n,\lambda}},\widehat{\Xi}_{t,t+h}^{n,\lambda})
       \cdot (\int_{t}^{t+h}\widehat{g}_r^n d\overleftarrow{\widehat{B}_r})\big]d\lambda\\
 =&\ \int_0^1\widehat{E}\big[(\partial_{\mu}G)(P_{\Xi_{t,t+h}^{n,\lambda}},\widehat{Y}_{t+h}^{n})
       \cdot (\int_{t}^{t+h}\widehat{g}_r^n d\overleftarrow{\widehat{B}_r})\big]d\lambda\\
 &\  -  \int_0^1(1-\lambda)\int_0^1\widehat{E}\big[\partial_{y}(\partial_{\mu}G)(P_{\Xi_{t,t+h}^{n,\lambda}},\widehat{\Gamma}_{t,t+h}^{n,\rho,\lambda})
       \cdot(\widehat{Y}_{t+h}^n-\widehat{Y}_t^n)\cdot (\int_{t}^{t+h}\widehat{g}_r^n d\overleftarrow{\widehat{B}_r})\big]d\rho d\lambda  \\
=&\   -  \int_0^1(1-\lambda)\int_0^1\widehat{E}\big[\partial_{y}(\partial_{\mu}G)(P_{\Xi_{t,t+h}^{n,\lambda}},\widehat{\Gamma}_{t,t+h}^{n,\rho,\lambda})
       \cdot(\widehat{Y}_{t+h}^n-\widehat{Y}_t^n)\cdot (\int_{t}^{t+h}\widehat{g}_r^n d\overleftarrow{\widehat{B}_r})\big]d\rho d\lambda   \\
=&\  I_{2,1}-I_{2,2},
  \end{split}
\end{equation}
where
\begin{equation*}
  \begin{split}
&I_{2,1}:=  \int_0^1(1-\lambda)\int_0^1\widehat{E}\big[\partial_{y}(\partial_{\mu}G)(P_{\Xi_{t,t+h}^{n,\lambda}},\widehat{\Gamma}_{t,t+h}^{n,\rho,\lambda})
       \cdot (\int_{t}^{t+h}\widehat{g}_r^n d\overleftarrow{\widehat{B}_r})^2\big]d\rho d\lambda, \\
&I_{2,2}:=  \int_0^1(1-\lambda)\int_0^1\widehat{E}\big[\partial_{y}(\partial_{\mu}G)(P_{\Xi_{t,t+h}^{n,\lambda}},\widehat{\Gamma}_{t,t+h}^{n,\rho,\lambda})
       \cdot (\int_{t}^{t+h}\widehat{g}_r^n d\overleftarrow{\widehat{B}_r})\cdot (\int_{t}^{t+h}\widehat{Z}_r^nd\widehat{W}_r)\big]d\rho d\lambda .
  \end{split}
\end{equation*}
A straight-forward computation for $I_{2,1}$ shows
\begin{equation}\label{eqA.3.25}
  \begin{split}
 I_{2,1}=&\  \int_0^1(1-\lambda)\int_0^1\widehat{E}\big[\partial_{y}(\partial_{\mu}G)(P_{\Xi_{t,t+h}^{n,\lambda}},\widehat{Y}_{t+h}^{n})
       \cdot (\int_{t}^{t+h}\widehat{g}_r^n d\overleftarrow{\widehat{B}_r})^2\big]d\rho d\lambda+R_{2,1}\\
 =&\  \int_0^1(1-\lambda)\int_0^1\widehat{E}\big[\partial_{y}(\partial_{\mu}G)(P_{\Xi_{t,t+h}^{n,\lambda}},\widehat{Y}_{t+h}^{n})
       \cdot (\int_{t}^{t+h}|\widehat{g}_r^n|^2 dr)\big]d\rho d\lambda+R_{2,1}\\
 =&\  \frac{1}{2}\widehat{E}\big[\partial_{y}(\partial_{\mu}G)(P_{Y_{t}^{n}},\widehat{Y}_{t}^{n})
       \cdot (\int_{t}^{t+h}|\widehat{g}_r^n|^2 dr)\big]+R_{2,1} +R_{2,2},
  \end{split}
\end{equation}
where
\begin{equation*}
  \begin{split}
&R_{2,1}:=  \int_0^1(1-\lambda)\int_0^1\widehat{E}\[\(\partial_{y}(\partial_{\mu}G)(P_{\Xi_{t,t+h}^{n,\lambda}},\widehat{\Gamma}_{t,t+h}^{n,\rho,\lambda})
       -\partial_{y}(\partial_{\mu}G)(P_{\Xi_{t,t+h}^{n,\lambda}},\widehat{Y}_{t+h}^{n})\)
        (\int_{t}^{t+h}\widehat{g}_r^n d\overleftarrow{\widehat{B}_r})^2\]d\rho d\lambda, \\
&R_{2,2}:=  \int_0^1(1-\lambda)\int_0^1\widehat{E}\[\(\partial_{y}(\partial_{\mu}G)(P_{\Xi_{t,t+h}^{n,\lambda}},\widehat{Y}_{t+h}^{n})
       -\partial_{y}(\partial_{\mu}G)(P_{Y_{t}^{n}},\widehat{Y}_{t}^{n})\)\cdot (\int_{t}^{t+h}|\widehat{g}_r^n|^2 dr)\]d\rho d\lambda,
  \end{split}
\end{equation*}
and with arguments similar to those for $R_1$ we see that
\begin{equation*}
 |R_{2,1}|+|R_{2,2}|\leq C\((E[\int_{t}^{t+h}\!\!|Z_r^n|^2dr] )^{\frac{1}{2}}+\sqrt{h}\)h.
\end{equation*}
As concerns $I_{2,2}$, we remark
\begin{equation}\label{eqA.3.26}
  \begin{split}
 I_{2,2}=&\  \int_0^1(1-\lambda)\int_0^1\widehat{E}\big[\partial_{y}(\partial_{\mu}G)(P_{\Xi_{t,t+h}^{n,\lambda}},\widehat{Y}_{t+h}^{n})
       \cdot (\int_{t}^{t+h}\widehat{g}_r^n d\overleftarrow{\widehat{B}_r})\cdot (\int_{t}^{t+h}\widehat{Z}_r^nd\widehat{W}_r)\big]d\rho d\lambda\\
&\ +  \int_0^1(1-\lambda)\int_0^1\widehat{E}\[\(\partial_{y}(\partial_{\mu}G)(P_{\Xi_{t,t+h}^{n,\lambda}},\widehat{\Gamma}_{t,t+h}^{n,\rho,\lambda})
       -\partial_{y}(\partial_{\mu}G)(P_{\Xi_{t,t+h}^{n,\lambda}},\widehat{Y}_{t+h}^{n})\) \\
 &\      \times(\int_{t}^{t+h}\widehat{g}_r^n d\overleftarrow{\widehat{B}_r})\cdot (\int_{t}^{t+h}\widehat{Z}_r^nd\widehat{W}_r)\]d\rho d\lambda,\\
  \end{split}
\end{equation}
where the first double integral at the right-hand side equals to zero.
We observe that the terms in $I_{2,2}$ are similar to those of $I_{1}$. Thus, using an argument analogous to that in the proof of $I_{1,1}$ and $I_{1,2}$,
and combining it with our result for $I_{2,1}$, we obtain
\begin{equation}\label{eqA.3.27}
  \begin{split}
 I_{2}= \frac{1}{2}\widehat{E}\[\partial_{y}(\partial_{\mu}G)(P_{Y_{t}^{n}},\widehat{Y}_{t}^{n})
           \cdot (\int_{t}^{t+h}|\widehat{g}_r^n|^2dr)\]+R_2,
  \end{split}
\end{equation}
where $R_2=R_{2,1}+ R_{2,2}-I_{2,2}$, and
\begin{equation*}
  |R_2|\leq C E\[\((\int_t^{t+h}\!\!|Z_r^n|^2 dr)^{\frac{1}{2}}+(E[\int_{t}^{t+h}\!\!|Z_r^n|^2dr] )^{\frac{1}{2}}+\sqrt{h}\)
  \cdot\(\int_{t}^{t+h}\!\!|Z_r^n|^2dr+h\)\],\ t_{k-1}^n\!\leq\! t\!<\!t+h\!\leq\! t_{k}^n.
\end{equation*}

Hence, from \eqref{eqA.3.8}, \eqref{eqA.3.23} and \eqref{eqA.3.27}, we have
\begin{equation}\label{eqA.3.28}
  \begin{split}
 &\ G(P_{Y_{t+h}^n})-G(P_{Y_t^n})= \frac{1}{2}\widehat{E}\[\partial_{y}(\partial_{\mu}G)(P_{Y_{t}^{n}},\widehat{Y}_{t}^{n})
           \cdot (\int_{t}^{t+h}(|\widehat{Z}_r^n|^2-|\widehat{g}_r^n|^2)dr)\]+R_{t,t+h}^n,
  \end{split}
\end{equation}
where $R_{t,t+h}^n:=R_{1}-R_{2}$, and
\begin{equation*}
  |R_{t,t+h}^n|\!\leq\! C E\[\((\int_t^{t+h}\!\!|Z_r^n|^2 dr)^{\frac{1}{2}}\!+\!(E[\int_{t}^{t+h}\!\!|Z_r^n|^2dr] )^{\frac{1}{2}}\!+\!\sqrt{h}\)
  \!\cdot\!\(\int_{t}^{t+h}\!\!|Z_r^n|^2dr\!+\!h\)\],\ t_{k-1}^n\!\leq\! t\!<\!t+h\!\leq\! t_{k}^n.
\end{equation*}
Now let $m\geq 1$, and put $s_i^m:=t+i2^{-m}h$, $0\leq i\leq 2^m$.
Then,
\begin{equation}\label{eqA.3.29}
  \begin{split}
 &\ G(P_{Y_{t+h}^n})-G(P_{Y_t^n})= \sum_{i=1}^{2^m}\(G(P_{Y_{s_{i}^m}})-G(P_{Y_{s_{i-1}^m}})\)\\
 =&\ \frac{1}{2}\widehat{E}\[\int_{t}^{t+h}\sum_{i=1}^{2^m}I_{[s_{i-1}^m,s_{i}^m]}(s)\partial_{y}(\partial_{\mu}G)(P_{Y_{s_{i-1}^m}^{n}},\widehat{Y}_{s_{i-1}^m}^{n})
           \cdot (|\widehat{Z}_s^n|^2-|\widehat{g}_s^n|^2)ds\]+\sum_{i=1}^{2^m}R_{s_{i-1}^m,s_{i}^m}^n,
  \end{split}
\end{equation}
with $\sum_{i=1}^{2^m}|R_{s_{i-1}^m,s_{i}^m}^n|\rightarrow 0$, as $m\rightarrow\infty$.
Indeed,
\begin{equation*}
  \begin{split}
  \sum_{i=1}^{2^m}|R_{s_{i-1}^m,s_{i}^m}^n|
\leq&\ E\[\max_{1\leq i\leq 2^m}\alpha(s_{i-1}^m,s_{i}^m)\cdot\(\int_{t}^{t+h}|Z_r^n|^2dr+h\)\]\\
\leq&\ C (E[(\max_{1\leq i\leq 2^m}\alpha(s_{i-1}^m,s_{i}^m))^2])^{\frac{1}{2}}=:J_m,
  \end{split}
\end{equation*}
where
\begin{equation*}
  \alpha(s,s')=(\int_s^{s'}|Z_r^n|^2 dr)^{\frac{1}{2}}+(E[\int_{s}^{s'}|Z_r^n|^2dr] )^{\frac{1}{2}}+\sqrt{s'-s},\ s\leq s'.
\end{equation*}
As $\displaystyle s\rightarrow\int_{0}^{s}|Z_r^n|^2dr$ is continuous, and bounded in $L^p$ ($p\geq1$),
it follows with the help of the dominated convergence theorem that $J_m\rightarrow0$.

Hence, letting $m\rightarrow\infty$ we obtain
\begin{equation}\label{eqA.3.30}
  \begin{split}
  G(P_{Y_{t+h}^n})-G(P_{Y_t^n})
 = \frac{1}{2}\widehat{E}\[\int_{t}^{t+h}\partial_{y}(\partial_{\mu}G)(P_{Y_{s}^{n}},\widehat{Y}_{s}^{n})
           \cdot (|\widehat{Z}_s^n|^2-|\widehat{g}_s^n|^2)ds\].
  \end{split}
\end{equation}
As this holds for all $ t_{k-1}^n\leq t< t+h\leq t_k^n$, $1\leq k\leq 2^n$, it holds for all $0\leq t<t+h\leq T$.
Moreover, from \eqref{eqA.3.5}, \eqref{eqA.3.6}, the boundedness and
the Lipschitz continuity of $G$ and $\partial_y(\partial_{\mu}G)$ and the dominated convergence theorem,
it follows that we can  take limit in \eqref{eqA.3.30}, as $n\rightarrow\infty$. This yields
\begin{equation}\label{eqA.3.31}
  \begin{split}
  G(P_{Y_{t+h}})-G(P_{Y_t})
 = \frac{1}{2}\widehat{E}\[\int_{t}^{t+h}\partial_{y}(\partial_{\mu}G)(P_{Y_{s}},\widehat{Y}_{s})
           \cdot (|\widehat{Z}_s|^2-|\widehat{g}_s|^2)ds\],\ 0\leq t<t+h\leq T.
  \end{split}
\end{equation}

\noindent\textbf{Step 2.} The extension of \eqref{eqA.3.31} to the multi-dimensional case is
straight-forward. For this reason we state only the result here.
Let $G\in C_b^{1,2}([0,T]\times\mathcal{P}_2(\mathbb{R}^d))$ and
$\xi\in L^{2}(\mathcal{F}_T;\mathbb{R}^d)$, $g\in\mathcal{H}^2_{\mathcal{F}}(0,T;\mathbb{R}^{d\times l})$
be bounded. Then, for the unique solution $(Y,Z)$ of the BDSDE
$\displaystyle Y_t=\xi+\int_t^Tg_sd\overleftarrow{B_s}-\int_t^TZ_sdW_s$, $t\in[0,T]$, we have, for $0\leq t<t+h\leq T$,
\begin{equation*}
  \begin{split}
G(t\!+\!h,P_{Y_{t\!+\!h}})\!-\!G(t,P_{Y_t})\!=\! \int_{t}^{t+h}\!\Big\{(\partial_tG)(s,P_{Y_{s}})
  \!+\!\frac{1}{2}\widehat{E}\[\mbox{tr}\big(\partial_{y}(\partial_{\mu}G)(s,P_{Y_{s}},\widehat{Y}_{s})
           \!\cdot\! (\widehat{Z}_s\widehat{Z}_s^{*}\!-\!\widehat{g}_s\widehat{g}_s^{*})\big)\]\Big\}ds.
  \end{split}
\end{equation*}

\noindent\textbf{Step 3.} We now consider $F\in C_b^{1,2}([0,T]\times\mathbb{R})$,
$l\in\mathcal{H}^2_{\mathcal{F}}(0,T;\mathbb{R})$ and $\xi\in L^{2}(\mathcal{F}_T;\mathbb{R})$,
where $l$ and $\xi$ are supposed to be bounded.
Our objective is to show that, for the unique solution $(U,V)\in \mathcal{S}_{\mathcal{F}}^2(0,T;\mathbb{R})\times\mathcal{H}_{\mathcal{F}}^2(0,T;\mathbb{R})$
of the BDSDE $\displaystyle U_t=\xi+\int_t^Tl_sd\overleftarrow{B_s}-\int_t^TV_sdW_s$, $t\in[0,T]$, it holds
\begin{equation*}
  \begin{split}
F(t+h,&U_{t+h})-F(t,U_t)= \int_{t}^{t+h}\!\Big\{(\partial_tF)(s,U_{s})
  \!+\!\frac{1}{2}\big((\partial_{xx}^2F)(s,U_{s})
           \!\cdot\! (|V_s|^2\!-\!|l_s|^{2})\big)\Big\}ds\\
 & - \int_{t}^{t+h} (\partial_{x}F)(s,U_{s})\cdot l_sd\overleftarrow{B_s}+  \int_{t}^{t+h} (\partial_{x}F)(s,U_{s})\cdot V_sdW_s,\
  t\in[0,T],\ P\mbox{-a.s.}
  \end{split}
\end{equation*}
For this we inspire from the proof of Theorem \ref{th9.1}. Let $\eta^1\in L^{\infty}(\mathcal{F}_{T};\mathbb{R})$
and $\eta^2\in L^{\infty}(\mathcal{F}_{0,T}^B;\mathbb{R})$ be such that, for $(\eta^1_t)\in\mathcal{H}^2_{\mathcal{F}^W}(0,T;\mathbb{R})$
and $(\eta^2_t)\in\mathcal{H}^2_{\mathcal{F}_{\cdot,T}^B}(0,T;\mathbb{R})$ with
\begin{equation*}
  \begin{split}
    & \eta^1(t):= E[ \eta^1\big|\mathcal{F}_t^W]= E[ \eta^1]+\int_0^t \eta^1_sdW_s,\ t\in[0,T],  \\
    & \eta^2(t):= E[ \eta^2\big|\mathcal{F}_{t,T}^B]= E[ \eta^2] +\int_t^T \eta^2_sd\overleftarrow{B_s},\ t\in[0,T],\  P\mbox{-a.s.},
  \end{split}
\end{equation*}
it holds that $(\eta^1_t)_{t\in[0,T]}$ and $(\eta^2_t)_{t\in[0,T]}$ are bounded. Note that, knowing $\mathcal{F}_{t}=\mathcal{F}_{t}^W\vee\mathcal{F}_{t,T}^B$,
$\eta^1$ and $\eta^2$ are conditionally independent, and, thus
\begin{equation*}
 E[F(t,U_t)\eta^1\eta^2]=E[F(t,U_t)E[\eta^1\eta^2|\mathcal{F}_{t}]]= E[F(t,U_t)\eta^1(t)\eta^2(t)],\ t\in[0,T].
\end{equation*}
Now, to apply Step 2, we put $\displaystyle G(t,\mu):=\int_{\mathbb{R}^3}F(t,x)yz\mu(dxdydz)$, $(t,\mu)\in[0,T]\times\mathcal{P}_2(\mathbb{R}^3)$, and
$Y_t:=(U_t,\eta^1(t),\eta^2(t))^{*}$, $Z_t:=(V_t,\eta^1_t,0)^{*}$, $t\in[0,T]$.
We observe that $(Y,Z)\in \mathcal{S}_{\mathcal{F}}^2(0,T;\mathbb{R}^3)\times\mathcal{H}_{\mathcal{F}}^2(0,T;\mathbb{R}^3)$
is the unique solution of the BDSDE
\begin{equation*}
  \left(   \begin{matrix}    U_t\\  \eta^1(t)\\ \eta^2(t) \end{matrix}\right)
  =\left(   \begin{matrix}   \xi\\  \eta^1\\ E[\eta^2] \end{matrix}\right)
  +\int_t^T\left(   \begin{matrix}    l_s\\  0\\ \eta^2_s \end{matrix}\right)d\overleftarrow{B_s}
  -\int_t^T\left(   \begin{matrix}    V_s\\  \eta^1_s\\ 0 \end{matrix}\right)dW_s,\ t\in[0,T].
\end{equation*}
This allows to apply the result of Step 2 to $G(t,P_{Y_t})=E[F(t,U_t)\eta^1(t)\eta^2(t)]$.
Note that, in particular, $(\partial_{\mu}G)(t,P_{Y_t},(x,y,z))=\big(\partial_xF(t,x)yz,F(t,x)z,F(t,x)y\big)$ and
put $g_t:=(l_t,0,\eta^2_t)^{*}$, $t\in[0,T]$. Consequently, a straight-forward computation yields
\begin{equation*}
  \begin{split}
  &\  E[F(t+h,U_{t+h})\eta^1\eta^2]-E[F(t,U_t)\eta^1\eta^2]\\
=  &\  E[F(t+h,U_{t+h})\eta^1(t+h)\eta^2(t+h)]-E[F(t,U_t)\eta^1(t)\eta^2(t)]=G(t+h,P_{Y_{t+h}})-G(t,P_{Y_t})\\
=  &\ \int_t^{t+h} \Big\{(\partial_tG)(s,P_{Y_{s}})
  \!+\!\frac{1}{2}\widehat{E}\[\mbox{tr}\Big(\partial_{(x,y,z)}(\partial_{\mu}G)(s,P_{Y_{s}},\widehat{Y}_{s})
           \!\cdot\! (\widehat{Z}_s\widehat{Z}_s^{*}\!-\!\widehat{g}_s\widehat{g}_s^{*})\Big)\]\Big\}ds\\
=  &\ \int_t^{t+h} \Bigg\{E[(\partial_tF)(s,U_{s})\eta^1(s)\eta^2(s)]
  \!+\!\frac{1}{2}\widehat{E}\[(\partial_{xx}^2F)(s,\widehat{U}_{s})\cdot(|\widehat{V}_s|^2-|\widehat{l}_s|^2)\widehat{\eta^1(s)}\widehat{\eta^2(s)}\]\\
&\   + \widehat{E}\[ (\partial_{x}F)(s,\widehat{U}_{s})\cdot(\widehat{V}_s\widehat{\eta^1_s}\widehat{\eta^2(s)}
                 -\widehat{l}_s\widehat{\eta^1(s)}\widehat{\eta^2_s} )    \]\Bigg\}ds,\ 0\leq t<t+h\leq T.
  \end{split}
\end{equation*}
As $\eta^2(s)= E[ \eta^2\big|\mathcal{F}_{s}]$, $s\in[0,T]$,
\begin{equation*}
  \begin{split}
&\ \int_t^{t+h}\widehat{E}\[ (\partial_{x}F)(s,\widehat{U}_{s})\cdot\widehat{V}_s\widehat{\eta^1_s}\widehat{\eta^2(s)} \]ds
=\int_t^{t+h}E\[ (\partial_{x}F)(s,U_{s})\cdot V_s\eta^1_s\eta^2 \]ds\\
=&\ E\[ E\[\int_t^{t+h}(\partial_{x}F)(s,U_{s})\cdot V_s\eta^1_s ds\big|\mathcal{F}_{0,T}^B\] \eta^2\]\\
=&\ E\[ E\[\big(\int_t^{t+h}(\partial_{x}F)(s,U_{s})\cdot V_sdW_s\big)\cdot\big(\int_0^{T} \eta^1_sdW_s\big)\big|\mathcal{F}_{0,T}^B\] \eta^2\]\\
=&\ E\[ E\[\big(\int_t^{t+h}(\partial_{x}F)(s,U_{s})\cdot V_sdW_s\big)\cdot\eta^1\big|\mathcal{F}_{0,T}^B\] \eta^2\]\\
=&\ E\[ \big(\int_t^{t+h}(\partial_{x}F)(s,U_{s})\cdot V_sdW_s\big)\cdot\eta^1 \eta^2\],\ 0\leq t<t+h\leq T.
  \end{split}
\end{equation*}
With a symmetric argument we see that
\begin{equation*}
  \begin{split}
&\ \int_t^{t+h}\widehat{E}\[ (\partial_{x}F)(s,\widehat{U}_{s})\cdot\widehat{l}_s\widehat{\eta^1(s)}\widehat{\eta^2_s} \]ds
=\int_t^{t+h}E\[ (\partial_{x}F)(s,U_{s})\cdot l_s\eta^1\eta^2_s \]ds\\
=&\ E\[ E\[\int_t^{t+h}(\partial_{x}F)(s,U_{s})\cdot l_s\eta^2_s ds\big|\mathcal{F}_{T}^W\] \eta^1\]\\
=&\ E\[ E\[\big(\int_t^{t+h}(\partial_{x}F)(s,U_{s})\cdot l_sd\overleftarrow{B_s}\big)\cdot\eta^2\big|\mathcal{F}_{T}^W\] \eta^1\]\\
=&\ E\[ \big(\int_t^{t+h}(\partial_{x}F)(s,U_{s})\cdot l_sd\overleftarrow{B_s}\big)\cdot\eta^1 \eta^2\],\ 0\leq t<t+h\leq T.
  \end{split}
\end{equation*}
Consequently, summarising the above computation, we obtain
\begin{equation*}
  \begin{split}
E\[(F(t+h,U_{t+h})-F(t,U_t))\eta^1&\eta^2\]= E\[\(\int_{t}^{t+h}\!\Big\{(\partial_tF)(s,U_{s})
  \!+\!\frac{1}{2}\big((\partial_{xx}^2F)(s,U_{s})
           \!\cdot\! (|V_s|^2\!-\!|l_s|^{2})\big)\Big\}ds\\
 & - \int_{t}^{t+h} (\partial_{x}F)(s,U_{s})\cdot l_sd\overleftarrow{B_s}+  \int_{t}^{t+h} (\partial_{x}F)(s,U_{s})\cdot V_sdW_s\)\eta^1\eta^2\],
  \end{split}
\end{equation*}
and recalling the arbitrariness of $\eta^1\in L^{\infty}(\mathcal{F}_{T};\mathbb{R})$, $\eta^2\in L^{\infty}(\mathcal{F}_{0,T}^B;\mathbb{R})$
with $(\eta^1_t)$, $(\eta^2_t)$ bounded, we conclude
\begin{equation}\label{eqA.3.32}
  \begin{split}
F(&t+h,U_{t+h})-F(t,U_t)= \int_{t}^{t+h}\!\Big\{(\partial_tF)(s,U_{s})
  \!+\!\frac{1}{2}\big((\partial_{xx}^2F)(s,U_{s})
           \!\cdot\! (|V_s|^2\!-\!|l_s|^{2})\big)\Big\}ds\\
 & - \int_{t}^{t+h} (\partial_{x}F)(s,U_{s})\cdot l_sd\overleftarrow{B_s}+  \int_{t}^{t+h} (\partial_{x}F)(s,U_{s})\cdot V_sdW_s,\
  0\leq t<t+h\leq T,\ P\mbox{-a.s.}
  \end{split}
\end{equation}
We remark that, standard approximation, \eqref{eqA.3.32} can be extended to the case
where $((\partial_tF)(s,U_{s}))\in\mathcal{H}^2_{\mathcal{F}}(0,T;\mathbb{R})$.

\noindent\textbf{Step 4.} We combine now the results of the steps 2 and 3 with keeping
the assumptions which have been made there. Let $H\in C_b^{1,2,2}([0,T]\times\mathbb{R}\times\mathcal{P}_{2}(\mathbb{R}))$.
Then, due to Step 2, for all $x\in\mathbb{R}$, $0\leq t<t+h\leq T$,
\begin{equation*}
  \begin{split}
H(t\!+\!h,x,\!P_{Y_{t\!+\!h}})\!-\!H(t,x,\!P_{Y_{t}})\!=\!\! \int_{t}^{t+h}\!\!\Big\{\!(\partial_tH)(s,x,P_{Y_{s}})
  \!+\!\frac{1}{2}\widehat{E}\[\partial_{y}(\partial_{\mu}H)(s,x,P_{Y_{s}},\widehat{Y}_{s})
            (|\widehat{Z}_s|^2\!-\!|\widehat{g}_s|^2)\!\]\!\Big\}ds,
  \end{split}
\end{equation*}
and from Step 3 for $F(t,x):=H(t,x,P_{Y_{t}})$, $(t,x)\in[0,T]\times\mathbb{R}$,
\begin{equation*}
  \begin{split}
 &\ H(t+h,U_{t+h},P_{Y_{t+h}})-H(t,U_t,P_{Y_t})=F(t+h,U_{t+h})-F(t,U_t)\\
 =&\ \int_{t}^{t+h}\!\Big\{(\partial_tF)(s,U_{s})
  +\frac{1}{2}\big((\partial_{xx}^2F)(s,U_{s})
           \cdot (|V_s|^2-|l_s|^{2})\big)\Big\}ds\\
 & - \int_{t}^{t+h} (\partial_{x}F)(s,U_{s})\cdot l_sd\overleftarrow{B_s}+  \int_{t}^{t+h} (\partial_{x}F)(s,U_{s})\cdot V_sdW_s\\
 =&\ \int_t^{t+h}\Bigg\{(\partial_tH)(s,U_s,P_{Y_s})+\frac{1}{2}(\partial_{xx}^2 H)(s,U_s,P_{Y_{s}})(|V_s|^2\!-\!|l_s|^2)\\
  &\ +\frac{1}{2}\widehat{E}\[\partial_{y}(\partial_{\mu}H)(s,U_{s},P_{Y_{s}},\widehat{Y}_{s})
           \cdot (|\widehat{Z}_s|^2-|\widehat{g}_s|^2)\]    \Bigg\}ds
           +\int_t^{t+h}   (\partial_{x} H)(s,U_s,P_{Y_{s}})\cdot V_sdW_s\\
  &\   -\int_t^{t+h}   (\partial_{x} H)(s,U_s,P_{Y_{s}})\cdot l_sd\overleftarrow{B_s},\ 0\leq t<t+h\leq T,\ P\mbox{-a.s.}
  \end{split}
\end{equation*}
\end{proof}

\subsection{Mean-field BDSDEs}

We first give two classical estimates for solutions of backward doubly stochastic differential equations (BDSDEs for short), the proof is standard,
the readers may refer to, e.g., \cite{PP1994} and \cite{LX2021}.
\begin{lemma} \label{leA.1}
Suppose $(Y^i,Z^i)$ is the unique solution of the following BDSDE with data $(f_i,g_i,\theta_i)$,
\begin{equation}\label{eqA.2.1}
  \left\{
   \begin{array}{l}
   dY_s^{i}=-f_i(s,Y_s^{i},Z_s^{i})ds-g_i(s,Y_s^{i},Z_s^{i})d\overleftarrow{B_s}+Z_s^{i}dW_s,\\
   Y_T^{i}=\theta_i,\\
   \end{array}
  \right.
\end{equation}
where the integral w.r.t. $B$ is the It\^{o} backward one, denoted by $d\overleftarrow{B}$,
$\theta_i\in L^2(\Omega,\mathcal{F}_T, P;\mathbb{R}^d)$, and the coefficients
 $f_i: [0,T]\times\Omega\times\mathbb{R}^k\times\mathbb{R}^{k\times d}\rightarrow\mathbb{R}^{k\times d}$
and $g_i:[0,T]\times\Omega\times\mathbb{R}^k\times\mathbb{R}^{k\times d}\rightarrow\mathbb{R}^{k\times l}$, $i=1,2$, are jointly measurable and satisfy:

\noindent\emph{\textbf{Assumption (H10.1)}} \emph{(i)} $(g_i(t,\cdot,0,0))_{t\in[0,T]}\in\mathcal{H}_{\mathcal{F}}^2(0,T;\mathbb{R}^{k\times l})$;\\
\emph{(ii)} $g_i$ is Lipschitz in $(y,z)$, i.e., there exist constants $C>0$, and $0<\alpha<1$ such that for all\\
 \mbox{ } \ \ \    $t\in[0,T]$, $y_1,y_2\in\mathbb{R}^k$, $z_1,z_2\in\mathbb{R}^{k\times d}$, $P$-a.s.,
 $$|g_i(t,y_1,z_1)-g_i(t,y_2,z_2)|^2\leq C|y_1-y_2|^2+\alpha|z_1-z_2|^2;$$
\emph{(iii)} $(f_i(t,\cdot,0,0))_{t\in[0,T]}\in\mathcal{H}_{\mathcal{F}}^2(0,T;\mathbb{R}^{k})$;\\
\emph{(iv)} $f_i$ is Lipschitz in $(y,z)$, i.e., there exists a constant $C>0$ such that for all  $t\in[0,T]$, $y_1,y_2\in\mathbb{R}^k$,\\
 \mbox{ } \ \ \  $z_1,z_2\in\mathbb{R}^{k\times d}$, $P$-a.s.,
      $$|f_i(t,y_1,z_1)-f_i(t,y_2,z_2)|\leq C(|y_1-y_2|+|z_1-z_2|).$$
Then, for $(\overline{Y},\overline{Z}):=(Y^1,Z^1)-(Y^2,Z^2)$, $\overline{f}:=f_1-f_2$, $\overline{g}:=g_1-g_2$, $\overline{\theta}:=\theta_1-\theta_2$, we have the
following estimates:

\noindent \emph{(1)} There exists a constant $C>0$ depending only on the Lipschitz constant of $f$ and $g$, such that, for $t\in[0,T]$, $P$-a.s.,
\begin{equation}\label{eqA.2.2}
  \begin{split}
    E[\sup_{s\in[0,T]}|\overline{Y}_{s}|^{2}]& +E[\int_{0}^{T}|\overline{Z}_{s}|^{2}ds]
    \leq CE\[|\overline{\xi}|^{2}+\int_{0}^{T}\(|\overline{f}(s,Y_s^1,Z_s^1)|^{2}
 + |\overline{g}(s,Y_s^1,Z_s^1)|^{2}\)ds\].
  \end{split}
\end{equation}
 \noindent \emph{(2)} We suppose that, for some $p\geq 2$, $\overline{C}_p\alpha^{\frac{p}{2}}<1$. Here $\overline{C}_p:=2^{p-1}C^{\ast}_{p}((\frac{p}{p-1})^p+1)C_p'$, $C^{\ast}_{p}:=2^{-p-2}3^{p}p^{3p}+2^{\frac{p}{2}}$, $C_p':=(\frac{p}{p-1})^p3^{p-1}\(2C^p5^{p-1}\vee (6p^3)^p 5^{\frac{p}{2}-1}\)$,
and $C$ is the Lipschitz constant in Assumption (H10.1). Then there exists $C_p\in\mathbb{R}_{+}$ only
depending on the bounds of the coefficients and on $p$, such that
\begin{equation}\label{eqA.2.3}
  \begin{split}
     E\[\sup_{s\in[0,T]}|\overline{Y}_{s}|^{p}\!+\!(\int_{0}^{T}|\overline{Z}_{s}|^{2}ds)^{\frac{p}{2}}\]
\leq  C_p E\[|\overline{\theta}|^{p}\!+\!(\int_{0}^{T}|\overline{f}(s,Y_s^1,Z_s^1)|ds)^{p}\!+\!(\int_{0}^{T}|\overline{g}(s,Y_s^1,Z_s^1)|^2ds)^{\frac{p}{2}}\!\].
  \end{split}
\end{equation}

\end{lemma}

We now consider a more general case of BDSDE \eqref{eqA.2.1}.
Let $f: [0,T]\times\Omega\times\mathbb{R}^k\times\mathbb{R}^{k\times d}\times\mathcal{P}_{2}(\mathbb{R}^{d+k+k\times d})\rightarrow\mathbb{R}^{k}$,
  $g: [0,T]\times\Omega\times\mathbb{R}^k\times\mathbb{R}^{k\times d}\times\mathcal{P}_{2}(\mathbb{R}^{d+k+k\times d})\rightarrow\mathbb{R}^{k\times l}$,
and $h: [0,T]\times\mathcal{P}_{2}(\mathbb{R}^{d+k+k\times d})\rightarrow\mathbb{R}^{k\times l}$ be jointly measurable and satisfy the following standard assumptions:

\noindent\textbf{Assumption (H10.2)}\\
(i) $(g(t,\cdot,0,0,\delta_0))_{t\in[0,T]}\in\mathcal{H}_{\mathcal{F}}^2(0,T;\mathbb{R}^{k\times l})$, where $\delta_0$ is the Dirac measure with mass at \\
\mbox{ } \ \ \ $0\in\mathbb{R}^{d+k+k\times d}$;\\
(ii) $g$ is Lipschitz in $(y,z,\mu)$ with some Lipschitz constants $C>0$, and $\alpha_1,\alpha_2>0$ with $0<\alpha_1+\alpha_2$\\
 \mbox{ } \ \ \  $<1$ such that, for all $t\in[0,T]$, $\mu,\mu'\in\mathcal{P}_{2}(\mathbb{R}^{d+k}\times\mathbb{R}^{k\times d})$,
 $y_1,y_2\in\mathbb{R}^k$, $z_1,z_2\in\mathbb{R}^{k\times d}$, $P$-a.s.,
 $$|g(t,y_1,z_1,\mu)-g(t,y_2,z_2,\mu')|^2\leq C|y_1-y_2|^2+\alpha_1|z_1-z_2|^2+W_{2,C,\alpha_2}(\mu,\mu')^2;$$
(iii) $(f(t,\cdot,0,0,\delta_0))_{t\in[0,T]}\in\mathcal{H}^2_{\mathcal{F}}(0,T;\mathbb{R}^k)$;\\
(iv) $f$ is Lipschitz in $(y,z,\mu)$, i.e., there exists a constant $C>0$ such that, for all
$y_1,y_2\in\mathbb{R}^k$,\\
\mbox{ } \ \ \ $z_1,z_2\in\mathbb{R}^{k\times d}$, $\mu,\mu'\in\mathcal{P}_{2}(\mathbb{R}^{d+k+k\times d})$, $t\in[0,T]$, $P$-a.s.,
      $$|f(t,y_1,z_1,\mu)-f(t,y_2,z_2,\mu')|\leq C(W_2(\mu,\mu')+|y_1-y_2|+|z_1-z_2|);$$
(v) $(h(t,\delta_0))_{t\in[0,T]}\in\mathcal{H}_{\mathcal{F}}^2(0,T;\mathbb{R}^{k\times l})$;\\
(vi) $h$ is Lipschitz in $\mu$, i.e., there exists a constant $C\!>\!0$ such that, for all $t\in[0,T]$, $\mu,\mu'\in$\\
\mbox{ } \ \ \ $\mathcal{P}_{2}(\mathbb{R}^{d+k+k\times d})$,
 $$|h(t,\mu)-h(t,\mu')|^2\leq CW_{2}(\mu,\mu')^2;$$
(vii) $\xi\in L^2(\Omega,\mathcal{F}_T, P;\mathbb{R}^k)$;\\
(viii) $X\in \mathcal{H}^2_{\mathcal{F}}(0,T;\mathbb{R}^d)$.

\begin{theorem} \label{thA.1}
Under Assumption (H10.2), the following mean-field BDSDE
\begin{equation}\label{eqA.2.4}
\begin{split}
  Y_t=\xi&+\int_t^Tf(s,Y_s,Z_s,P_{(X_s,Y_s,Z_s)})ds+\int_t^Tg(s,Y_s,Z_s,P_{(X_s,Y_s,Z_s)})d\overleftarrow{B_s}\\
  &+\int_t^Th(s,P_{(X_s,Y_s,Z_s)})d\overleftarrow{B_s}-\int_t^TZ_sdW_s,\ 0\leq t\leq T,
  \end{split}
\end{equation}
has a unique solution $(Y,Z)\in\mathcal{S}_{\mathcal{F}}^2(0,T;\mathbb{R}^k)\times\mathcal{H}_{\mathcal{F}}^2(0,T;\mathbb{R}^{k\times d})$.
\end{theorem}
For the proof the reader is referred to \cite{LX2021}. Similar to Lemma \ref{leA.1}, we also have the following estimates for mean-field BDSDE.
In particular, we consider dimension $k=1$.

\begin{theorem} \label{thA.2}
Let $(Y^i,Z^i)$ be the unique solution of the following BDSDE with data $(f^i,g^i,h^i,\xi^i)$,
\begin{equation}\label{eqA.2.5}
\begin{split}
  Y_t^i=\xi^i&+\int_t^Tf^i(s,Y_s^i,Z_s^i,P_{(X_s^i,Y_s^i,Z_s^i)})ds+\int_t^Tg^i(s,Y_s^i,Z_s^i,P_{(X_s^i,Y_s^i,Z_s^i)})d\overleftarrow{B_s}\\
  &+\int_t^Th^i(s,P_{(X_s^i,Y_s^i,Z_s^i)})d\overleftarrow{B_s}-\int_t^TZ_s^idW_s,\ 0\leq t\leq T,
  \end{split}
\end{equation}
where $(f^i,g^i,h^i,\xi^i,X^i)$ satisfies Assumption (H10.2), $i=1,2$.

We denote $(\overline{X},\overline{Y},\overline{Z}):=(X^1,Y^1,Z^1)-(X^2,Y^2,Z^2)$, $\overline{f}:=f^1-f^2$, $\overline{g}:=g^1-g^2$,
$\overline{h}:=h^1-h^2$, $\overline{\xi}:=\xi^1-\xi^2$.
Then, there exists a constant $C>0$ only depending on the Lipschitz constants of the coefficients such that, $P$-a.s., $t\in[0,T]$,
\begin{equation}\label{eqA.2.6}
  \begin{split}
    E[&\sup_{s\in[0,T]}|\overline{Y}_{s}|^{2}] +E[\int_{0}^{T}|\overline{Z}_{s}|^{2}ds]
    \leq CE\[|\overline{\xi}|^{2}
+\int_{0}^{T}|\overline{f}(s,Y_s^1,Z_s^1,P_{(X_s^1,Y_s^1,Z_s^1)})|^{2}ds\\
&+\!\int_{0}^{T}|\overline{g}(s,Y_s^1,Z_s^1,P_{(X_s^1,Y_s^1,Z_s^1)})|^{2}ds\!+\!\int_{0}^{T}|\overline{h}(s,P_{(X_s^1,Y_s^1,Z_s^1)})|^{2}ds
\!+\!\int_{0}^{T}E[|\overline{X}_s|^2]ds\].\\
  \end{split}
\end{equation}
\end{theorem}
For the proof the reader is referred to Proposition 3.2 in \cite{LX2021}.

\begin{corollary} \label{corA.0}
For $i=1,2$, let $(Y^i,Z^i)\in\mathcal{S}_{\mathcal{F}}^2(0,T;\mathbb{R})\times\mathcal{H}_{\mathcal{F}}^2(0,T;\mathbb{R}^{ d})$ be the unique solution of BDSDE \eqref{eqA.2.5} with data $(f^i,g^i,h^i,\xi^i)$, $i=1,2$. We suppose that $f^i$, $g^i$ and $h^i$ are of the form   $f^i(s,y,z,P_{(Y_s,Z_s)})=u^i(s,y,z)+\widehat{E}[v^i(s,\widehat{Y}_s,\widehat{Z}_s)]$,
$g^i(s,y,z,P_{(Y_s,Z_s)})=a^i(s,y,z)+\widehat{E}[b^i(s,\widehat{Y}_s,\widehat{Z}_s)]$,
$h^i(s,P_{(Y_s,Z_s)})=\widehat{E}[q^i(s,\widehat{Y}_s,\widehat{Z}_s)]$,
where $u^i,v^i,a^i,b^i$ defined over $[0,T]\times\Omega\times\mathbb{R}\times\mathbb{R}^{ d}$,
and $q^i$ defined over $[0,T]\times\mathbb{R}\times\mathbb{R}^{ d}$ satisfy Assumption (H10.1).
Then, we have, for a constant $C$ only depending on the Lipschitz constants in (H10.1), $P$-a.s.,
\begin{equation}\label{eqA.2.7}
  \begin{split}
    E[&\!\sup_{s\in[0,T]}|\overline{Y}_{s}|^{2}]\!+\!E[\int_{0}^{T}|\overline{Z}_{s}|^{2}ds]
\!\leq\! CE\[|\overline{\xi}|^{2}\!+\!\int_{0}^{T}\!\(|(u^1-u^2)(s,Y_s^1,Z_s^1)|^{2}\!+\!|(a^1-a^2)(s,Y_s^1,Z_s^1)|^{2}\)ds\\
& \!+\!\int_{0}^{T}\!\(|\widehat{E}[(v^1-v^2)(s,Y_s^1,Z_s^1)]|^{2}ds\!+\!|\widehat{E}[(b^1-b^2)(s,Y_s^1,Z_s^1)]|^{2}ds
 \!+\!|\widehat{E}[(q^1-q^2)(s,Y_s^1,Z_s^1)]|^{2}\)ds\].\\
  \end{split}
\end{equation}
\end{corollary}

\begin{theorem} \label{thA.3}
Let Assumption (H10.2) be satisfied. We suppose that,
for some $p\geq 2$, $\alpha_1,\alpha_2>0$ are small enough, such that $\overline{C}_p(\alpha_1+\alpha_2)^{\frac{p}{2}}<1$. Here $\overline{C}_p:=2^{p-1}C^{\ast}_{p}((\frac{p}{p-1})^p+1)C_p'$, where $C_p':=(\frac{p}{p-1})^p3^{p-1}\(2C^p5^{p-1}\vee (6p^3)^p 5^{\frac{p}{2}-1}\)$,
 $C^{\ast}_{p}:=2^{-p-2}3^{p}p^{3p}+2^{\frac{p}{2}}$,
and $C$ is the Lipschitz constant in Assumption (H10.2). Let $(Y,Z)$ denote the unique solution
of the mean-field BDSDE \eqref{eqA.2.4} with data $(f,g,h,\xi)$.
Then there exists $C_p\in\mathbb{R}_{+}$ only
depending on the Lipschitz constants in (H10.1) and on $ p$, such that
\begin{equation}\label{eqA.2.9}
  \begin{split}
    E[\sup_{s\in[0,T]}&|Y_{s}|^{p}]+E\[\(\int_{0}^{T}|Z_{s}|^{2}ds\)^{\frac{p}{2}}\]\leq C_p E\[|\xi|^{p}+\(\int_{0}^{T}|X_{s}|^{2}ds\)^{\frac{p}{2}}\\
& +(\int_{0}^{T}|f(s,0,0,\delta_0)|ds)^{p}
 +(\int_{0}^{T}|g(s,0,0,\delta_0)|^2ds)^{\frac{p}{2}}+(\int_{0}^{T}|h(s,\delta_0)|^2ds)^{\frac{p}{2}}\].
  \end{split}
\end{equation}
\end{theorem}
For the proof the reader is referred to Proposition 3.3 in \cite{LX2021}.

Now we give the estimates for a special type of mean-field BDSDEs, which are used frequently in our work. We suppose that

\noindent\textbf{Assumption (H10.3)} (i) $\xi\in L^2(\Omega,\mathcal{F}_T, P;\mathbb{R})$;\\
(ii) $\lambda=(\lambda_s)$, $\beta=(\beta_s)$, $\delta=(\delta_s)$, $\gamma=(\gamma_s)$ are bounded $\{\mathcal{F}_s\}$-adapted measurable processes;\\
(iii) $\zeta=(\zeta_s)$, $\theta=(\theta_s)$, $\eta=(\eta_s)$, $\rho=(\rho_s)$ are bounded $\{\mathcal{A}_s\}$-adapted measurable processes with\\
 \mbox{}\ \ \ \ \ \  $\mathcal{A}_s=\mathcal{F}_s\otimes\widehat{\mathcal{F}}$,
 where $(\widehat{\Omega},\widehat{\mathcal{F}},\widehat{P})$ is a copy of $(\Omega,\mathcal{F},P)$;\\
(iv) There exist constants $\alpha_1,\alpha_2>0$ with $0<\alpha_1+\alpha_2<1$ such that $|\delta|^2\leq\alpha_1$, $|\rho|^2\leq\alpha_2$;\\
(v) $R=(R(s))$ is $\{\mathcal{F}_s\}$-adapted and measurable with $E[\int_0^T|R(r)|^2dr]<+\infty$;\\
(vi) $H=(H(s))$ is $\{\mathcal{F}_s\}$-adapted and measurable with $E[\int_0^T|H(r)|^2dr]<+\infty$.

With similar arguments as for the Theorems \ref{thA.2} and \ref{thA.3}, we get the following corollaries.

\begin{corollary} \label{corA.1}
Suppose Assumption (H10.3) is satisfied. Let $(Y,Z)\in\mathcal{S}_{\mathcal{F}}^2(0,T;\mathbb{R})\times\mathcal{H}_{\mathcal{F}}^2(0,T;\mathbb{R}^{ d})$
be the unique solution of the linear BDSDE
\begin{equation}\label{eqA.2.10}
  \begin{split}
    Y_s=& \xi+\int_s^T\(R(r)+\lambda_rY_r+\widehat{E}[\zeta_r\widehat{Y}_r]+\gamma_rZ_r+\widehat{E}[\theta_r\widehat{Z}_r]\)dr\\
    & +\!\int_s^T\!\(H(r)+\beta_rY_r\!+\!\widehat{E}[\eta_r\widehat{Y}_r]\!+\!\delta_rZ_r\!+\!\widehat{E}[\rho_r\widehat{Z}_r]\)d\overleftarrow{B_r}
    \!-\!\int_s^T\!Z_rdW_r\ s\in[0,T],\ P\mbox{-a.s.},
  \end{split}
\end{equation}
where $(\widehat{Y},\widehat{Z})$ denotes a copy of $(Y,Z)$ on $(\widehat{\Omega},\widehat{\mathcal{F}},\widehat{P})$
(i.e., $P_{(Y,Z)}=\widehat{P}_{(\widehat{Y},\widehat{Z})}$).
Then there exists a constant $C\in\mathbb{R}_{+}$ only depending on the bounds of the coeffcients such that
\begin{equation}\label{eqA.2.11}
  \begin{split}
    &\ E[\sup_{s\in[0,T]}|Y_{s}|^{2}]+E[\int_{0}^{T}|Z_{s}|^{2}ds]\leq CE\[|\xi|^{2}+\int_{0}^{T}|R(t)|^{2}dt+\int_{0}^{T}|H(t)|^{2}dt\].\\
  \end{split}
\end{equation}
\end{corollary}

\begin{corollary} \label{corA.2}
Let Assumption (H10.3) hold. Moreover, we suppose that,
for some $p\geq 2$, $\overline{C}_p(\alpha_1+\alpha_2)^{\frac{p}{2}}<1$. Here $\overline{C}_p:=2^{p-1}C^{\ast}_{p}((\frac{p}{p-1})^p+1)C_p'$, $C_p':=(\frac{p}{p-1})^p3^{p-1}\(2C^p5^{p-1}\vee (6p^3)^p 5^{\frac{p}{2}-1}\)$, and $C^{\ast}_{p}:=2^{-p-2}3^{p}p^{3p}+2^{\frac{p}{2}}$,
where $C$ is the bound of the processes $\lambda$, $\beta$, $\delta$, $\gamma$,
$\zeta$, $\theta$, $\eta$ and $\rho$ in Assumption (H10.3).
Then there exists $C_p\in\mathbb{R}_{+}$ only
depending on the bounds of the coefficients and on $p$, such that, for the unique solution $(Y,Z)$
of the BDSDE \eqref{eqA.2.10}, we have
\begin{equation}\label{eqA.2.12}
  \begin{split}
    & E[\sup_{s\in[0,T]}|Y_{s}|^{p}]+E\[\(\int_{0}^{T}|Z_{s}|^{2}ds\)^{\frac{p}{2}}\]
\leq C_pE\[|\xi|^{p}+\(\int_{0}^{T}|R(t)|dt\)^{p}+\(\int_{0}^{T}|H(t)|^{2}dt\)^{\frac{p}{2}}\].
  \end{split}
\end{equation}
\end{corollary}

Let us, finally, briefly discuss
\subsection{ The special case of $g$ affine with respect to $z$}
The results are mainly used in the proof of Proposition \ref{prop9.2}.

Let $f: \mathbb{R}^d\times\mathbb{R}\times\mathbb{R}^d\times\mathcal{P}_{2}(\mathbb{R}^d\times\mathbb{R}\times\mathbb{R}^d)\rightarrow\mathbb{R}$,
$g: \mathbb{R}^d\times\mathbb{R}\times\mathbb{R}^d\times\mathcal{P}_{2}(\mathbb{R}^d\times\mathbb{R}\times\mathbb{R}^d)\rightarrow\mathbb{R}^{l}$,
$h: \mathcal{P}_{2}(\mathbb{R}^d\times\mathbb{R}\times\mathbb{R}^d)\rightarrow\mathbb{R}^{l}$
and $\Phi :\mathbb{R}^d\times\mathcal{P}_{2}(\mathbb{R}^d)\rightarrow\mathbb{R}$ be deterministic functions and satisfying:\\
\textbf{Assumption (H10.4)}\\
(i) Let $\Phi\!\in\! C_b^{1,1}(\mathbb{R}^{d}\times\mathcal{P}_{2}(\mathbb{R}^{d}))$,
$f\!\in\! C_b^{1,1}(\mathbb{R}^{d+1+d}\times\mathcal{P}_{2}(\mathbb{R}^{d+1+d}))$
and $h\!\in\! C_b^{1}(\mathcal{P}_{2}(\mathbb{R}^{d+1+d});\mathbb{R}^{l})$.\\
(ii) The coefficient $g$ is affine in $z$: for all $x\in\mathbb{R}^{d}$, $y\in\mathbb{R}$, $z\in\mathbb{R}^{d}$, $\mu\in\mathcal{P}_{2}(\mathbb{R}^{d+1}\times\mathbb{R}^d)$,
$$g(x,y,z,\mu)=g^1(x,y,\mu)+g^2(\mu(\cdot\times\mathbb{R}\times\mathbb{R}^{d}))z,$$
where $g^1\in\! C_b^{1,1}(\mathbb{R}^{d+1}\times\mathcal{P}_{2}(\mathbb{R}^{d+1+d});\mathbb{R}^{l})$ and
$g^2\in\! C_b^{1}(\mathcal{P}_{2}(\mathbb{R}^{d});\mathbb{R}^{l\times d})$.
In addition we suppose $|g^2|^2\leq \alpha_1$,
$\sum_{k=1}^{d}\sum_{i=1}^l|(\partial_{\mu}g_{i}^1)_{d+1+k}|^2\leq \alpha_2$,
for constants $\alpha_1,\alpha_2>0$ with $0\!<\!\alpha_1\!+\!\alpha_2\!<\!1$.

Given $x\in\mathbb{R}^{d}$ and $\xi\in L^{2}(\mathcal{G}_t;\mathbb{R}^d)$ we consider the following both BDSDEs:
\begin{equation}\label{eqA.4.1}
  \left\{
   \begin{array}{l}
   dY_s^{t,\xi}=-f(\Pi_s^{t,\xi},P_{\Pi_s^{t,\xi}})ds-
   \big(g(\Pi_s^{t,\xi},P_{\Pi_s^{t,\xi}})+h(P_{\Pi_s^{t,\xi}})\big)d\overleftarrow{B_s}+Z_s^{t,\xi}dW_s,\ s\in[t,T],\\
   Y_T^{t,\xi}=\Phi (X_T^{t,\xi},P_{X_T^{t,\xi}}),\\
   \end{array}
  \right.
\end{equation}
and
\begin{equation}\label{eqA.4.2}
  \left\{
   \begin{array}{l}
   dY_s^{t,x,\xi}=-f(\Pi_s^{t,x,\xi},P_{\Pi_s^{t,\xi}})ds-
   \big(g(\Pi_s^{t,x,\xi},P_{\Pi_s^{t,\xi}})+h(P_{\Pi_s^{t,\xi}})\big)d\overleftarrow{B_s}+Z_s^{t,x,\xi}dW_s,\ s\in[t,T],\\
   Y_T^{t,x,\xi}=\Phi (X_T^{t,x,\xi},P_{X_T^{t,\xi}}),\\
   \end{array}
  \right.
\end{equation}
where $\Pi_s^{t,\xi}:=(X_s^{t,\xi},Y_s^{t,\xi},Z_s^{t,\xi})$, $\Pi_s^{t,x,\xi}:=(X_s^{t,x,\xi},Y_s^{t,x,\xi},Z_s^{t,x,\xi})$.
Recall that the processes $X^{t,\xi}$ and $X^{t,x,\xi}$ are the solution of SDEs \eqref{eq3.1} and \eqref{eq3.2}, respectively.

By using similar techniques as in the proofs of the Theorems \ref{th6.1} and \ref{th6.2} and the Propositions \ref{prop4.2}, \ref{prop6.1} and \ref{prop6.2},
we obtain the following Propositions \ref{prop6.1+1}, \ref{prop6.1+2} and \ref{prop6.6}.
For simplicity of redaction but without loss of generality, we restrict to the dimensions $d = 1$, $l=1$
and to $f(\Pi_s^{t,x,\xi},P_{\Pi_s^{t,\xi}})=f(Z_s^{t,x,P_{\xi}})$,
$g(\Pi_s^{t,x,\xi},P_{\Pi_s^{t,\xi}})=g^2(P_{X_s^{t,\xi}})Z_s^{t,x,P_{\xi}}$, $h(P_{\Pi_s^{t,\xi}})=h(P_{(Y_s^{t,\xi},Z_s^{t,\xi})})$
and $\Phi(X_T^{t,x,P_{\xi}},P_{X_T^{t,\xi}})=\Phi(X_T^{t,x,P_{\xi}})$.

\begin{proposition} \label{prop6.1+1}
Let the Assumptions (H5.1) and (H10.4) hold true.
Then for all $(t,x)\in[0,T]\times \mathbb{R}$, $\xi\in L^{2}(\mathcal{G}_t;\mathbb{R})$,
$(Y_s^{t,x,P_{\xi}},Z_s^{t,x,P_{\xi}})\!\in\! L^2(t,T;(\mathbb{D}^{1,2})^{2})$
and a version of $\{D_{\theta}Y_s^{t,x,P_{\xi}},$
$D_{\theta}Z_s^{t,x,P_{\xi}}:\theta,s\in[t,T]\}$ is given by:\\
\indent \emph{(i)} $D_{\theta}Y_s^{t,x,P_{\xi}}=0$, $D_{\theta}Z_s^{t,x,P_{\xi}}=0$, $t\leq s<\theta\leq T$;\\
\indent \emph{(ii)} $\{D_{\theta}Y^{t,x,P_{\xi}},D_{\theta}Z^{t,x,P_{\xi}}:s\in[\theta,T]\}$
is the unique solution of the linear BDSDE: $s\in[t,T]$,
\begin{equation}\label{eq6.2+1}
  \begin{split}
D_{\theta} Y_s^{t,x,P_{\xi}}=&\ \partial_{x}\Phi (X_T^{t,x,P_{\xi}})D_{\theta}X_T^{t,x,P_{\xi}}
+\int_s^T\partial_{z}f(Z_r^{t,x,P_{\xi}})D_{\theta}Z_r^{t,x,P_{\xi}}dr\\
 &\ +\int_s^Tg^2(P_{X_r^{t,\xi}})D_{\theta}Z_r^{t,x,P_{\xi}}d\overleftarrow{B_r}
 \!-\!\int_s^T\!D_{\theta}Z_r^{t,x,P_{\xi}}dW_r,\ d \theta dP\mbox{-a.e.},\ t\leq\theta\leq s.
  \end{split}
\end{equation}
Moreover, $Z_s^{t,x,P_{\xi}}=P\mbox{-}\displaystyle{\lim_{s<u\downarrow s}}D_{s}Y_u^{t,x,P_{\xi}}$, $d s dP$-a.e.
Furthermore, if in addition Assumption (H4.2) is satisfied, then there exists a constant $C_p>0$ only depending on the Lipschitz constants of the coefficients,
 such that for all $t\in[0,T]$, $x,x'\in\mathbb{R}$,
$\xi,\xi'\in L^{2}(\mathcal{G}_t;\mathbb{R})$, $P$-a.s.,
\begin{equation}\label{eq6.2+2}
  \begin{split}
    &\ \mbox{\emph{(i)}}\ E\[\sup_{s\in[t,T]}|D_{\theta}Y_s^{t,x,P_{\xi}}|^p+(\int_t^T|D_{\theta}Z_s^{t,x,P_{\xi}}|^2ds)^{\frac{p}{2}}\]\leq C_p,\
    \mbox{for all}\ p\in[2,p_0];\\
  &\ \mbox{\emph{(ii)}}\ E\[\sup_{s\in[t,T]}|D_{\theta}Y_s^{t,x,P_{\xi}}-D_{\theta}Y_s^{t,x',P_{\xi'}}|^p
  +(\int_t^T|D_{\theta}Z_s^{t,x,P_{\xi}}-D_{\theta}Z_s^{t,x',P_{\xi'}}|^2ds)^{\frac{p}{2}}\]\\
 &\ \hspace{30pt} \leq C_p\(|x-x'|^p+W_2(P_{\xi},P_{\xi'})^p\),\ \mbox{for all}\ p\in[2,\frac{p_0}{2}] .
  \end{split}
\end{equation}
In particular, there exists a constant $C_p>0$ only depending on the Lipschitz constants of the coefficients,
 such that for all $x,x'\in\mathbb{R}$,
$\xi,\xi'\in L^{2}(\mathcal{G}_t;\mathbb{R})$, $d s dP$-a.e., $s\in[t,T]$,
\begin{equation}\label{eq6.2+3}
  \begin{split}
    &\ \mbox{\emph{(i)}}\ E[|Z_s^{t,x,P_{\xi}}|^p]\leq C_p,\ \mbox{for all}\ p\in[2,p_0];\\
  &\ \mbox{\emph{(ii)}}\ E[|Z_s^{t,x,P_{\xi}}-Z_s^{t,x',P_{\xi'}}|^p] \leq C_p\(|x-x'|^p+W_2(P_{\xi},P_{\xi'})^p\),\ \mbox{for all}\ p\in[2,\frac{p_0}{2}] .
  \end{split}
\end{equation}
\end{proposition}
\begin{proof}
It is standard to prove that $Y_s^{t,x,P_{\xi}}$ and $Z_s^{t,x,P_{\xi}}$ are Malliavin differentiable under our assumptions,
and $Z_r^{t,x,P_{\xi}}=P\mbox{-}\displaystyle{\lim_{r<s\downarrow r}}D_rY_s^{t,x,P_{\xi}}$, so we omit proving this here.
Thus, it suffices to prove \eqref{eq6.2+2}. For \eqref{eq6.2+2}, from Lemma \ref{leA.1}-(2) we get
$ \displaystyle E\[\sup_{s\in[t,T]}|D_{\theta}Y_s^{t,x,P_{\xi}}|^{p_0}+(\int_t^T|D_{\theta}Z_s^{t,x,P_{\xi}}|^2ds)^{\frac{p_0}{2}}\]\leq C_{p_0}.$
Then from H\"{o}lder's inequality, for all $p\in[2,p_0]$ we have
$ \displaystyle  E\[\sup_{s\in[t,T]}|D_{\theta}Y_s^{t,x,P_{\xi}}|^p+(\int_t^T|D_{\theta}Z_s^{t,x,P_{\xi}}|^2ds)^{\frac{p}{2}}\]$ $\leq C_p.$
 Now we prove \eqref{eq6.2+2}-(ii).

For all $0\leq t\leq T$, $t\leq\theta\leq s\leq T$, $x,x'\in\mathbb{R}$, $\xi,\xi'\in L^{2}(\mathcal{G}_t;\mathbb{R})$, from \eqref{eq6.2+1} we get the following BDSDE:
\begin{equation}\label{eq6.2+4}
  \begin{split}
 & D_{\theta}Y_s^{t,x,P_{\xi}}\!\!-\!D_{\theta} Y_s^{t,x',P_{\xi'}}
 \!=\!  I(t,x,P_{\xi},x',P_{\xi'})\!+\!\int_s^T\!\!\!R(r,x,P_{\xi},x',P_{\xi'})dr\!+\!\int_s^T\!\!\!H(r,x,P_{\xi},x',P_{\xi'})d\overleftarrow{B_r}\\
&\quad \!+\!\int_s^T\!\partial_{z}f(Z_r^{t,x,P_{\xi}})\(D_{\theta}Z_r^{t,x,P_{\xi}}\!-\!D_{\theta}Z_r^{t,x',P_{\xi'}}\)dr
 \!+\!\int_s^T\!g^2(P_{X_r^{t,\xi}})\(D_{\theta}Z_r^{t,x,P_{\xi}}\!-\!D_{\theta}Z_r^{t,x',P_{\xi'}}\)d\overleftarrow{B_r}\\
&\quad  -\int_s^T\(D_{\theta}Z_r^{t,x,P_{\xi}}-D_{\theta}Z_r^{t,x',P_{\xi'}}\)dW_r,
  \end{split}
\end{equation}
where
\begin{equation*}
  \begin{split}
 &\ I(t,x,P_{\xi},x',P_{\xi'}):= \partial_{x}\Phi (X_T^{t,x,P_{\xi}})D_{\theta}X_T^{t,x,P_{\xi}}-\partial_{x}\Phi (X_T^{t,x',P_{\xi'}})D_{\theta}X_T^{t,x',P_{\xi'}};\\
 &\ R(r,x,P_{\xi},x',P_{\xi'}):=\(\partial_{z}f(Z_r^{t,x,P_{\xi}})-\partial_{z}f(Z_r^{t,x',P_{\xi'}})\)D_{\theta}Z_r^{t,x',P_{\xi'}};\\
 &\ H(r,x,P_{\xi},x',P_{\xi'}):=\(g^2(P_{X_r^{t,\xi}})-g^2(P_{X_r^{t,\xi'}})\)D_{\theta}Z_r^{t,x',P_{\xi'}}.
  \end{split}
\end{equation*}
From \eqref{eq6.2+2}-(i), the Lemmas \ref{leA.1} and \ref{le3.1}, Proposition \ref{prop3.1} and Remark \ref{re4.1} as well as from our assumptions it follows that
\begin{equation}\label{eq6.2+5}
  \begin{split}
  &\  E\[\sup_{s\in[\theta,T]}|D_{\theta}Y_s^{t,x,P_{\xi}}-D_{\theta}Y_s^{t,x',P_{\xi'}}|^{\frac{p_0}{2}}
  +(\int_{\theta}^T|D_{\theta}Z_s^{t,x,P_{\xi}}-D_{\theta}Z_s^{t,x',P_{\xi'}}|^2ds)^{\frac{p_0}{4}}\]\\
 \leq&\  C_{\frac{p_0}{2}} E\[|I(t,x,P_{\xi},x',P_{\xi'})|^{\frac{p_0}{2}}\!+\!(\int_{\theta}^{T}\!|R(r,x,P_{\xi},x',P_{\xi'})|dr)^{\frac{p_0}{2}}\!
 +\!(\int_{\theta}^{T}\!|H(r,x,P_{\xi},x',P_{\xi'})|^2dr)^{\frac{p_0}{4}}\]\\
 \leq&\  C_{\frac{p_0}{2}} \!\(\!E\[\(\int_{\theta}^T|(\partial_{z}f)(Z_r^{t,x,P_{\xi}})\!-\!(\partial_{z}f)(Z_r^{t,x',P_{\xi'}})|^2dr\)^{\frac{p_0}{2}}\]\)^{\frac{1}{2}}
 \!\cdot\! \(E\[\(\int_{\theta}^T|D_{\theta}Z_r^{t,x',P_{\xi'}}|^2dr\)^{\frac{p_0}{2}}\]\)^{\frac{1}{2}}\\
 &\ +\!C_{\frac{p_0}{2}} E\[\(\int_{\theta}^T|g^2(P_{X_r^{t,\xi}})\!-\!g^2(P_{X_r^{t,\xi'}})|^2\cdot|D_{\theta}Z_r^{t,x',P_{\xi'}}|^2dr\)^{\frac{p_0}{4}}\]
 \!+\! C_{\frac{p_0}{2}}\(|x\!-\!x'|^{\frac{p_0}{2}}\!+\!W_2(P_{\xi},P_{\xi'})^{\frac{p_0}{2}}\)\\
\leq &\ C_{\frac{p_0}{2}}\(|x-x'|^{\frac{p_0}{2}}+W_2(P_{\xi},P_{\xi'})^{\frac{p_0}{2}}\)
     +C_{\frac{p_0}{2}} \(E\[\(\int_{\theta}^T|Z_r^{t,x,P_{\xi}}-Z_r^{t,x',P_{\xi'}}|^2dr\)^{p_0}\]\)^{\frac{1}{2}} \\
&\ +C_{\frac{p_0}{2}}\sup_{r\in[\theta,T]}W_2(P_{X_r^{t,\xi}},P_{X_r^{t,\xi'}})^{\frac{p_0}{2}}
      E\[\!\(\!\int_{\theta}^T \!|D_{\theta}Z_r^{t,x',P_{\xi'}}|^2dr\)^{\frac{p_0}{4}}\]\\
\leq&\ C_{\frac{p_0}{2}}\(|x-x'|^{\frac{p_0}{2}}+W_2(P_{\xi},P_{\xi'})^{\frac{p_0}{2}}\).
  \end{split}
\end{equation}
Then from H\"{o}lder's inequality, for all $p\in[2,\frac{p_0}{2}]$ we have \eqref{eq6.2+2}-(ii).
\end{proof}

\begin{remark} \label{reA.1}
In analogy to the proof in Proposition \ref{prop4.0}, it can easily be checked that
under the Assumptions (H5.1) and (H10.4), there exists a constant $C>0$ only depending on the Lipschitz
constants of the coefficients, such that for all $t\in[0,T]$, $x,x'\in\mathbb{R}$, $\xi,\xi'\in L^{2}(\mathcal{G}_t;\mathbb{R})$,
 $P$-a.s.,
\begin{equation}\nonumber
  \begin{split}
  E\[\sup_{s\in[t,T]}|D_{\theta}Y_s^{t,x,P_{\xi}}|^2+\int_t^T|D_{\theta}Z_s^{t,x,P_{\xi}}|^2ds\]\leq C,\
   \mbox{and}\ E[|Z_s^{t,x,P_{\xi}}|^2]\leq C.
  \end{split}
\end{equation}
\end{remark}

\begin{proposition} \label{prop6.1+2}
Under the Assumptions (H5.1) and (H10.4), the $L^2$-derivative, $(\partial_xY^{t,x,P_{\xi}},$
$\partial_xZ^{t,x,P_{\xi}})$, of the solution of Eq. \eqref{eq4.2} with respect to $x$ exists and is the unique solution of
the following BDSDE:
\begin{equation}\label{eq6.2+6}
  \begin{split}
 \partial_{x}Y_s^{t,x,P_{\xi}}=&\ \partial_{x}\Phi (X_T^{t,x,P_{\xi}})\partial_{x}X_T^{t,x,P_{\xi}}
+\int_s^T\partial_{z}f(Z_r^{t,x,P_{\xi}})\partial_{x}Z_r^{t,x,P_{\xi}}dr\\
 &\ +\int_s^Tg^2(P_{X_r^{t,\xi}})\partial_{x}Z_r^{t,x,P_{\xi}}d\overleftarrow{B_r}
 -\int_s^T\partial_{x}Z_r^{t,x,P_{\xi}}dW_r,\  s\in[t,T].
  \end{split}
\end{equation}
Furthermore, if in addition Assumption (H4.2), then there exists a constant $C_p>0$ only depending on the Lipschitz
constants of the coefficients, such that for all $t\in[0,T]$, $x,x'\in\mathbb{R}$, $\xi,\xi'\in L^{2}(\mathcal{G}_t;\mathbb{R})$,
 $P$-a.s.,
\begin{equation}\label{eq6.2+7}
  \begin{split}
  &\ \emph{(i)}\ E\[\sup_{s\in[t,T]}|\partial_xY_s^{t,x,P_{\xi}}|^p+(\int_t^T|\partial_xZ_s^{t,x,P_{\xi}}|^2ds)^{\frac{p}{2}}ds\big|\mathcal{G}_t\]\leq C_p,\
   \mbox{for all}\ p\in[2,p_0],\\
  &\ \emph{(ii)}\ E\[\sup_{s\in[t,T]}|\partial_xY_s^{t,x,P_{\xi}}-\partial_xY_s^{t,x',P_{\xi'}}|^p
  +(\int_t^T|\partial_xZ_s^{t,x,P_{\xi}}-\partial_xZ_s^{t,x',P_{\xi'}}|^2 ds)^{\frac{p}{2}}\big|\mathcal{G}_t\]\\
  &\hspace{1cm}     \leq C_p\(|x-x'|^p+W_2(P_{\xi},P_{\xi'})^p\),\ \mbox{for all}\ p\in[2,\frac{p_0}{2}].
  \end{split}
\end{equation}
\end{proposition}
Indeed, using \eqref{eq6.2+3}, the proof is similar to that of Theorem \ref{th6.1} and Proposition \ref{prop6.1}, and so we omit it here.

\begin{remark} \label{reA.2}
In analogy to the proof in Proposition \ref{prop4.0}, it can easily be checked that
under the Assumptions (H5.1) and (H10.4), there exists a constant $C>0$ only depending on the Lipschitz
constants of the coefficients, such that for all $t\in[0,T]$, $x,x'\in\mathbb{R}$, $\xi,\xi'\in L^{2}(\mathcal{G}_t;\mathbb{R})$,
 $P$-a.s.,
\begin{equation}\nonumber
  \begin{split}
E\[\sup_{s\in[t,T]}|\partial_xY_s^{t,x,P_{\xi}}|^2+\int_t^T|\partial_xZ_s^{t,x,P_{\xi}}|^2dsds\big|\mathcal{G}_t\]\leq C.
  \end{split}
\end{equation}
\end{remark}

\begin{proposition} \label{prop6.6}
Assume (H5.1) and (H10.4). Then, for all $0\leq t\leq s \leq T$, $x\in\mathbb{R}$,
the lifted processes $L^{2}(\mathcal{G}_t;\mathbb{R})\ni\xi\rightarrow Y_s^{t,x,\xi}:=Y_s^{t,x,P_{\xi}}\in L^{2}(\mathcal{F}_s;\mathbb{R})$,
and $L^{2}(\mathcal{G}_t;\mathbb{R})\ni\xi\rightarrow (Z_s^{t,x,\xi}:=Z_s^{t,x,P_{\xi}})\in\mathcal{H}_{\mathcal{F}}^2(t,T;\mathbb{R})$
are Fr\'{e}chet differentiable, with the Fr\'{e}chet derivatives
\begin{equation}\label{eq6.39}
 \begin{split}
   & DY_s^{t,x,\xi}(\eta)=\overline{E}\[\partial_{\mu}Y_s^{t,x,P_{\xi}}(\overline{\xi})\overline{\eta}\],\ s\in[t,T],\ P\mbox{-}a.s.,\\
  & DZ_s^{t,x,\xi}(\eta)=\overline{E}\[\partial_{\mu}Z_s^{t,x,P_{\xi}}(\overline{\xi})\overline{\eta}\],\  dsdP\mbox{-}a.e.,
 \end{split}
\end{equation}
for $\eta\in L^{2}(\mathcal{G}_t;\mathbb{R})$, where for all $y\in\mathbb{R}$,
$(\partial_{\mu}Y_r^{t,x,P_{\xi}}(y),\partial_{\mu}Z_r^{t,x,P_{\xi}}(y))\in\mathcal{S}^2_{\mathcal{F}}(t,T;\mathbb{R})\times \mathcal{H}^2_{\mathcal{F}}(t,T;\mathbb{R})$
is the unique solution of the following BDSDE:
\begin{equation}\label{eq6.40}
  \begin{split}
  \partial_{\mu}&Y_s^{t,x,P_{\xi}}(y)\!\!=\!\! \partial_{x}\Phi (X_T^{t,x,P_{\xi}})\partial_{\mu}X_{T}^{t,x,P_{\xi}}(y)
  \!\!+\!\! \int_s^T\!\!\partial_{z}f(Z_r^{t,x,P_{\xi}})\partial_{\mu}Z_r^{t,x,P_{\xi}}(y)  dr
  \!\!+\!\!\int_s^T\!\!g^2(P_{X_r^{t,\xi}})\partial_{\mu}Z_r^{t,x,P_{\xi}}(y) d\overleftarrow{B_r}\\
&\   +\int_s^T\widehat{E}\[\!(\partial_{\mu}g^2)(P_{X_r^{t,\xi}},\widehat{X}_r^{t,y,P_{\xi}}) \partial_{x}\widehat{X}_r^{t,y,P_{\xi}}
 \!+\!(\partial_{\mu}g^2)(P_{X_r^{t,\xi}},\widehat{X}_r^{t,\widehat{\xi}}) \partial_{\mu}\widehat{X}_r^{t,\widehat{\xi},P_{\xi}}(y) \]Z_r^{t,x,P_{\xi}}d\overleftarrow{B_r}\\
&\   +\int_s^T\widehat{E}\[\Big\langle(\partial_{\mu}h)(P_{(Y_r^{t,\xi},Z_r^{t,\xi})},(\widehat{Y}_r^{t,y,P_{\xi}},\widehat{Z}_r^{t,y,P_{\xi}})),
\left(   \begin{matrix}    \partial_{x}\widehat{Y}_r^{t,y,P_{\xi}}\\ \partial_{x}\widehat{Z}_r^{t,y,P_{\xi}} \end{matrix}\right)\Big\rangle\]d\overleftarrow{B_r}\\
&\ +\int_s^T\widehat{E}\[\Big\langle(\partial_{\mu}h)(P_{(Y_r^{t,\xi},Z_r^{t,\xi})},(\widehat{Y}_r^{t,\widehat{\xi}},\widehat{Z}_r^{t,\widehat{\xi}})),
 \left(   \begin{matrix}   \partial_{\mu}\widehat{Y}_r^{t,\widehat{\xi},P_{\xi}}(y)\\ \partial_{\mu}\widehat{Z}_r^{t,\widehat{\xi},P_{\xi}}(y) \end{matrix}\right) \Big\rangle\]d\overleftarrow{B_r}
  -\int_s^T\partial_{\mu}Z_r^{t,x,P_{\xi}}(y)dW_r.
  \end{split}
\end{equation}
Furthermore, if in addition Assumption (H4.2) is satisfied, then
there is a constant $C_p>0$ only depending on the Lipschitz constants of the coefficients,
such that, for all $t\in[0,T]$, $x,x',y,y'\in\mathbb{R}$, $\xi,\xi'\in L^{2}(\mathcal{G}_t;\mathbb{R})$,
$P$-a.s.,
\begin{equation}\label{eq6.41}
  \begin{split}
\emph{(i)}&\ E\[\sup_{s\in[t,T]}|\partial_{\mu}Y_s^{t,x,P_{\xi}}(y)|^p+(\int_t^T|\partial_{\mu}Z_s^{t,x,P_{\xi}}(y)|^2ds)^{\frac{p}{2}}\]\leq C_p,\
     \mbox{for all}\ p\in[2,p_0],\\
\emph{(ii)}&\ E\[\sup_{s\in[t,T]}|\partial_{\mu}Y_s^{t,x,P_{\xi}}(y)-\partial_{\mu}Y_s^{t,x',P_{\xi'}}(y')|^p
   +(\int_t^T|\partial_{\mu}Z_s^{t,x,P_{\xi}}(y))-\partial_{\mu}Z_s^{t,x',P_{\xi'}}(y')|^2ds)^{\frac{p}{2}}\]\\
&\ \leq C_p\(|x-x'|^p+|y-y'|^p+W_2(P_{\xi},P_{\xi'})^p\),\ \mbox{for all}\ p\in[2,\frac{p_0}{2}].
  \end{split}
\end{equation}
\end{proposition}
Using \eqref{eq6.2+3} and \eqref{eq6.2+7}, the proof is similar to that of Theorem \ref{th6.2} and Proposition \ref{prop6.2}, and so we omit it here.

\begin{remark} \label{reA.3}
In analogy to the proof in Proposition \ref{prop4.0}, it can easily be checked that
under the Assumptions (H5.1) and (H10.4) but without (H4.2), there exists a constant $C>0$ only depending on the Lipschitz
constants of the coefficients, such that for all $t\in[0,T]$, $x,x',y,y'\in\mathbb{R}^d$, $\xi,\xi'\in L^{2}(\mathcal{G}_t;\mathbb{R}^d)$,
 $P$-a.s.,
\begin{equation}\nonumber
  \begin{split}
E\[\sup_{s\in[t,T]}|\partial_{\mu}Y_s^{t,x,P_{\xi}}(y)|^2+\int_t^T|\partial_{\mu}Z_s^{t,x,P_{\xi}}(y)|^2ds\]\leq C.
  \end{split}
\end{equation}
\end{remark}

\end{document}